\documentclass[12pt]{article}
\usepackage{defs}
\usepackage[normalem]{ulem}
\usepackage{url}

\newcommand{\Spc}{\mathbf{Spc}}
\newcommand{\CAlg}{\mathbf{CAlg}}
\newcommand{\cU}{\mathcal{U}}
\newcommand{\Rings}{\mathbf{Rings}}

\newcommand{\Ho}{\mathrm{Ho}}
\newcommand{\Free}{\mathrm{Free}}
\newcommand{\bs}{\mathbf{s}}
\newcommand{\Subm}{\mathbf{RelProp}}\newcommand{\Mor}{\mathbf{Mor}}\newcommand{\Fib}{\mathrm{Fib}}
\newcommand{\Cone}{\mathrm{Cone}}
\newcommand{\Submor}{\mathbf{Subm}^{or}}
\title{Differential function spectra, the differential Becker-Gottlieb transfer, and applications to differential algebraic $K$-theory}
\author{Ulrich Bunke\thanks{NWF I - Mathematik,
Universit{\"a}t Regensburg,
93040 Regensburg,
GERMANY, ulrich.bunke@mathematik.uni-regensburg.de}
\, and David Gepner\thanks{Department of Mathematics, Purdue University, 150 N. University Street, West Lafayette, IN 47907 USA, dgepner@math.purdue.edu}
}

\newcommand{\cZ}{\mathcal{Z}}
\begin{document}
\maketitle
\begin{abstract}
{We develop differential algebraic $K$-theory for rings of integers in number fields and we construct a cycle map from geometrized bundles of modules over such a ring to the differential algebraic $K$-theory.  We also treat some of the foundational aspects of differential cohomology, including differential function spectra and the differential Becker-Gottlieb transfer.
We then state a transfer index conjecture about the equality of the Becker-Gottlieb transfer
and the analytic transfer defined by Lott.
In support of this conjecture, we derive some non-trivial consequences which are provable by independent means.}
\end{abstract}

\tableofcontents

\section{Introduction}\label{nullekhfkfweewfw}

The study of differential extensions of generalized cohomology theories merge the fields of homotopy theory and  differential geometry. Roughly speaking, a differential cohomology class combines the information on the underlying homotopy theoretic cohomology class and related characteristic forms with secondary data. It is this secondary information which makes differential cohomology theory relevant in applications.

\bigskip

Historically, the  first example of a differential cohomology theory was the differential extension of integral cohomology defined in terms of Cheeger-Simons characters \cite{MR827262}. For  a complex vector bundle   with connection on a smooth manifold Cheeger and Simons constructed   Chern classes in differential integral cohomology which combine the information on the underlying topological Chern classes with the characteristic forms given by Chern-Weil theory. The differential Chern classes contain non-trivial secondary invariants for flat bundles.

Much later, motivated by developments in mathematical physics, 
the differential extension of topological complex $K$-theory attracted a lot of attention
and was popularized among mathematicians in particular by the work of Freed and  Hopkins \cite{MR1769477}, \cite{MR1919425}, \cite{MR2664467}, \cite{MR2732065}, \cite{MR2231056}, \cite{MR2286784}. Differential topological complex $K$-theory captures the information on the underlying topological $K$-theory class of a vector bundle with connection in combination with the Chern character form, again given by Chern-Weil theory.

Differential extensions of bordism theories and Landweber exact cohomology theories have been constructed in \cite{MR2550094}.
The construction of the differential extension in all these examples starts of from a geometric cycle/relation model of the underlying generalized cohomology theory. A  homotopy theoretic construction of a differential extension for an arbitrary cohomology theory {was first} given by Hopkins and Singer in \cite{MR2192936}.

\bigskip

The main example of the present paper is the differential extension $\widehat{KR}$ of the cohomology theory $KR^{*}$ represented by the connective algebraic $K$-theory spectrum $KR$ of a number ring $R$. This differential extension  will be defined by a homotopy theoretic construction. Let us explain $\widehat{KR}^{0}$ in greater detail. Given a manifold $M$ we can consider the cohomology group $KR^{0}(M)$.  
 A differential cohomology  class $x\in \widehat{KR}^{0}(M)$  combines the information of an underlying cohomology class $I(x)\in KR^{0}(M)$ with a closed differential form $R(x) \in Z^{0}(\Omega A (M))$. Here $\Omega A$ denotes the differential forms with coefficients in the graded vector space $A:=K_{*}(R)\otimes \R$. The latter has been  calculated  explicitly by Borel \cite{MR0387496}, see Theorem   \ref{borel}. 
 The cohomology class of {the} form $R(x)$ 
  is equal to the image of $I(x)$ under the natural map  $KR^{0}(M)\to H^{0}(M;A)$ (see Definition \ref{jun216} for $I$ and $R$).

\bigskip

The construction which associates a differential cohomology class to a geometric object (e.g. differential Chern or $K$-theory class to a complex vector bundle with connection in the above examples) will be called a cycle map.  Usually, a cycle map is {a}  natural transformation of set-valued contravariant functors on the category of smooth manifolds
$$\{\mbox{geometric objects on $M$}\}\stackrel{\cycl}{\to} \{\mbox{differential cohomology  classes on $M$}\}\ .$$
{Often the geometric objects form monoids (e.g. bundles with direct sum), and in this case {one additionally} requires that the cycle map is additive.}
In most cases where the differential cohomology theory is build from a geometric model the construction of a cycle map is easy, sometimes even tautological. For example, in the model  \cite{MR2664467} of differential topological $K$-theory the cycle of complex vector bundle $V$ with connection $\nabla^{V}$ is simply the class $[V,\nabla^{V}]$ represented by the pair $(V,\nabla^{V})$. For geometric models of differential extensions of generalized cohomology theories it turned out to be difficult to obtain the functorial properties, e.g. differential extensions of natural transformations between different cohomology theories.

In this respect, the homotopy theoretic constructions  like \cite{MR2192936} or the one presented in the present paper are much better behaved. On the other hand,
it is turned out to be  difficult to construct cycle maps into the homotopy theoretic version of differential cohomology theory. This also applies to {the} case considered in the present paper. We are going to invest a great effort to construct a cycle map (Section \ref{aug0505}) which associates a differential algebraic $K$-theory class \begin{equation}\label{ccasjhaksjhcakashcsacas} \cycl(\cV,h^{\cV})\in \widehat{KR}^{0}(M)\end{equation}  to a locally constant sheaf $\cV$ of finitely generated
 $R$-modules with a geometry $h^{\cV}$ (Definition \ref{defpre}) on a smooth manifold $M$.

This discrepancy raises the question of comparison of different models.
Partial answers are given by the uniqueness theorems \cite{MR2608479}. The main assumption
there is that the coefficients of the cohomology theory are torsion groups in odd degree.
Therefore these uniqueness results apply well to ordinary cohomology, topological $K$-theory and bordism theories (see \cite{MR2740650} and \cite{MR2734150} for applications)
but not to the algebraic $K$-theory of a number rings where the non-torsion coefficients all live in degree zero and   odd degrees.

We now explain the cycle map in greater detail. Let $\cV$ be a locally constant sheaf of finitely generated $R$-modules on a manifold $M$. Then we can define a class
\begin{equation}\label{dqwdqwdqwqdqw3211}[\cV]\in KR^{0}(M)\ . \end{equation}  
The geometry $h^{\cV}$ is given by a collection of   metrics on the flat vector bundles on $M$
induced by $\cV$ via the various embeddings $R\hookrightarrow \R$ and $R\hookrightarrow \C$.
Using the characteristic forms for flat vector bundles with metrics introduced by Bismut-Lott in \cite{MR1303026} together with the explicit description of the graded vector space $A=K_{*}(R)\otimes \R$ by Borel \cite{MR0387496}, we define a closed form $\omega(\cV,h^{\cV})\in Z^{0}(\Omega A(M))$, see \eqref{cfr}. The  cycle class \eqref{ccasjhaksjhcakashcsacas} differentially refines \eqref{dqwdqwdqwqdqw3211} in the sense that
$$R(\cycl(\cV,h^{\cV}))=\omega(\cV,h^{\cV})\ , \qquad I(\cycl(\cV,h^{\cV}))=[\cV]\ .$$

\bigskip

A differential cohomology theory is in particular  a contravariant functor from the category of smooth manifolds to $\Z$-graded abelian groups. Its topological counterpart has wrong-way (Umkehr, integration) maps for suitably oriented fibre bundles $\pi:W\to B$. The notion of a differential orientation and the refinement of the wrong-way map to differential cohomology has been discussed in general in \cite{MR2192936}, {\cite{skript}} and precise versions for geometric models have been given in \cite{MR2664467},\cite{MR2550094}, \cite{MR2674652}, \cite{MR2602854}.

If one also has a descent of geometric objects along $\pi$ one can consider the compatibility with the cycle map. A (differential)  index theorem is the statement that the following diagram commutes
\begin{equation}\label{sep2401}\xymatrix{\{\mbox{(geometric) objects on $W$}\}\ar[d]^{descent}\ar[r]^(.4){\cycl}&\{\mbox{(differential) cohomology classes on $W$}\}\ar[d]^{wrong \:way \:map}\\\{\mbox{(geometric) objects on $B$}\}\ar[r]^(.4){\cycl}&\{\mbox{(differential) cohomology classes on $B$}\}}\ .\end{equation}
The most prominent example of an index theorem is the index theorem of Atiyah-Singer which 
compares the descent of vector bundles along $K$-oriented maps with the wrong-way map in topological $K$-theory. Its differential refinement has been given in  \cite{MR2602854}.

All wrong-way maps for generalized cohomology theories in topology   use the suspension isomorphisms
(i.e. the fact that a cohomology theory is represented by a spectrum) in an essential way.
From this it is clear that for a differential refinement of the wrong-way map one has to take into account properly  the spectrum aspect of  a differential cohomology theory.  A first construction of a differential cohomology spectrum has been given in \cite{MR2192936}.  In the present paper we  take the
spectrum aspect of differential cohomology as a starting point  and give in Definition \ref{may271} a direct and straight-forward homotopy theoretic construction of a differential function spectrum. The differential cohomology will then be defined in terms of homotopy groups of the differential function spectrum{, see Definition \ref{difdef1}.}
We consider our construction here is a short-cut to the construction of \cite{MR2192936}.

\bigskip

A particularly simple example of a wrong-way map is the Becker-Gottlieb transfer \cite{MR0377873} which does not require  any kind of orientation. It is defined for every smooth fibre bundle $$\pi:W\to B$$  with closed fibres and thus applies to
the cohomology theory $KR^{*}$ represented by the algebraic $K$-theory spectrum of a number ring  $R$. 
The Becker-Gottlieb transfer provides a homomorphism {(see \eqref{jun182})} $$\tr^{*}:KR^{0}(W)\to KR^{0}(B)\ .$$  {In the following we review  some constructions explained in detail in Subsection \ref{kljfelwkfwefwefwefewffewfef}.} Using \eqref{dqwdqwdqwqdqw3211} we define  the topological index  $$\ind^{top}(\cV):=\tr^{*} [\cV]\in KR^{0}(B)$$  of the sheaf $\cV$ of finitely generated $R$-modules on $W$.
On the geometric side one can descend   $\cV$   from $W$
to $B$ using the sheaf theoretic higher-derived images $R^{i}\pi_{*}(\cV)$. 
We define the analytic index by
$$\ind^{an}(\cV):=\sum_{i}(-1)^{i}[R^{i}\pi_{*}(\cV)]\in KR^{0}(B)\ .$$

In this case the corresponding index theorem is the Dwyer-Weiss-Williams index theorem stating    that
$$\ind^{top}(\cV)=\ind^{an}(\cV)\ .$$ It
 has been proven in this form in \cite{MR1982793} and refines  an earlier  characteristic class version of Bismut-Lott
\cite{MR1303026}. In order to define a differential Becker-Gottlieb transfer we must choose a pair {consisting} of a vertical Riemannian metric and a horizontal distribution for the bundle $\pi:W\to B$. In Definition \ref{may272}  such a pair is  
called a Riemannian structure.
One {of our} main results (presented in Section \ref{sep2201}, {Theorem \ref{may3102}}) is the construction of a differential Becker-Gottlieb transfer 
(and the verification of all its expected functorial properties) for smooth fibre bundles with closed fibres equipped with Riemannian structures. {In the case of differential algebraic $K$-theory it is a homomorphism
$$\hat \tr:\widehat{KR}^{0}(W)\to \widehat{KR}^{0}(B)\ .$$ We define the differential topological index by $$\hind^{top}(\cV,h^{\cV}):=\hat \tr(\cycl(\cV,h^{\cV}))\in \widehat{KR}^{0}(B)\ .$$}{As mentioned above, if we fix an embedding of the number ring  $R$ into $\R$ or $\C$, then $\cV$ induces   a flat vector bundle.    The $i$'th cohomology of the  fibrewise de Rham complex twisted    with such a vector bundle can be identified with  the corresponding vector bundle 
 for the sheaf $R^{i}\pi_{*}(\cV)$ on $B$.   Using the Riemannian structure on $\pi$, by  fibrewise Hodge theory, we can represent the cohomology by fibrewise harmonic forms. In particular, we obtain    an induced  $L^{2}$-metric and  eventually a   geometry
 $h^{R^{i}\pi_{*}(\cV)}$.  
 In analogy to the topological case we define the naive differential analytical index by $$\hind^{an}_{0}(\cV,h^{\cV}):=\sum_{i} (-1)^{i} \cycl(R^{i}\pi_{*}(\cV),h^{R^{i}\pi_{*}(\cV)})\in \widehat{KR}^{0}(B)\ .$$
 It is called naive since we {cannot} expect that it is equal to the differential topological index. 
 Indeed, 
 by the result of   \cite{MR1303026} we have
 $$R(\hind^{top}(\cV,h^{\cV}))- R(\hind^{an}_{0}(\cV,h^{\cV}))=d\cT\ ,$$
 where $R$ is as in Definition \ref{jun216}, and $\cT\in \Omega A^{-1}(B)/\im(d)$
 is a version of the Bismut-Lott higher torsion form. So following  \cite{MR1724894} we define
  the differential analytic index by $$\hind^{an}(\cV,h^{\cV}):=\hind^{an}_{0}(\cV,h^{\cV})+a(\cT)\in \widehat{KR}^{0}(B)\ ,$$
  where $a$ is  again one  of the structure maps of the differential cohomology theory $\widehat{KR}^{*}$ explained in Definition \ref{jun216}.}


{The transfer index conjecture (TIC), Conjecture \ref{may3103}, is the now the natural differential refinement of the Dwyer-Weiss-Williams index theorem and states 
that $$\hind^{top}(\cV,h^{\cV})\stackrel{?}{=}\hind^{an}(\cV,h^{\cV})\ .$$
Our main result here  is the development of all ingredients needed to state the TIC.  }
We are far from having a proof of this transfer index conjecture in differential algebraic $K$-theory. 
Instead, in Subsection \ref{sep0901} we are going to discuss special cases and consequences which can be proved by different means. Some of these consequences are deep theorems in different fields like arithmetic geometry or global analysis. Therefore we expect that a proof of Conjecture  \ref{may3103} will be quite complicated.

To finish this introduction let us explain the structure of the present paper and its main achievements.
\begin{enumerate}
\item {Section \ref{nullekhfkfweewfw} {consists only of} this introduction.}
\item  In Section \ref{sep2202} we review the axioms for a differential extension of a generalized cohomology theory. Our first main result is the {construction of the differential function spectrum}.
We define differential cohomology in terms of the homotopy groups of the differential function spectrum
and verify that this definition satisfies the required axioms. We further discuss the dependence on data.
\item In Section \ref{aug0505} we review the definition of the algebraic $K$-theory {of} number rings $R$ and its rational calculation due to Borel. {In Corollary \ref{efewklfjwlfwlfjelwkfwfewfewfw43242424324}, to   we associate differential data to number rings in a functorial way and define 
  differential algebraic $K$-theory $\widehat{KR}^{0}$ by specialization the general 
  Definition \ref{difdef1}.  The main result of Section \ref{aug0505}} is {the construction of a cycle map} which associates classes $\cycl(\cV,h^{\cV})\in \widehat{KR}^0(M)$ to 
locally constant sheaves of {finitely generated}   $R$-modules with geometry on $M$.
\item In Section \ref{sep2201} we construct, {for an {arbitrary} generalized cohomology theory},  {the differential Becker-Gottlieb transfer} and verify its functorial properties.
\item In Section \ref{sep2205} we { state the {T}ransfer {I}ndex {C}onjecture  \ref{may3103} }.
Our main contribution here is to work out  {links between this conjecture  and interesting deep results in other branches of mathematics}.
\item Our homotopy theoretic construction of the differential function spectrum takes place in an $\infty$-category of sheaves of spectra on smooth manifolds. {We introduce and motivate the notation and language at the place where it first appears in the text. In Section \ref{sep2210}    we give a comprehensive account of the theory and prove
most of the technical facts.  The first five subsections are basic for understanding the present paper, and the less experienced reader should  throw a glance to these subsections at an early stage. The remaining subsubsections
 starting with Subsubsection \ref{jun148} are extremely technical and should be read only if one wants to understands the details
 of the proofs where there material is used. 
}

\end{enumerate}

{The theory of differential cohomology developed in this paper can be refined to theories of multiplicative and twisted differential cohomology, and will be treated in future work \cite{buta,Bunke:2013kx,Bunke:2013uq,Bunke:2014fk}.}

\section{Differential function spectra}\label{sep2202}

\subsection{Differential cohomology -- the axioms}\label{aug0701}

{In this subsection we introduce the basic features of a differential cohomology theory in an axiomatic way along the lines of  \cite{MR2608479}, \cite{MR2365651}. The actual construction will be  pursued in Subsection \ref{fjwelkfjwelkfwfwfwf}. } 

\bigskip

A generalized cohomology theory is a homotopy invariant contravariant functor 
$$E^{*}:\Top^{op} \to \Ab_{\Z-gr}$$
from the category of topological spaces $\Top$ to the category of $\Z$-graded abelian groups  $\Ab_{\Z-gr}$
with some additional structures like suspension isomorphisms or Mayer-Vietoris sequences.
 Examples which are relevant for the present paper are ordinary integral cohomology $H\Z^{*}$, topological complex $K$-theory $KU^{*}$, and the algebraic $K$-theory $KR^{*}$ of a ring $R$.
 
Under suitable hypotheses a generalized cohomology theory $E^{*}$ can be represented by a spectrum $E$. {For a smooth manifold {(with corners)}  $M$ we can consider the underlying topological space $M_{top}$ and form the function spectrum $E^{M_{top}}$. Then the $E$-cohomology groups of $M$ can be expressed in terms of the function spectrum as \begin{equation}\label{dqwhdjwqhdwjkdwqhdkqhdk}
E^{k}(M)\cong \pi_{-k}(E^{M_{top}})
\end{equation} for all $k\in \Z$.}
{At this point we can work with any of the known models for the homotopy theory of spectra so we refrain from specifying one. We refer to {Example \ref{dqwjdklqwdwqwqdqdwq} for a description of such a model and
to Example 
  \ref{dkqwjdlkqwdjldjldjqwdoi2u1oe21e12e} for   the $\infty$-category of spectra}. For the Brown representability theorem we refer to \cite[Ch.9 ]{MR1886843}}. 
 
In the present paper we will always assume that a generalized cohomology theory is represented by a spectrum. Given a cohomology theory 
$E^{*}$ represented by a spectrum $E$ we have the notion of a differential extension $\widehat E^{*}$ of $E^{*}$. {It is a contravariant functor  \begin{equation}\label{hjkewhfwjhfewhfkjwhfjewhwhek}
\widehat E^{*}:\Mf^{op}\to  \Ab_{\Z-gr}
\end{equation}
 from the category of smooth manifolds with values again in $\Z$-graded abelian groups.}
The differential $E$-cohomology $\widehat E^*$ extends the restriction of {the functor}  $E^{*}$ from the category of topological spaces to the category of smooth manifolds by differential forms.
Note that a generalized cohomology {theory} $E^{*}$ is completely determined by its restriction to smooth manifolds {\cite{MR2775352}}, which implies that we can recover the cohomology theory $E^{{*}}$ from { the functor  $\widehat{E}^{*}$ together with the structure maps $R,I,a$ of differential cohomology, see Definition \ref{jun216}.} 


We start with a description of the data and basic constructions on which the  differential cohomology will depend. {For every abelian group $A$ we can define an Eilenberg-MacLane spectrum $H(A)$ such that
$$\pi_{*}(H(A))\cong \left\{\begin{array}{cc}A&{*}=0\\ 0&{*}\not=0\end{array}\right. \ ,$$
for example we have $H\Z:=H(\Z)$. More generally, to  a chain complex $A\in \Ch$ of abelian groups  we can naturally associate  an    Eilenberg-MacLane spectrum $H(A)$ with  $$\pi_{*}(H(A))\cong H_{*}(A)\ .$$ 
We refer to {Example \ref{jcsdkjcdksjhcdkschkdscer}} for details about the Eilenberg-MacLane spectrum construction.} 


For an abelian group $G$ we can define a Moore spectrum $MG$ {(not to be confused with $H(G)$)} which is characterized uniquely up to equivalence by the property that
$$H\Z_{*}(MG)\cong\left\{\begin{array}{cc}G&*=0\\0&*\not=0\end{array}\right\} \ .$$
We write $E{G}:=E\wedge M{G}$ for the smash product {of a spectrum} $E$
with the Moore spectrum $M{G}$ of {the abelian group} ${G}$. {Since the additive group  $\R$ is flat over $\Z$,} the effect of smashing with $M\R$ on the level of homotopy groups is tensoring with $\R$
$$\pi_{*}(E\R)\cong \pi_{*}(E)\otimes \R\ ,$$ see \cite[{Eq. (2.1)}]{MR551009}. 
There is a natural map $\epsilon_{\R}:E\to E\R$ of spectra which induces the canonical map
$$\pi_{*}(E)\to \pi_{*}(E\R)\cong \pi_{*}(E)\otimes \R\ , \quad x\mapsto x\otimes 1$$
on the level of homotopy groups.

\bigskip

In order to define the notion of a  differential extension of $E$ we
choose a chain complex $A\in \Ch_{{\R}}$ of real vector spaces  together with an equivalence of spectra $$c:E\R\stackrel{{\simeq}}{\to} H(A)\ .$$  
\begin{ddd}\label{aug1001}
The triple $(E,A,c)$ will be called the {differential data}
 for  the differential extension of $E$.
\end{ddd}
The complex $A$ is, of course, not arbitrary but necessarily satisfies
$$\pi_{*}(E)\otimes \R\cong H_{*}(A)\ .$$
From the data $(E,A,c)$ we {obtain the composition}
\begin{equation}\label{apr040113}E^{*}\stackrel{\epsilon_{\R}}{\to} E\R^{*}\stackrel{c}{\to} H(A)^{*}\end{equation}
{of natural  transformations between  generalized cohomology theories.}

 For a smooth manifold $M$ we  consider the    complex $\Omega A(M)$ of smooth $A$-valued differential forms. {It is defined (see  Definition \ref{deRham}) as the evaluation of the tensor product of $\Ch_{\R }$-valued sheaves  $\Omega\otimes_{\R}{L}(\underline{A})$, where $\Omega$ is the sheaf of de Rham complexes and  ${L}(\underline{A})$ is the locally constant sheaf generated by $A$. Here we will consider $A$ as a cohomological chain complex  using the convention $A^{n}:=A_{-n}$.
} {Note that if  $M$ is compact, then we can avoid this sheaf {theoretic} description and simply  have an isomorphism of chain complexes $$\Omega A(M)\cong \Omega(M)\otimes_{\R} A\ ,$$ see Remark \ref{jkdhwqjkdhkwqdwqhdqwdkjhwqdqwdwqdwqdwd}.}  
 {The cohomology of the complex of $A$-valued forms $\Omega A(M)$ is related with the cohomology $H(A)^{*}(M)$ defined topologically by 
 the} de Rham isomorphism.    \begin{equation}\label{jul0315}j:H^{*}(\Omega A(M))\stackrel{\cong }{\to} H(A)^{*}(M) \ .\end{equation} It is induced from {a map between smooth spectra}
  \eqref{hgjhgdjhgqwd873iu2kjd31}   by evaluation at $M$ and applying $\pi_{-*}$.   The de Rham isomorphism is natural in the manifold $M$ (see Remark \ref{lkqwdjqwkld89uoij3d3d})   and the chain complex $A$  by Corollary \ref{lqejhdlwdjqwodu321oidwqd}. 
\begin{rem}{\rm {The de Rham isomorphism is induced by a natural map of spectra \begin{equation}\label{fwefewfefewfewfwfef}   j:H(\Omega A(M))\to H(A)^{M_{top}}\ , \end{equation} see Proposition \ref{canequ}. This will be needed in the constrution of differential cohomology in Subsection \ref{fjwelkfjwelkfwfwfwf}.   The effective treatment of  this map uses elements of sheaf theory with values in $\infty$-categories which will be discussed later. Since the details are not needed to state the axioms for differential cohomology we will not present them at this place.} 
 }\hB\end{rem}
 If $C$ is a {cohomological} chain complex, then we let $$Z^{k}(C):=\ker(d:C^{k}\to C^{k+1})\subseteq C^{k}$$ denote the {group} of cycles in degree $k$. We have a natural map $$[\dots]:Z^{k}(C)\to H^{k}(C)$$ which associates to a cycle the cohomology class it represents.
 The de Rham isomorphism thus  induces a natural transformation
 \begin{equation}\label{apr040213}Z^{*}(\Omega A)\xrightarrow{[\dots ]} H^{*}(\Omega A) \stackrel{j}{\to}H(A)^{*}\end{equation}
 between contravariant functors from manifolds to $\Z$-graded abelian groups.
  
 \bigskip
 
We now fix  differential  data $(E,A,c)$ for $E$.
 
  \begin{ddd}\label{jun216}
 A differential extension  $\widehat E^{*}$  of the cohomology theory $E^{*}$ is a tuple $(\widehat E^{*},I,R,a)$ consisting  of
\begin{enumerate}
\item a functor $\widehat E^{*}:\Mf^{op}\to \Ab_{\Z-gr}$
from smooth manifolds to $\Z$-graded abelian groups,
\item a natural transformation $I:\widehat E^{*}\to E^{*}$ {({the} underlying class)},
\item a natural transformation $R:\widehat E^{*}\to Z^{*}(\Omega A)$, {({the} curvature)}
\item and a natural transformation $a:\Omega A^{*-1}/\im(d)\to \widehat E^{*}$.
\end{enumerate}
These objects must satisfy the following axioms:
\begin{enumerate}
\item[i.] $R\circ a=d$.
\item[ii.] The diagram
$$\xymatrix{\widehat E^{*}\ar[r]^{I}\ar[d]^{R}&E^{*}\ar[d]^{{\eqref{apr040113}}}\\Z^{*}(\Omega A)\ar[r]^{{\eqref{apr040213}}}&H(A)^{*} }$$
commutes.
\item[iii.] The sequence
\begin{equation}\label{apr140213}E^{*-1}\xrightarrow{j^{-1}\circ c\circ \epsilon_{\R}}\Omega A^{*-1}/\im(d)\stackrel{a}{\to} \widehat E^{*}\stackrel{I}{\to} E^{*}\to 0\end{equation}
is exact.
\end{enumerate}
The kernel $$\widehat E^{*}_{flat}:=\ker\left(R:\widehat E^{*}\to Z^{*}(\Omega A)\right)\subseteq \widehat E^{*}$$ {of the curvature}
is called the flat part of $\hat E^{{*}}$.
\end{ddd}

\begin{rem}{\rm 
As explained in \cite{MR2608479}, in general these axioms do not characterize the differential extension of $E^{*}$ uniquely. For uniqueness one needs additional structures like an integration. We will not discuss these structures in the present paper. {But note that the differential extensions constructed in the present paper using the Hopkins-Singer method admit these additional structure naturally, see \cite{skript}.}
} \hB \end{rem}

\bigskip

{The differential cohomology functor $\widehat{E}^{*}$ is not homotopy invariant. Its deviation from homotopy invariance is measured by the homotopy formula \cite[(1)]{MR2608479} {which is a formal consequence of the axioms.} Let $$x\in \widehat{E}^{*}([0,1]\times M)$$ and $f_{i}:M\to [0,1]\times M$ be the inclusions corresponding to the endpoints of the interval. Then we have
\begin{equation}\label{feb2710}f_1^{*}(x)-f_0^{*}(x)=a\left(\int_{I\times M/M} R(x)\right)\ .\end{equation}}

\begin{rem}{\rm

The first example of a differential extension of a cohomology theory was the differential extension
of integral cohomology represented by the Eilenberg-MacLane spectrum $H\Z$.
In this example constructed by Cheeger-Simons \cite{MR827262}  the {groups $\widehat{H\Z}^{*}(M)$  were  realized as   groups} of differential characters.  {See Example \ref{vsdsdlkvkjlvjew09iu0wef} for an alternative {construction}.}
We refer to \cite{MR2664467}, \cite{MR2550094} and the literature cited therein for further examples constructed using geometric cycles.
A general construction of a differential extension for an arbitrary cohomology theory has been given by Hopkins-Singer in 
\cite{MR2192936}.

In the present paper we give a general construction of a differential extension $\widehat E^{*}$ in terms of the differential function spectrum, Definition  \ref{may271}. The evaluation on $M$ of the differential function spectrum for $E$ is a spectrum $\Diff(E)(M)$ which
refines the function spectrum {$E^{M_{top}}$} by Eilenberg-MacLane spectra associated with
$A$-valued differential forms on $M$ such that
$$ \widehat E^{0}(M):=\pi_{0}(\Diff(E)(M))$$ is the underlying functor of the differential extension. 
In order to get the differential cohomology in all degrees we use shifts
$$\widehat E^{n}(M):=\pi_{0}(\Diff(\Sigma^{n}E)(M))\ .$$
Using a slightly different
language, the differential function spectrum has first been constructed in \cite{MR2192936}.}
\hB\end{rem}

\subsection{The construction of the differential function spectrum}\label{fjwelkfjwelkfwfwfwf}

In this section we construct for every choice of data $(E,A,c)$ as in Definition  \ref{aug1001} a differential function spectrum $\Diff(E,A,c)$. We further define the differential cohomology   in terms of homotopy groups of $\Diff(E,A,c)$.
 {In order to formulate the construction properly we will use the language of $\infty$-categories. We will explain and motivate  our usage of $\infty$-categories in a non-technical way.   \bigskip
 
 The following toy example should serve as a motivation.
 
 \begin{ex}\label{vsdsdlkvkjlvjew09iu0wef}{\rm
  Assume that we want to construct the differential extension $\widehat{H\Z}^{p}$ of integral cohomology for some $p\in \nat$. We choose this example since the graded vector space $A$ is just $\R$ in degree zero and the  spectra $H\Z$ and $H\R$ are  well-known. In order to avoid shifts at this place on the one hand,  but to define something interesting on the other hand,  we consider a cut-off at degree $p\in \nat$. Later we will  concentrate on the case $p=0$.
   
   So 
  the corresponding data is $(H\Z,\R,c:H\Z\to H\R)$. Implicitly we fix some point-set model of spectra.

     For a manifold $M$  we can then define the function spectrum $H\R^{M_{top}}$ whose homotopy groups are the real cohomology groups of $M$, see \eqref{dqwhdjwqhdwjkdwqhdkqhdk}.
We further let $\sigma^{\ge p}\Omega(M)$ be the truncation of the de Rham complex to degrees $\ge p$. 
We want to define a spectrum $\Diff^{p}(H\Z,\R,c)(M)$ so that it fits into a 
  homotopy pull-back diagram of spectra \begin{equation}\label{gdhghjewgdjwedgewdiuzi}
\xymatrix{\Diff^{p}(H\Z,\R,c)(M)\ar[r]\ar[d]&H(\sigma^{\ge p}\Omega(M))\ar[d]\\
H\R^{ M_{top}}\ar[r]^{j^{-1}}&H(\Omega(M))}\ .
\end{equation}
It then has 
  the correct homotopy type  so that \begin{equation}\label{t354t343t43t34t3t3t34t}
\pi_{-p}(\Diff^{p}(H\Z,\R,c)(M))\cong \widehat{H\Z}^{p}(M)\ .
\end{equation} The corresponding calculations are similar to the ones performed in Subsection \ref{kjefklwefwefwefwefewfwefwef}. The right vertical map in \eqref{gdhghjewgdjwedgewdiuzi}  is induced by the canonical inclusion $\sigma^{\ge p }\Omega(M)\to \Omega(M)$ provided we  understood the Eilenberg-MacLane spectrum construction $H$ as a functor from chain complexes to our chosen point-set model of  spectra. The lower horizontal arrow is an inverse to the de Rham equivalence \eqref{fwefewfefewfewfwfef}.

   The point set level construction of the diagram  \eqref{gdhghjewgdjwedgewdiuzi} involves various choices. First of all we must choose  the spectrum level Rham equivalence $j:H(\Omega(M))\to H\R^{ M_{top}}$ and its inverse. We can only expect to fix these maps up to contractible choice. In order to construct the homotopy pull-back we must replace the pull-back diagram in our point-set model of spectra fibrantly. Finally, a point-set level functoriality of the  Eilenberg-MacLane construction $H$ is problematic. But note that the pull-back-diagram \eqref{gdhghjewgdjwedgewdiuzi}  is fixed uniquely up to contractible choice in an appropriate space of such diagrams. So the differential function spectrum $\Diff^{p}(H\Z,\R,c)(M)$  can be defined uniquely up to contractible choice. 

Now note that we actually want to consider the differential function spectrum as a contravariant functor which associates to a manifold
$M$ a spectrum $\Diff^{p}(H\Z,\R,c)(M)$. Because of the choices involved in the construction we can not expect such a functoriality on the point-set level. Our construction gives the following: 
 
\begin{enumerate}
\item To each manifold we associate a spectrum $\Diff^{p}(H\Z,\R,c)(M)$.
\item To each smooth map $f:M\to M^{\prime}$ we can associate a map between spectra $f^{*}:\Diff^{p}(H\Z,\R,c)(M^{\prime})\to \Diff^{p}(H\Z,\R,c)(M)$.
\item To each pair of composable 
morphisms $f:M\to M^{\prime}$ and $f^{\prime}:M^{\prime}\to M^{\prime\prime}$ we can associate a homotopy between $f^{*}\circ   f^{\prime,*} $ and $ (f^{\prime}\circ f)^{*}$.
\item  We have higher homotopies for  chains of composable morphisms. 
\end{enumerate} 
The dependence of $ \Diff^{p}(H\Z,\R,c)$ on the choices is not problematic since at the end we are interested in its homotopy groups. The functor $M\mapsto \widehat{H\Z}^{p}(M)$ given by \eqref{t354t343t43t34t3t3t34t} is a well-defined contravariant functor from manifolds to abelian groups.
\hB
}
\end{ex}

{The language of $\infty$-categories provides the right framework to formulate such  functoriality up to contractible choice in a rigorous and efficient way}.  In the following we explain roughly how this works.  
We model $\infty$-categories by quasi-categories (a standard reference is \cite{MR2522659}) which are simplicial sets
satisfying the inner horn filling conditions.  Functors between $\infty$-categories are    just  maps between simplicial sets. 
If $\bC$ is a small ordinary category, then its nerve $\Nerve(\bC)$ is first of all a simplicial set, but   actually turns out to   be an  $\infty$-category.

For the domain of our functors,
we replace the opposite of the category of manifolds $\Mf$ (more precisely a small skeleton consisting of submanifolds of some $\R^{\infty}$) by its nerve $\Nerve(\Mf^{op})$. The target of the functors is the stable infinity category (a reference for this notion is Ch.1 of \cite{highalg}) of spectra $\Sp$, see Definition \ref{dkqwjdlkqwdjldjldjqwdoi2u1oe21e12e}. It is characterized by the property that its homotopy category is the usual stable homotopy category.
\begin{rem}{\rm We define $\Sp$ as the stabilization of the $\infty$-category of pointed spaces. An alternative way   to construct $\Sp$ is to start  with some point-set model for the category of spectra $ \check{ \Sp}$, e.g. the one of \cite{MR513569} or \cite{MR1417719}. We will need such a model when we consider the Becker-Gottlieb transfer. Some details will be given in Remark \ref{dqwjdklqwdwqwqdqdwq}. As a first step one takes the nerve $\Nerve(\check{\Sp})$. Then, in $\infty$-categories, one inverts 
  the appropriate class of stable weak equivalences $W$ in order to get an equivalence $$\Sp\simeq \Nerve(\check{ \Sp})[W^{-1}]\ .$$ Here for any $\infty$-category $\bC$ and a set of morphism  $W$ the $\infty$-category $\bC[W^{-1}]$ is characterized by a universal property. We refer to  Subsection 
  \ref{jul0501}  for details of this construction. 
  }\hB\end{rem}

The   $\infty$-category of functors between two $\infty$-categories $\bC,\bD$ is defined as  the simplicial mapping space
$$\Fun(\bC,\bD):=\Map(\bC,\bD)\ .$$
So the $\infty$-category of contravariant  functors from manifolds to spectra  is the  mapping space 
$$\Fun(\Nerve(\Mf^{op}),\Sp):=\Map(\Nerve(\Mf^{op}),\Sp)\ .$$ 
This definition makes the   data which we  roughly described  at the end of the Example \ref{vsdsdlkvkjlvjew09iu0wef} precise. 
 
We will construct the differential function spectrum functor $\Diff(E,A,c)$ as a point in the simplicial set
$\Map(\Nerve(\Mf^{op}),\Sp)$ which is well-defined up to contractible choice (we also say essentially unique).
 } 

\bigskip

{
In the present paper we will  not only  consider contravariant functors from the category of manifolds with values in spectra $\Sp$, but also in  various other $\infty$-categories like chain complexes, spaces, commutative monoids etc. 
Let   $\bC$ denote some $\infty$-category.  
\begin{ddd} 
We will call   functors from $\Nerve(\Mf^{op})$ to $\bC$ smooth objects in $\bC$ or $\bC$-valued presheaves, and will denote the $\infty$-category  of them   by
$$\Sm(\bC):=\Fun(\Nerve(\Mf^{op}),\bC)\ .$$
\end{ddd}

{A typical aspect of a cohomology theory is the   Mayer-Vietoris sequence. It will be  encoded in the present set-up as a descent or sheaf condition. {To this end we equip the category $\Mf$ with a Grothendieck pre-topology  which associates to   every manifold $M\in \Mf$ the collection  of covering families of $M$. Here a covering family of $M$ is an at most countable familiy
$\cU:=(U_{\alpha})_{\alpha\in A}$ of open subsets of $M$ such that $\bigcup_{\alpha\in A}U_{\alpha}=M$.
}
{If $\cU$ is a  covering family of $M$, then we can form the manifold $U:=\coprod_{\alpha\in A}U_{\alpha}$
which has a natural map $U\to M$. We can further define   the simplicial   manifold $U^{\bullet}_{M}$ called the \v{C}ech nerve. {In simplicial degree $n\in \nat$  it is given by}    \begin{equation}\label{f4uhfjkhhfk3hf3fmnj334f487fz3487fz34ff34f3f}
\underbrace{U_{M}^{n}:=U\times_{M}\dots\times_{M}U}_{n+1\:\: factors}\ .
\end{equation}


\begin{ddd}\footnote{{This definition corrects the  definitions given in \cite{Bunke:2013uq} and \cite{Bunke:2013kx}, where the first condition was forgotten.}}\label{jun215}
A smooth object $F\in \Sm(\bC)$ satisfies descent (or is called a sheaf), if the following two conditions are satisfied:
{\begin{enumerate}
\item For every at most countable  family of manifolds $(M_{\alpha})_{\alpha\in A}$ the natural map
$$F(\coprod_{\alpha} M_{\alpha})\to \prod_{\alpha\in A} F(M_{\alpha})$$ is an equivalence.
\item
For every   covering family $\cU  $ of a manifold $M$   the natural morphism {\begin{equation}\label{gergrekgjreglkrejg09i5034534543535435345}
F(M)\to \lim_{\Nerve(\Delta)} F (U_M^{\bullet} )
\end{equation}} is an equivalence.
\end{enumerate}}
\end{ddd}

Implicitly, we assume that the limits appearing the definition above exist in $\bC$.  

\begin{rem}\label{ljfewfewfwfewfewfewfwf67678463842344234234}{\rm
In the technical part of the present paper every object or diagram belongs to some specified $\infty$-category. So it should be always clear in which $\infty$-category a limit or colimit is taken. For example, the limit in \eqref{gergrekgjreglkrejg09i5034534543535435345} is understood in $\bC$. The limit or colimit  is an object fitting into an appropriate conical extension of the diagram satisfying a universal property. In the case of ordinary categories, the limit object is determined  up to isomorphism while for proper $\infty$-categories it is fixed up to equivalence.
For details we refer to  \cite[Ch. 4]{MR2522659}.

With the exception of some motivating remarks we will avoid {talking} about homotopy limits or homotopy colimits. But note the following. Consider an  ordinary category $\bC$      with a set of {arrows} $W$ and a diagram $X\in \bC^{I}$. Then using the identification \eqref{lkdjlkqwjdldwqdqwd} we can consider it as an object of $\Fun(\Nerve(I),\Nerve(\bC))$. Let $\iota:\Nerve(\bC)\to \Nerve(\bC)[W^{-1}]$ be the localization. One can assume that $ \Nerve(\bC)[W^{-1}]$ has the same objects as $\Nerve(\bC)$.
Then the objects $\lim_{\Nerve(I)} (\iota X)$ and $\colim_{\Nerve(I)}(\iota X)$  can be considered  as the homotopy limit $\holim_{I}X$ or homotopy colimit $\hocolim_{I}X$ of the diagram. Explicit examples of calculations of limits and colimits of this kind are discussed in Examples \ref{kdjqwdkljqlwdjlqwjdlwqdwqdwqdqwdwqdwqddqwdqd} and \ref{fjfwlefjlwefjlewfewfewffewfewfef}.
}\hB\end{rem}
We let 
  $$\Sm^{desc}(\bC)\subseteq \Sm(\bC)$$  be the  full $\infty$-subcategory of sheaves. We refer to  Subsection \ref{smoothobjects} for more details about sheaves.}  
  
  \bigskip
  
  Thus a reformulation of the goal of   the present subsection is that we want to  construct  
an essentially unique object 
$$\Diff(E,A,c)\in \Sm^{desc}(\Sp)\ .$$}Often we will
use the shorter notation $\Diff(E)$ instead of $\Diff(E,A,c)$.

\bigskip

{Our goal is now to set up the pull-back diagram \eqref{difdef} below. 
It is the functorial analog of \eqref{gdhghjewgdjwedgewdiuzi}.  
We start with the lower left corner.

\bigskip

There is a  particularly well-behaved class of $\infty$-categories which are called presentable  \cite[Def. 5.5.0.1]{MR2522659},  see Example \eqref{dlkqjwdlkjqwljqwdqwdqwdqw}. The basic example of such a category is the category of spaces $\Spc$. One way to construct this $\infty$-category  is to start with the nerve of the category of simplicial sets $\Nerve(\sSet)$.  If $W$ denotes  the weak equivalences, then we define the $\infty$-category of spaces by  
$$\Spc:= \Nerve(\sSet)[W^{-1}]\ ,$$ see Definition \ref{lkdjlkqwjwqlkdjwqdwqdwqd}.
A presentable $\infty$-category $\bC$ is tensored and cotensored over the $\infty$-category of spaces. In particular we have a functor \begin{equation}\label{jchjsachakscsacscaca}
\Spc^{op}\times \bC\to \bC\ , \quad (X,C)\mapsto C^{X}\ .
\end{equation}

  We can use topological spaces, i.e. objects of $\Top$, in order to present spaces. Indeed, the adjunction \begin{equation}\label{f2fhfjkh23f789f32ff32f}
|-|:\sSet\leftrightarrows \Top:\sing
\end{equation}
 between geometric realization and the singular complex functor induces an equivalence of $\infty$-categories \begin{equation}\label{hdqhjwdqdwqddqdqdd134123}
\Nerve(\Top)[W^{-1}]\simeq\Spc\ ,
\end{equation} see Example \ref{lkejfewlkfewfwffewfwe34}.
We let $X_{\infty}\in \Spc$ denote the image of a topological space $X\in \Top$  {under the composition  \begin{equation}\label{qwdqwdwqdqdu3hkj321e12e}
\Nerve(\Top)\to \Nerve(\Top)[W^{-1}]\stackrel{ \eqref{hdqhjwdqdwqddqdqdd134123}}{\simeq}\Spc\ .
\end{equation}     So the functor
$\Mf\to \Top$ which sends   a manifold $M$ to its underlying topological space $M_{top}$ induces a functor
\begin{equation}\label{grgegg434353535345} \Nerve(\Mf)\to \Nerve(\Top)\stackrel{\eqref{qwdqwdwqdqdu3hkj321e12e}}{\to}   \Spc\ , \quad M\mapsto M_{top,\infty}\ .\end{equation}
Together with \eqref{jchjsachakscsacscaca} it provides a functor
$\Nerve(\Mf^{op})\times \bC\to \bC$.

\begin{ddd}\label{kjdlkqwwqddwqdqdq} The 
  adjoint  of this functor will be denoted by 
$\Funk:\bC\to \Sm(\bC)$. It 
associates to an object $C\in \bC$ the smooth function object $\Funk(C)\in  \Sm(\bC)$.
\end{ddd}
By definition, for $C\in \bC$ and $M\in \Mf$ we have
$$\Funk(C)(M)\simeq C^{M_{top,\infty}}\ .$$
\begin{rem}\label{klfjwklfewlfewfewlfewfew23423454943}{\rm {We will often simplify the notation and write
$C^{M}$ for $C^{M_{top,\infty}}$.}}\hB \end{rem}
{By Lemma \ref{kdqqwdwqjldwqdqdqwdq343224342342242}} the smooth function object  $\Funk(C) $ satisfies descent, i.e. we actually have defined a functor
 \begin{equation}\label{}
\Funk:\bC\to \Sm^{desc}(\bC)\ .
\end{equation}
For the details we refer to Subsection  \ref{smoothobjects}.  

\bigskip

Note that the diagram \eqref{gdhghjewgdjwedgewdiuzi} involves the   functor $M\mapsto H\R^{M_{top}}$. In order to make this precise we observe that  the $\infty$-category $ \Sp $ is an example of a presentable $\infty$-category to which we can apply the above construction.  For a spectrum $E$ and a manifold $M$, by definition, we have an equivalence
$\Funk(E)(M)\simeq E^{ M_{top,\infty}}$
in $\Sp$. This finishes the description of the lower left corner in \eqref{difdef}.
}

\begin{rem}{\rm  Recall that  $E^{*}$ denote the cohomology theory represented by a spectrum $E$.
It can be expressed in terms of the smooth function spectrum as follows. Let $M$ be a manifold and $n\in \Z$. Then we have a canonical isomorphism
\begin{equation}\label{hgr1}\pi_{n}(\Funk(E)(M))\cong E^{-n}(M)\ .\end{equation}}

\end{rem}


\bigskip

{We now set-up the right part of \eqref{difdef}.  
For a chain complex $A\in \Ch_{\R}$ of real vector spaces we {introduce the de Rham complex $$\Omega A\in \Sm(\Nerve(\Ch)) $$ of smooth forms with coefficients in $A$  in Definition \ref{deRham}.}} If $$C: \quad \cdots \to C^{-1}\to C^{0} \to C^{1}\to \cdots$$ is a {cohomological} chain complex, then 
we let \begin{equation}\label{gerggrklgrjeglregjelrgjerlre}
\sigma C :\quad \cdots\to 0\to C^{0}\to C^{1}\to \cdots
\end{equation}
 be its truncation to non-negative degrees.
Note that we have a natural inclusion
\begin{equation}\label{frfwefwewfwf234}\sigma C\to C\ .\end{equation}
{This truncation functor $\sigma:\Ch\to \Ch$ extends {objectwise} to
smooth chain  complexes, so we may form the truncated de Rham complex  with coefficients in $A$ $$\sigma \Omega A\in  \Sm(\Nerve(\Ch))\ .$$} The inclusion morphism of   smooth chain complexes \begin{equation}\label{frfwefwewfwf2kjkjlkjljljlkjaljldw34} \sigma \Omega A\to \Omega A  \end{equation}  will eventually give the right vertical arrow in the square (\ref{difdef}) below.
{
For the category of chain complexes $\Ch$ we consider the set $W$ of quasi-isomorphisms, see Example \ref{lkjklfjqwlfjqlwkjfwqfwqfq}. The $\infty$-category
$\Nerve(\Ch)[W^{-1}]$ is a presentable stable $\infty$-category whose homotopy category models the derived category of
chain complexes of abelian groups \cite[Ch. 1.3]{highalg}. The functor \begin{equation}\label{cljlljlkjljke9} \Nerve(\Ch)\to \Nerve(\Ch)[W^{-1}]\ , \quad {C\mapsto C_{\infty}}  \end{equation} extends objectwise to a functor 
$$\Sm(\Nerve(\Ch))\to \Sm(\Nerve(\Ch)[W^{-1}])$$ between the $\infty$-categories of smooth objects.
We will use the notation $$ \Omega A_{\infty}\ , \sigma \Omega A_{\infty}\in \Sm(\Nerve(\Ch)[W^{-1}])$$
for the images of $\Omega A$ and $ \sigma \Omega A$ under this morphism. 
By Lemma \ref{feb1601},   $ \Omega A_{\infty}$ and $ \sigma \Omega A_{\infty}$
statisfy descent, i.e. we have
\begin{equation}\label{fwfwffewkjhkjhrk23r2r} \Omega A_{\infty}\ ,\sigma \Omega A_{\infty}\in \Sm^{desc}(\Nerve(\Ch)[W^{-1}])\ .\end{equation}
 }


{We consider  the   Eilenberg-MacLane spectrum functor} introduced in  Definition   \ref{jdjwqldjwqldjwqljdwqopiop12e}.  
\begin{equation}\label{fweffefewfewfewfewfwf}
H:\Nerve(\Ch)[W^{-1}] \to { \Sp}\ , \quad C\mapsto \map(\Z_{\infty},C) \ .
\end{equation} 
Here $\Z$ is considered as a chain complex sitting in degree zero, and $\map(-,-)$ is the mapping spectrum functor  \eqref{kjdkqwjdlqkwdjlwqdqwdqwd} $$\Nerve(\Ch)[W^{-1}]^{op}\times \Nerve(\Ch)[W^{-1}]\to \Sp\ .$$

{
\begin{rem}{\rm  Similar to the simplification explained in   Remark \ref{klfjwklfewlfewfewlfewfew23423454943}, in order to simplify the notation we write $$H(C):=H(C_{\infty})$$ for a chain complex $C\in \Ch$.}\hB \end{rem}}

{The objectwise application of $H$  provides}  a functor ({also denoted by $H$})
\begin{equation}\label{dqwwqdwqdq3244324242}
H:\Sm (\Nerve(\Ch)[W^{-1}])\to \Sm ({ \Sp})\ .
\end{equation}

  Since
 \eqref{fweffefewfewfewfewfwf}   preserves limits  the functor \eqref{dqwwqdwqdq3244324242} preserves   sheaves and thus restricts to
$$H:\Sm^{desc}(\Nerve(\Ch)[W^{-1}])\to \Sm^{desc}({ \Sp})\ .$$
 The morphism \eqref{frfwefwewfwf2kjkjlkjljljlkjaljldw34}   induces,   after application of $H$, the right vertical map in \eqref{difdef}.}

\bigskip

{In order to set-up  \eqref{difdef} it remains to define the lower horizontal arrow.}  Given a chain complex of real vector spaces $A\in \Ch$ we can   form two smooth spectra with descent, namely  $$\Funk(H(A))\ ,\: H(\Omega A) \in    \Sm^{desc}({ \Sp})\ .$$

 They are related by a  canonical equivalence called the de Rham equivalence.

\begin{prop}\label{canequ} There is a canonical equivalence 
\begin{equation}\label{hgjhgdjhgqwd873iu2kjd31}
j:H(\Omega A)\stackrel{\simeq}{\to}\Funk(H(A))
\end{equation}
 of smooth spectra.
 \end{prop}
 \proof
 The de Rham equivalence is given by the composition 
 \eqref{jhekjhekjhkd12983ud128oidkl1d1d}
below.

\begin{rem}\label{lkqwdjqwkld89uoij3d3d}{\rm  By considering the de Rham equivalence as a morphism between smooth spectra we 
 encode the naturality in the manifold $M$ of the  spectrum level de Rham equivalence \eqref{fwefewfefewfewfwfef}.  
The construction of  the de Rham equivalence as the composition \eqref{jhekjhekjhkd12983ud128oidkl1d1d} below does not yield a preferred direction of this map. We choose the direction of the map as in \eqref{hgjhgdjhgqwd873iu2kjd31}
since in a preliminary version of this paper we had a different construction    via integration of forms which had this preferred direction.
  }\hB \end{rem}

The construction of the Rham equivalence is a simple consequence of the homotopy invariance  of the de Rham complex $\Omega A_{\infty}\in \Sm^{desc}(\Nerve(\Ch)[W^{-1}])$.  We start with introducing the necessary notions {and refer {to Subsection \ref{kjdkljkdljlqwdwqdqwdwqddqwdqdwqdwqd}} for more details}.

\bigskip

Let $\bC$ be a presentable $\infty$-category and $I:=[0,1]$ denote the unit interval. \begin{ddd}\label{jdqwdjqwldqdwqdwqdqd}
A smooth object $F\in \Sm(\bC)$   is homotopy invariant if the projection $  M\times I\to M$ induces an equivalence $F(M)\stackrel{\simeq}{\to} F(  M\times I)$ for all manifolds $M$.  \end{ddd}
We let $$\Sm^{h}(\bC)\subseteq \Sm(\bC)\ , \quad \Sm^{desc,h}(\bC)\subseteq \Sm^{desc}(\bC)$$ denote the full subcategories of homotopy invariant smooth objects of $\bC$, and homotopy invariant smooth objects of $\bC$ which in addition satisfy descent.  
 
{
We now assume that $\bC\in \{\Spc,\Sp, {\Nerve}(\Ch)[W^{-1}]\}$ or $\bC=\Nerve(\cC)$ for some complete and cocomplete one category $\cC$. This covers all cases needed in the present paper. 
More generally we can assume \ref{wijfwijeoifjeowifjw}.}
\begin{lem}\label{djqwjkqwdqwdqdqd213}
The functors  $$\Funk:\bC\to \Sm^{desc,h}(\bC)\ ,\quad \ev_{*}:\Sm^{desc,h}(\bC)\to \bC$$ (evaluation at $*$) are inverse to each other equivalences of $\infty$-categories. In particular, for $F\in \Sm^{desc,h}(\bC)$ there is a natural equivalence
 \begin{equation}\label{dkqwkdqwdjdljdlkwd9879}
\Funk(F(*))\stackrel{\simeq}{\to} F\ . \end{equation}
\end{lem}
\proof  This will be shown as Lemma \ref{djqwjkqwdqwdqdqd2131} in the appendix. \hB

In order to construct the de Rham equivalence we now use that $\Omega A_{\infty}$   and hence $ H(\Omega A)$ are homotopy invariant by Lemma  \ref{djqwkdjqwlwsqdqw}.   If we further use    the canonical identification \begin{equation}\label{fiwjeiofjwefewofewfewfewf234}
\Omega A(*)\cong A\ ,
\end{equation} then we get   the de Rham equivalence  as the following composition \begin{equation}\label{jhekjhekjhkd12983ud128oidkl1d1d}
\Funk(H(A))\stackrel{\eqref{fiwjeiofjwefewofewfewfewf234}}{\simeq}  \Funk(H(\Omega A)(*))\stackrel{\eqref{dkqwkdqwdjdljdlkwd9879}}{\simeq} H(\Omega A)\ .
\end{equation}
This finishes the proof of Proposition \ref{canequ}.
\hB 

An inspection of the proof of Proposition \ref{canequ} yields that the de Rham equivalence is natural in the chain complex $A\in \Ch_{\R}$.

\begin{kor}\label{lqejhdlwdjqwodu321oidwqd}
We have constructed an equivalence   
$$H(\Omega-)\simeq \Funk(H(-))$$  between functors $ \Nerve(\Ch_{{\R}})\to \Sm(\Sp)$.
In particular, if $A\to A^{\prime}$ is a morphism between chain complexes {of real vector spaces}, then we have a commuting
square in $\Sm(\Sp)$
$$\xymatrix{H(\Omega A)\ar[d]^{\simeq}\ar[r]&H(\Omega A^{\prime})\ar[d]^{\simeq}\\
\Funk(H(A))\ar[r]&\Funk(H(A^{\prime}))}\ .$$\end{kor}

{\begin{rem}{\rm 
In Proposition  \ref{flkijweflwefewfewfewfwe23442434} we will improve the de Rham equivalence further to an equivalence between  symmetric monoidal functors.}
\end{rem}
}

\bigskip

  The input for the construction of the differential function spectrum is the datum
$(E,A,c)$ as in Definition  \ref{aug1001}. {This datum together with    de Rham equivalence gives rise to a} map of smooth spectra
\begin{equation}\label{jun147}
\rat:\Funk(E)\xrightarrow{e_{\R}} \Funk(E\R)\xrightarrow{c} \Funk(H(A)) \stackrel{\eqref{jhekjhekjhkd12983ud128oidkl1d1d}}{\simeq} H(\Omega A)\ .\end{equation}


   \begin{ddd}\label{may271}
We define the differential function spectrum $$\Diff(E,A,c)\in \Sm^{desc}( {\Sp} )$$ as the pull-back 
\begin{equation}\label{difdef}\xymatrix{\Diff(E,A,c)\ar[d]\ar[r]&H(\sigma \Omega A)\ar[d] \\\Funk(E)\ar[r]^{\rat}& H(\Omega A) }.\end{equation}
\end{ddd}

\begin{rem}{\rm
{In this definition, {according to the first part of Remark \ref{ljfewfewfwfewfewfewfwf67678463842344234234},} the  pull-back defining $\Diff(E,A,c)$ is understood as a limit in the $\infty$-category  of smooth spectra $\Sm^{desc}( {\Sp} )$. In order to ensure its existence note that  $\Sp$, and hence  $\Sm^{desc}(\Sp)$, are presentable $\infty$-categories and therefore  admit all small limits. 
 Introducing the differential function spectrum in this way we encode the functoriality on the manifold $M$ automatically. By evaluation at a manifold $M$ we get a pull-back 
 \begin{equation}\label{difder23r2r32r2r32rf}\xymatrix{\Diff(E,A,c)(M)\ar[d]\ar[r]&H(\sigma \Omega A(M))\ar[d]\\ E^{M}\ar[r] & H(\Omega A(M)) }\end{equation}
in the $\infty$-category $\Sp$ of spectra. Even more concretely, if we realize  this diagram in a point-set  model $\check{\Sp}$, then   $\Diff(E,A,c)(M)$ is equivalent to the homotopy pull-back of a point-set version of the diagram
$$E^{M_{top}}\rightarrow  H(\Omega A(M))\leftarrow  H(\sigma \Omega A(M))\ .$$
We conclude that Definition \ref{may271} is really  the $\infty$-categorical version of the construction explained in the motivating Example \ref{may271}.
}
}\hB \end{rem}

\bigskip

{Finally we define the differential cohomology functors by   taking homotopy groups.
For $n\in \Z$ we let $\Sigma^{n}$ denote the shift functor on the categories $\Sp$, $\Ch$ and 
$\Nerve(\Ch)[W^{-1}]$.
 The Eilenberg-MacLane spectrum functor  is compatible with the shift, i.e} we have an equivalence $\Sigma^{n} H(A){\simeq} H(\Sigma^{n}A)$.
The datum $(E,A,c)$ thus induces a datum
$$\Sigma^{n}(E,A,c):=(\Sigma^{n}E,\Sigma^{n}A,c_{n})\ ,$$ where
$$c_{n}:\Sigma^{n}E\R\xrightarrow{\Sigma^{n}c}\Sigma^{n} H(A){\simeq} H(\Sigma^{n}A)\ .$$

The functor {$\pi_{0}: {\Sp}\to \Ab$ induces a functor
$$\pi_{0}:\Sm({\Sp})\to \Sm(\Nerve(\Ab))$$ {from smooth spectra to smooth abelian groups.}}  {Note that {by \eqref{fhewfhwefhewfiuizuii3u3244234234ewddewd} a smooth abelian group is nothing else than a functor $\Mf^{op}\to \Ab$.}
\begin{ddd}\label{difdef1}
The differential $E$-cohomology in degree $n\in\Z$ is the smooth {abelian group
$$\widehat E^{n}:=\pi_{0}(\Diff(\Sigma^{n}(E,A,c)))\in \Sm(\Nerve(\Ab))\ .$$}
\end{ddd}
 
{This finishes the construction of the fuctor \eqref{hjkewhfwjhfewhfkjwhfjewhwhek}.}
In the next subsection we will construct the natural transformations $R$, $I$, and $a$ and verify the 
axioms stated in Definition \ref{jun216}.
 
\subsection{Homotopy groups and long exact sequences}\label{kjefklwefwefwefwefewfwefwef}

In this subsection we calculate the smooth abelian groups $\pi_{n}(\Diff(E))\in \Sm(\Nerve(\Ab))${, provide the structure maps of the differential cohomology theory $\widehat E^{*}$, and verify the axioms stated in Definition \ref{jun216}. We further discuss the Mayer-Vietoris sequence.}

\bigskip

 By $Z^{0}(\Omega A)\in  \Sm(\Nerve(\Ab))$ we denote the smooth abelian group which associates to a manifold $M$  the abelian group of cycles of degree zero in the chain complex  $\Omega A(M)$. We  consider $Z^{0}(\Omega A)$ as a smooth  chain complex $Z^{0}(\Omega A)^{\bullet}\in \Sm(\Nerve(\Ch))$ which is concentrated in degree zero. {We then have an isomorphism of smooth abelian groups
 \begin{equation}\label{jul0602}Z^{0}(\Omega A)\cong  \pi_{0}(H(Z^{0}(\Omega A)^{\bullet})) \ .\end{equation}
 Furthermore, the inclusion map of smooth chain complexes $$i:Z^{0}(\Omega A)^{\bullet} \to {\sigma}\Omega A$$  
 induces a map of smooth spectra} \begin{equation}\label{eqq2}H(Z^{0}(\Omega A)^{\bullet})\to H(\sigma\Omega A) \to H(\Omega A)\stackrel{\eqref{jhekjhekjhkd12983ud128oidkl1d1d}}{\to}\Funk(H(A))\ .\end{equation}
These maps give the first two maps in the composition 
$$Z^{0}(\Omega A)\stackrel{(\ref{jul0602})}{\cong} \pi_{0}(H(Z^{0}(\Omega A)^{{\bullet}}))\xrightarrow{(\ref{eqq2})} \pi_{0}(\Funk(H(A)))\xrightarrow{(\ref{hgr1}) \:and\: {c} } E\R^{0}$$
used to define the fibre product in equation \eqref{exats} of the following proposition.
 \begin{prop}\label{hgr}
We have  isomorphisms {of smooth abelian groups}
\begin{equation}\label{jfwe89f8923rkjhnkwjfwef98u32f23f}
\pi_{n}(\Diff(E))\cong \left\{\begin{array}{cc}
E^{-n}
&n\le -1\\
E\R/\Z^{-n-1}&n\ge 1\\ \widehat E^{0}\:\:
\mbox{see\:below}&n=0
\end{array}\right\}\ .
\end{equation}
The smooth abelian group $\widehat E^{0}$ fits into the following exact sequences
\begin{equation}\label{exats}E^{-1}\stackrel{c}{\to} E\R^{-1}\xrightarrow{a} \widehat E^{0}\xrightarrow{(I,R)} E^{0}\times_{E\R^{0}} Z^{0}(\Omega A)\to 0\end{equation}
{and
\begin{equation}\label{exats111}E^{-1}\stackrel{c}{\to} \Omega A^{-1}/\im(d)\xrightarrow{a} \widehat E^{0}\xrightarrow{I} E^{0} \to 0\ ,\end{equation}}where
the maps $I,R$ and $a$ will be constructed in the course of the proof.
\end{prop}
\proof
Let $\tau^{\ge 1}: \Sp \to  \Sp $ be functor which takes $(0)$-connective covers. 
It induces a functor denoted by the same symbol $$\tau^{\ge 1}:\Sm( \Sp )\to \Sm( \Sp) \ .$$ If     we define the intermediate smooth spectrum 
$X$   by the pull-back in $ \Sm( \Sp  )$
\begin{equation}\label{eqq4}\xymatrix{X\ar[d]\ar[r]&\tau^{\ge 1}H(\sigma \Omega A)\ar[d]\\\tau^{\ge 1}
\Funk(E)\ar[r]&\tau^{\ge 1}H(\Omega A)}\ ,\end{equation}
then we get a map $X\to \Diff(E)$ induced by the natural map of pull-back diagrams $(\ref{eqq4})\to (\ref{difdef})$ coming from the  
 morphism
$\tau^{\ge 1}\to \id$ . 
By definition the spectrum $E\R/\Z$ fits into a pull-back square  
$$\xymatrix{\Sigma^{-1}E\R/\Z\ar[r]\ar[d]&0\ar[d]\\E\ar[r]^{e_{\R}}&E\R}$$
in $ \Sp $.
{The} {smooth} function spectrum construction $$\Funk: \Sp \to \Sm( \Sp )$$ {preserves all limits, so in particular}   pull-back squares. {It further}   commutes with shifts. {Therefore we get}
 a pull-back square of smooth spectra
\begin{equation}\label{jul1711}\xymatrix{\Sigma^{-1}\Funk(E\R/\Z)\ar[r]\ar[d]&0\ar[d]\\\Funk(E)\ar[r]^{e_{\R}}&\Funk(E\R)}\ .\end{equation}
This can be extended to a diagram of squares\\[0.3cm]
\begin{equation}\label{eqq4111}\xymatrix{\Sigma^{-1}\Funk(E\R/\Z)\ar@/^{1cm}/@{.>}[rr]\ar[r]\ar[d]&0\ar@{=}[r]\ar[d]&0\ar[d]\\\Funk(E)\ar@/_{1cm}/@{.>}[rr]^{\rat}\ar[r]^{e_{\R}}&\Funk(E\R)\ar[r]^{j^{-1}\circ c}& H(\Omega A)}\ .\end{equation}

We have $\pi_{k}(\tau^{\ge 1}H(\sigma \Omega A))\cong 0$ for $k\le 0$. Since also
 $H^{-k}(\sigma \Omega A)\cong 0$ for $k>0$ and $$\pi_{k}(\tau^{\ge 1}H(\sigma \Omega A))\cong H^{-k}(\tau^{\ge 1}H(\sigma \Omega A))$$ we conclude that $\pi_{k}(\tau^{\ge 1}H(\sigma \Omega A))\cong 0$ for all $k\in \Z$ and hence 
$\tau^{\ge 1}H(\sigma \Omega A){\simeq} 0$.
Thus we also have a map
$X\to \Sigma^{-1}\Funk(E\R/\Z)$ induced by the natural map from the pull-back diagram (\ref{eqq4}) to the outer square of (\ref{eqq4111}).
Each of these pull-back diagrams provides a long exact sequence of homotopy groups (similar to (\ref{eq104})), and these sequences are related by corresponding maps.
  
{Using  the Five Lemma we get the following isomorphisms
  of smooth abelian groups
$$\pi_{n}(\Diff(E))\stackrel{{\cong}}{\leftarrow} \pi_{n}(X)\stackrel{{\cong}}{\to} \pi_{n+1} (\Funk(E\R/\Z))$$
 for $n\ge 1$.}
Using (\ref{hgr1}) in order to identify the group on the right with $E\R/\Z^{-n-1}(M)$ we get the  {isomorphism \eqref{jfwe89f8923rkjhnkwjfwef98u32f23f}} in the case $n\ge 1$.

We let $$\tau^{\le -1}:\Sm( \Sp )\to \Sm( \Sp )$$ be the cofibre of the natural transformation $\tau^{\ge 0}\to \id$.
We get a map
$\Diff(E)\to Y$, where the intermediate smooth spectrum $Y\in \Sm( \Sp )$ is defined by the pull-back
$$\xymatrix{Y\ar[d]\ar[r]&\tau^{\le -1}H(\sigma \Omega A)\ar[d]\\ \tau^{\le - 1}
\Funk(E)\ar[r]&\tau^{\le -1}H(\Omega A)}\ .$$
Since we have an equivalence
$$\tau^{\le -1}H(\sigma \Omega A )\stackrel{{\simeq}}{\to}\tau^{\le -1}H(\Omega A)$$
we conclude that
$Y\to  \tau^{\le - 1}
\Funk(E)$ is an equivalence. It again follows from the Five-Lemma that
$$\pi_{n}(\Diff(E))\stackrel{\cong}{\to}\pi_{n}(Y)$$ for
$n\le -1$. This gives the {isomorphism \eqref{jfwe89f8923rkjhnkwjfwef98u32f23f}} for the case $n\le -1
$.

{
 Finally, in order to see the exact sequences in the case $n=0$, we use the long exact sequences
associated to the fibre sequences in $\Sm(\Sp)$
\begin{equation}\label{eq104}\cdots\to \Sigma^{-1} H(\Omega A)\stackrel{a}{\to} \Diff(E)\stackrel{(I,R)}{\to} \Funk(E)\times H(\sigma \Omega A)\to H(\Omega A)\to \cdots\ .\end{equation}
and \begin{equation}\label{fweflkfkellklklwlkejflkwejfklweouoi}
\Funk(\Sigma^{-1} E)\to \Sigma^{-1}H(\Omega A/\sigma \Omega A)\stackrel{{a}}{\to} \Diff(E)\stackrel{I}{\to} \Funk(E)\end{equation}
which follow immediately from the definition of $\Diff(E)$ as a pull-back \eqref{difdef}.}
The homotopy groups of all  entries different from  $\Diff(E)$ are known.
In particular,
$$\pi_{0}(H(\sigma \Omega A))\stackrel{ }{\cong}   Z^{0}(\Omega A)\ , {\quad \pi_{0}(\Sigma^{-1}H(\Omega A/\sigma \Omega A))\cong \Omega A^{-1}/\im(d)}\ .$$
This gives the desired result.  
\hB 

The following proposition justifies to call the smooth abelian group $\widehat E^{n}=\pi_{0}(\Diff(\Sigma^{n}E))$ the $n^\mathrm{th}$  differential $E$-cohomology. 
We verify that it satisfies the axioms for a differential extension Definition \ref{jun216} (see \cite{MR2365651} and \cite{MR2608479}). It suffices to consider the case $n=0$.

\begin{prop}\label{jun0305} {
The smooth abelian groups $E^{0}$ and $\widehat{E}^0$ fit into} the following diagram 
\begin{equation}\label{eqq144}\xymatrix{&E\R/\Z^{-1}\ar[dr]\ar[rr]^{\partial}&&E^{0}\ar[dr]&\\
E\R^{-1}\ar[dr]\ar[ur]&&\widehat E^{0}\ar[ur]^{I}\ar[dr]^{R}&&E\R^{0}\\
&\Omega A^{-1}/\im(d)\ar[ur]^{a}\ar[rr]^{d}&&Z^{0}(\Omega A)\ar[ur]&}\ .
\end{equation}
Moreover, we have an isomorphism
 \begin{equation}\label{eqq14}E\R/\Z^{-1}\cong \widehat E^{0}_{flat} \end{equation}
 identifying the flat part (Definition \ref{jun216} ) with a homotopy theoretic object.
\end{prop}
\proof
The maps $R$ and $I$ have been constructed in Proposition \ref{hgr}. The map $a$ is induced by the map in {\eqref{fweflkfkellklklwlkejflkwejfklweouoi}} denoted by the same symbol.
There is a natural map of pull-back diagrams ${\eqref{jul1711}}\to (\ref{difdef})$
which induces the maps of smoth spectra and their zeroth homotopy groups
$$\Sigma^{-1}\Funk(E\R/\Z)\to \Diff(E) ,\quad E\R/\Z^{-1}\to \widehat E^{0}\ .$$
We get the commutativity of (\ref{eqq144})  and {the isomorphism} (\ref{eqq14}) by an analysis of the induced maps of long-exact sequences of homotopy groups.
\hB

{Up to this point we could have worked in $\Sm(\Sp)$. The subcategory $\Sm^{desc}(\Sp)\subseteq \Sm(\Sp)$ of smooth spectra satisfying descent is closed under limits. Since the right and lower corners of the pull-back diagram \eqref{difdef} all  belong to  $\Sm^{desc}(\Sp) $ we could conlcude in Definition \ref{may271} that $\Diff(E)\in \Sm^{desc}(\Sp)$. This additional  descent property leads to a Mayer-Vietoris sequence for differential cohomology which we now formulate {in its} simplest form.}} We consider a smooth manifold $M$ and  assume that $U,V\subseteq M$ are open subsets such that $M=U\cup V$.
\begin{prop}\label{efjwlfkwefwefwef}
We have a long exact sequence of groups
$$\cdots\to E\R/\Z^{-{2}}(U\cap V)\to\widehat E^{0}(M)\to \widehat E^{0}(U)\oplus \widehat E^{0}(V)\to \widehat E^{0}(U\cap V)\to \ E^{1}(M)\to\cdots$$
which extends to the left and right by the long exact {Mayer-Vietoris sequences for the cohomology theories} $E\R/\Z^{*}$ and $E^{*}$.
\end{prop}
\proof
Since $\Diff(E)$ satisfies descent we have
a pull-back square in $\Sp$
$$\xymatrix{\Diff(E)(M)\ar[d]\ar[r]&\Diff(E)(U)\ar[d]\\\Diff(E)(V)\ar[r]&\Diff(E)(U\cap V)}\ .$$
It induces the Mayer-Vietoris sequence in view of the calculation of homotopy groups of Proposition \ref{hgr}.
\hB 

\subsection{{Differential} Data and Transformations}\label{aug1006}
In this subsection we explain how $\Diff(E,A,c)$ depends on the data $(E,A,c)$. 
 {We start with introducing the}  $\infty$-category  {$\widehat{\Sp}$} of such data and describe $\Diff(E,A,c)$ as a functor {$\widehat{\Sp}\to \Sm^{desc}(\Sp)$. Then we introduce the notion of canonical data and discuss the possibility of choosing functorial canonical data.}  
 
 \bigskip
 
{Note that $\infty$-categories themselves {form} an $\infty$-category. The pull-back in the following  Definition is understood in this $\infty$-category of $\infty$-categories. We will describe some aspects of the construction more concretely below.}

\begin{ddd}
We define the $\infty$-category of data $\widehat{\Sp}$ as the pullback
 \begin{equation}\label{ferfre43r23rr2r23r32r}\xymatrix{
\widehat{\Sp}\ar[r]\ar[d] & \Nerve(\Ch)\ar[d]^{H}\\
 \Sp \ar[r]^{\dots\wedge M\R} &  \Sp }\ .
\end{equation}
 \end{ddd}

More concretely,  {an object of $\widehat{\Sp}$} consists of a spectrum $E$, a chain complex $A$, and a specified  equivalence $c:E\mathbb{R}\simeq H(A)$.
Similarly, a map of data $$(f,g,\phi):(E,A,c)\to (E^\prime,A^\prime,c^\prime)$$ consists of a map of spectra $f:E\to E'$, a map of chain complexes $g:A\to A'$, and a homotopy $\phi$ from $H(g)\circ c$ to $c'\circ f\land M\mathbb{R}$.

Let {$\Lambda^2_2$}  be the nerve of the category of the shape
$$\bullet \rightarrow \bullet \leftarrow \bullet\ .$$
{Using Corollary \ref{lqejhdlwdjqwodu321oidwqd}   we can consider the} functor
$$P:\widehat{\Sp}\to \Fun(\Lambda^2_2,\Sm( \Sp))$$
which maps the triple $(E,A,c)$ to the diagram
$$\Funk(E)\stackrel{rat}{\to}  H(\Omega A_{\infty})  \leftarrow  H(\sigma \Omega A_{\infty})\ .$$
{At this point it is important that left right upper corner in \eqref{ferfre43r23rr2r23r32r}  {is} $\Nerve(\Ch)$ and not $\Nerve(\Ch)[W^{-1}]$ since the cut-off $\sigma$ is not well-defined  on $\Nerve(\Ch)[W^{-1}]$.}
If we compose {this functor} with the limit $$\lim_{\Lambda^2_2}:\Fun(\Lambda^2_2,\Sm( \Sp ))\to \Sm( \Sp )\ ,$$ then we obtain the differential function spectrum functor (Definition \ref{may271}) $$\Diff=\lim_{\Lambda^2_2}\circ P:\widehat{\Sp}\to \Sm( \Sp )\ .$$ It actually takes values in $ \Sm^{desc}( \Sp )$.
This construction {finally} makes clear how the differential function spectrum depends on the {differential data} $(E,A,c)$.

\begin{kor}
We have constructed a functor
$$\Diff:\widehat{\Sp}\to \Sm^{desc}(\Sp)\ , \quad (E,A,c)\mapsto \Diff(E,A,c)$$
\end{kor}

\bigskip

We now discuss the possibility of the choice of canonical {differential} data.
There is a version of the Eilenberg-MacLane equivalence ({see Example  \ref{kdjqwkldjqwdljlwqwqd1341323213213123}})
$$H_{\Q}:\Nerve(\Ch_{\Q})[W^{-1}]\stackrel{\simeq}{\to} \Mod(H\Q)\ , \quad C\mapsto \map(\Q_{\infty},C)$$
where $\Ch_{\Q}$ denote the category of chain complexes over $\Q$.
We call a spectrum $E$ rational if it is in the essential  image of the forgetful functor
$$\Mod(H\Q)\to \Sp \ .$$
For a chain complex of rational vector spaces $A_{\Q}\in \Nerve(\Ch_{\Q})[W^{-1}]$ we can choose an equivalence
$A_{\Q}\stackrel{\simeq}{\to} H_{*}(A_{\Q})$, where we consider the homology $H_{*}(A_{\Q})$ of $A_{\Q }$ as a complex with zero differential.
After application of $H_{\Q}$ we get an equivalence
$$H_{\Q}(A_{\Q})\stackrel{\simeq}{\to} H_{\Q}(H_{*}(A_{\Q}))$$
such that
$$\xymatrix{\pi_{*}(H_{\Q}(A_{\Q}))\ar[d]^{\cong}\ar[r]^{\cong}&\pi_{*}(H_{\Q}(H_{*}(A_{\Q})))\ar[d]^{\cong}\\
H_{*}(A_{\Q})\ar@{=}[r]&H_{*}(A_{\Q})}$$ 
commutes.
 
Since  we have an equivalence  $M\Q{\simeq} H\Q$   of the Moore and the Eilenberg-MacLane spectra associated to $\Q$ the spectrum $E\Q{:=}E\wedge M\Q$ is rational for any spectrum  $E\in  \Sp $.
If we {define the complex $A_{\Q}$ by  $A_{\Q}:=\pi_{*}(E)\otimes \Q\in \Nerve(\Ch_{\Q})[W^{-1}]$ with trivial differential},  then there exists an equivalence
$$c_{\Q}:E\Q\stackrel{\simeq}{\to} H_{\Q}(A_{\Q})$$ which induces the canonical identification in homotopy groups.
We define $A:=\pi_{*}(E)\otimes \R$. Using the  canonical equivalences ${(E\Q)}\R\simeq E\R$ {and 
$\cF(H_{\Q}(A_{\Q}))\R\simeq H(A)$ (where $\cF$ forget the $H\Q$-module structure)}
the equivalence
$c_{\Q}$ induces an equivalence $$c:E\R\stackrel{\simeq}{\to} H(A)\ .$$

\begin{ddd}\label{candef}
{Any differential data $(E,A,c)\in \widehat{\Sp}$ which arises in the manner described above will be referred to as canonical differential data of $E$}. 
\end{ddd}

For $A_{\Q},B_{\Q}\in \Nerve(\Ch_{\Q})[W^{-1}]$ we know that a map $A_{\Q}\to B_{\Q}$ which induces the zero map in homology  is equivalent to the zero map. This easily follows from the fact that
$A_{\Q}$ and $B_{\Q}$ are equivalent to their homology complexes with trivial differential. This fact implies via the equivalence $H$ that a map between rational spectra
is equivalent to the zero map if it induces the zero map in homotopy. It follows that
the canonical {differential data} $(E,A,c)\in \widehat{\Sp}$ of $E$ is unique up to homotopy.
But one should be careful as the {construction} which associates to $E$ the canonical {differential data} $(E,A,c)$
can not be turned into a functor even  when considered with values in the homotopy category $\mathrm{Ho}(\widehat{\Sp})$ of $\widehat{\Sp}$.

Let $f:E\to E^{\prime}$ be a map of spectra and
$(E,A,c)$ and $(E^{\prime},A^{\prime},c^{\prime})$ be the canonical data.
Then we get an induced  map $f_{*}:A\to A^{\prime}$. By the above there exists a homotopy
$H(f_{*})\circ c_{\Q}\sim c_{\Q}^{\prime}\circ (f\wedge M\Q)$. A choice of this homotopy $\phi_{\Q}$
induces a homotopy $\phi:H(f_{*})\circ c\sim c^{\prime}\circ (f\wedge M\R)$ and therefore
 produces a map
$$(f,f_{*},\phi):(E,A,c)\to (E^{\prime},A^{\prime},c^{\prime})$$
in $\mathrm{Ho}(\widehat{\Sp})$.

But in general there is no way to choose this homotopy $\phi_{\Q}$ naturally. In general the group
$$\pi_{1}(\map(E\Q,H(A_{\Q}^{\prime})))\cong \hom(A_{\Q},A_{\Q}^{\prime}[-1])\cong \prod_{n\in \Z}\hom(A_{\Q,n},A^{\prime}_{\Q,n+1})\ .$$ acts simply transitively on the set of choices of this homotopy. If we apply this to $\id:E\to E$
we see that we can define the canonical datum only uniquely up to the automorphism group
$$\prod_{n\in \Z}\hom(A_{\Q,n},A_{\Q,n+1})\ .$$ After fixing the choice of a canonical {differential data}
this group acts by automorphisms on the functor
$\widehat E^{0}$. In the case of differential algebraic $K$-theory we describe the action in (\ref{nov0102}).

Let $ \Sp^{\Q-ev} \subset \Sp $ be the full subcategory of rationally
even spectra, i.e. spectra satisfying $\pi_{i}(E)\otimes \Q=0$ for all odd $i\in \Z$. If $E$ is rationally even, then $$\prod_{n\in \Z}\hom(\pi_{n}(E\Q),\pi_{n+1}(E\Q))=0$$ and therefore
the canonical {differential data} is well-defined in $\mathrm{Ho}(\widehat{\Sp})$   up to unique isomorphism.

\begin{kor}\label{fkwwlewfwefewfewfewfeff}
We have a functor $${\mathrm{Ho}}( \Sp^{\Q-ev}) \to \mathrm{Ho}(\widehat{\Sp})$$ which associates to a rationally even spectrum the class of its canonical {differential} data. A similar statement holds true for rationally odd spectra.
 \end{kor}
This corollary complies with the results of \cite{MR2608479}.

\begin{ex}{\rm 
Examples of rationally even spectra are $H\Z$, the bordism spectra $MBO\langle n \rangle$ for $n\ge 0$, $MU$ and the topological $K$-theory spectra $KU$ and $KO$. The connected cover 
$KR\langle 1 \rangle$ of the algebraic $K$-theory spectrum of a number ring is rationally odd by Theorem \ref{borel}. Since $\pi_{0}(KR\Q)=\Q$ and, for example, $\dim_{\Q}\pi_{5}(KR\Q)>0$ (Theorem \ref{borel}),
the  algebraic $K$-theory spectrum of a number ring itself 
is neither rationally even nor odd. The  group
$$\hom(\pi_{0}(KR\Q),\pi_{1}(KR\Q))\cong  K_{1}(R)\otimes \Q$$
acts on $\widehat{KR}^{0}$ by automorphisms. An explicit formula for this action will be given in (\ref{nov0102}).
}\hB\end{ex}
{\begin{ex}{\rm
We consider the complex $K$-theory spectrum $KU$ and the periodic rational cohomology spectrum
$HP\Q:=H(\Q[b,b^{-1}])$ with $\deg(b)=2$. The Chern character defines a map of spectra $\ch_{\Q}:KU\to  HP\Q$.
Since both spectra are rationally even by Corollary \ref{fkwwlewfwefewfewfewfeff} we get a map between differential data
$$\hat \ch_{\Q}:(KU,\pi_{*}(KU)\otimes \R,c_{KU})\to (HP\Q, \Q[b,b^{-1}],c_{HP\Q})$$
which unique up to homotopy. 
It induces a well-defined
  transformation between differential cohomology theories
$$\hat \ch_{\Q}:\widehat{KU}^{*}\to \widehat{HP\Q}^{*}\ .$$
This is the differential refinement of the Chern character whose existence and uniqueness was first shown in  \cite[Thm. 6.2]{MR2664467}. The comparison of  our present construction of $\hat \ch_{\Q}$ with the efforts needed
in  \cite[Sec. 6]{MR2664467} demonstrates the usefulness of the homotopy-theoretic approach to the construction of natural transformations between differential cohomology theories.
}\hB  
\end{ex}
 }




\section{Cycle maps}\label{aug0505}

\subsection{Introduction}

It often occurs  in topology or differential topology that {various structures} are classified by invariants in a cohomology theory $E^{*}$. In this case one motivation to consider the differential extension $\widehat E^{*}$ is because it may capture invariants of additional geometric structures.

Complex vector bundles present a typical example. They are classified by their $KU$-theory classes, whereas differential $KU$-theory can be used to capture additional information about connections.
As our definition of differential $KU$-theory in the present paper is homotopy theoretic, it requires some work to construct the refined geometric invariants. The result will be encoded into a cycle map. {We will discuss this case as a warm-up.}
We then move on to one of the main results of this paper; namely, the construction of a cycle map for locally constant sheaves of finitely generated projective $R$-modules with geometry taking values in the differential version of the algebraic $K$-theory of $R$, where $R$ is  the ring of integers in a number field (Definition \ref{sep0701} and Theorem \ref{thecycleproj1}).

\subsection{Complex $K$-theory -- a warm-up}\label{jul0311}

{In this subsection we discuss the cycle map for differential complex $K$-theory. Cycles for differential complex $K$-theory classes are complex vector bundles with hermitean metrics and metric connections.
Since complex $K$-theory is rationally even and admits geometric models this case is easy enough in order to serve as an introduction to the much more complicated case of algebraic $K$-theory.}

\bigskip

We let $KU$ denote the complex $K$-theory spectrum. {Since it is rationally even  we can choose the canonical  differential data $(KU,A,c)$ functorially (Definition \ref{candef} and Corollary \ref{fkwwlewfwefewfewfewfeff}). The $\Z$-graded $\R$-vector space $A$ is the underlying {graded} vector space of the ring  of Laurent polynomials
\begin{equation}\label{apr060113}A:=   \R[b,b^{-1}] \ ,\end{equation}} {where $b$ has homological degree $-2$ and is called the Bott element.} {Using the isomorphism \begin{equation}\label{fewffwefewfefewffefewf4234324}
A\cong \pi_{*}(KU)\otimes \R
\end{equation} given by the Chern character we construct the differential  function spectrum $$\Diff(KU):=\Diff(KU,A,c)$$ with \begin{equation}\label{ewfwfewfewfwf8978923432424}
c:KU\stackrel{can}{\to} H( \pi_{*}(KU)\otimes \R)\stackrel{\eqref{fewffwefewfefewffefewf4234324}}{\simeq} H(A)
\end{equation}  as in   Definition \ref{may271}. According  to   Definition \ref{difdef1} the  differential $K$-theory group  $\widehat{KU}^{0}(M)$ of a smooth manifold $M$ is  then defined by evaluation of  the differential  function spectrum at $M$ and applying $\pi_{0}$. The group $\widehat{KU}^{0}(M)$ will be the target of the cycle map for $M$.
 }

\bigskip 

{The domain of the cycle map is the commutative monoid $\bvect_{h,\nabla}(M)$ of isomorphism classes of geometric vector bundles on $M$.}
A geometric vector bundle $\bV=(V,h^{V},\nabla^{V})$ over a manifold $M$ is by definition a complex vector bundle $V\to M$ together with a hermitean metric $h^{V}$ and a hermitean (i.e. metric preserving, see (\ref{aug1002})) connection $\nabla^{V}$.  For simplicity we will also use the notation $\bV$ for the isomorphism class of the geometric vector bundle $\bV$. The monoid structure on $\bvect_{h,\nabla}(M)$ is given by the direct sum.

 {If $f:M^{\prime}\to M $ is a smooth map and $\bV $ is a geometric vector bundle on $M $, then we can define the pull-back $f^{*}\bV $ which is a geometric vector bundle on $M^{\prime}$.} Hence we can define the smooth commutative monoid
$$\bvect_{h,\nabla}\in  \Sm(\CommMon(\Set))$$ (i.e. a contravariant functor from manifolds to commutative monoids)     which associates to a manifold $M$ the commutative monoid of isomorphism classes of geometric vector bundles, and to $f:M^{\prime}\to M $ the pull-back $f^{*}$.  

\bigskip

We have maps of smooth commutative monoids {(i.e. natural transformations between functors) 
$$\hat I:\bvect_{h,\nabla}\to KU^{0}\ , \quad \bV\to [V]$$ and $$  \hat R:\bvect_{h,\nabla}\to Z^{0}(\Omega A)\ , \quad \bV\mapsto \ch(\nabla^{V})\ .$$ The first  sends} {the isomorphism class $\bV  {\in\bvect_{h,\nabla}(M)}$ of a} geometric vector bundle   on $M$ to the $KU$-theory class
$$[V]\in KU^{0}(M)$$ of the underlying complex vector bundle. The second associates to  $\bV$   the {Chern character 
form}
\begin{equation}\label{jul0312}\ch(\nabla^{V}):=\Tr \exp\left(-\frac{R^{\nabla^{V}}}{ b\ 2\pi i}\right) \in Z^{0}(\Omega A(M))\ ,\end{equation}  {where $R^{\nabla^{V}}\in \Omega^{2}(M,\End(V))$ denotes the curvature of the connection $\nabla^{V}$.} 

\begin{ddd}\label{jul0310}
A cycle map is a map of smooth commutative monoids
$$\cycl:\bvect_{h,\nabla}\to \widehat{KU}^{0}$$
such that the following diagram   commutes:
\begin{equation}\label{rh2urhi2uri2hr32r32r32r32fwf}
\xymatrix{&&{KU}^{0}\\\bvect_{h,\nabla}\ar[urr]^{\hat I}\ar[drr]_{\hat R}\ar[r]^-{\cycl}&\widehat{KU}^{0}\ar[ur]_{I}\ar[dr]^{R}\\
&&Z^{0}(\Omega A)}\ .
\end{equation}
\end{ddd}

\begin{rem}\label{jkdljwqdlwqdqwdqwdwdqwdq}{\rm
{A cycle map associates to a geometric vector bundle $\bV$ on a manifold $M$ a differential $K$-theory class $\cycl(\bV)\in \widehat{KU}^{0}(M)$. 
The diagram \eqref{rh2urhi2uri2hr32r32r32r32fwf} ensures the compatibility of this association with the underlying $K$-theory class and the Chern character form.  The definition of the cycle map as a map between smooth commutative monoids furthermore encodes the additivity and naturality in the manifold. So for bundles $\bV_{0},\bV_{1}\in \bvect_{h,\nabla}(M)$ and with $\bV:=\bV_{0}\oplus \bV_{1}$ we require
$$\cycl(\bV)=\cycl(\bV_{0})+\cycl(\bV_{1})\ .$$ Furthermore,  for a smooth map $f:M^{\prime}\to M$  and a bundle $\bV\in  \bvect_{h,\nabla}(M)$ we require the relation
$$f^{*}\cycl(\bV)=\cycl(f^{*}\bV)\ .$$
}\hB}
\end{rem}
{In the following we will show that a cycle map exists and is unique.}
{In order to connect with the geometric approach to differential $KU$-theory developed in \cite{MR2664467} we assume for the remainder of the present subsection that all smooth objects are defined on, or are restricted to the subcategory of compact manifolds $\Mf_{comp}\subseteq \Mf$. We will not write the restrictions explicitly.}
We have the following theorem. 

\begin{theorem}\label{dkdlqwdwqdwqddwqd}
There exists a unique cycle map
$\cycl:\bvect_{h,\nabla}\to \widehat{KU}^{0}$.
\end{theorem}
\proof
For the moment let $\widehat{KU}^{0}_{BS}$ be the differential extension of complex $K$-theory defined by the model \cite{MR2664467}. By construction, in this model differential $KU$-theory classes are represented by geometric families. In particular a geometric vector bundle,  {as}  a special kind of geometric family, represents a differential $KU$-theory class.

We now use the result of \cite{MR2608479} which states that the functor
$\widehat{KU}^{0}$ is uniquely characterized by the axioms for a differential extension
$(KU,I,R,a)$. 
 We thus have a unique comparison isomorphism \begin{equation}\label{jhfjkwehfkwefewfwefewfwef}
\comp:\widehat{KU}^{0}\stackrel{{\cong}}{\to} \widehat{KU}^{0}_{BS}\ .
\end{equation}
 For the model $\widehat{KU}^{0}_{BS}$ the cycle map is {tautological}. It sends the isomorphism class of a geometric vector bundle $\bV$ on a manifold $M$ to the {differential $K$-theory} class
$\cycl_{BS}(\bV)\in \widehat{KU}_{BS}^{0}(M)$ it represents. {Note that the analog of the diagram \eqref{rh2urhi2uri2hr32r32r32r32fwf} {is a built-in feature of} the construction of $\widehat{KU}^{0}_{BS}$.}
Therefore $$\cycl:=\comp^{-1}\circ \cycl_{BS}:\bvect_{h,\nabla}\to \widehat{KU}^{0}$$
is a cycle map. 

We now show uniqueness.
If $\cycl^{\prime}$ is a second cycle map, then in view of the exact sequence (\ref{exats}) the difference $\cycl-\cycl^{\prime}$ factorizes as the composite
$$\delta: \bvect_{h,\nabla}\to  KU\R^{-1}/KU^{-1}\to\widehat{KU}^{0}\ .$$ For every geometric vector bundle $\bV$ on a manifold $M$
there exists a map $f:M\to M^{\prime}$ and a geometric vector bundle $\bV^{\prime}$ on $M^{\prime}$ such that
$f^{*}\bV^{\prime}\cong \bV$ and the map  $f^{*}:H^{odd}(M^{\prime};\Q)\to H^{odd}(M;\Q)$ vanishes.
This implies that
$f^{*}:KU\R^{-1}(M^{\prime})\to KU\R^{-1}(M)$ vanishes and therefore
$$\delta (\bV)=\delta(f^{*} \bV^{\prime})=f^{*}\delta(\bV^{\prime})=0\ .$$ 
For $M^{\prime}$ one can take an approximation of the $\dim(M)+1$-skeleton
of the classifying space $BU$ such that a classifying map of $V$ factorizes over
an embedding $f:M\to M^{\prime}$.
We refer to \cite{MR2740650} for similar arguments.
 \hB


{{Since the cycle map is additive it can naturally be extended to virtual vector bundles.} It follows from  \cite{MR2732065} that this extension   is surjective.}


{
\begin{ex}{\rm Here we illustrate the cycle map by a simple example.
For $u\in U(1)$ we consider the flat line bundle $\bV(u)$ on $S^{1}$ with holonomy $u$. We want to calculate $$\cycl(\bV(u))\in \widehat{KU}^{0}(S^{1})\ .$$
{We have isomorphisms
 $$KU^{-1}(S^{1})\cong \Z\ , \quad KU^{0}(S^{1})\stackrel{\dim}{\cong} \Z\ , \quad H(A)^{-1}(S^{1})\cong \R \ ,$$
 where the inverse of the first one sends the integer $n\in \Z$ to the $K$-theory class represented by the map $S^{1}\to U(1)$ of degree $n$, and the third isomorphism is given by $b\int_{S^{1}}$.} By \eqref{exats} we have an exact sequence
$$0\to  \R/\Z\stackrel{{\mathbf{a}}}{\to} \widehat{KU}^{0}(S^{1})\stackrel{\dim}{\to} \Z \to 0\ .$$
From \cite{MR2664467} we know that $\widehat{KU}^{0}_{BS}(S^{1})$ takes values in commutative rings, and that the trivial line bundle represents the multiplicative unit $\beins$. We use the notation $\beins\in \widehat{KU}^{0}(S^{1})$ also for image of this unit under the isomorphism \eqref{jhfjkwehfkwefewfwefewfwef}. Hence we have $\beins=\cycl(\bV(1))$. We claim that
$${\cycl(\bV(u))}=\beins+{\mathbf{a}}(\frac{1}{2\pi i}\log(u))\ .$$
One way to see this is to use the homotopy formula. To this end we construct a geometric bundle
$\tilde \bV$ over $I\times S^{1}$ which restricts to
$\bV(1)$ on $\{0\}\times S^{1}$ and to $\bV(u)$ on $\{1\}\times S^{1}$.
Then by \eqref{feb2710} and \eqref{jul0312} we have  
$${\cycl(\bV(u))}=\beins+a(\int_{I\times S^{1} } \ch(\tilde \bV))=a(-\int_{I\times S^{1} }  \frac{1}{{b}2\pi i} R^{\nabla^{\tilde V}})  \ .$$
On the other hand, we can express the holonomy of $\tilde \bV_{|\{1\}\times S^{1}}$ along the circle in terms of a curvature integral
$$u=\exp(- \int_{I\times S^{1}}  R^{\nabla^{\tilde V}})\ .$$ 
Alternatively, one can refer to \cite[Lemma 5.3]{MR2664467}.

}\hB
\end{ex}
}

\subsection{The spectrum $KR$}\label{jun27333}
 
 {
In this subsection we    provide a construction of the {connective} algebraic $K$-theory spectrum $KR\in \Sp$ for a unital ring $R$.}


\bigskip

{  This construction involves   some theory of symmetric monoidal $\infty$-categories as developed in \cite{highalg}. We first explain the necesseary background.

If $\bC$ is a symmetric monoidal $\infty$-category, then we let $\CAlg(\bC)$ denote the $\infty$-category of commutative algebra objects in $\bC$. 
Commutative algebras in sets and spaces will be called commutative monoids, and we use the corresponding notation $$\CommMon( \Set):=\CAlg(\Nerve(\Set)) \ , \quad \CommMon(\Spc):=\CAlg(\Nerve(\Spc))\ .$$
  The latter is the $\infty$-categorical model of the homotopy theory of $E_{\infty}$-spaces.

\bigskip

We consider the symmetric monoidal one-category of small categories $\Cat$ with respect to the cartesian product. 
A commutative algebra object $\bC\in \CAlg(\Nerve(\Cat))$  is a very rigid form of a symmetric monoidal category were
the associators and the symmetries are given by identities. 

\begin{rem}\label{vjsdklvjlsdvsdvdsvdsvsdvsd}{\rm 
In this remark we give two reasons to consider the symmetric monoidal $\infty$-category $\Nerve(\Cat)[W^{-1}]$ of categories.

For a   ring $R$ we can consider the symmetric monoidal category $\Mod^{fin}(R)$  of finitely generated $R$-modules, where the symmetric monoidal structure is induced by the direct sum. After secretely passing to a small skeleton we can consider $\Mod^{fin}(R)$ as an object of $\Cat$. But since its associators and the symmetries are non-identities in general, it is not a commutative algebra in $\Nerve(\Cat)$. With  appropriate  definitions it would 
 be a  commutative algebra in the symmetric monoidal two-category of categories. But we prefer to go directly to $\infty$-categories. We will consider $\Mod^{fin}(R)$ as a commutative algebra in a symmetric monoidal $\infty$-category
 $\Nerve(\Cat)[W^{-1}]$ which we are going to explain below.

The functoriality in the ring provides the second reason.
  If $\phi:R\to S$ is a map between unital rings, then we obtain a symmetric monoidal functor
 $$\phi_{*}:\Mod^{fin}(R)\to \Mod^{fin}(S)\ , \quad M\mapsto M\otimes_{R}S\ .$$
Forgetting the symmetric monoidal structures for a moment  observe  that this construction does not give a functor $\Mod^{fin}:\Rings\to \Cat$. Indeed, if $\psi:S\to T$ is a second morphism between rings, then $\psi_{*}\circ \phi_{*}$ is not equal to $(\psi\circ \phi)_{*}$ but only naturally isomorphic. These isomorphisms satisfy higher coherence relations  which are encoded in the statement that the construction provides  functor between $\infty$-categories $$\Mod^{fin}:\Nerve(\Rings)\to \Nerve(\Cat)[W^{-1}]\ .$$  
This extends to the symmetric monoidal case.

  }\hB
 \end{rem}
 
 In the symmetric monoidal one-category  $\Cat$ we consider the class $W^{-1}$ of equivalences and form the symmetric monoidal $\infty$-category $\Nerve(\Cat)[W^{-1}]$. If $\bC\in \Cat$ is a small symmetric monoidal category, then we will use the same symbol in order to view it as an object of $\bC\in \CAlg(\Nerve(\Cat)[W^{-1}])$. To use $\bC_{\infty}$ instead  would be more consistent but our intention is to simplify the notation in the following.

  We now observe that the nerve functor $\Nerve: \Cat\to \sSet$ preserves cartesian products and maps
 equivalences between categories to weak equivalences between simplicial sets. It therefore descends to a symmetric monoidal functor between symmetric monoidal $\infty$-categories 
 $$\Nerve:\Nerve(\Cat)[W^{-1}] \to \Spc$$
 and finally induces a functor between the categories of commutative algebras
 $$\CAlg(\Nerve(\Cat)[W^{-1}] )\to  \CommMon(\Spc)\ .$$
 
  }

 {The functor $\pi_{0}:\Spc\to \Set$ is symmetric monoidal and therefore preserves commutative monoids. The $\infty$-category of groups in spaces 
 $$\CommGroups(\Spc)\subset   \CommMon(\Spc)$$ is characterized  in Definition  \ref{djqwldjlqwdqwd324234234} as the full subcategory 
of commutative  monoids in spaces which are mapped by $\pi_{0}$ to groups. The $\infty$-category $\CommGroups(\Spc)$ models the homotopy theory of group-like $E_{\infty}$-spaces.}

The inclusion of the full subcategory of commutative groups  into commutative monoids has as left-adjoint  which is called the group completion functor \eqref{khfkqhwfkhqwkfhqwkfhkhwqfqw34124124}
\begin{equation}\label{okt2001}\Omega B:\CommMon(\Spc)\leftrightarrows \CommGroups(\Spc):{inclusion}\ .\end{equation}
Finally, we have the well-known equivalence   \eqref{feb1703}
 \begin{equation}\label{okt2002}\spp:\
 \CommGroups(\Spc)\leftrightarrows \Sp_{\ge 0} :\Omega^{\infty}\end{equation}
between the $\infty$-categories of commutative groups in 
 {spaces} and connective spectra. If we postcompose $\spp$ with the  inclusion $\Sp_{\ge 0}\to \Sp$ of connective spectra into all spectra   we obtain a map  $\CommGroups(\Spc)\to \Sp$ which we will also be denote  by $\spp$.
  
\bigskip

{Let $R$ be a unital ring. Then} we consider the symmetric monoidal subcategory $$i\bP(R)\subseteq \Mod^{fin}(R)\ , \quad  {i\bP(R)\in \CAlg(\Nerve(\Cat)[W^{-1}])}$$  of finitely generated projective $R$-modules and isomorphisms ({this is indicated by the letter $i$ in front of $\bP$}).  {Its} nerve is a commutative monoid in the $\infty$-category of  spaces: $$
\Nerve(i\bP(R) )\in \CommMon(\Spc) \ .$$
 
\begin{ddd}\label{nov0101}
The algebraic $K$-theory spectrum of the ring $R$ is defined as
\begin{equation}\label{okt2003}KR:=\spp(\Omega B(\Nerve(i\bP(R) )))\in \Sp\ .\end{equation}
\end{ddd}
{We refer to 
\cite[Sec. IV]{weibel} for a verification  that $\pi_{*}(KR)$ is
 Quillen's higher $K$-theory of the ring $R$.}
 
 \bigskip
 
We actually have   a functor
$$i\bP:\Nerve(\Rings)\to \CAlg(\Nerve(\Cat)[W^{-1}])\ .$$
Consequently we  have defined an algebraic $K$-theory functor
\begin{equation}\label{jkjfl23jf32fjl2opip3r2r}K(-):= \spp(\Omega B(\Nerve(i\bP(-) ))):\Nerve(\Rings)\to \Sp\ .\end{equation}

 \subsection{The topological cycle map}\label{ewjflwefjeljlfkwfwefewf234}
 
{
 In this subsection we describe the topological cycle map which associates a class $[\cV]\in KR^{0}(M)$ to a locally constant sheaf of finitely generated projective $R$-modules $\cV$.   Later in Subsection \ref{may121270} and for $R$ a number ring we will drop the assumption ``projective''. The main technical result is the construction of the map \eqref{efjkhewjfhwkefhewfejwfhewkfhewjkfhw}.
 }
 
 \bigskip

{ Let $R$ be a unital ring.}

 \begin{rem}{\rm 
{The algebraic $K$-theory spectrum $KR$ represents a cohomology theory. For completeness of the presentation and later reference we recall various equivalent expressions  of  $KR$-cohomology groups of a manifold $M$. We use the {cotensor} \eqref{jchjsachakscsacscaca} and {the underlying topological space functor} \eqref{grgegg434353535345} in order to define the function spectrum $KR^{M}\in \Sp$ {(see Remark \ref{klfjwklfewlfewfewlfewfew23423454943} for the simplified notation)}. The algebraic $K$-theory groups of the manifold $M$  are then derived from the function spectrum by 
  \begin{equation}\label{dilqjwdlkqwdwqdwqdqwdwq}
KR^{*}(M) :=  \pi_{-*}(KR^{M})\ .
\end{equation}
 In particular we have following expressions for the zeroth cohomology:
 \begin{equation}\label{hhefkhwekfhwekfefewfewf234234}
KR^{0}(M)\stackrel{def}{=}\pi_{0}(KR^{M})\cong \pi_{0}(\Omega^{\infty}(KR^{M}))\cong  \pi_{0}((\Omega^{\infty}KR)^{M})\ .
\end{equation} }
\hB }\end{rem} 

 \bigskip

  {In order to describe the domain of the topological cycle map we use the language of stacks. A  prestack on the category of manifolds $
\Mf$ is an object of
$ \Sm( \Nerve(\Cat)[W^{-1}])$. It is called a stack if it belongs to the subcategory 
$$\Sm^{desc}( \Nerve(\Cat)[W^{-1}])\subset \Sm( \Nerve(\Cat)[W^{-1}])\ .$$ Similarly, a symmetric monoidal prestack on  $\Mf$ is an object of  $\Sm(\CAlg(\Nerve(\Cat)[W^{-1}]))$. We again consider the subcategory $$\Sm^{desc}(\CAlg(\Nerve(\Cat)[W^{-1}]))\subset \Sm (\CAlg(\Nerve(\Cat)[W^{-1}]))$$ of symmetric monoidal stacks. Note that  a symmetric monoidal prestack is a symmetric monoidal stack if and only if the underlying prestack is a stack. }

\begin{rem}\label{kldlqwdqwdwqdqwdqwd}{\rm 
A contravariant functor from the category of manfolds $\Mf$ to the one-category of categories $\Cat$ is{  called a strict   prestack. In the case of a strict prestack, for a} pair of smooth maps $f:M\to M^{\prime}$ and $f^{\prime}:M^{\prime}\to M^{\prime\prime}$ we have an equality $(f^{\prime}\circ f)^{*} =f^{*}\circ f^{\prime,*}$. If we associate to a manifold $M$ the category of locally constant sheaves of finitely generated $R$-modules, then we do not have this equality but just a natural isomorphism. The definition of a prestack given above encodes these isomorphisms and the higher coherence relations properly. See also the second part of Remark \ref{vjsdklvjlsdvsdvdsvdsvsdvsd}. Typically, prestacks of bundle-like objects are stacks. 
}\hB \end{rem}

\begin{rem}\label{kldlqwdqwdwqdqwdqwd1}{\rm 
In the following we let $\bC_{\infty}\in \Nerve(\Cat)[W^{-1}]$ be the image of $\bC\in \Nerve(\Cat)$ under the localization 
$\Nerve(\Cat)\to \Nerve(\Cat)[W^{-1}]$. We apply this notation also to smooth objects. 
We say that the strict prestack
$X\in \Sm(\Nerve(\Cat))$ presents a prestack $Y\in \Sm(\Nerve(\Cat)[W^{-1}])$, if there is an equivalence
$X_{\infty}\simeq  Y$.
In the proof of Lemma \ref{ejfwelfjwelfefwefewf} we will use the fact that one can present a stack of bundles by a  
strict prestack. {See e.g. \cite{Kim:2014fk} for a similar statement.}
 
 }\hB \end{rem}
   
   {

{We consider a locally constant sheaf $\cV$ of finitely generated $R$-modules on a manifold $M$. 
  \begin{rem}{\rm Let us recall what this means. First of all $\cV$ is a sheaf of $R$-modules on the site of open subsets of $M$ with respect to the topology given by open coverings. Locally constant means, that for every
  $x\in M$ and contractible neighbourhood $U$  of $x$   the restriction 
  $\cV(U)\to \cV_{x}$   is an isomorphism, where $\cV_{x}$ denotes the stalk of $\cV$ at $x$. The adjective ``finitely generated''  then means that $\cV_{x}$ is a finitly generated $R$-module for every $x\in M$. A similar convention applies to conditions like ``free'' or ``projective''. Note that $\cV(U)$ for arbitrary open subsets $U$ of $M$ may not be finitely generated. The happens e.g. for the locally constant sheaf obtained by sheafifying the constant presheaf with values $R$ if $U$ has infinitely many components.
    }\end{rem}\hB}
   
   If $f:M^{\prime}\to M $ is a smooth map between manifolds and
$\cV $ is a  locally constant sheaf of finitely generated $R$-modules on $M $, then we can define the  locally constant sheaf $f^{*}\cV$  of finitely generated   $R$-modules on $M^{\prime} $. 
As mentioned in Remark \ref{kldlqwdqwdwqdqwdqwd} we have isomorphisms for iterated pull-backs which  satisfy higher coherence relations. The pull-back furthermore preserves direct sums and sheaves of   projective modules. \begin{ddd}\label{fjewjflwefjlwfoiiouoiuoifwefwefwef} We define the 
symmetric monoidal stacks
$$ \loc^{proj}(R)\subset \loc(R)\in \Sm^{desc}(\CAlg(\Nerve(\Cat)[W^{-1}]))\ .$$ such that
the stack $\loc(R)$ (or $\loc^{proj}(R)$, respectively) associates to a manifold $M$ the symmetric monoidal category $\loc(R)(M)$ (or $\loc^{proj}(R){(M)}$, respectively)  of locally constant sheaves of finitely generated (projective) $R$-modules on $M$ and isomorphisms.
\end{ddd}

If $\phi:R\to R^{\prime}$ is a homomorphism of   rings, then we have an  induced map of symmetric monoidal stacks
\begin{equation}\label{kjldjlkjlkjdl1kdj1o2idop12d12d12d1d}\phi_{*}:\loc(R)\to \loc(R^{\prime})\end{equation} which  sends $\cV\in \loc(R)(M)$ to $\cV\otimes_{R}R^{\prime}\in \loc(R^{\prime})(M)$. This construction preserves projectivity.
We {conclude} that the symmetric monoidal stacks $\loc^{proj}(R)$ and $\loc(R)$ depend functorially on the ring, i.e. we actually have  functors between $\infty$-categories
\begin{equation}\label{gerggerg43r79834534534543fgerge}
\loc^{proj},\  \loc:\Nerve(\Rings)\to \Sm^{desc}(\CAlg(\Nerve(\Cat)[W^{-1}]))\ .
\end{equation}

\bigskip

Note that $\pi_{0}:\Spc\to \Nerve(\Set)$ is a symmetric monoidal functor and hence preserves monoids.
  \begin{ddd} \label{feb1720} We define the smooth monoid $$\bloc^{proj}(R) :=  \pi_{0}(\Nerve(\loc^{proj}(R)){)}\in \Sm(\CommMon(\Set)) $$
which associates to a manifold $M$ the monoid of isomorphism classes of locally constant sheaves of finitely generated projective $R$-modules.
\end{ddd}
 
 In order to relate the stack $\loc^{proj}(R)$ with $i\bP(R)$ we use the theory of homotopy invariant smooth objects.
 Recall the Definition \ref{jdqwdjqwldqdwqdwqdqd} of  homotopy invariance of a smooth object $F\in \Sm(\bC)$ of some $\infty$-category $\bC$. 
Equivalently, $F$ is homotopy invariant, if the inclusion $\{0\}\to I$ induces an equivalence
$F(I\times M)\stackrel{\simeq}{\to} F(M)$ for all manifolds $M$.


\bigskip

We are going to apply this to the stack  $\loc^{proj}(R)$.

\begin{lem}\label{jkljljdqlwdkwqdqdqd}
The stack  $\loc^{proj}(R)$ is homotopy invariant.
\end{lem}
 \proof For every manifold $M$ the
 restriction $$\loc^{proj}(R)(I\times M)\to   \loc^{proj}(M)$$ induced by the inclusion $\{0\}\to I$ is an equivalence of categories. This follows immediately from the fact that the objects of the stack are locally constant sheaves. 
\hB
 
It follows that $\Nerve(\loc^{proj}(R))$ is homotopy invariant, too. We thus have the following corollary of Lemma \ref{djqwjkqwdqwdqdqd213}.
\begin{kor}\label{fkljlwkefewfewwfweffwe}
We have an equivalence
$$\Funk( \Nerve(\loc^{proj}(R)(*)) )\stackrel{{\simeq}}{ \to}  \Nerve(\loc^{proj}(R)) $$ in $ \Sm^{desc}(\CAlg(\Nerve(\Cat)[W^{-1}]))$ which is natural in the ring $R$.
\end{kor}
 }

 \bigskip
 
 \begin{rem}{\rm 
There are analogs of Definition \ref{feb1720}, Lemma \ref{jkljljdqlwdkwqdqdqd} and Corollary \ref{fkljlwkefewfewwfweffwe} for  $\loc(R)$.}\hB\end{rem}
{
Essentially by definition we have an equivalence
$$  \loc^{proj}(R)(*)\simeq i\bP(R)\ .$$ Applying  Corollary \ref{fkljlwkefewfewwfweffwe}   we get an equivalence
\begin{equation}\label{jwqjejqejkwqe989898i234jkkjkjerqkjwr}
\Funk( \Nerve(i\bP(R)))\stackrel{\simeq}{\to} \Nerve(\loc^{proj}(R))
\end{equation}
in $\Sm^{desc}(\CommMon(\Spc))$.

\bigskip

 The composition of the units of the adjunctions (\ref{okt2001}) and (\ref{okt2002}) in view of the definition (\ref{okt2003}) 
provides a map of commutative monoids in spaces
\begin{equation}\label{okt2005}\Nerve(i\bP(R) )\to \Omega B (\Nerve(i\bP(R) ))\to \Omega^{\infty} (KR)\ .\end{equation}
We get an induced map in $\Sm^{desc}(\CommMon(\Spc))$:
 \begin{equation}\label{efjkhewjfhwkefhewfejwfhewkfhewjkfhw}
\Nerve(\loc^{proj}(R))\stackrel{\eqref{jwqjejqejkwqe989898i234jkkjkjerqkjwr}, \simeq}{\leftarrow} \Funk( \Nerve(i\bP(R)))
 \stackrel{\eqref{okt2005}}{\to} \Funk(\Omega^{\infty} (KR))\ .
\end{equation}
All our constructions are functorial in the ring. \begin{kor}\label{kldjqlkwjdlqdjqdqwdjqlwkd}We actually  have constructed a map
$$\Nerve(\loc^{proj}(-))\to \Funk(\Omega^{\infty} (K(-)))$$
in
$$\Fun(\Nerve(\Rings), \Sm^{desc}(\CommMon(\Spc)))\ .$$\end{kor}
 }
 
%
%
%
{

 \begin{ddd}\label{jul1702} The topological cycle map is the 
 transformation of smooth {monoids}
\begin{equation}\label{okt2006} {\hat I}:\bloc^{proj}(R)=\pi_{0}(\Nerve(\loc^{proj}(R)))\stackrel{\eqref{efjkhewjfhwkefhewfejwfhewkfhewjkfhw}}{\to} \pi_0(\Funk_\infty(\Omega^{\infty} (KR))){\stackrel{\eqref{hhefkhwekfhwekfefewfewf234234}}{\cong} KR^{0}}\ .\end{equation}
\end{ddd}
{Often} we write $$[\cV]:=\hat I(\cV)\ .$$
{In Definition \ref{may121311},  {under additional assumptions on $R$,} we will extend the transformation $\hat I$ to all of $\bloc(R)$.}

\bigskip

The functoriality in the ring $R$ is encoded in the statement that we  actually  have a  natural transformation
$$\hat I: \bloc^{proj}(-)\to  K(-)^{0}$$
in $$\Fun(\Nerve(\Rings),\Sm(\CommMon(\Set)))\ .$$
}

Under certain conditions   the topological cycle map splits short exact sequences. This property refines its additivity.  
\begin{prop}\label{jwlkjclkjwlkwjclcwewfefeklw9889234}
We assume that $R$ is a regular commutative ring. If $$\cV:0\to \cV_{0}\to \cV_{1}\to \cV_{2}\to 0$$
is a short exact sequence of  objects  $\cV_{i}\in \loc^{proj}(R)(M)$, then we have the equality   
$$[\cV_{1}]=[\cV_{0}]+[\cV_{2}]$$
in $KR^{0}(M)$.
\end{prop}
\proof
Since $R$ is a regular ring, {by the fundamental theorem (see e.g. \cite[Thm. 6.3]{weibel}),} the inclusion $c:R\to R[x]$ induces an   equivalence $c:KR\to KR[x]$ of $K$-theory spectra. For $i=0,1$ let  $s_i:R[x]\to R$ be the specializations at $x=i$. Since $s_{0}\circ c=s_{1}\circ c $ we have the equality
$$s_{0}^{*}=s_{1}^{*}:KR[x]^{0}(M)\to KR^{0}(M)\ .$$
Let us assume that $M$ is connected and choose a base point $m\in M$.
The sheaves $\cV_{j}$ for  $j=0,1,2$ are then determined by the corresponding holonomy representations  $\rho_{j}$ of the
fundamental group $\pi_{1}(M,m)$ on the fibres $\cV_{i,m}$. Since the fibres are projective $R$-modules   the sequence of fibres splits and we have $\cV_{1,m}\cong \cV_{0,m}\oplus \cV_{2,m}$. The representation $\rho_{1}$ is then related with $\rho_{0}$ and $\rho_{2}$ by
$$\rho_{1}=\left(\begin{array}{cc}\rho_{0}&\sigma\\0&\rho_{2} \end{array} \right)\ ,$$
where $\sigma:\pi_{1}(M,m)\to \Hom_{R}(\cV_{2,m},\cV_{0,m})$ is some map.
The representation
$$\tilde \rho_{1}=\left(\begin{array}{cc}\rho_{0}\otimes 1&\sigma\otimes x\\0&\rho_{2} \otimes 1\end{array} \right)\ .$$
determines an exact sequence of locally constant  sheaves of finitely generated projective $R[x]$-modules
$$\tilde \cV:0\to \cV_{0}\otimes_{R} R[x]\to \tilde \cV_{1}\to \cV_{2}\otimes_{R}R[x]\to 0$$
such that $s_1^{*}\tilde \cV_{1}\cong \cV_{1}$ and
$s_{0}^{*}\tilde \cV_{1} \cong  \cV_{0}\oplus \cV_{2} $.
Since the topological cycle map is additive this implies
$$[\cV_{1}]=[s^{*}_{1}\tilde \cV]=[s_{0}^{*}\tilde \cV]= [\cV_{0}\oplus \cV_{2}]= [\cV_{0}]+[\cV_{2}]\ .$$ \hB

\subsection{Kamber-Tondeur forms }\label{klejflkewfjwlfjlwefwefwef}

{In this subsection we recall the construction of Kamber-Tondeur forms. Its main result is the   construction of the  maps \eqref{ofewjfljwelkfjwewefiwejfoiwefef} for all $j\ge 0$. They represent   cohomology classes  
  which are the building blocks of the Borel regulator whose construction we will recall in Subsection \ref{wejfwlef}.}

\bigskip

{
A reference for the following constructions with differential forms is
\cite{MR1303026}. We consider a flat vector bundle with a hermitean metric $\bV=(V,\nabla^{V},h^{V})$  on a smooth manifold $M$. {If $M$ has boundary faces, then we require that the metric   $h^{V}$ is normally parallel near the boundary faces (we say shortly: ``with product structure'').}
Using the metric we define the adjoint connection $\nabla^{V,*_{h^{V}}}$ on $V$. It is characterized by the identity
$$dh^{V}(\phi,\psi)=h^{V}(\nabla^{V}\phi,\psi)+h^{V}(\phi,\nabla^{V,*_{h^{V}}}\psi)$$
for sections $\phi,\psi$ of $V$.
The metric $h^{V}$ is called flat (or equivalently, $\nabla^{V}$ is hermitean with respect to $h^{V}$) if \begin{equation}\label{aug1002}\nabla^{V,*_{h^{V}}}=\nabla^{V}\ .\end{equation}
 In general it is not possible to choose a flat metric  on $V$.
  The deviation from flatness of $h^{V}$ is measured by the difference
 \begin{equation}
 \label{jul0594} 
\omega(h^{V}):=\nabla^{V,*_{h^{V}}}-\nabla^{V}\in \Omega^{1}(M,\End(V))\ . \end{equation}
This form can be used to define the
 Kamber-Tondeur forms (see \cite[(0.2)]{MR1303026}) by
\begin{equation}\label{jun102}\omega_{k}(\bV):=\frac{1}{(2\pi i)^{j}2^{k}}\Tr\, \omega(h^{V})^{k}\in Z^{k}(\Omega(M))\end{equation}
for all odd $k=2j+1\ge 1$. {A priori,} {the right-hand side of \eqref{jun102} is just a complex-valued  form}. But one checks that
the Kamber-Tondeur forms are  indeed real and closed (the latter property follows from the flatness of $\nabla^{V}$).} 
The cohomology class $$[\omega_{k}(\bV)]\in H^{k}(M;\R)$$
does not depend on the choice of the metric $h^{V}$. 
\begin{rem}{\rm The formula (\ref{jun102}) defines Kamber-Tondeur forms in the Bismut-Lott normalization
which is dictated by the local index theory methods of \cite{MR1303026}
culminating in the construction of higher analytic torsion forms. 
{The relation of the Kamber-Tondeur forms with the Borel regulator will be explained in}  Subsection \ref{jun274333}. {In Subsection \ref{jul1705} we compare various different normalizations appearing in the literature}.
}\hB\end{rem}}

{
We will also consider the form
$$\omega_{0}(\bV):=\dim(V)\in Z^{0}(\Omega(M))$$
as a Kamber-Tondeur form. For even $k>0$ we set $\omega_{k}:=0$.

\bigskip

In the following we record some properties of the Kamber-Tondeur forms  which follow immediately from the definition.
\begin{prop} 
\label{feb1704}
 \begin{enumerate}
\item {\rm Naturality:} If $f:M^{\prime}\to M$ is a smooth map and $\bV:=(V,\nabla^{V},h^{V})$ is a flat bundle with metric on $M$, then we can define the pull-back $$f^{*}\bV:=(f^{*}V,\nabla^{f*V},h^{f^{*}V})\ .$$ In this situation  we have the equality
$$f^{*}\omega_{k}(\bV)=\omega_{k}( f^{*}\bV)\ .$$
\item {\em Additivity:} For two flat bundles with metrics  $\bV_{i}$, $i=0,1$ on $M$   we can form the sum $\bV=\bV_{0}\oplus \bV_{1}$. In this case we have   the equality
 \begin{equation}\label{add}\omega_{k}(\bV)=\omega_{k}(\bV_{0})+\omega_{k}(\bV_{1})\ .\end{equation} 
\item {\em Conjugation:} Let $\bar \bV$ be the conjugate bundle. Then
for all $2j+1=k>0$ we have
$$\omega_{k}(\bar \bV):=(-1)^{j}\omega_{k}(\bV)\ .$$ 
\item {\em Vanishing for flat metrics:} If $\bV$ is a flat bundle with  a flat metric, then \begin{equation}\label{okt2403}\omega_{k}(\bV)=0\end{equation} for all $k>0 $. 
\item {\em Tensor products and determinants:} \begin{equation}\label{jbjkqwdkqwdkqdqkwdqwdwd}
\omega_{1}(\bV_{0}\otimes \bV_{1})=\omega_{1}(\bV_{0})\wedge \omega_{0}(\bV_{1})+\omega_{0}(\bV_{0})\wedge \omega_{1}(\bV_{1})\ , \quad  
\omega_{1} (\bV)=\omega_{1}(\det(\bV))\ .\end{equation}\end{enumerate}
\end{prop}
}

\begin{rem}{\rm 
{
Our next goal is to explain how the Kamber-Tondeur form $\omega_{k}$ induces a primitive cohomology class
\begin{equation}\label{erhgkleghlkrgrgregergergerge}
\omega_{k}\in H^{k}(\Nerve(i\bP(\C) ;\R)\ .
\end{equation} We will actually not just construct such a class but a map \eqref{ofewjfljwelkfjwewefiwejfoiweffefeef} of commutative monoids in $\Spc$ representing this class whose relation with the geometric construction is encoded in the diagram \eqref{dzqwduwqidqwdwqdwqdqdqwdqwdd} below. This additional precision will be the key to the construction of the differential cycle map.  We use the complex case to demonstrate the method. In Subsection \ref{jun273} we will employ   the same method in the case of number rings.
The construction will again use some $\infty$-categorical technology we will explain in passing.  }}\hB \end{rem}

\bigskip

{
 
{Using  \eqref{qwdqwdwqdqdu3hkj321e12e},  a commutative monoid $M$  in topological spaces (or in sets considered as discrete topological spaces) gives rise to a {commutative} monoid $M_{\infty}$ {in spaces.} If $\bC$ is a symmetric monoidal category, then $\pi_{0}(\Nerve(\bC))$ is a {commutative} monoid in sets. {We 
 can thus} form the commutative monoid $\pi_{0}(\Nerve(\bC))_{\infty}$ in spaces.
\begin{lem}
If   $\bC$ is a symmetric monoidal category, then we have a  {canonical}  map  \begin{equation}\label{cfwefwefewfewfwfwfwf4r23r23r1} \Nerve(\bC)\to \pi_{0}(\Nerve(\bC))_{\infty}\end{equation}  in $ \CommMon(\Spc)\ .$ \end{lem}
\proof 
The natural transformation $\id\to\pi_{0}(-)_{\infty}$ between endofunctors on $\Spc$  is the unit of the adjunction
$$\pi_{0}(-) :\Spc\leftrightarrows  \Nerve(\Set):(-)_{\infty}$$
 between symmetric monoidal functors. {It} therefore induces a natural transformation  on the level of commutative monoids.
 \hB
 
 Note that the transformation \eqref{cfwefwefewfewfwfwfwf4r23r23r1} is natural in $\bC$ by construction.

\bigskip

For a chain complex $C\in \Ch$ let $Z^{0}(C)$ denote the abelian group of zero cycles.
\begin{lem}
If $C$ is a chain complex, then we have a natural  map  \begin{equation}\label{cfwefwefewfewfwfwfwf4r23r23r}  Z^{0}(C)_{\infty}\to \Omega^{\infty} (H(\sigma C))\ .
\end{equation} 
 in $ \CommMon(\Spc) $.
 \end{lem}
 \proof 
 We have a natural inclusion of chain complexes $Z^{0}(C)\to \sigma C$. Applying the Eilenberg-MacLane correspondence and $\Omega^{\infty}$ we get the map
\begin{equation}\label{jkhdkjdhkqehdqwdwqdqwd23423432432edwq}
\Omega^{\infty} H(Z^{0}(C))\to \Omega^{\infty} (H(\sigma C))
\end{equation}  in $ \CommMon(\Spc) $. 
 
 In general, if 
  $A$ is an abelian group, then {by Lemma \ref{fjelfjwkfjekjfwlipo31}} we have a natural equivalence $$H(A)\simeq \spp(A_{\infty})\ .$$ In view of the equivalence \eqref{okt2002}, by an application of $\Omega^{\infty}$ we get the equivalence
$$\Omega^{\infty}(H(A))\simeq A_{\infty}\ .$$ We apply this to $A:=Z^{0}(C)$ in order to get the first {equivalence} in the composition
$$Z^{0}(C)_{\infty}\simeq \Omega^{\infty} (H(Z^{0}(C)))\stackrel{\eqref{jkhdkjdhkqehdqwdwqdqwd23423432432edwq}}{\to} \Omega^{\infty}(H(\sigma C))\ .$$ \hB

Note again, that by construction the map \eqref{cfwefwefewfewfwfwfwf4r23r23r} is natural in the chain complex $C\in \Ch$.

 \bigskip

 {
   We consider the symmetric monoidal stack  $$\loc_{geom}(\C)\in \Sm^{desc}(\CAlg(\Nerve(\Cat)[W^{-1}]))$$ which  associates to a manifold $M$ the symmetric monoidal category of flat complex vector bundles  with metric  {(with product structure)}
$(V,\nabla^{V},h^{V})$ and isomorphisms,  and to a smooth map  $f:M^{\prime}\to M$ the pull-back. The symmetric monoidal structure is   given by the direct sum. We further define the smooth monoid
$$\bloc_{geom}(\C):=\pi_{0}(\Nerve(\loc_{geom}(\C)))\in \Sm(\CommMon(\Set))\ .$$
 
}

   By Proposition \ref{feb1704} the Kamber-Tondeur form $\bV\mapsto \omega_{k}(\bV)$ provides a 
 map of smooth monoids
 \begin{equation}\label{ergergelkjlk2r32398ro32r23r23r}
 \bloc_{geom}(\C) \to Z^{0}( \Omega \R[k])\ .
\end{equation}   
  }

{
If we combine   {the three constructions explained above}, 
then we get a map    \begin{equation}\label{fwekjhfwjkfkwefefwefewfewfewffwefwefewf} \Nerve(\loc_{geom}(\C))\stackrel{\eqref{cfwefwefewfewfwfwfwf4r23r23r1}}{\to} \bloc_{geom}(\C)_{\infty}\stackrel{\eqref{ergergelkjlk2r32398ro32r23r23r}}{\to} Z^{0}(\Omega \R[k])_{\infty}\stackrel{\eqref{cfwefwefewfewfwfwfwf4r23r23r}}{\to} \Omega^{\infty}(H(\sigma \R[k]))
\end{equation}
in $\Sm^{desc}(\CommMon(\Spc))$.
}

\bigskip

{
We have a forgetful morphism  $$\loc_{geom}(\C)\to \loc(\C)$$ between symmetric monoidal stacks. It induces a map \begin{equation}\label{fewffewfwfwe234423424324324}
\Nerve(\loc_{geom}(\C))\to \Nerve(\loc(\C))
\end{equation} in 
$\Sm^{desc}(\CommMon(\Spc))$. By Lemma \ref{jkljljdqlwdkwqdqdqd}  the object $ \Nerve(\loc(\C))$ is homotopy invariant.{ {Below, in Lemma \ref{ejfwelfjwelfefwefewf},}
  we will give a categorical interpretation of the map \eqref{fewffewfwfwe234423424324324}. Let $\bC$ be a presentable $\infty$-category.
 We let $$\Sm^{h}(\bC)\subseteq \Sm(\bC)$$ denote the full subcategory of homotopy invariant smooth objects of $\bC$.  This  inclusion is the right adjoint of an adjunction
$$\cH^{pre}:\Sm(\bC)\leftrightarrows \Sm^{h}(\bC):inclusion$$
whose left adjoint is called homotopification.
We let furthermore $$\Sm^{desc,h}(\bC)\subseteq \Sm^{desc}(\bC)$$ denote the full subcategory of homotopy invariant smooth objects of $\bC$ which in addition satisfy descent. We again have an adjunction \begin{equation}\label{jfkjfljkl4f24f2f24r}
\cH:\Sm^{desc}(\bC)\leftrightarrows \Sm^{desc,h}(\bC):inclusion\ .
\end{equation}
 These facts are shown in Lemma \ref{oiwioewoifwoifjwfoj245239852345}.

\bigskip

 Let $E\to F$ be a morphism in $\Sm^{desc}(\bC)$. \begin{ddd}\label{dlqwkjqlwkdjwqldjwqdou231ildejodidqw}
We say that this morphism presents $F$ as a homotopification of 
$E$, {if $F$ is homotopy invariant and the canonical morphism
$\cH(E)\to F$ given by the universal property of the homotopification}
 is an equivalence.
\end{ddd}

%

\begin{lem}\label{ejfwelfjwelfefwefewf}
The  
map $$\Nerve_{geom}(\loc(\C))\stackrel{\eqref{fewffewfwfwe234423424324324}}{\to}\Nerve(\loc(\C))$$ presents $ \Nerve(\loc(\C))$ as a homotopification of $ \Nerve(\loc_{geom}(\C))$.
\end{lem}
The proof will be given at the end of this subsection.}}
  We further use that by Lemma \ref{djqlwdjqwldjwqldwqdwqdwqdwqdq}
$$\Omega^{\infty}(H(\sigma \Omega \R[k]))\stackrel{\eqref{frfwefwewfwf2kjkjlkjljljlkjaljldw34}}{\to} \Omega^{\infty}(H(\Omega \R[k]))$$
presents $\Omega^{\infty}(H(\Omega \R [k]))$ as the homotopification of $\Omega^{\infty}(H(\sigma \Omega \R[k]))$.
If we apply the homotopification functor \eqref{jfkjfljkl4f24f2f24r} to the map \eqref{fwekjhfwjkfkwefefwefewfewfewffwefwefewf}, then we get the middle square of
the following diagram 
  \begin{equation}\label{dzqwduwqidqwdwqdwqdqdqwdqwdd}
\xymatrix{&\Nerve(\loc_{geom}(\C))\ar[d]\ar[r]^{\eqref{fwekjhfwjkfkwefefwefewfewfewffwefwefewf}}&\Omega^{\infty}(H(\sigma \Omega \R[k]))\ar[d]    &  \\\Funk(\Nerve(i\bP(\C)) \ar[r]_{\eqref{jwqjejqejkwqe989898i234jkkjkjerqkjwr}}^{\simeq}\ar@/_1cm/[rrr]^{\omega_{k}}&
\Nerve(\loc(\C))\ar[r]&\Omega^{\infty}(H(\Omega\R[k]))\ar[r]_{\eqref{jhekjhekjhkd12983ud128oidkl1d1d}}^{\simeq}&\Funk(\Omega^{\infty}(H\R[k]))} 
\end{equation}  in $ \Sm^{desc}(\CommMon(\Spc))$. The  lower horizontal composition defines $\omega_{k}$.
} 
  
  \bigskip
  
  
  \bigskip
  
 {The evaluation of $\omega_{k}$ at the one point manifold $*$ defines the  map \begin{equation}\label{ofewjfljwelkfjwewefiwejfoiweffefeef}
   \omega_{k}(*):\Nerve(i\bP(\C)) \to  \Omega^{\infty}(H\R[k])\end{equation}
 in $\CommMon(\Spc)$. 
 This map represents the primitive cohomology class \eqref{erhgkleghlkrgrgregergergerge}.} It extends to a map, also denoted by $\omega_{k}$,  between spectra as follows. {The target of \eqref{ofewjfljwelkfjwewefiwejfoiweffefeef} is a group in $\Spc$.
By the universal property of the group completion we get  the dotted arrow and the upper triangle in
\begin{equation}\label{ldkjqwdkljlqwdjqlwdqwdqwd}\xymatrix{ \Nerve(i\bP(\C))\ar[r]^{\omega_{k}(*)}\ar[d]& \Omega^{\infty}(H\R[k])\\\Omega B( \Nerve(i\bP(\C))) \ar@{..>}[ur]\ar[r]^{\simeq}&\Omega^{\infty} (K\C)\ar[u]}\ .\end{equation}
The lower horizontal equivalence reflects the definition of $K\C$ and the right vertical map and the lower triangle are defined so that the latter commutes.}
 {{The adjoint of the right vertical map} with respect to \eqref{okt2002} is {the desired  map}  \begin{equation}\label{ofewjfljwelkfjwewefiwejfoiwefef}\omega_{k}:K\C\to H\R[k]\end{equation} in $\Sp$. 
}
 

 \bigskip 
 
{\proof[of Lemma \ref{ejfwelfjwelfefwefewf}] 
We must show that the natural map
$$\cH( \Nerve(\loc_{geom}(\C)))\to  \Nerve(\loc(\C))$$ is an equivalence in $\Sm(\CommMon(\Spc))$.
As explained in the proof of \cite[Lemma 6.3]{Bunke:2013uq}
 it suffices to show that the corresponding map of the underlying smooth spaces  is an equivalence, i.e. we can neglect the commutative monoid structure.  We now use the criterion
Lemma \ref{dkhkqlwdjqlwdwqdwdwqdqwdw}.
We must show that \begin{equation}\label{jdhqwkjdhkqwdwqdwqdwqd}
\bs(\Nerve(\loc_{geom}(\C)))\to \bs(\Nerve(\loc(\C)))
\end{equation}
is an equivalence in $\Sm(\Spc)$, where $\bs$ is defined by \eqref{dqwhdkjqwhdwqdqdhwqkd89ue12e12ee1e21e1e12e12}. If we forget the symmetric monoidal structure, as mentioned in  Remark \ref{kldlqwdqwdwqdqwdqwd1}, we can choose a strict model $\loc^{\flat}(\C)\in \Sm(\Nerve(\Cat))$ of the stack $\loc(\C)\in \Sm(\Nerve(\Cat)[W^{-1}])$.  We can further assume that
$\loc_{geom}(\C)$ is represented by a strict model
$\loc^{\flat}_{geom}(\C)\in \Sm(\Nerve(\Cat))$ such that
$\loc^{\flat}_{geom}(\C)(M)$ is the category of pairs $((V,\nabla^{V}),h^{V})$  of flat vector bundles   $(V,\nabla^{V})\in \loc^{\flat}(\C)(M)$ and a metric $h^{V}$ on $V$.
Then the map of simplicial spaces $$\Nerve(\loc_{geom}(\C)^{Y(\Delta^{\bullet})}(M))\to \Nerve(\loc(\C)^{Y(\Delta^{\bullet})}(M))$$
is represented by the map of bi-simplicial sets
$$\Nerve(\loc^{\flat}_{geom}(\C)(M\times \Delta^{\bullet}))\to \Nerve(\loc^{\flat}(\C)(M\times \Delta^{\bullet}))\ .$$
We now fix the simplicial degree $p\in \nat$ in the direction of the nerve and obtain a map of simplicial sets 
\begin{equation}\label{efwefkjfkjewfhewkfhk32iouo23e32e32e23e32e3}
  \Nerve(\loc^{\flat}_{geom}(\C)(M\times \Delta^{\bullet}))[p]\to \Nerve(\loc^{\flat}(\C)(M\times \Delta^{\bullet}))[p]\ . 
\end{equation}
We claim that this map is an acyclic Kan fibration.
An  element in $ \Nerve(\loc^{\flat}(\C)(M\times \Delta^{\bullet}))[p]$ is a $p$-tuple of isomorphisms of flat bundles
$$\cV:=(V_{0}\stackrel{\cong}{\to} \dots\stackrel{\cong}{\to} V_{p})\ .$$ The fibre of \eqref{efwefkjfkjewfhewkfhk32iouo23e32e32e23e32e3} over $\cV$ is in bijection to  the set of all metrics on $V_{0}$.
We now observe that a metric prescribed on  
$M\times \partial \Delta^{n}$ (or on $M\times \Lambda$ for some horn $\Lambda\subset \Delta^{n}$)
can be extended to $M\times \Delta^{n}$. 
At this point  the assumption that the metric has a product structure is crucial in order to  ensure the possibililty of choosing an  extension which is smooth at the corners.
{In order to construct such an extension from $M\times \partial \Delta^{n}$ to $M\times \Delta^{n}$ we first extend the metric to a tubular neighourhood of the boundary  using the parallel transport in the normal direction. We then glue with an arbitrary metric over the rest of the simplex using a partition of unity. In the case of a horn we
first use this construction in order to extend the metric over the remaining face and then proceed as above.}
We conclude that the diagonal
 $$\diag(\Nerve(\loc^{\flat}_{geom}(\C)(M\times \Delta^{\bullet})))\to \diag(\Nerve(\loc^{\flat} (\C)(M\times \Delta^{\bullet})))
$$ is a weak equivalence of simplicial sets.
Since the diagonal is a concrete model for the colimit in $\Spc$ appearing in the definition \eqref{dqwhdkjqwhdwqdqdhwqkd89ue12e12ee1e21e1e12e12} of $\bs$ we conclude that 
\eqref{jdhqwkjdhkqwdwqdwqdwqd} is an equivalence. \hB}

 \subsection{Borel's regulator}\label{wejfwlef}

{In this subsection we recall the calculation of $K_{*}(R)\otimes\R$ due to  Borel. This gives an explicit description {(see Definition \ref{kjqwldqwdqwdwqdqwdqd})} of 
the graded vector space $A$ appearing in the differential data $(KR,A,c)$ for $KR$.
{The missing piece, the map $c$, will be constructed in Subsection \ref{jun273}.}

}

\bigskip

%
%

Recall that
the algebraic $K$-groups of a  ring $R$ are   defined as the homotopy groups
$$K_{i}(R):=\pi_{i}(KR)\ ,\quad i\ge 0$$
of the algebraic $K$-theory spectrum introduced in Definition \ref{nov0101}.  The first two groups,  $K_{0}(R)$ and $K_{1}(R)$,  have a simple algebraic description which is explained e.g. in the classical book of Milnor  \cite{MR0349811}. The higher algebraic $K$-groups of a ring were first defined by Quillen \cite{MR0338129}.  
Since we consider connective $K$-theory we have $K_{i}(R)=0$ for $i<0$ by definition.

\bigskip

We now consider a number field $k$, i.e. a finite  field extension $\Q\subseteq k$ of the field of  rational numbers.
We let $R\subset k$  denote the ring of integers in $k$,  defined as the integral closure of $\Z$ in $k$.

\begin{ex}\label{fkjwelfwefewfwf}{\rm 
A typical example, relevant in Subsection \ref{jul0110}, is the cyclotomic field $k:=\Q(\xi)$ for an $r^\textrm{th}$  root of unity $\xi$, {where $r\ge 2$}. {Its}  ring of integers
is given by $$R\cong \Z[\xi]/(1+\xi+\dots+\xi^{r-1})\ .$$  
}\hB \end{ex}

\bigskip 

{For a number ring $R$}  we have isomorphisms
$$K_{0}(R)\cong \Z\oplus \Cl(R)\ ,\quad K_{1}(R)\cong R^{*}\ ,$$
where $\Cl(R)$ denotes the ideal class group of $R$ and $R^{*}\subset R$ is the group of multiplicative units. The ideal class group $\Cl(R)$ is finite {(see \cite[Thm. 6.3]{MR1697859})}, and the summand $\Z$ measures the rank of free $R$-modules. Dirichlet's unit theorem describes the group of units $R^{*}\subset R$ as an extension of 
$\Z^{r_{\R}+r_{\C}-1}$ by the finite group of roots of unity in $R$, where the integers $r_{\R}$ and $2r_{\C}$ denote the {number} of real and complex embeddings of $k$, {respectively}.

\begin{ex}\label{fkjwelfwefewfwf1}{\rm 
In our Example \ref{fkjwelfwefewfwf} of    the cyclotomic field $\Q(\xi)$ with an $r^\textrm{th}$ root of unity $\xi$, if $r\ge {3}$ is a prime, we have
$r_{\R}=0$ and $r_{\C}=\frac{r-1}{2}$.
In contrast if
{
 $r=2$, then  we have $r_{\R}=1$ and $r_{\C}=0$.}
 }\hB \end{ex}

The higher algebraic $K$-groups are  much more complicated.  For a number ring $R$
the ranks of the abelian groups $K_{i}(R)$ are finite and have been calculated by Borel
 \cite[Prop. 12.2]{MR0387496}. Recall that
$\pi_{i}(KR\R)\cong \pi_{i}(KR)\otimes \R$.
Together with the classical calculations for $i=0,1$ \cite{MR0349811} we have: \begin{theorem}[Borel]\label{borel}
If $R$ is a number ring, then $\dim_{\R}(\pi_{i}(KR\R))$ is four-periodic for $i\ge 2$ and  given by the following table
$$\begin{array}{|c||c|c|c|c|}\hline i \:\:(j\ge 1)&4j-2&4j-1&4j&4j+1\\\hline \dim_{\R}(\pi_{i}(KR\R))&0&r_{\C}&0&r_{\R}+r_{\C}\\\hline\end{array}$$
where $r_{\R}$ and $2r_{\C}$ are  the numbers of real and complex  embeddings.
Furthermore, we have
$\pi_{0}(KR\R)\cong \R$
 and
 $\pi_{1}(KR\R)\cong R^{*}\otimes \R$, where $R^{*}$ is the finitely generated group of units of $R$ which has rank $r_{\R}+r_{\C}-1$ by Dirichlet's unit theorem.
\end{theorem}
A more precise statement {which includes the choice of a dual basis} will be given in Proposition \ref{jun261}.
It is furthermore known by a result of Quillen that the groups $K_{*}(R)$ for the ring of integers of a number field are finitely generated.

%
%

\bigskip

{
 Let $\Sigma(R)$ denote the set of embeddings $R\hookrightarrow \C$. If the ring is clear, then we simplify the notation and write $\Sigma:=\Sigma(R)$.
 The group
$\Z/2\Z$ acts on $\Sigma$  by complex conjugation.
The set $\Sigma$ equivariantly decomposes as the disjoint  union 
$\Sigma=\Sigma_{\R}\sqcup \Sigma_{\C}$ of the subsets of real and complex embeddings,
where $\Sigma_{\R}=\Sigma^{\Z/2\Z}\subseteq \Sigma$ is the subset of fixed points. Note that
$$r_{\R}=|\Sigma_{\R}|\ ,\quad r_{\C}=|\Sigma_{\C}/(\Z/2\Z)|\ .$$ 
 
 \bigskip
 
By the functoriality of the $K$-theory spectrum \eqref{jkjfl23jf32fjl2opip3r2r}, for an embedding $\sigma\in \Sigma$ we get a map of spectra $KR\to K\C$. 
We compose this map with   \eqref{ofewjfljwelkfjwewefiwejfoiwefef} in order to define  a map of spectra \begin{equation}\label{kqedlkjdlqkwjlwqkjqo32jle}
\omega_{k}(\sigma):KR\to H\R[k]\ .\end{equation}
We will use the same symbol for the cohomology class represented by this map. 
 A more precise version of the calculation of Borel stated in  Theorem \ref{borel}
is now given by the following proposition.
\begin{prop}\label{jun261} Let $R$   a number  ring and $\Sigma:=\Sigma(R)$.
 \begin{enumerate}
 \item For $j\ge 0$ and $\sigma\in \Sigma$ we have the relation
$\omega_{2j+1}(\sigma)=(-1)^{j}\omega_{2j+1}(\bar \sigma)$. 
  \item  If $j\ge 1$, then the family of classes
$(\omega_{2j+1}(\sigma))_{\sigma\in \Sigma}$ generate    $H^{2j+1} (KR;\R)$ subject to the relation 1.  \item
The  cohomology classes $(\omega_{1}(\sigma))_{\sigma\in  \Sigma}$  generate
$H^{1}(KR;\R)$ subject to the relation given in 1. and the additional relation $\sum_{\sigma\in \Sigma} \omega_{1}(\sigma)=0$.  
\item For $\sigma\in \Sigma$ the cohomology class $\omega_{0}:=\omega_{0}(\sigma)$    is a basis of  $H^{0}(KR;\R)$ and independent of the choice of $\sigma$.
\end{enumerate}
 \end{prop}
 \proof
 In Subsection \ref{jun274333} we will verify that the cohomology classes used by Borel  to detect the homotopy of $KR$ coincide, up to normalization,
with the Kamber-Tondeur classes \eqref{kqedlkjdlqkwjlwqkjqo32jle}. {The main observation is that the matrix $(\psi_{\sigma,\sigma^{\prime}})_{\sigma,\sigma^{\prime}\in \Sigma}$ studied in Proposition \ref{jul0231} is diagonal with non-vanishing entries.}
\hB

\begin{rem}{\rm 
Observe that the relation 1. implies that $\omega_{2j+1}(\sigma)=0$ if $\sigma$ is real and
$j$ is odd. It is easy to see that the table of dimensions given in Theorem \ref{borel} is compatible with the statement of Proposition \ref{jun261}.}\hB \end{rem}

\begin{ddd}\label{kjqwldqwdqwdwqdqwdqd}
We define the graded vector space $\tilde A(R)$ with an action of $\Z/2\Z$ as follows:
\begin{enumerate}
\item $\tilde A(R)_{k}:=\R^{\Sigma(R)}$ for $k=0$ or $k\ge 0$ and odd.
\item $\tilde A(R)_{k}:=0$ for $k> 1$ even or $k<0$.
\end{enumerate}
For $j\in \nat$ the action of the non-trivial element of $\Z/2\Z$ on $\tilde A_{2j+1}(R)$ by  is given by $$f\mapsto \bar f(\sigma):=(-1)^{j}f(\bar \sigma)\ .$$ The action of $\Z/2\Z$ on $\tilde A(R)_{0}$ is induced by the action on $\Sigma(R)$. 

\bigskip

We furthermore define the graded vector subspace $$A(R)\subset \tilde A(R)$$ by
\begin{enumerate}
\item $A(R)_{0}:=\{f\in  \tilde A(R)_{0}\:|\: \mbox{$f$ is constant}\}$
\item $A(R)_{1}:=\{f \in  \tilde A(R)_{1}^{\Z/2\Z}\:|\: \sum_{\sigma\in \Sigma(R)}f(\sigma)=0\} $  \item $A(R)_{2j+1}:= \tilde A(R)_{2j+1}^{\Z/2\Z}$ for $j\ge 1$
\item $A(R)_{k}:=0$ for $k<0$, or  $k>0$ and even.
\end{enumerate}
 \end{ddd}
 
 If the ring is clear, then we shorten the notation and write $A:=A(R)$.

\bigskip

By Proposition \ref{jun261} the graded vector space $A(R)$ is an explicit model for the dual of the cohomology $H\R^{*}(KR)$, or equivalently, for $K_{*}(R)\otimes \R$.  
For $\sigma\in \Sigma$ we let $\delta_{\sigma}\in \R^{\Sigma}$ be the function which takes the zero value everywhere except at the point $\sigma$ where the value is $1$. For $k=0$ or $k>0$ and odd we can consider this function as an element $\delta_{k,\sigma}\in \tilde A(R)_{k}$. 
We further consider the cohomology classes $\omega_{k}(\sigma)$ as homomorphisms
\begin{equation}\label{dkqwdjqwkldlwqkwdwqdooiuo123213213}
\omega_{k}(\sigma):K_{k}(R)\otimes \R \to \R\ .
\end{equation}
\begin{kor}\label{jlefjewljfewlfjewfu23oilrkf}
We have an isomorphism (the Borel regulator) \begin{equation}\label{kdqwkdjlwqkdjwqkdjq987}
c:\sum_{k} \sum_{\sigma\in \Sigma} 
 \omega_{k}(\sigma)\delta_{k,\sigma}:K_{*}(R)\otimes \R\stackrel{\cong}{\to} A(R)_{*}\ ,\end{equation}
  where  the outer sum is taken over $k=0$ and $k\ge 1$, odd.
\end{kor}
Apriori the map \eqref{kdqwkdjlwqkdjwqkdjq987} takes values in $\tilde A(R)$, but the relations between the cohomology classes $\omega_{k}(\sigma)$ listed in Proposition \ref{kjqwldqwdqwdwqdqwdqd} imply that the map takes values in the subspace $A(R)$.
 
 \bigskip

}

Let $\phi:R\to R^{\prime}$ be a homomorphism of number rings.   If $\sigma:R^{\prime}\to \C$ is an embedding of $R^{\prime}$,  then
$\phi^{*}\sigma:=\sigma\circ \phi:R\to \C$ is the induced embedding of $R$ . We thus get an induced map
$$\phi^{*}:\Sigma(R^{\prime})\to\Sigma(R)$$
which is $\Z/2\Z$-equivariant. The pull-back of functions with $\phi^{*}$ provides a $\Z/2\Z$-equivariant map of graded vector spaces
$$\phi_{*}:\tilde A(R)\to \tilde A(R^{\prime})\ .$$ It obviously restricts to a map
of subspaces
$$\phi_{*}:A(R)_{k}\to A(R^{\prime})_{k}$$ for all $k\not=1$.  
In order to see that it restricts in the case $ k=1$ as well,  we use that 
$|  \phi^{*,-1}(\sigma)|=[k^{\prime}:k]$ is independent of $\sigma\in \Sigma(R)$ and therefore for $f\in \R^{\Sigma(R)}$
$$\sum_{\sigma^{\prime}} f(\phi^{*}(\sigma^{\prime}))=\sum_{\sigma\in \Sigma(R)} \sum_{\{\sigma^{\prime}\in \Sigma(R^{\prime})\:|\: \phi^{*}(\sigma^{\prime})=\sigma\}} f(\sigma)  =[k^{\prime}:k]\sum_{\sigma\in \Sigma(R)} f(\sigma)\ .$$
We therefore have defined a functor
$$A:\Rings\to \Ab_{\Z-gr}\ , \quad R\mapsto A(R)\ .$$
As a direct consequence of the definitions the Borel regulator is natural in the number ring $R$.
\begin{kor}
If $\phi:R\to R^{\prime}$ is a morphism between number rings, then the following diagram commutes:
$$\xymatrix{K_{*}(R)\ar[r]^{c}\ar[d]^{\phi_{*}}&A(R)_{*}\ar[d]^{\phi_{*}}\\K_{*}(R^{\prime})\ar[r]^{c}&A(R^{\prime})_{*}}\ .$$
\end{kor}

\subsection{Characteristic forms} \label{jun273}

In this subsection we introduce a characteristic form \eqref{cfr} for locally constant sheaves of finitely generated  $R$-modules with geometry. We show how it can be used to refine Borel's isomorphism  given in Corollary   \ref{jlefjewljfewlfjewfu23oilrkf}  to a spectrum level  map $c:KR\to H(A)$. This map will be the third entry   in the  data $(KR,A,c)$ for differential algebraic $K$-theory of a number ring $R$.

\bigskip

Let $M$ be a smooth manifold and  $R$ be {a} number ring. We consider a locally constant sheaf $\cV$ of finitely generated  $R$-modules on $M$. 
For an embedding $\sigma\in \Sigma$ we can form the locally constant sheaf of finite-dimensional complex vector spaces $\cV_{\sigma}:=\cV\otimes_{\sigma} \C$. This sheaf is the sheaf of parallel sections of a uniquely determined complex vector bundle $V_{\sigma}\to M$ with a flat connection $\nabla^{V_{\sigma}}$. This vector bundle will be referred to as the complex vector bundle associated to $\cV$ and $\sigma$.

  \begin{ddd}\label{defpre} A geometry on a locally constant sheaf $\cV$ of finitely generated $R$-modules is the choice of a family
$h^{\cV}:=(h^{V_{\sigma}})_{\sigma\in \Sigma}$ of hermitean metrics ({with product structure}) on the bundles $V_{\sigma}$
such that 
\begin{equation}\label{jun101}h^{V_{\bar \sigma}}=\bar h^{V_{\sigma}}\end{equation} for 
the complex places $\sigma$ of $R$. We call the pair $(\cV,h^{\cV})$ a geometric $R$-bundle.
\end{ddd}
 
For $\sigma\in \Sigma$ we let $\bV_{\sigma}:=(V_{\sigma},\nabla^{V_{\sigma}},h^{V_{\sigma}})$ denote the associated flat complex vector bundle with metric.

\bigskip

\bigskip

Recall the definition \eqref{jun102} of the Kamber-Tondeur forms. 
 {Let $(\cV,h^{\cV})$ be a geometric $R$-bundle.}  If $(\sigma,\bar \sigma)$ is a conjugate pair of  complex places of $k$, then $\bV_{\bar \sigma}\cong \bar \bV_{\sigma}$ and consequently by Proposition \ref{feb1704}, 3.,   $$\omega_{k}(\bV_{\sigma})= (-1)^{j}\omega_{k}(\bV_{\bar \sigma})$$
 for all $2j+1=k>0$. This implies that the relation described in Proposition \ref{jun261}, 1. holds true on the form level. Unfortunately, the additional relation Proposition \ref{jun261}, 3. in the degree-one case is not true on the form level. But this defect can be corrected. By Lemma 5.21 of \cite{buta} there exists a natural choice of a  function $\lambda(\cV,h^{\cV})\in C^{\infty}(M)$  such that
 $$\sum_{\sigma\in \Sigma} \omega_{1}(\bV_{\sigma})=d\lambda(\cV,h^{\cV})\ .$$
\begin{rem}{\rm For completeness let us recall some details of its construction. We form the complex line bundle $$\bW:=\bigotimes_{\sigma\in \Sigma} \det(\bV_{\sigma})\ .$$  It has a natural    $\Z$-structure which can be used to define a normalized flat metric. Its   associated Kamber-Tondeur form $\omega_{1}$ vanishes by Proposition \ref{feb1704}, 4. 
 We then let $\lambda(\cV,h^{\cV})$
be the corresponding  transgression of $\omega_{1}$ such that
$$d\lambda(\cV,h^{\cV})=\omega_{1}(\bW)\stackrel{\eqref{jbjkqwdkqwdkqdqkwdqwdwd}}{=}\sum_{\sigma\in \Sigma} \omega_{1}(\bV_{\sigma})\ .$$
As a   consequence of the construction  the function $\lambda(\cV,h^{\cV})$ depends on the data only via the determinant. In particular  we have the identity
\begin{equation}\label{fkejfhwefhwkef}
\lambda(\cV,h^{\cV})=\lambda(\det(\cV),h^{\det(\cV)})\ .
\end{equation}
}\hB\end{rem}
 It is clear that the map $(\cV,h^{\cV})\mapsto d\lambda(\cV,h^{\cV})$ is additive   and natural with respect to pull-back. For $k\in \nat$, $k=0$ or odd,  we consider the elements \begin{equation}\label{gdqjwhqjwduiiuiuiuiui}\kappa_{k}:=\sum_{\sigma\in \Sigma} \delta_{k,\sigma}\in \tilde A(R)_{k} \ .\end{equation} 
 The following definition is now 
guided by 
\eqref{kdqwkdjlwqkdjwqkdjq987}. 
{\begin{ddd}\label{defpre1}
We define the characteristic form of a geometric $R$-bundle $(\cV, h^{\cV})$ by
\begin{equation}\label{cfr}\omega (\cV,h^{\cV}):= \sum_{k}  \sum_{\sigma\in \Sigma}     \delta_{k,\sigma} \omega_{2j+1}(\bV_{\sigma})-d\lambda(\cV,h^{\cV})\kappa_{1}    \in Z^{0}(\Omega A(M))\ ,\end{equation}   where  the outer sum is taken over $k=0$ and $k\ge 1$, odd. \end{ddd}
The right-hand side of \eqref{cfr} apriori belongs to $Z^{0}(\Omega \tilde A(R)(M))$. But the relations between the Kamber-Tondeur forms  discussed above and the addition of the term involving $\kappa_{1}$ ensure that we define an element of the subgroup $Z^{0}(\Omega A(M))$.

 }

 \bigskip
 
%

 In order to formalize the properties of the characteristic form for geometric $R$-bundles we introduce the stack of geometric $R$-bundles
 $$\loc_{geom}(R)\in \Sm^{desc}(\CAlg(\Nerve(\Cat)[W^{-1}]))\ .$$
 It associates to a manifold $M$ the symmetric monoidal category of geometric $R$-bundles and isomorphisms, and to a map $f:M^{\prime}\to M$ the pull-back. The symmetric monoidal structure is given by the direct sum. We further consider  the symmetric monoidal substack $$\loc_{geom}^{proj}(R)\subseteq \loc_{geom}(R) $$
 of geometric $R$-bundles with projective fibres.   
 In view of the definition  of a geometry these stacks fit into pull-back squares   \begin{equation}\label{knknlknldnqwldnqldnqwd}\xymatrix{\loc^{proj}_{geom}(R)\ar[r]\ar[d]&\loc_{geom}(R)\ar[d]\ar[r]&\prod_{\sigma\in \Sigma^{*}} \loc_{geom}(\C)\ar[d]\\\loc^{proj}(R)\ar[r]&\loc(R)\ar[r]&\prod_{\sigma\in \Sigma^{*}}\loc(\C) } ,\end{equation}
 where the  vertical maps forget geometry, and the right horizontal maps are given by the collection of functors $\cV\mapsto V_{\sigma}$ and $(\cV,h^{\cV})\mapsto \bV_{\sigma}$ for $\sigma\in \Sigma^{*}$, where $\Sigma^{*}\subset \Sigma$ is a choice of representatives of the orbit space  $\Sigma/(\Z/2\Z)$. We further define the smooth monoids
 $$\bloc_{geom}(R):=\pi_{0}(\Nerve(\loc_{geom}(R)))\ , \quad  \bloc^{proj}_{geom}(R):=\pi_{0}(\Nerve(\loc^{proj}_{geom}(R)))$$ in $\Sm(\CommMon(\Set))$.
 As an immediate consequence of the construction and Proposition \ref{feb1704} we get:
 \begin{kor}
The characteristic form \eqref{cfr} induces a map between smooth monoids\footnote{{The notation $\hat R$ resembles the notation $R$ for the curvature map of a differential cohomology theory and is not related to the ring $R$.}}

 \begin{equation}\label{kdqkjdqldjqldjqldqwdqw312}\hat R:\bloc_{geom}(R)\to Z^{0}(\Omega A)\ .\end{equation}\end{kor}

%

%
%

{We now proceed step by step as in Subsection \ref{klejflkewfjwlfjlwefwefwef} in order to obtain a map of spectra
$$c:KR\to H(A)\ ,$$ the analog of \eqref{ofewjfljwelkfjwewefiwejfoiwefef}. }
{In analogy to \eqref{fwekjhfwjkfkwefefwefewfewfewffwefwefewf}
we get a map  in $\Sm^{desc}(\CommMon(\Spc))$ \begin{equation}\label{fwekjhfwjkfkwefefwefewfewfewffwefwefewf1} \Nerve(\loc^{proj}_{geom}(R))\stackrel{\eqref{cfwefwefewfewfwfwfwf4r23r23r1}}{\to} \bloc^{proj}_{geom}(R)_{\infty}\stackrel{\eqref{kdqkjdqldjqldjqldqwdqw312}}{\to} Z^{0}(\Omega A)_{\infty}\stackrel{\eqref{cfwefwefewfewfwfwfwf4r23r23r}}{\to} \Omega^{\infty}(H(\sigma \Omega A))
\end{equation}
}

{
The left vertical map in \eqref{knknlknldnqwldnqldnqwd}    induces a map \begin{equation}\label{fewffewfwfwe2344234243243241}
\Nerve(\loc^{proj}_{geom}(R))\to \Nerve(\loc^{proj}(R))
\end{equation} in 
$\Sm^{desc}(\CommMon(\Spc))$. By Lemma \ref{jkljljdqlwdkwqdqdqd}  the object $ \Nerve(\loc^{proj}(R))$ is homotopy invariant. Recall Definition \ref{dlqwkjqlwkdjwqldjwqdou231ildejodidqw}.
\begin{lem}
The  
map \eqref{fewffewfwfwe2344234243243241} presents $ \Nerve(\loc^{proj}(R))$ as a homotopification of $ \Nerve(\loc^{proj}_{geom}(R))$.
\end{lem}
\proof
This follows from Lemma \ref{ejfwelfjwelfefwefewf} if we use the presentation
\eqref{knknlknldnqwldnqldnqwd} and the fact that homotopification preserves pull-backs over homotopy invariant objects (Lemma \ref{ehfwkejhfewkfwdwqwqdqdqd}).
Alternatively one can show this lemma directly. \hB

We further use that by Lemma \ref{djqlwdjqwldjwqldwqdwqdwqdwqdq}
$$\Omega^{\infty}(H(\sigma \Omega A))\stackrel{\eqref{frfwefwewfwf2kjkjlkjljljlkjaljldw34}}{\to} \Omega^{\infty}(H(\Omega A))$$
presents $\Omega^{\infty}(H(\Omega A))$ as the homotopification of $\Omega^{\infty}(H(\sigma \Omega A))$.
If we apply the homotopification functor to the map \eqref{fwekjhfwjkfkwefefwefewfewfewffwefwefewf1}, then we get the middle square of
the following diagram 
  \begin{equation}\label{dzqwduwqidqwdwqdwqdqdqwdqwdd1}
\xymatrix{&\Nerve(\loc^{proj}_{geom}(R))\ar[d]\ar[r]^{\eqref{fwekjhfwjkfkwefefwefewfewfewffwefwefewf}}&\Omega^{\infty}(H(\sigma \Omega A))\ar[d]^{\eqref{frfwefwewfwf2kjkjlkjljljlkjaljldw34}}    &  \\\Funk(\Nerve(i\bP( R)) \ar[r]_{\eqref{jwqjejqejkwqe989898i234jkkjkjerqkjwr}}^{\simeq}\ar@/_1cm/[rrr]^{r}&
\Nerve(\loc^{proj}(R))\ar[r]&\Omega^{\infty}(H(\Omega A))\ar[r]_{\eqref{jhekjhekjhkd12983ud128oidkl1d1d}}^{\simeq}&\Funk(\Omega^{\infty}(H (A)))} 
\end{equation}  in $ \Sm^{desc}(\CommMon(\Spc))$. The lower horizontal composition defines the map $r$.
} 
  
  \bigskip
  
  
  \bigskip
  
 {Evaluation of $r$ at $*$ defines the  map \begin{equation}\label{ofewjfljwelkfjwewefiwejfoiweffefeef1}
   r(*):\Nerve(i\bP(R)) \to  \Omega^{\infty}(H(A))\end{equation}
 in $\CommMon(\Spc)$. 
By the universal property of the group completion we get  the dotted arrow and the upper triangle in
\begin{equation}\label{ldkjqwdkljlqwdjqlwdqwdqwd1}\xymatrix{ \Nerve(i\bP(R))\ar[r]^{r(*)}\ar[d]& \Omega^{\infty}(H (A))\\\Omega B( \Nerve(i\bP(R)))) \ar@{..>}[ur]\ar[r]^{\simeq}&\Omega^{\infty} (KR)\ar[u]}\ .\end{equation}
The lower horizontal equivalence reflects the definition of $K R$ and the   right vertical map and the lower triangle are defined so that the latter commutes.}
 { 
Its adjoint with respect to \eqref{okt2002} is the desired  map  \begin{equation}\label{ofewjfljwelkfjwewefiwejfoiwefef1}c:KR\to H (A)\end{equation} in $\Sp$. 
}
 
For later use we note the following fact. If we combine the diagram \eqref{dzqwduwqidqwdwqdwqdqdqwdqwdd1} with $\Funk(\eqref{ldkjqwdkljlqwdjqlwdqwdqwd1})$ and \eqref{fwekjhfwjkfkwefefwefewfewfewffwefwefewf1}, then we get a commuting diagram
\begin{equation}\label{dwqdqwdqwkdqiejio12je12e1}\xymatrix{\Nerve(\loc^{proj}_{geom}(R))\ar[dr]\ar[rr]&&\Omega^{\infty} (H(\sigma \Omega A))\ar[dd])\\\Funk(\Nerve(i\bP(R)))\ar[dr]^{r(*)}\ar[d]\ar[r]^{\simeq}&\Nerve(\loc^{proj}(R))\ar[rd]&\\
\Funk(\Omega^{\infty}(K R))\ar[r]^{c}&\Funk(\Omega^{\infty}(H (A)))&\Omega^{\infty} (H(\Omega A))  \ar[l]^{\eqref{jhekjhekjhkd12983ud128oidkl1d1d}}_{\simeq}}
\end{equation}
in $\Sm(\CommMon(\Spc))$.

\bigskip

Let now $\sigma\in \Sigma$ be an embedding $\sigma: R\to \C$ and $k\in \nat$ be zero or $k\ge 1$ and odd. Using that $ A_{k}\subset \R^{\Sigma}$ by Definition \ref{kjqwldqwdqwdwqdqwdqd} we have a projection $e_{\sigma,k}:A\to \R[k]$ given by evaluation at $\sigma$.

\begin{prop}\label{jkjkfekjfw78782378423874ef}
The following diagram commutes:
\begin{equation}\label{caskcaskljclkjlkwqe}\xymatrix{KR\ar[r]^{c}\ar[d]^{\sigma}&H(A)\ar[d]^{e_{\sigma,k}}\\
K\C\ar[r]^{\omega_{k}(\sigma)}&H\R[k]} .\end{equation}
\end{prop}
\proof
We give the argument for $k\not=0$. 
We have the equality
$$e_{\sigma,k} (\omega(\cV,h^{\cV}))=\omega_{k}(\bV_{\sigma})$$ in $Z^{k}(\Omega(M))$.
Now every map used in the construction of $c$ can be compared by a square with a corresponding map in the construction of $\omega_{k}(\sigma)$. Putting all these squares together we get the square \eqref{caskcaskljclkjlkwqe}. 
In the case $k=1$ we must modify the argument slightly in order to take the correction term involving $d\lambda(\cV,h^{\cV})$ into account.  \hB 
 
{The} Corollary \ref{jlefjewljfewlfjewfu23oilrkf} and Proposition  \ref{jkjkfekjfw78782378423874ef} together have the following consequence:
\begin{kor}\label{fewjflwefwelfewklfjewlflewf23u298759285234}
The map $c:KR\to H(A)$ constructed in \eqref{ofewjfljwelkfjwewefiwejfoiwefef1} induces the isomorphism \eqref{kdqwkdjlwqkdjwqkdjq987} and hence an equivalence $KR\R\to H(A)$.
\end{kor}

By Corollary \ref{fewjflwefwelfewklfjewlflewf23u298759285234} the map $c$ in \eqref{ofewjfljwelkfjwewefiwejfoiwefef1}  qualifies as a third entry in a differential data $(KR,A,c)$, see Definition \ref{aug1001}. 
\begin{ddd}\label{jhdkjdqwdqwdqwdqwdqwdwqd}
The differential data for the algebraic $K$-theory spectrum of a number ring $R$ is the triple $(KR,A,c)$, where $c$ is given by \eqref{ofewjfljwelkfjwewefiwejfoiwefef1}. 
 \end{ddd}

 \bigskip
 
 We now discuss the functoriality of the construction in the number ring $R$. Let $\phi:R\to R^{\prime}$ be a homomorphism of number rings. We then have a morphism of stacks
 $$\phi_{*}:\loc_{geom}(R)\to \loc_{geom}(R^{\prime})$$  which  refines the map \eqref{kjldjlkjlkjdl1kdj1o2idop12d12d12d1d}
   from the topological to the geometric case. The geometric version 
sends  an object $(\cV,h^{\cV})\in \loc_{geom}(R)(M)$ to $(\cV^{\prime},h^{\cV^{\prime}})\in \loc_{geom}(R^{\prime})(M)$. Its underlying sheaf  of $R^{\prime}$-modules is defined by   $\cV^{\prime}:=\cV\otimes_{R}R^{\prime}$. If $\sigma^{\prime}\in \Sigma(R^{\prime})$ and $\sigma:=\phi^{*}\sigma^{\prime}$, then we have a canonical isomorphism of flat vector bundles $V_{\sigma}\cong V_{\sigma^{\prime}}^{\prime}$. We use these isomorphisms in order to define the geometry $h^{\cV^{\prime}}$ induced by $h^{\cV}$.
  
It now immediately follows from  the Definition \ref{defpre1} that the diagram of smooth monoids
$$\xymatrix{\bloc_{geom}(R)\ar[d]^{\phi_{*}}\ar[r]^{\hat R}&Z^{0}(\Omega(A(R))\ar[d]^{\phi_{*}}\\
\bloc_{geom}(R^{\prime})\ar[r]^{\hat R}&Z^{0}(\Omega A(R^{\prime}))}$$
commutes.

An inspection of the arguments leading to the diagram \eqref{dwqdqwdqwkdqiejio12je12e1} now shows that its construction is functorial in the ring $R$. Thus let $\Box$ be the category
of the shape
\begin{equation}\label{ewfwfewfwj2lkr23iori32rio2ri322r2r2}
\xymatrix{\bullet\ar[r]\ar[d]&\bullet\ar[d]\\\bullet \ar[r]&\bullet} .
\end{equation}
\begin{kor} \label{eilwfjwelkfewfewf322r23lk23r2347249834234234}
We have constructed a functor between $\infty$-categories
$$\Nerve(\Rings) \to \Sm(\CommMon(\Spc))^{{\Nerve}(\Box)}$$
which sends the ring $R$ to the square \begin{equation}\label{fhjhfjkwefhjkeiiowef4354654656234}
\xymatrix{\Nerve(\loc^{proj}_{geom}(R))\ar[d]\ar[r]&\Omega^{\infty}(H(\sigma\Omega(A(R))))\ar[d]\\
\Funk(\Omega^{\infty}(KR))\ar[r]&\Omega^{\infty}(H(\Omega(A(R))))}\ .\end{equation}
\end{kor}

Analysing further the construction of the differential data for $KR$ we conclude:

\begin{kor}\label{efewklfjwlfwlfjelwkfwfewfewfw43242424324}
We have a functor
$\Nerve(\Rings)\to \widehat{\Sp}$
which sends a ring $R$ to the differential data 
$(KR,A(R),c)$ for $KR$. In particular, for a homomorphism $\phi:R\to R^{\prime}$ of number rings we have a natural transformation $$\phi_{*}:\widehat{KR}^{0}\to \widehat{KR^{\prime}}^{0}$$
between differential algebraic $K$-theory functors.
\end{kor}

\subsection{The cycle map for geometric locally constant sheaves}

The goal of the present {s}ubsection is to construct the differential cycle map for  geometric $R$-bundles with projective fibres, where $R$ is a number ring.   

\bigskip

We have maps of smooth monoids  
$$\hat{I}:\bloc^{proj}_{geom}(R)\to KR^{0}\ ,\quad \hat{R}:\bloc^{proj}_{geom}(R)\to Z^{0}(\Omega A)$$
which associate to the isomorphism class of a geometric $R$-bundle
$(\cV,h^{\cV})\in \bloc^{proj}_{geom}(M)$ the underlying class ${\hat I(V)}\in KR^{0}(M)$ (see Definition {\ref{jul1702} )} 
and the characteristic form
$\omega(\cV,h^{\cV})\in Z^{0}(\Omega A(M))$ given in   \eqref{cfr}{, see also \eqref{kdqkjdqldjqldjqldqwdqw312} for $\hat R$}.
\begin{ddd} \label{sep0701}
A differential cycle map   is a map of smooth monoids
$$\cycl:\bloc^{proj}_{geom}(R)\to \widehat{KR}^{0}$$ such that the diagram
\begin{equation}\label{m2159}\xymatrix{&&KR^{0}\\\bloc^{proj}_{geom}(R)\ar[urr]^{\hat I}\ar[drr]_{\hat R}\ar[r]^{\cycl}&\widehat{KR}^{0}\ar[ur]_{I}\ar[dr]^{R}\\
&&Z^{0}(\Omega A)}  
\end{equation}
commutes.
\end{ddd}
This definition should be compared with Definition \ref{jul0310}. The analog of Remark \ref{jkdljwqdlwqdqwdqwdwdqwdq} applies here. 

\begin{prop}\label{thecycle}
There exists a differential cycle map.
     \end{prop}
     \proof
     We will define a map between smooth commutative monoids \begin{equation}\label{kqdkqjddlkwqdjqwldkjqwdq98798}
     \Nerve(\loc_{geom}^{proj}(R))\to \Omega^{\infty}( \Diff(KR,A,c))\end{equation} which will induce the cycle map after application of $\pi_{0}$.
     As a right-adjoint the functor $\Omega^{\infty}$ preserves pull-backs. In view of  the {d}efinition of the differential function spectrum \eqref{difdef} we have a pull-back \begin{equation}\label{difdj2h32j3h2kj3h2khk23hef}\xymatrix{\Omega^{\infty}(\Diff(KR,A,c))\ar[d]\ar[r]&\Omega^{\infty}(H(\sigma \Omega A))\ar[d]^{{\eqref{frfwefwewfwf234}}}\\\Funk(\Omega^{\infty}(KR))\ar[r]^{\rat}& \Omega^{\infty}(H(\Omega A)) }.\end{equation} in $\Sm(\CommMon(\Spc))$.
    In order to construct the map \eqref{kqdkqjddlkwqdjqwldkjqwdq98798} we must therefore construct a commuting diagram \begin{equation}\label{jdjkdkqjwdhqdqwkdqwdkqwd}
   \xymatrix{ \Nerve(\loc_{geom}^{proj}(R)) \ar[d]\ar[r]&\Omega^{\infty}(H(\sigma \Omega A))\ar[d]^{{\eqref{frfwefwewfwf234}}}\\\Funk(\Omega^{\infty}(KR))\ar[r]^{\rat}& \Omega^{\infty}(H(\Omega A)) }.\end{equation} 
This is exactly {\eqref{fhjhfjkwefhjkeiiowef4354654656234}.}  \hB

If we take into account the addition{al} functoriality of the topological cycle map (Corollary \ref{kldjqlkwjdlqdjqdqwdjqlwkd}), the differential data (Corollary \ref{efewklfjwlfwlfjelwkfwfewfewfw43242424324}) and {of} the square \eqref{jdjkdkqjwdhqdqwkdqwdkqwd} (Corollary \ref{eilwfjwelkfewfewf322r23lk23r2347249834234234}) we conclude that the cycle map is natural in the ring.

\begin{kor}
The {map  \eqref{kqdkqjddlkwqdjqwldkjqwdq98798} is an evaluation} of a map between $\infty$-categories $$\Nerve(\Rings)\to \Sm(\CommMon(\Spc))^{{\Nerve(\bullet \to\bullet)}} \ .$$ In particular, if $\phi:R\to R^{\prime}$ is a homomorphism between number rings, then
the following diagram commutes
$$\xymatrix{\bloc^{proj}_{geom}(R)\ar[r]^{cycl}\ar[d]^{\phi_{*}}&\widehat{KR}^{0}\ar[d]^{\phi_{*}}\\
\bloc^{proj}_{geom}(R^{\prime})\ar[r]^{cycl}&\widehat{KR^{\prime}}^{0}}\ .$$
\end{kor}

\hB

In Definition \ref{uliapr1210} we will  remove the condition ``projective''. This  construction will then finish the proof of  the following theorem. 
  \begin{theorem}\label{thecycleproj1}
There exists a  differential cycle map 
$$\bloc_{geom}(R)\to \widehat{KR}^{0}$$ \end{theorem}

{
\begin{rem}{\rm 
 Note that Theorem \ref{thecycleproj1} only asserts the existence of a cycle map.
  It is clear that the  cycle map can be altered using the action of $K_{1}(R)\otimes \Q$ by automorphisms of $\widehat{KR}^{0}$ (see {Lemma \ref{feb210112}}).     }\hB \end{rem}
} 



\subsection{Some calculations with $\widehat{KR}^{0}(*)$}\label{jun274} 
{
In this subsection we discuss some calculations with the  differential algebraic $K$-theory {$\widehat{KR}^{0}(*)$ of the number ring} $R$ evaluated on a  point.  
The results serve as pedagogical examples, but will also be used  {at} various places below.
}

\subsubsection{Calculation of $\widehat{KR}^{0}(*)$}


\begin{lem}
The group $\widehat{KR}^{0}(*)$ fits into an exact sequence 
\begin{equation}\label{jun202}0\to K_{1}(R)\otimes \R/\Z \to \widehat{KR}^{0}(*)\xrightarrow{I} K_{0}(R)\to 0\ .\end{equation}\end{lem}
\proof We have canonical isomorphisms
$$KR^{-1}(*)\cong K_{1}(R)\ ,\quad \Omega A^{-1}(*)/\im(d)\cong A_{1}\cong K_{1}(R)\otimes \R\ .$$ With
these identifications the map  
$c:KR^{-1}(*)\to \Omega A^{-1}(*)/\im(d)$ in \eqref{exats111} identifies with the canonical map
$K_{1}(R)\to K_{1}(R)\otimes \R$ given by $x\mapsto x\otimes 1$. 
In particular,
$$\Omega A^{-1}(*)/(\im(d)+\im(c))\cong K_{1}(R)\otimes \R/\Z\ .$$
The exact sequence (\ref{exats111})
{specializes to} the exact  sequence \eqref{jun202} \hB 
 
\subsubsection{The action of $\Aut(KR,A,c)$.}
We will describe the action of the  
 group
$$\hom(\pi_{0}(KR\Q),\pi_{1}(KR\Q))\cong  K_{1}(R)\otimes \Q$$ on
the functor on $\widehat{KR}^{0}$. 

\bigskip

 {The automorphism corresponding to $x\in K_{1}(R)\otimes \Q$
will be denoted by $\Phi_{x}:\widehat{KR}^{0}\to \widehat{KR}^{0}$.  
The} element $x\in 
K_{1}(R)\otimes \Q$ gives naturally a class $$[x]\in K_{1}(R)\otimes \R/\Z\subset \widehat{KR}^{0}(*) .$$
We consider the set of integers $\Z$ as a discrete manifold.  Using the fact that $\Diff(KR,A,c)$ satisfies descent
 we can define the class $u_{x}\in \widehat{KR}^{0}(\Z)$ such that it restricts to $n[x]$ at the point $n\in \Z$.
We have a natural homomorphism {
$d:\widehat{KR}^{0}\to \Funk(\Z)$} of smooth groups whose evaluation at $M$ is 
$$d:\widehat{KR}^{0}(M)\stackrel{I}{\to} KR^{0}(M)\stackrel{\dim}{\to} \Hom_{\Mf}(M,\Z) \ .$$ {
\begin{lem}\label{feb210112}
For $x\in K_{1}(R)\otimes \Q$ the homomorphism  $\Phi_{x}:\widehat{KR}^{0}(M)\to \widehat{KR}^{0}(M)$ is given by 
 \begin{equation}\label{nov0102} \Phi_{x}(y)=y+d(y)^{*} u_{x}\ .\end{equation}
\end{lem}}

 \bigskip
 
\subsubsection{Determinants}\label{fewjfwlfwefweff}
 
 {We will 
  discuss the compatibility of the cycle map with forming determinants of $R$-modules.}
   The result will be important when we consider the torsion of chain complexes in the framework of differential algebraic $K$-theory.  
  
  \bigskip

 The one-dimensional free $R$-module $R$ has a canonical metric $h^{R}$ given by the standard euclidean metric $h^{R_{\sigma}}:=h^{\C}$ on $R_{\sigma}= R\otimes_{\sigma}\C\cong \C$ for all places $\sigma\in \Sigma$. We let
\begin{equation}\label{jun201}\beins:=\cycl(R,h^{R})\in \widehat{KR}^{0}(*)\ .\end{equation}
For a finitely generated projective $R$-module $V$ we consider its top exterior power $\det(V)$.
{As a consequence of the    Steinitz theorem, see \cite[Thm. 1.6]{MR0349811}}, in $K_{0}(R)$  we have the relation
$$[V]-\dim(V)=[\det(V)]-1\ .$$
We are going to refine this relation to $\widehat{KR}^{0}(*)$.
We consider the finitely generated {projective} $R$-module $V$ as a locally constant sheaf on the point $*$. In this {way} we have a well-defined  notion of a geometry, Definition \ref{defpre}.
For each complex embedding $\sigma\in \Sigma$ of $R$ we have a canonical isomorphism
$$\det(V)_{\sigma}\cong \det(V_{\sigma})\ .$$
Therefore, if we are given a geometry $h^{V}$ on $V$, then we have an induced geometry $h^{\det(V)}$
such that $h^{\det(V)_{\sigma}}$ is the metric on $\det(V_{\sigma})$ induced by $h^{V_\sigma}$ for each complex embedding $\sigma\in \Sigma$.
\begin{lem}\label{jun1820}
In $\widehat{KR}^{0}(*)$ we have the relation
\begin{equation}\label{feb2101}\cycl(V,h^{V})-\dim(V)\ \beins=\cycl(\det(V),h^{\det(V)})-\beins\ .\end{equation}
\end{lem}
\proof
We first assume that
\begin{equation}\label{jun1813}(V,h^{V})\cong \bigoplus_{i=1}^{\dim(V)-1} (R,h^{R})\oplus (\det(V),h^{\det(V)})\ .\end{equation}
The additivity of the cycle map  gives
$$\cycl(V,h^{V})\cong (\dim(V)-1)\ \beins+\cycl(\det(V),h^{\det(V)})\ .$$
Hence the assertion of the Lemma holds true in this special case.
We now use the fact (the Steinitz theorem, see \cite[Thm. 1.6]{MR0349811}) that every finitely generated projective $R$-module can be written in the form
$$V\cong R^{\dim(V)-1}\oplus \det(V)\ .$$ Hence we must compare the variations of {both sides of (\ref{feb2101})} with respect to the choice of the geometry.

We consider now a finitely generated projective $R$-module $V$ with two geometries
$h^{V}_{i}$, $i=0,1$.  Then we can form the locally constant sheaf $\cV$ on the unit interval with fibre $V$ and choose a geometry $h^{\cV}$ which interpolates between the geometries $h^{V}_{i}$ on the ends.
By the homotopy formula \eqref{feb2710}  we have
$$\cycl(V,h^{V}_{1})-\cycl(V,h^{V}_{0})=a\left(\int_{0}^{1} R(\cycl(\cV,h^{\cV}))\right) $$
and $$\cycl(\det(V),h^{\det(V)}_{1})-\cycl(\det(V),h^{\det(V)}_{0})=a\left(\int_{0}^{1} R(\cycl(\det(\cV),h^{\det(\cV)}))\right)\ .$$
We now use
\begin{eqnarray*}\label{}R(\cycl(\cV,h^{\cV}))&\stackrel{\eqref{cfr}}{=}&\sum_{\sigma\in  \Sigma} \delta_{1,\sigma} \omega_{1}(\bV_{\sigma})-d\lambda(\cV,h^{\cV})\kappa_{1}\\&\stackrel{\eqref{fkejfhwefhwkef}, \eqref{jbjkqwdkqwdkqdqkwdqwdwd}}{=}&\sum_{\sigma\in  \Sigma} \delta_{1,\sigma} \omega_{1}(\det(\bV_{\sigma}))-d\lambda(\det(\cV),h^{\det(\cV)})\kappa_{1}\\&\stackrel{\eqref{cfr}}{=}&R(\cycl(\det(\cV),h^{\det(\cV)}))\ ,\end{eqnarray*}
where $\kappa_{1}$ is defined in \eqref{gdqjwhqjwduiiuiuiuiui}. 
We get
$$\cycl(V,h^{V}_{1})-\cycl(V,h^{V}_{0})=\cycl(\det(V),h^{\det(V)}_{1})-\cycl(\det(V),h^{\det(V)}_{0})\ .$$
This implies the desired result. 
\hB

\subsubsection{Rescaling the metric}\label{efjwefkefewffewffewfefewf}

 Here we analyse the effect of scaling the geometry.
 
 \bigskip
 
  Let $(V,h^{V})\in \loc_{geom}^{proj}(R)(*)$ be a projective $R$-module of rank one with geometry $h^{V}=(h^{V_{\sigma}})_{\sigma\in \Sigma}$. Given a collection $\lambda =(\lambda_{\sigma})_{\sigma\in \Sigma}$
  of non-zero real numbers such that $\lambda_{\bar \sigma}=\lambda_{\sigma}$ we can define the rescaled metric
$\lambda^{2} h^{V}:=(\lambda_{\sigma}^{2} h^{V_{\sigma}})_{\sigma\in \Sigma}$.
In order to simplify formulas we introduce the projections of $\delta_{1,\sigma}\in \tilde A^{-1}$ to $A^{-1}$ given by 
\begin{equation}\label{dh12j1l2jlkjdljlj21ldjl12jd21d}\bar \delta_{1,\sigma}:=\delta_{1,\sigma}-\frac{1}{|\Sigma|} \kappa_{1}\in A^{-1}\ .\end{equation}
The formula of the following lemma is closely related with \cite[Eq. (28)]{MR1724894}, see also Subsubsection \ref{jul0481}).
\begin{lem}\label{jul0420}
In $\widehat{KR}^{0}(*)$ we have the following identity:
$$\cycl(V,\lambda^{2} h^{V})-\cycl(V,h^{V})=a\left( \sum_{\sigma\in \Sigma}  \ln( |\lambda_{\sigma}|) \ \bar \delta_{1,\sigma} \right)\ .$$

\end{lem}
\proof
We consider the locally constant sheaf of $R$-modules $\cV$ on the unit interval $[0,1]$ generated by $V$.
We define a geometry $h^{\cV}$ on $\cV$ by linearly interpolating between the geometries $h$ and $\lambda^{2} h$.
Hence $$h^{V_{\sigma}}(t)=(1+t(\lambda^{2}_{\sigma}-1)) h^{V_{\sigma}}\ .$$
This gives by (\ref{jun102})
$$\omega_{1}(\bV_{\sigma})(t)= \frac{1}{2}\frac{(\lambda^{2}_{\sigma}-1)dt}{1+t(\lambda_{\sigma}^{2}-1)}\ .$$
We have 
$$R(\cycl(\cV,h^{\cV}))=\kappa_{0} +\frac{1}{{2}}\sum_{\sigma\in \Sigma}\frac{(\lambda_{\sigma}^{2}-1)dt}{1+t(\lambda_{\sigma}^{2}-1)} \bar \delta_{1,\sigma} \ .$$
 By the homotopy formula  \cite[Equation (1)]{MR2608479} we have 
\begin{eqnarray*}\cycl(V,\lambda^{2} h^{V})-\cycl(V,h^{V})&=&a\left(\int_{0}^{1}   R(\cycl(\cV,h^{\cV}))\right)\\
&=& a\left({\sum_{\sigma\in \Sigma}\frac{1}{{2}}\int_{0}^{1}\frac{(\lambda_{\sigma}^{2}-1)}{1+t(\lambda_{\sigma}^{2}-1)} dt \ \bar \delta_{1,\sigma}} \right)\\&=&a\left( \sum_{\sigma\in \Sigma} \ln(|\lambda_{\sigma}| ) \ \bar \delta_{1,\sigma} \right)\ .\end{eqnarray*}
 \hB

\subsection{Extension of the cycle maps from projective to finitely generated $R$-bundles} \label{may121270}

{In general, the cohomology $H^{i}(M;\cV)$ of a compact manifold $M$ with coefficients in a {locally constant} sheaf $\cV$ of finitely generated projective $R$-modules is a {finitely generated $R$-module   {(see Lemma \ref{kjhdqkwdqwdqwdqwdqwdqwdwqd})}, but not necessarily projective, $R$-module.  In order to define the  analytic index in Subsection \ref{jun183}}   it is therefore essential to extend the topological and differential cycle maps
{(Definition \ref{jul1702} and Definition \ref{sep0701})} to all geometric $R$-bundles.
 This is the main goal of the present subsection.   
 
 \subsubsection{Extension of the topological cycle map}

 \begin{rem}{\rm The natural approach would be to consider differential $G$-theory and to use the equivalence of $K$ and $G$-theory for number rings. Since our approach is based on group completion,  we  {cannot} apply it to $G$-theory directly. }
\end{rem}
}

\bigskip

{Let $R$ be some integral domain. If $V$ is an $R$-module, then an element $v\in V$ is called a torsion element if
there exists a {nonzero} element $a\in R$ such that $av=0$. The subset of torsion elements is a submodule
$\Tors(V)\subseteq V$. The quotient $V/\Tors(V)$ does not contain torsion elements. An $R$-module with $\Tors(V)\cong \{0\}$ is called torsion-free.
}

{
\begin{ass}\label{feklwflewewewfewfewfw78z8234234324324234}
 In the following we assume that $R$ is a regular noetherian integral domain with the property that every finitely generated torsion-free $R$-module is projective. \end{ass}
Note that a number ring satisfies Assumption \ref{feklwflewewewfewfewfw78z8234234324324234}.}
{See \cite{Bass62} for some related results.}

\bigskip

Let $M$ be a manifold and $T$ be a locally constant sheaf of finitely generated torsion $R$-modules.
Then we can construct a canonical resolution \begin{equation}\label{jwejflejewfkllkjweflkewf8767887}
0\to E(T)\to F(T)\to T\to 0\ , \end{equation}
where $E(T)$ and $F(T)$ are locally constant sheaves  of  finitely generated 
projective $R$-modules. The locally constant sheaf $F(T)$ is determined by the condition that on every contractible open subset $U\subset M$ we have 
$F(T)(U)=R[T(U)]$, the free $R$-module generated by $R[T(U)]$.
There is a canonical augmentation $F(T)\to T$, and $E(T)$ is defined as its kernel.
As a subsheaf of a locally constant sheaf of finitely generated projective $R$-modules it is again
of this type.

We consider the difference
\begin{equation}\label{may121206}\cZ(T):=[F(T)]-[E(T)]\in KR^{0}(M)\ .\end{equation}
In the following we show that this class does not depend on the choice of the resolution of $T$.
Let $$0\to A\to B\to T\to 0$$ be
any resolution of $T$ by locally constant sheaves of finitely generated projective $R$-modules.
\begin{lem}\label{may121001}
In $KR^{0}(M)$ we have
$$[B]-[A]=\cZ(T)\ .$$
\end{lem}
\proof
We consider a degree-wise surjective morphism of resolutions. It extends to a diagram
\begin{equation}\label{may121002}\xymatrix{&0\ar[d]&0\ar[d]&&\\&K\ar[d]\ar[r]^{\cong}&L\ar[d]&&\\
0\ar[r]&U\ar[r]\ar[d]&V\ar[d]\ar[r]&T\ar@{=}[d]\ar[r]&0\\
0\ar[r]&A\ar[r]\ar[d]&B\ar[r]\ar[d]&T\ar[r]&0\\
&0&0&&}\ .\end{equation}
{Since we assume that $R$ is regular we can apply Proposition \ref{jwlkjclkjwlkwjclcwewfefeklw9889234}  to}  conclude that
$$[V]-[U]=[B]-[A]$$
in $KR^{0}(M)$. 
In order to finish the argument we consider a second resolution $$0\to A^{\prime} \to B^{\prime}\to T\to 0$$ of the torsion sheaf $T$.
We then must find a resolution  which maps degree-wise surjectively to both.
To this end we consider the sum
$$0\to A\oplus A^{\prime}\to B\oplus B^{\prime}\to T\oplus T\to 0$$
and form the pull-back along the diagonal $T\to T\oplus T$. We get the resolution
$$0\to A\oplus A^{\prime}\to V\to T\to 0\ .$$
The natural projections now induce the required degree-wise surjective maps. 
\hB

Let $\cV$ be a locally constant sheaf of finitely generated $R$-modules on a manifold $M$. {Then we can define the torsion subsheaf $\Tors(\cV)\subseteq \cV$ as follows. For $U\subseteq M$ a section $a\in \cV(U)$ belongs to $\Tors(\cV)(U)$
if for every $u\in U$ the stalk $a_{u}\in \cV_{u}$ is a torsion element. We define the projective quotient as the quotient $\Proj(\cV):=\cV/\Tors(\cV)$ in sheaves of $R$-modules on $M$. By our assumptions on $R$ we know that   $\Proj(\cV)$ is a sheaf of finitely generated
projective $R$-modules.

\begin{rem}{\rm It is not true in general, that $\Tors(\cV)(U)$ consists of torsion elements. For example,  consider the sheaf $\cV$ of $\Z$-modules
on the zero-dimensional manifold $\nat$ with $\cV_{m}\cong\Z/m\Z$.
Then $\Tors(\cV)=\cV$, but the section $\nat \ni m\mapsto [1]\in \Z/m\Z$ is not a torsion element in
$\cV(\nat)\cong \prod_{m\in \Z} \Z/m\Z$.
}\hB
\end{rem}

If $\cV$ is a sheaf of $R$-modules on $M$, then 
there exists a natural exact sequence
\begin{equation}\label{apr123013}0\to \Tors(\cV)\to \cV\to \Proj(\cV)\to 0\ .\end{equation}
{It is clear that the functors  $\Tors$ and $\Proj$ preserve direct sums.}
}

\bigskip

{
By Definition \ref{jul1702} we have a map of smooth commutative monoids
$$\hat I:\bloc^{proj}(R)\to KR^{0}\ .$$
We define the smooth commutative monoid
\begin{equation}\label{iwjiejfwofwefewfewflkewfjlfefwef}
\bloc(R):=\pi_{0}(\loc(R))\in \Sm(\CommMon(\Set))
\end{equation} (see \eqref{gerggerg43r79834534534543fgerge} for $\loc(R)$),
{which associates to a manifold $M$ the commutative monoid of isomorphism classes of locally constant sheaves of finitely generated $R$-modules on $M$.} Then we have an inclusion
$\bloc^{proj}(R)\subseteq \bloc(R)$, and the topological cycle map on $\bloc^{proj}(R)$  is already defined in Definition \ref{jul1702}.

\begin{ddd}\label{may121311} {Suppose that $R$ satisfies Assumption \ref{feklwflewewewfewfewfw78z8234234324324234}.} Then
we define the extension  of the topological cycle  map
 \begin{equation}\label{may121250}
[...]= \hat I:\bloc(R)\to KR^{0}\end{equation}
by the formula
$$[\cV]:=[\Proj(\cV)]+\cZ(\Tors(\cV))\ , \quad  \cV\in \bloc(M)\ .$$
\end{ddd}

\begin{rem}\label{jqgdqhjwqqdwqwqdqw}{\rm
{At the moment, the extension $\hat I$ is well-defined as a transformation between set-valued functors since the sequence \eqref{apr123013} is canonical and preserved under pull-back along smooth maps. 
The additivity of $\hat I$ requires an additional justification. It follows from the more general result Lemma \ref{fjweklfefwefwef} below. Alternatively, one can verify additivity directly using Lemma \ref{may121001} and the fact that a sum of two torsion sheaves can be resolved by the sum of resolutions of the summands.}
}\hB \end{rem}
}

\subsubsection{Generalized additivity}

{
We keep the Assumption {\ref{feklwflewewewfewfewfw78z8234234324324234}}   that $R$ is a {regular noetherian integral domain}   such that any submodule of a finitely generated projective $R$-module is again projective.}
{We consider a short exact sequence
\begin{equation}\label{kwejfkejwflefewfwf8979832424}\cV: 0\to \cV_{0}\to \cV_{1}\to \cV_{2}\to 0\end{equation}
of objects in $\loc(R)(M)$. The following Lemma generalizes Proposition \ref{jwlkjclkjwlkwjclcwewfefeklw9889234}.
\begin{lem}\label{fjweklfefwefwef}
In $KR^{0}(M)$ we have $$0=\sum_{i=0}^{2} (-1)^{i}[\cV_{i}]\ .$$ 
\end{lem}
\proof
By Definition \ref{may121311} we have
$$\sum_{i=0}^{2}(-1)^{i} [\cV_{i}]=\sum_{i=0}^{2}(-1)^{i} [\Proj(\cV_{i})]+ \sum_{i=0}^{2}(-1)^{i} [\Tors(\cV_{i})]\ .$$
We now use that \eqref{kwejfkejwflefewfwf8979832424}   induces a diagram
\begin{equation}\label{kbelkjlerjekjrleglkrejglregregregergwfewff}
\xymatrix{
 &&0\ar[d]&0\ar[d]&0\ar[d]&\\\Tors(\cV):&0\ar[r]&\Tors(\cV_{0})\ar[r]\ar[d]&\Tors(\cV_{1})\ar[r]^{!}\ar[d]&\Tors(\cV_{2})\ar[r]\ar[d]&0\\\cV:&0\ar[r]&\cV_{0}\ar[r]\ar[d]^{!!}&\cV_{1}\ar[r]\ar[d]&\cV_{2}\ar[r]\ar[d]&0\\\Proj(\cV):&0\ar[r]&\Proj(\cV_{0})\ar[r]\ar[d]&\Proj(\cV_{1})\ar[r]\ar[d]&\Proj(\cV_{2})\ar[r]\ar[d]&0\\
&&0&0&0&}
\end{equation}
   in $\loc(R)(M)$ of vertical and horizontal exact (in the sense of sheaves) sequences.  Indeed, the arrow {named} $!$ is surjective since the arrow {named} $!!$ locally splits.}
By   Proposition \ref{jwlkjclkjwlkwjclcwewfefeklw9889234} we have
$$0=\sum_{i=0}^{2}(-1)^{i} [\Proj(\cV_{i})]\ .$$
In order to conclude the assertion of the lemma we must show it in the case of a short exact sequence of torsion sheaves which is then applied to the upper horizontal sequence.

  \bigskip 
  
  From  now on we assume  that the  $\cV_{i}$ are torsion sheaves for $i=0,1,2$.
In this case we consider the diagram
\begin{equation}\label{vwvwnenwvewvnewvnwvnw23849823u43242343242}
\xymatrix{
 &&0\ar[d]&0\ar[d]&0\ar[d]&\\&0\ar[r]&E(\cV_{0})\ar[r]^{!!}\ar[d]&E(\cV_{1})\ar[r]\ar[d]&\cB_{2}\ar[r]\ar[d]^{!!!}&0\\ &0\ar[r]&F(\cV_{0})\ar[r]^{!}\ar[d]&F(\cV_{1})\ar[r]\ar[d]&\cA_{2}\ar[r]\ar[d]^{!!!!}&0\\&0\ar[r]&\cV_{0}\ar[r]\ar[d]&\cV_{1}\ar[r]\ar[d]&\cV_{2}\ar[r]\ar[d]&0\\
&&0&0&0&}
\end{equation}  where $F(-)$ and $E(-)$ are as in \eqref{jwejflejewfkllkjweflkewf8767887}.
The arrow marked by $!$ is injective, since $F(-)$ preserves injections.
The arrow marked by $!!$ is obtained by the universal property of the kernel and is injective, too. The objects $\cA_{2}$ and $\cB_{2}$ are defined as quotients such that the horizontal sequences are exact. The map  marked by $!!!$  is again obtained by universal properties of the quotient and is injective. The map $!!!!$ is surjective. Finally one checks that the right column is exact in the middle. We now claim that $\cA_{2}$ is projective. This can be checked on fibres. We fix a point $m\in M$. Then
the fibre of $\cA_{2,m}$ at $m$ is  freely generated by  the complement of the image of the inclusion $\cV_{0,m}\to \cV_{1,m}$. It follows that $\cB_{2}$ is projective, too.

We
apply  Proposition \ref{jwlkjclkjwlkwjclcwewfefeklw9889234} to the two upper horizontal sequences  and Lemma \ref{may121001} to the three vertical sequences in order to conclude that
$$0=\sum_{i=0}^{2}(-1)^{i} [\cV_{i}]\ .$$
\hB

Let now $$\cV: 0\to \cV^{0}\to \cV^{1}\to \dots \to \cV^{n}\to 0$$
be a   complex of objects $\cV^{i}\in \loc(R)(M)$. It has cohomology objects
  $H^{i}(\cV)\in \loc(R)(M)$ for $i=0,\dots,n$.  The following proposition generalizes Proposition \ref{jwlkjclkjwlkwjclcwewfefeklw9889234}.

 \begin{prop}\label{may121101}
 In $KR^{0}(M)$ we have the identity
 $$\sum_{i=0}^{n}(-1)^{n} [\cV^{i}]=\sum_{i=0}^{n} (-1)^{n}[H^{i}(\cV)]\ .$$
\end{prop}
\proof
We consider the difference of the left- and right hand sides
$$\delta:=\sum_{i=0}^{n}(-1)^{n} [\cV^{i}]-\sum_{i=0}^{n} (-1)^{n}[H^{i}(\cV)] \ .$$
We argue by induction on $n$ that $\delta=0$. The lemma is obvious in the case {$n=0$}.

\bigskip

{Assume now that $n\ge 1$. Then we} have the exact sequences
$$0\to  H^{0}(\cV)\to \cV^{0}\to \cV^{0}/Z^{0}(\cV)\to 0\ ,$$
$$0\to \cV^{0}/Z^{0}(\cV)\to Z^{1}(\cV)\to H^{1}(\cV)\to 0\ , $$
 and
$$0\to Z^{1}(\cV) \to \cV^{1} \to \cV^{1}/Z^{1}(\cV)\to 0\ .$$
We now use Lemma \ref{fjweklfefwefwef} three times in order to conclude that
$$[H^{0}(\cV)]-[H^{1}(\cV)]=[\cV^{0}]-[\cV^{1}]+[\cV^{1}/Z^{1}(\cV)]\ .$$
We  can therefore replace
$\cV$ by 
$$0\to \cV^{1}/Z^{1}(\cV)\to \dots \to \cV^{n}\to 0$$
without changing $\delta$. By induction we now conclude that $\delta=0$. \hB

\subsubsection{Extension of the differential cycle map}

{We now assume that $R$ is a number ring. In particular, it satisfies Assumption \ref{feklwflewewewfewfewfw78z8234234324324234}.}
{Let $T$ be again a locally constant sheaf of finitely generated torsion $R$-modules on a manifold $M$. The canonical geometry on $R$ (see Subsubsection \ref{fewjfwlfwefweff}) induces a geometry on the   sheaf of free $R$-modules  $F(T)$ generated by $T$.  
 Note that the induced maps  $E(T)_{\sigma}\to F(T)_{\sigma}$ are isomorphisms of complex vector bundles for all embeddings $\sigma\in \Sigma$.  Using these isomorphisms we obtain a geometry $h^{E(T)}$ on $E(T)$ by restriction. We consider the difference
\begin{equation}\label{may121207}\hat \cZ(T):=\cycl(F(T),h^{F(T)})-\cycl(E(T),h^{E(T)})\in \widehat{KR}^{0}(M)\ .\end{equation}
We again show that this class is actually independent of the choice of the resolution if the geometries are related as above. So
let  $$0\to A\to B\to T\to 0$$ be any resolution and $h^{B}$ a metric on $B$. Then again we define a metric $h^{A}$ by restriction. 
\begin{lem}\label{may121201}
In $\widehat{KR}_{flat}^{0}(M)$ we have
$$\hat \cZ(T)=\cycl(B,h^{B})-\cycl(A,h^{A})\ .$$
\end{lem}
\proof
Since $(A_{\sigma},h^{A_{\sigma}})$ and $(B_{\sigma},h^{B_{\sigma}})$ are unitarily equivalent  for all $\sigma\in \Sigma$ we conclude that
$\cycl(B,h^{B})-\cycl(A,h^{A})$ is flat. 
It follows from a homotopy argument that this difference does not depend on the choice of $h^{B}$.
We now consider a diagram \eqref{may121002}. A metric $h^{V}$ determines metrics on $U,K,L$ by restriction, and a quotient metric on $B$. The quotient metric on $A$ coincides with the metric obtained by restriction from $B$.
Furthermore, $K$ and $L$ become isometric, and the complexifications of the  two vertical sequences in  \eqref{may121002} are isometric for every $\sigma\in \Sigma$.
We now apply Lott's relation \eqref{jul2023} to the two vertical sequences {of sheaves of projective $R$-modules} which have equal analytic torsions  in order to conclude that
$$\cycl(B,h^{B})-\cycl(A,h^{A})=\cycl(V,h^{V})-\cycl(U,h^{U})\ .$$
We now finish the argument for the independence of the choice of the resolution similarly as in the proof of 
Lemma \ref{may121001}.
\hB 
}

\begin{rem}{\rm 
{In order to avoid a cicular argument note that in the proof of Lemma \ref{may121201} we use Lott's relation for complexes sheaves of projective $R$-modules which is a Theorem by \cite{buta}. 
The Lemma will be employed to verify Lott's relation in general (without projectivity assumption) in the proof of Lemma \ref{uliapr1501}.
 }}\hB \end{rem}

\begin{ex}{\rm 
{Here is an example which shows that $\hat \cZ(T)$ can be  non-trivial.
We consider the ring $$R:=\Z[a]/(a^{2}-2)$$ and  the torsion module $T$ defined by
$$0\to R\stackrel{5+a}{\to} R\to T\to 0\ .$$ We calculate
$\cZ(T)\in \widehat{KR}^{0}(*)$. We first observe that $I(\hat \cZ(T))=[R]-[R]=0$.
Hence there exists $f\in \Omega A^{-1}(*)\cong A^{-1}$ such that $a(f)=\hat \cZ(T)$.
We have
$$\hat \cZ(T)=\cycl(R,h_{can})-\cycl(R,h)\ ,$$
where $h_{can}$ is the canonical geometry (see Subsection \ref{fewjfwlfwefweff}), and the geometry $h$ has the following description.
We have two real embeddings $\sigma_{\pm}:R\to \C$ given by $a\mapsto \pm \sqrt{2}$. Then
$h_{\sigma_{\pm}}=\frac{1}{|5\pm \sqrt{2}|^{2}} h_{can,\sigma_{\pm}}$.
We conclude with  Lemma \ref{jul0420}  that
$$f= -\ln  (5+ \sqrt{2})  \:\bar \delta_{1,\sigma_{+}} - \ln (5-\sqrt{2})\: \bar \delta_{1,\sigma_{-}} 
\ .$$
The units of $R$ are spanned by 
$-1$ and  $1+a$ (note that $(a+1)(a-1)=1$). We have by a similar calculation as above
$$c(R^{*})= \Z\left[ \ln  (1+\sqrt{2}) \bar \delta_{1,\sigma_{+}} +  \ln (1-\sqrt{2}) \: \bar \delta_{1,\sigma_{-}} \right]\subset A^{-1}\ .$$
We see that $f\not\in c(R^{*})$ and hence
$\hat \cZ(T)=a(f)\not=0$.
}
}\hB \end{ex}

{

{
If $\cV$ is a locally constant sheaf of finitely generated $R$-modules, then for every place $\sigma$ of $R$ the canonical map
$V_{\sigma}\to \Proj(\cV)_{\sigma}$  
of complex vector bundles is an isomorphism which is used to transfer a geometry $h^{ \cV}$ on $\cV$ to a geometry $h^{\Proj(\cV)}$ on $\Proj(\cV)$. Here $\Proj(\cV)_{\sigma}$ denotes the complex vector bundle associated to the locally constant sheaf $\Proj(\cV)\otimes_{R,\sigma} \C$}

\begin{ddd}\label{uliapr1210} We  extend the differential cycle
 map Definition \ref{sep0701} to
a map between {smooth commutative monoids}
\begin{equation}\label{may121260}\cycl:\bloc_{geom}(R)\to \widehat{KR}^{0}\  , \quad \cycl(\cV,h^{\cV}):=\cycl(\Proj(\cV),h^{\Proj(\cV)})+\hat \cZ(\Tors(\cV))\ . \end{equation}
\end{ddd}

\begin{rem}{\rm
{In order to justify this definition note that for a smooth map $f:M^{\prime}\to M $ we obviously have the equalities \begin{eqnarray*}
 \cycl(\Proj(f^{*}\cV),h^{\Proj(f^{*}\cV)})&=&f^{*}\cycl(\Proj(\cV),h^{\Proj(\cV)})\ ,\\ \quad \hat \cZ( \Tors(f^{*}\cV))&=&f^{*}(\hat \cZ(f^{*}\Tors(\cV)))\ .\end{eqnarray*} Furthermore
the map $\cV\mapsto  \hat \cZ( \Tors( \cV))$ is additive as a consequence of  Lemma \ref{may121201}, compare with Remark \ref{jqgdqhjwqqdwqwqdqw}. The extended cycle map fits into an obvious analog of the commutative diagram \eqref{m2159}.
}}\hB \end{rem}

\subsubsection{Reidemeister torsion}
{ 
We finish this subsection with a geometric analog of Proposition \ref{may121101} in the case of $M=*$. This calculation  will play a role in Subsection \ref{may121280}.

\bigskip

{Let $R$ be a number ring.}
We consider a complex (not necessarily exact)
$$\cV:0\to V^{0} \to V^{1}\to \dots \to V^{n}\to 0$$
of finitely generated $R$-modules and a collection
$h^{\cV}:=(h^{V^{i}})$ of geometries. We call $h^{\cV}$ a geometry on the complex $\cV$. For each embedding
$\sigma\in \Sigma$ we define {the  one-dimensional complex vector spaces}
$$\det(\cV_{\sigma}):=\bigotimes_{i=0}^{n} \det(V^{i}_{\sigma})^{(-1)^{i}}\ , \quad \det(H(\cV_{\sigma})):=\bigotimes_{i=0}^{n} \det(H^{i}(\cV_{\sigma}))^{(-1)^{i}}\ ,$$ {where $\cV_{\sigma}$ is the chain complex of complex vector spaces obtained by base-change of $\cV$ along $\sigma$.}
{For every $\sigma\in \Sigma$ there} exists a canonical isomorphism (compare \cite[Sec. 1]{MR1189689})
\begin{equation}\label{ji3jfio2fj23987f9798f23f23f}\phi_{\sigma}(\cV):\det(\cV_{\sigma})\stackrel{\cong}{\to}\det(H(\cV_{\sigma}))\ .\end{equation}
The geometry $h^{\cV}$ induces for every $\sigma\in \Sigma$ a collection of metrics  $$h^{\cV_{\sigma}}=(h^{V^{i}_{\sigma}})_{i=0,\dots,n}$$ which in turn induces a metric on
$\det(\cV_{\sigma})$. \begin{ddd}We define the Reidemeister metric $h^{\det(H(\cV_{\sigma}))}_{RM}$  on $\det(H(\cV_{\sigma}))$ by the condition
 that the canonical isomorphism \eqref{ji3jfio2fj23987f9798f23f23f}
 becomes is an isometry.\end{ddd}
 }
{
  If  $$h^{\cH}=(h^{H^{i}(\cV)})_{i=1,\dots,n}$$ is  a choice of  geometries on the cohomology modules, then we get an induced geometry $h^{\det(H(\cV))}$.
\begin{ddd}\label{kewjkweelkfjflewfewfefewfeefewfewf}
For every $\sigma\in \Sigma$
the Reidemeister torsion $\tau_{\sigma}(\cV,h^{\cV},h^{\cH})\in \R^{>0}$   is {the unique {positive} real number} such that
\begin{equation}\label{may121210} {\tau_{\sigma}(\cV,h^{\cV},h^{\cH})^{2}\: h^{\det(H(\cV))_{\sigma} }=  \:h^{\det(H(\cV_{\sigma}))}_{RM}}\ .\end{equation} {We further define  
\begin{equation}\label{hkhkfhkefhfhehfekjwhfkewfhewfwfwefwefwfwe}
\ln \tau(\cV,h^{\cV},h^{\cH}):= \sum_{\sigma\in \Sigma} \ln \tau_{\sigma}(\cV,h^{\cV},h^{\cH})\:  \bar \delta_{1,\sigma}\in A^{-1}\ .
\end{equation}}
\end{ddd} We refer to  \eqref{dh12j1l2jlkjdljlj21ldjl12jd21d} for the definition of the elements $\bar \delta_{1,\sigma}\in A^{-1}$.
}

{

{
\begin{rem}\label{ewlfkjwelfjewlfewfuo3234}{\rm If we fix the geometry $h^{\cV}$, then we can define an induced geometry $h_{Hodge}^{\cH}$ on the cohomology $ H(\cV)$
by applying Hodge theory to the complexes $\cV_{\sigma}$. With this choice the Reidemeister torsion
coincides with the analytic torsion, i.e. we have the equality 
\begin{equation}\label{jhkfhkejfhekfhekfewfewfewfewfewf}
\sqrt{\prod_{i=0}^{n} \det (\Delta^{\prime}_{i,\sigma})^{i(-1)^{i}}}\stackrel{def}{=}\tau_{an}(\cV_{\sigma},h^{\cV_{\sigma}})=\tau_{\sigma}(\cV,h^{\cV},h^{\cH}_{Hodge})\ .
\end{equation}
In this formula $\Delta^{\prime}_{i,\sigma}$ is the Laplacian in degree $i$ of the complex $(\cV_{\sigma},h^{\cV_{\sigma}})$ restricted to the orthogonal complement of its kernel. We will also need the alternative formula
\begin{equation}\label{glerjglkjrelgreljrelkgjrelgregregregregreg}
\tau_{an}(\cV_{\sigma},h^{\cV_{\sigma}})=\prod_{i} |\det(d^{\prime}_{i,\sigma})|^{(-1)^{i+1}}\ ,
\end{equation}
where $d_{i,\sigma}$ is the $i$th differential of $\cV_{\sigma}$ and the superscript $\prime$ means that we consider it as a map from the orthogonal complement of its kernel to its image. Note that the definition of $|\det(-)|$ involves the metrics $h^{\cV_{\sigma}}$.

}\hB\end{rem}
}

{
\begin{ex}{\rm 
The following example is a useful tool to check  normalizations and signs.
We consider the complex
$$\cV:0\to R\stackrel{d_{0}}{\to} R\to 0\to 0$$
starting in degree $0$ with $d_{0}:=\id_{R}$. We consider the geometry $h^{\cV}$ given by the canonical metric
$h^{R}$ on $V^{0}$ and the scaled canonical metric $\lambda^{2} h^{R}$ on $V^{1}$ for a collection $\lambda=(\lambda_{\sigma})_{\sigma\in \Sigma}$ of positive real numbers such that $\lambda_{\sigma}=\lambda_{\bar \sigma}$ (see \ref{efjwefkefewffewffewfefewf}).
Then $V^{i}_{\sigma}\cong \C$ for all $\sigma\in \Sigma$ and $i\in \{0,1\}$.  
Under these identifications
the adjoint of $d_{0,\sigma}$ is given by the multiplication   $d_{0,\sigma}^{*}=\lambda^{2}_{\sigma} $. Consequently, we have the formula $\Delta_{1,\sigma}=\lambda_{\sigma}^{2}$ for the Laplacian and hence by \eqref{jhkfhkejfhekfhekfewfewfewfewfewf}
\begin{equation}\label{lklkjlkjfdlewjflewfewfewfewfewfewfewfewfewf}
\tau_{an}(\cV_{\sigma},h^{\cV_{\sigma}})=\frac{1}{\lambda_{\sigma}}\ .
\end{equation}
We can also calculate the absolute value of the determinant of the differential
$d_{0,\sigma}$. We observe that the multiplication  $\lambda_{\sigma}:(\C,\lambda^{2}_{\sigma}\|.\|^{2}) \to (\C,\|.\|^{2})$ is an isometry. Hence we must form the absolute value of the  determinant of   
 the composition
$$\C\stackrel{d_{0,\sigma}}{\to} \C\stackrel{\lambda_{\sigma}}{\to} \C$$
which is given by $|\det(d_{0,\sigma})|=\lambda_{\sigma}$. In view of \eqref{glerjglkjrelgreljrelkgjrelgregregregregreg} we get again \eqref{lklkjlkjfdlewjflewfewfewfewfewfewfewfewfewf}.

The complex $\cV$ is acyclic so that  $h^{\cH}$ is canonically determined. We have
\begin{eqnarray*}\sum_{i=0}^{n}(-1)^{n} \cycl(V^{i},h^{V^{i}})&=&\cycl(R,h^{R})-\cycl(R,\lambda^{2}h^{R})\\&\stackrel{Lemma \:\eqref{jul0420}}{=}& -a(\sum_{\sigma\in \Sigma} \ln(\lambda_{\sigma})\bar \delta_{1,\sigma})\\&\stackrel{\eqref{lklkjlkjfdlewjflewfewfewfewfewfewfewfewfewf}}{=}&a(\sum_{\sigma\in \Sigma}\ln\:\tau_{an}(\cV_{\sigma},h^{\cV_{\sigma}}  )\bar \delta_{1,\sigma})\\&\stackrel{\eqref{jhkfhkejfhekfhekfewfewfewfewfewf}}{=}&a(\ln\:\tau (\cV,h^{\cV},h^{\cH}  ))
\ .\end{eqnarray*}
A more general case of this relation will be shown below in Proposition \ref{may121220}.

We have   canonical isomorphisms $\det(\cV_{\sigma}):=\C\otimes \C^{-1}\cong \C$ and
$\det(H(\cV_{\sigma}))\cong \C$. Under these isomorphisms $\phi_{\sigma}(\cV)$ is the identity.
Furthermore, the first isomorphism sends
$h^{\det(\cV_{\sigma})}$ to $\lambda_{\sigma}^{-2}\|.\|^{2}$.
It follows that $$h_{RM}^{\det(H(\cV_{\sigma}))}=\frac{1}{\lambda_{\sigma}^{2}}\|.\|^{2}\ , \quad h^{\det (H(\cV))_{\sigma}}=\|.\|^{2}\ .$$
By comparison with  \eqref{may121210}  we get
$$\tau(\cV_{\sigma},h^{\cV_{\sigma}},h^{\cH_{\sigma}})=\frac{1}{\lambda_{\sigma}}\ .$$
This is again compatible with  \eqref{lklkjlkjfdlewjflewfewfewfewfewfewfewfewfewf} and \eqref{jhkfhkejfhekfhekfewfewfewfewfewf}.

}\hB 
\end{ex}
}

 \begin{prop}\label{may121220}
In $\widehat{KR}^{0}(*)$ we have
$$ \sum_{i=0}^{n}(-1)^{n} \cycl(V^{i},h^{V^{i}})=\sum_{i=0}^{n}(-1)^{n} \cycl(H^{i}(\cV),h^{H^{i}(\cV)}) +a( \ln \tau(\cV,h^{\cV},h^{\cH}))\ .$$
\end{prop}
\proof
}We define
\begin{eqnarray*}\delta&:=&\sum_{i=0}^{n}(-1)^{n} \cycl(V^{i},h^{V^{i}})-\sum_{i=0}^{n}(-1)^{n} \cycl(H^{i}(\cV),h^{ H^{i}(\cV)})\\&& -a(   \ln \tau (\cV,h^{\cV},h^{\cH})  ) \ .\end{eqnarray*}
We first show that $\delta$ is independent of the geometries $h^{\cV}$ and $h^{\cH}$. To this end we consider  families of such geometries $\tilde h^{\cV}$ and $\tilde h^{\cH}$ parametrized by $\R$ and let $\tilde \delta \in \widehat{KR}^{0}(\R)$ be defined by \begin{eqnarray*}\tilde \delta&:=&\sum_{i=0}^{n}(-1)^{n} \cycl( \pr^{*}V_{i},\tilde h^{V_{i}})-\sum_{i=0}^{n}(-1)^{n} \cycl(\pr^{*} H^{i}(\cV),\tilde h^{ H^{i}(\cV)})\\&&  - a(  \ln \tau(\cV,\tilde h^{\cV},\tilde h^{\cH}) )\ ,\end{eqnarray*} where $\pr:\R \to *$ is the projection.
We now observe  that $R(\tilde \delta)=0$. First of all note that the torsion part of the cohomology does not contribute since classes of the form $\hat \cZ(T)$ are flat (see the proof of Lemma \ref{may121201}).  Let us now write
$$R\left(\sum_{i=0}^{n}(-1)^{n} \cycl(\pr^{*}V^{i},\tilde h^{V_{i}})-\sum_{i=0}^{n}(-1)^{n} \cycl(\pr^{*}\Proj(H^{i}(\cV)),\tilde h^{H^{i}(\cV)})\right)=a\:dt $$
for some $a\in C^{\infty}(\R,A^{-1}) $ (note that there is no contribution of $A^{0}$).
Using Lemma \ref{jun1820} we reduce the modules $V^{i}$ and $\Proj(H^{i}(\cV))$ to their determinants and observe, {using a similar calculation as in the proof of Lemma \ref{jul0420}}, that $$a \:  dt=d\ln \tau(\cV,\tilde h^{\cV},\tilde h^{\cH})\ .$$ {This implies $R(\tilde \delta)=0$}.

\bigskip

Since $R(\tilde \delta)= 0$
it now follows from a homotopy argument that $\delta$ does not depend on the geometry.

\bigskip

In order to show that $\delta=0$ we can  choose a very special geometry.
We essentially repeat the proof of  Lemma \ref{may121101} and  argue by induction on $n$ starting at the trivial case $n=0$.   We have the exact sequences
\begin{equation}\label{dhwjkjkqwdhwhkddkwqdqwdq} \cC:0\to  H^{0}(\cV)\to V^{0}\stackrel{\pi}{\to} V^{0}/Z^{0}(\cV)\to 0\ ,\end{equation}
\begin{equation}\label{dhwjkjkqwdhwhkddkwqdqwdq1}0\to V^{0}/Z^{0}(\cV)\to Z^{1}(\cV)\to H^{1}(\cV)\to 0\ , \end{equation}
 and
\begin{equation}\label{dhwjkjkqwdhwhkddkwqdqwdq2}0\to Z^{1}(\cV) \to V^{1} \to V^{1}/Z^{1}(\cV)\to 0\ .\end{equation}
Assume that we have chosen a geometry $h^{V}$. We then take $h^{\cH}:=h^{\cH}_{Hodge}$, {see Remark \ref{ewlfkjwelfjewlfewfuo3234}.}  On $V^{0}_{\sigma}/Z^{0}(\cV_{\sigma})$ we consider the metric induced by the inclusion of this quotient into $V^{1}_{\sigma}$. We further choose the metric on $Z^{1}(\cV_{\sigma})$   induced by the   inclusion into $V^{1}_{\sigma}$. Finally we choose the metric on   $V^{1}_{\sigma}/Z^{1}(\cV_{\sigma})$ to be  the quotient metric.  
 With these choices the only map in the three sequences which is not a partial isometry is $\pi$.
 We furthermore have defined a geometry on the complex $$\cV^{\prime}:0\to \cV^{1}/Z^{1}(\cV)\to \dots \to \cV^{n}\to 0$$ of length $n-1$ and starting in degree $1$,  and on its cohomology by the condition that it  is isometric to the cohomology of $\cV$ in degrees $\ge 2$. 
In order to proceed with the induction we must show that the difference
 \begin{eqnarray}
 \delta(\cV)-\delta(\cV^{\prime})
&=&\cycl(V^{0},h^{V^{0}})-\cycl(V^{1},h^{V^{1}})+\cycl(V^{1}/Z^{1}(\cV),h^{V^{1}/Z^{1}(\cV)}) \nonumber
\\&&-\cycl(H^{0}(\cV),h^{H^{0}(\cV)})+\cycl(H^{1}(\cV),h^{H^{1}(\cV)})\label{flejklqfjwkfjqlkdjlqqwdwqdqwdwd}\\
&&-a( \ln\: \tau(\cV,h^{\cV},h^{\cH}))+a( \ln\: \tau(\cV^{\prime},h^{\cV^{\prime}},h^{\cH^{\prime}}))
\nonumber\end{eqnarray} 
vanishes.
Since we equip the cohomology with the metrics induced by Hodge theory we have the relations  {(see Remark \ref{ewlfkjwelfjewlfewfuo3234})} $$\tau(\cV,h^{\cV},h^{\cH})=\tau_{an}(\cV,h^{\cV})\ , \quad  \tau(\cV^{\prime},h^{\cV^{\prime}},h^{\cH^{\prime}})=\tau_{an}(\cV^{\prime},h^{\cV^{\prime}})\ .$$
 The analytic torsions of the sequences \eqref{dhwjkjkqwdhwhkddkwqdqwdq1} and \eqref{dhwjkjkqwdhwhkddkwqdqwdq2} are trivial {since their differentials are partial isometries and we have the formula \eqref{glerjglkjrelgreljrelkgjrelgregregregregreg}.}
 We can therefore apply Lott's relation  \eqref{jul2023} twice
 in order to reduce \eqref{flejklqfjwkfjqlkdjlqqwdwqdqwdwd} to
\begin{eqnarray*}
 \delta(\cV)-\delta(\cV^{\prime})
&=&\cycl(V^{0},h^{V^{0}})-\cycl(H^{0}(\cV),h^{H^{0}(\cV)})-\cycl(V^{0}/Z^{0}(\cV),h^{V^{0}/Z^{0}(\cV)})\\&&-a( \ln\: \tau_{an}(\cV,h^{\cV}))+a( \ln\: \tau_{an}(\cV^{\prime},h^{\cV^{\prime}}))
\end{eqnarray*}
 Again using \eqref{glerjglkjrelgreljrelkgjrelgregregregregreg}, the analytic torsion
of the complex \eqref{dhwjkjkqwdhwhkddkwqdqwdq}    is given by  $$\tau_{an}(\cC_{\sigma},h^{\cC_{\sigma}})=|\det(\pi_{\sigma})|  \ .$$
Let us define
$$
\ln \tau_{an}(\cC,h^{\cC}) :=\sum_{\sigma} \ln |\det \pi_{\sigma}|\:\bar \delta_{1,\sigma}\ .$$
 Lott's relation  \eqref{jul2023} for $\cC$ thus gives
 $$\cycl(V^{0},h^{V^{0}})-\cycl(H^{0}(\cV),h^{H^{0}(\cV)})-\cycl(V^{0}/Z^{0}(\cV),h^{V^{0}/Z^{0}(\cV)})=-a( \ln \tau_{an}(\cC,h^{\cC}))\ .$$

 {Here we use the fact that the degree-zero component of the torsion form $\cT_{\sigma}$   of $(\cC_{\sigma},h^{\cC_{\sigma}})$ defined in \cite[A.2]{MR1724894} 
 is by \cite[Thm. 2.25]{MR1724894} equal to
$ \ln \tau_{an}(\cC_{\sigma},h^{\cC_{\sigma}})$. We get 
 $$\delta(\cV)-\delta(\cV^{\prime})=-a( \ln \tau_{an}(C,h^{\cC}))-a( \ln\: \tau_{an}(\cV,h^{\cV}))+a( \ln\: \tau_{an}(\cV^{\prime},h^{\cV^{\prime}}))\ .$$
Using \eqref{glerjglkjrelgreljrelkgjrelgregregregregreg} we express the analytic torsions of the complexes $\cV$ and $\cV^{\prime}$  in terms of the absolute values of the determinants of their differentials.  {Since the differentials of $\cV$  match with those of $\cV^{\prime}$ for $i\ge 1$ and $\pi$ is the zeroth differential of $\cV$}
we see that
$$-\ln \:\tau_{an}(\cC,h^{\cC})- \ln\:\tau_{an}(\cV,h^{\cV})+ \ln\:\tau_{an}(\cV^{\prime},h^{\cV^{\prime}})=0\ .$$
This implies 
$\delta(\cV)-\delta(\cV^{\prime})=0$\ . \hB}

\section{Transfers in differential cohomology}\label{sep2201}
 
\subsection{Introduction}

We consider a proper submersion $\pi:W\to B$ between smooth manifolds.
{Equivalently, $\pi$ is a locally trivial fibre bundle of manifolds with closed fibres.} The Becker-Gottlieb transfer associated to $\pi$ is a   map of spectra in topological spaces (well-defined up to contractible choice) 
\begin{equation}\label{sep0904}\tr :\Sigma^{\infty}_{+} B_{top}\to \Sigma^{\infty}_{+} W_{top}\end{equation}
first constructed in \cite{MR0377873} (this construction will be recalled below, see (\ref{may279})).
It induces for any spectrum $E$ the cohomological  Becker-Gottlieb transfer
\begin{equation}\label{jun182}\tr ^{*}:E^{*}(W)\to E^{*}(B)\ .\end{equation}
The cohomological Becker-Gottlieb transfer  behaves naturally with respect to cartesian diagrams of the form (\ref{jun105}), {iterated fibre bundles $W\to B\to Z$,} and spectrum maps $E\to E^{\prime}$.
{Moreover, it can be characterized in an axiomatic way, see \cite{MR1621939}.}

The goal of the present section is to refine the cohomological  Becker-Gottlieb transfer to a transfer in differential cohomology
$$\hat\tr: \widehat E^{*}(W)\to \widehat E^{*}(B)\ .$$
This differential transfer depends on the choice of a Riemannian structure on the submersion
$\pi$.
In Subsection \ref{sep0902} we introduce this notion and formulate our statements about {the} existence and functorial properties of 
the differential Becker-Gottlieb transfer. Proofs will be deferred to subsequent subsections.

\subsection{Differential Becker-Gottlieb transfer} \label{sep0902}

In this subsection we state the main theorem about the existence and the properties of the differential transfer. We further lay down the plan of the proof.

\bigskip

 We start with the definition of the additional geometric structures needed to refine the topological Becker-Gottlieb transfer to its differential version.

\begin{ddd}\label{may272}
A Riemanian structure $g:=(g^{T^{v}\pi},T^{h}\pi)$ on a submersion $\pi:W\to B$ consists of
\begin{enumerate}
\item a metric $g^{T^{v}\pi}$ on the vertical {tangent} bundle $T^{v}\pi:=\ker(d\pi)\subseteq TW$,
\item a horizontal {subbundle} $T^{h}\pi\subseteq TW$ (i.e. {vector} subbundle  such that $T^{v}\pi\oplus T^{h}\pi=TW$).
\end{enumerate}
\end{ddd}
 
The following Lemma is a well-known fact about the Riemannian geometry of bundles.
\begin{lem}\label{kldlqwdwqdwqd}
A Riemannian structure $g$ on a submersion $\pi$ induces a connection $\nabla^{T^{v}\pi}$ on the vertical {tangent} bundle $T^{v}\pi$.
\end{lem} \proof 
  Note that \begin{equation}\label{fwefkfjkjewfjklkjlkwefewf}
d\pi_{|T^{h}\pi} :T^{h}\pi\to \pi^{*} TB
\end{equation} is an isomorphism of vector  bundles.
We choose a Riemannian metric $g^{TB}$ on the base $B$. It induces a metric $h^{T^{h}\pi}$ on the horizontal {sub}bundle $T^{h}\pi$ such that \eqref{fwefkfjkjewfjklkjlkwefewf}  becomes an isometry.
The orthogonal sum
$g^{TW}:=g^{T^{v}\pi}\oplus g^{T^{h}\pi}$ is a Riemannian metric on the total space $W$ of the submersion and induces a Levi-Civita connection $\nabla^{LC,TW}$ on the tangent bundle $TW$. Its projection
$\nabla^{T^{v}\pi}$ to the vertical {tangent} bundle turns out to be  independent of the choice of the metric $g^{TB}$ (see \cite[Ch. 9]{MR1215720}).  The connection $\nabla^{T^{v}\pi}$ is the desired connection on the vertical {tangent} bundle $T^{v}\pi$ induced by the Riemannian structure $g$. \hB

For a submersion $\pi:W\to B$ we let $\Lambda\to W$ denote the orientation bundle associated to the real vector bundle $T^{v}\pi\to W$, see (\ref{jul1601}). The Chern-Weil representative of  the Euler class of $T^{v}\pi$ associated to the connection $\nabla^{T^{v}\pi}$ is the  Euler form  (see (\ref{jul1603}) for a definition) $$e(\nabla^{T^{v}\pi})\in Z^{n}(\Omega(W, \Lambda))\ , $$ where  $n:=\dim(W)-\dim(B)$  is the fibre dimension. In view of Lemma \ref{kldlqwdwqdwqd} it depends on the Riemannian {structure} so that  it is natural to  use the notation  \begin{equation}\label{feb2401}e(g):=e(\nabla^{T^{v}\pi})\ .\end{equation}  

We now assume that the submersion $\pi:W\to B$ is proper {with closed fibres}. 
{Equivalently, $W\to B$ is a locally trivial fibre bundle with closed fibres.} Then we have an integration map of chain complexes
$$\int_{W/B}:\Omega  (W,\Lambda)[n]\to \Omega(B)\ .$$
{
\begin{rem}\label{efjwekffjelfeofuo343442}{\rm Since the domain consists of forms with twisted coefficients in the orientation bundle this integration does not require to choose a fibrewise orientation of the bundle.
}\hB
\end{rem}
}
Since the Euler form $e(g)$ is closed we obtain for any chain complex $A$ an induced morphism of chain complexes
 $$\int_{W/B}\cdots\wedge e(g):\Omega A(W)\to \Omega A(B)$$
 which preserves degree.
For {differential data} $(E,A,c)$ consisting of a  spectrum $E\in \Sp$, a chain  complex $A\in \Ch$ of real vector spaces, and an equivalence $c:E\R\stackrel{\simeq}{\to} H(A)$ let $\widehat E^{*}$ denote the associated
differential extension of $E$ introduced in Definition \ref{difdef1}.
The main goal of the present section is to prove: \begin{theorem}\label{may3102}
For every proper submersion $W\to B$ with a Riemannian structure $g$, {there exists a} canonical differential transfer $\hat \tr :\widehat E^{0}(W)\to \widehat E^{0}(B)$
such that
\begin{equation}\label{may3103222}\xymatrix{\Omega A^{-1}(W)\ar[d]_{\int_{W/B}\cdots\wedge e(g)}\ar[r]^{a}&\widehat E^{0}(W)\ar[d]^{\hat \tr}\ar[r]^{I}\ar@/^{1cm}/[rr]^{R}& E^{0}(W)\ar[d]^{\tr^{*}}&Z^{0}(\Omega A(W))\ar[d]^{\int_{W/B}\cdots\wedge e(g)}\\\Omega A^{-1}(B)\ar[r]^{a}&\widehat E^{0}(B)\ar[r]^{I}\ar@/_{1cm}/[rr]^{R}& E^{0}(B)&Z^{0}(\Omega A(B))
}\end{equation}
commutes. Furthermore:
\begin{enumerate}
\item  The differential transfer behaves naturally {with respect to} pull-back in cartesian diagrams (\ref{jun105}).
\item The transfer is functorial for iterated bundles with Riemannian structures (see Proposition \ref{jul2015} for a precise formulation).
\item The differential transfer $\hat \tr$ preserves the flat subgroups and the  diagram \begin{equation}\label{jun0304}
\xymatrix{\widehat E_{flat}^{0}(W)\ar[d]^{\hat \tr}\ar[r]^{\cong}&E\R/\Z^{-1}(W)\ar[d]^{\tr^{*}}\\
\widehat E^{0}_{flat}(B)\ar[r]^{\cong}&E\R/\Z^{-1}(B) \ .}
\end{equation}
  commutes.
\item If  $\phi:(E,A,c)\to (E^{\prime},A^{\prime},c^{\prime})$   is a morphism between differential data, 
then 
\begin{equation}\label{aug1012}\xymatrix{\widehat E^{0}(W)\ar[r]^{\phi}\ar[d]^{\hat \tr}&\widehat E^{\prime, 0}(W)\ar[d]^{\hat \tr}\\
\widehat E^{0}(B)\ar[r]^{\phi}&\widehat E^{\prime,0}(B) \ .}\end{equation}
 commutes. \end{enumerate}\end{theorem}
 
\begin{rem}{\rm 
The differential Becker-Gottlieb transfer $\hat \tr$ depends non-trivially on the choice of the Riemannian structure $g$. If we want to  emphasize the dependence on this choice, then we will use the notation $\hat \tr(g)$.}\hB
\end{rem}

 \begin{rem}{\rm
By proving Theorem \ref{may3102} we construct a particular choice of a differential transfer map. We do not know whether a differential transfer map is characterized uniquely by naturality and the diagram
(\ref{may3103222}).\hB
}\end{rem}

{
The differential transfer will be induced by a map of differential function spectra $$\Diff(E)(W)\to \Diff(E)(B)\ .$$ In view of the {Definition \ref{may271}} of the differential function spectrum as a pull-back, in order to obtain   such a map}
we must construct  a commutative diagram
\begin{equation}\label{may278}\xymatrix{\Funk(E)(W)\ar[d]^{\tr}\ar[r]^{\rat}&H(\Omega A(W))\ar[d]^{H(\int_{W/B}\cdots\wedge e(g))}&H(\sigma\Omega A(W))\ar[l]\ar[d]^{H(\int_{W/B}\cdots\wedge e(g))}\\\Funk(E)(B)\ar[r]^{\rat}&H(\Omega A(B))&H(\sigma \Omega A(B))\ar[l]}\end{equation}
in the $\infty$-category of smooth spectra $ \Sp$. For Theorem \ref{may3102}, 1. this construction must be  natural in the bundle $\pi:W\to B$} and its Riemannian structure. For \eqref{aug1012} the construction should be natural in the data $(E,A,c)$, too.
 The main effort lies in the left square.

\begin{rem}{\rm 
In order to formulate the naturality of the diagram (\ref{may278}) properly  we consider it as a diagram 
in the $\infty$-category of $\Bundle_{geom,trans}$-spectra,
where $\Bundle_{geom,trans}$ is the category of triples $(\pi,g,b)$ consisting of a proper submersion $\pi$ equipped with Riemannian and  transfer structures $g$ and $b$ (precise definitions will be given below). 
In Subsection \ref{jun143} we construct the right square. In Subsection \ref{jun144} we construct the left vertical arrow. Finally, in Subsection \ref{jun145} we provide the left square.  
The end of the proof of Theorem \ref{may3102}  now goes as follows.
Given $\pi:W\to B$ with Riemannian structure $g$ we choose a transfer datum $b$.
If we evaluate (\ref{may278}) at $(\pi,g,b)$, then we get an induced map of differential function spectra
$\Diff(E)(W)\to \Diff(E)(B)$ which in turn induces the differential transfer
$$\hat \tr(g,b):\widehat E^{0}(W)\to \widehat E^{0}(B)$$
on the level of zeroth   homotopy groups.
We must finally show:
\begin{lem}\label{sep2104}
The {differential} transfer $\hat \tr(b,g)$ does not depend on the choice of the transfer datum $b$.
\end{lem} 
 The proof of this Lemma will be given at the end of Subsection \ref{jun145}. } \hB \end{rem}
%
%
%

\subsection{Geometric bundles and integration of forms}\label{jun143}
 
In this subsection we make the right square of (\ref{may278}) precise. On the way we introduce some notation used in order to formulate the functoriality of the construction in the bundle $\pi:W\to B$ properly.

\bigskip

We consider the category $\Bundle$ whose objects are proper submersions 
$\pi:W\to B$ between manifolds {with closed fibres  (also called bundles)}, and whose morphisms $$(f,F):(\pi^{\prime}:W^{\prime}\to B^{\prime})\to (\pi:W\to B)$$ are cartesian squares  \begin{equation}\label{jun105}\xymatrix{W^{\prime}\ar[d]^{\pi^{\prime}}\ar[r]^{F}&W\ar[d]^{\pi}\\B^{\prime}\ar[r]^{f}&B}\ .\end{equation} Adopting a similar convention as for smooth objects (see Subsection \ref{smoothobjects}), for an $\infty$-category $\bC$ we write
 $$\Bundle(\bC):=\Fun(\Nerve(\Bundle^{op}),\bC)$$ for the $\infty$-category of bundle objects in $\bC$.
Similar conventions apply to $\Bundle_{geom}$-objects below.

\bigskip

We  introduce the bundle set
$$\Riem\in \Bundle(\Nerve(\Set))$$
which associates to each bundle $\pi:W\to B$ the set of Riemannian structures
$\Riem(\pi)$ on $\pi$ as introduced in Definition \ref{may272}. 
\begin{rem}{\rm Here are the details for compleneteness of the presentation. Given a morphism $(f,F):\pi^{\prime}\to \pi$ as in diagram (\ref{jun105}) and and a Riemannian structure  $$g=(g^{T^{v}\pi},T^{h}\pi)\in \Riem(\pi)$$ on $\pi$ we define the induced Riemannian structure
$$(f,F)^{*}g=(g^{T^{v}\pi^{\prime}},T^{h}\pi^{\prime})\in \Riem(\pi^{\prime})$$ on $\pi^{\prime}$ as follows.
 Since the square \eqref{jun105} is cartesian  the differential $dF:TW^{\prime}\to F^{*}TW $ of $F$ restricts to an isomorphism
$dF_{|T^{v}\pi^{\prime}}:T^{v}\pi^{\prime}\stackrel{\cong}{\to} F^{*}T^{v}\pi$. We define the induced vertical metric as a pull-back
$$g^{T^{v}\pi^{\prime}}:=(dF_{|T^{v}\pi^{\prime}})^{*}g^{F^{*}T^{v}\pi}\ .$$ Furthermore, we define the induced horizontal {subbundle} by  
$$T^{h}\pi^{\prime}:=dF^{-1}( F^{*}T^{h}\pi)\ .$$ } \hB\end{rem}

We let
\begin{equation}\label{may291}q:\Bundle_{geom}\to \Bundle\end{equation} denote the Grothendieck construction of 
$\Riem$. {\begin{rem}{\rm The objects of the category $\Bundle_{geom}$ are pairs
$(\pi:W\to B,g)$ consisting of a bundle together with a Riemannian structure $g\in \Riem(\pi)$. Morphisms $(\pi^{\prime},g^{\prime})\to (\pi,g)$ are morphisms $(f,F):\pi^{\prime}\to \pi$ in $\Bundle$ such that
$(f,F)^{*}g=g^{\prime}$.}\hB \end{rem}}
{We refer to the objects of $\Bundle_{geom} $ as bundles with geometry.}
The functor  $q$ forgets the Riemannian structure.

\bigskip

We have two functors $$\dom,\ran:\Bundle\to \Mf$$ given by the evaluation at the domain and range
$$\dom(\pi:W\to B):=W\ ,\quad \ran(\pi:W\to B):=B\  .$$
They induce corresponding pull-back functors between the categories of smooth objects, $\Bundle_{geom}$-objects, and their $\infty$-{categorical} analogs. 
Given a smooth object  $X$, we define (in order to shorten the notation) the 
$\Bundle_{geom}$-objects  
\begin{equation}\label{may281}Xd:=X\circ \dom\circ q\ ,\quad Xr:=X\circ \ran\circ q\ .\end{equation}
\begin{rem}{\rm  The usage of this notation is as follows. Assume that $X$ is a smooth set and we have a natural way to  define a push-forward $X(W)\to X(B)$ for a bundle $\pi:W\to B$ which depends on the choice of a Riemannian structure on $\pi$. Using the notation above we can formalize this by saying that we define
a map of $\Bundle_{geom}$-sets $Xd\to Xr$. } \hB \end{rem}

{As a first instance of this reasoning, we are going to} define morphisms of $\Bundle_{geom}$-chain complexes (i.e. $\Bundle_{geom}$-objects in $\Ch$) \begin{equation}\label{hdjwdkqjwdhkjhkhdkjhkqwdoiuoiqdqwdqwdqdd}
\int:\Omega A d\to \Omega A r\ ,\quad \int_{\sigma}:\sigma \Omega A d\to \sigma \Omega A r
\end{equation}
as follows. If $(\pi:W\to B,g)\in \Bundle_{geom}$, then {we define the morphism of chain complexes by }
\begin{equation}\label{fewlfjewlkjlkewfewf45454353454354}
\int:\Omega A (W)\to \Omega A(B)
\end{equation}
$$\Omega A(W)\ni \alpha\mapsto \int_{W/B} \alpha\wedge e(g)\in \Omega A(B)\ .$$

{
\begin{rem}{\rm 
It is clear how to define the integral for $\alpha\in \Omega(W)\otimes_{\R} A$ (see also Remark \ref{efjwekffjelfeofuo343442}). 
 If $B$ is not compact, then
$\Omega (W)\otimes_{\R} A$ is a proper subset of $\Omega A(W)$, see Remark \ref{jkdhwqjkdhkwqdwqhdqwdkjhwqdqwdwqdwqdwd}. In order to extend the definition
of the integral {to} all of $\Omega A(W)$ we use that $B$ is locally compact and work locally in $B$. We refer to Subsubsection \ref{kwjfweklfjewlfjelwfewopipi3ir23r20pr32r32r32r} for more details.
}\hB
\end{rem}}

{{The integral \ref{fewlfjewlkjlkewfewf45454353454354}  depends} on the Riemannian structure on $\pi$ via the Euler form
 $e(g)\in Z^{n}(\Omega(W,\Lambda))$ defined in \eqref{feb2401}.}
The transformation $\int_{\sigma}$ is obtained from $\int$ by restriction.  We get a  commutative diagram in $\Bundle_{geom}(\Nerve(\Ch))$:  \begin{equation}\label{sep0903}\xymatrix{\Omega A d\ar[d]^{\int}&\sigma\Omega A d\ar[l]\ar[d]^{\int_{\sigma}}\\\Omega A r&\sigma\Omega A r\ar[l]}\ .\end{equation}
Applying the Eilenberg-MacLane functor objectwise 
 we obtain the commutative diagram 
 \begin{equation}\label{feb2701}\xymatrix{H(\Omega Ad)\ar[d]^{H(\int)}&H(\sigma\Omega Ad)\ar[l]\ar[d]^{H(\int_{\sigma})}\\H(\Omega A r)&H(\sigma\Omega A r)\ar[l]}\end{equation} in the $\infty$-category  
$\Bundle_{geom} (\Sp )$.
This is the right square of (\ref{may278}).

\bigskip

An inspection of the construction of the diagram \eqref{feb2701} shows that it is functorial in the chain complex $A$. Let $\Box$ denote the category of the shape \eqref{ewfwfewfwj2lkr23iori32rio2ri322r2r2}.
\begin{kor}\label{lfjwlekfjeluo32ieifhqkf}
{There exists} a map
$$\Nerve(\Ch_{\R})\to {\Bundle_{geom}}(\Sp)^{{\Nerve}(\Box)}$$
which associates to a chain complex $A$ the commuting diagram   \eqref{feb2701}.
\end{kor}

\subsection{Transfer structures and the Becker-Gottlieb transfer}\label{jun144}

In this subsection we discuss the details of the construction of the Becker-Gottlieb transfer \cite{MR0377873}.
Its main goal is the construction of the left vertical arrow in (\ref{may278}) stated here as Lemma \ref{kjdklqwdqwdqdwqdwqdqwdwqdqdqdqdq}.

\bigskip

We first observe that the map (\ref{sep0904}) depends on choices. In order to formalize {this} we {will} subsume these choices under the notion of transfer {data}.
\begin{ddd}\label{jul2001}
Transfer data on a proper submersion {with closed fibres} $\pi:W\to B$ consists of 
\begin{enumerate}
\item an integer $k\ge 0$ (the dimension of the transfer data),
\item a fibrewise embedding
$$\xymatrix{W \ar[rr] \ar[dr]&& B \times \R^{k}\ar[dl]\\&B &}\ ,$$
\item an extension of the embedding to an open embedding $\emb:N\to B\times \R^{k}$ of the fibrewise normal bundle
$N\to W$ of this inclusion. 
\end{enumerate}
\end{ddd}
{We will simply write $b$ for a given choice of transfer data on $\pi:W\to B$.}
There is an obvious construction of a stabilization $b\mapsto S^{l}(b)$ which increases the dimension of $b$ from $k$ to $k+l$. Stabilization introduces an equivalence relation amongst transfer data.
The equivalence classes are referred to as stable transfer data.

\bigskip

If $V\to X$ is a vector bundle on a {topological} space $X$, then $X^{V}$ denotes its Thom space, a pointed topological space.

\begin{rem}\label{djqwdkjqwdkdjqwldjqwdqwdqwdwqd}{\rm For completeness we recall the details of the definition of the Thom space. In order to construct it we first consider the fibrewise compactification $V\subset \bar V$ obtained by attaching the  bundle $S_{\infty}(V)\to V$ of   spheres at $\infty$. The Thom space is then defined as a quotient $X^{V}:=\bar V/S_{\infty}(V)$.  
We have a natural inclusion \begin{equation}\label{ejg12e1g2hje21e2iu1e2i1oeu21oie21e21e21e}
V\to X^{V}\ .
\end{equation}
 A trivialization $V\cong X\times \R^{k}$ of vector bundles induces {a}  homeomorphism  $X^{V}\cong \Sigma^{k}_{+}X$. 
 Here $X_{+}$ denotes the space obtained from $X$ by attaching a disjoint base point and we write $\Sigma^{k}_{+}X:=\Sigma^{k}X_{+}$.
 Furthermore, an inclusion
$V\to W$ of bundles  induces a map between Thom spaces $X^{V}\to X^{W}$.}\hB \end{rem}

\begin{rem}{\rm 
In the following we work with the underlying topological spaces of  manifolds. In order to simplify the notation we will frequently omit the subscript $(-)_{top}$ and simply write $B$ for $B_{top}$.
For example, if $V\to B$ is a real vector bundle we let $B^{V}\in \Top_{*}$ denote the Thom space.
This abbreviates the longer symbol $(B_{top})^{V_{top}}$. Furthermore, the symbol $\Sigma^{k}_{+}B$ denotes
$\Sigma^{k}_{+}B_{top}\in \Top_{*}$.
}\hB\end{rem}
A transfer datum $b$ {of dimension $k$} gives rise to a map between pointed topological spaces \begin{equation}\label{wqdqjwdqwjdwqjkdwjqkdkwqqwdwqdwqdw}
 \Sigma^{k}_{+}B \cong B^{B\times \R^{k}}\xrightarrow{\clps} W^{N}\xrightarrow{z} W^{T^{v}\pi\oplus N}\cong W^{W\times \R^{k}}\cong  \Sigma^{k}_{+} W \ .
\end{equation}
Here  $$\clps:B^{B\times \R^{k}}\to W^{N}$$ is the collapse map.
In detail, 
the collapse map is given as follows: \begin{enumerate} \item  On the image of the embedding $\emb:N\to B\times \R^{k}$ 
the collapse map is defined as $\emb^{-1}$. \item   On the complement of  the image of   $\emb$ it is the constant map with value the base point of the Thom space 
 $W^{N}$. \end{enumerate} 
The map $z:W^{N}\to W^{T^{v}\pi\oplus N}$ is induced by the  zero section of the {vertical} tangent bundle $T^{v}\pi$.  Finally, the homeomorphism $W^{T^{v}\pi\oplus N}\cong W^{W\times \R^{k}}$ is induced by the canonical trivialization of vector bundes  $T^{v}\pi\oplus N\cong W\times \R^{k}$.

\bigskip

For the geometric construction of the Becker-Gottlieb transfer we fix the usual point-set model of the category of spectra $\check{\Sp}$ {in topological spaces}.

\begin{rem}\label{dqwjdklqwdwqwqdqdwq}{\rm Here are some details about the one-category $\check{\Sp}$ and how it is related with $\Sp$. {It is the analog of the version of spectra in simplicial sets considered in  \cite{MR513569}.}  An object  $$\check{E}=((E_{n}),(\sigma_{n}))\in \check{\Sp}$$ is a collection of pointed {topological} spaces and {structure} maps  indexed by $\nat$, where $\sigma_{n}:\Sigma E_{n} \to E_{n+1}$. Morphisms $f:\check{E}\to \check{E}^{\prime}$ in $\check{\Sp}$ are families of maps $f=(f_{n}:E_{n}\to E_{n}^{\prime})$ which are compatible with the structure maps. 

\bigskip

 If $X$ is a based topological space, then we let $\check{\Sigma}^{\infty} X\in \check{\Sp}$ be the 
suspension spectrum  of $X$ whose $n$'th space is given by $\Sigma^{n} X$. For an unbased topological space $X$ we
 set $$\check{\Sigma}^{\infty}_{+}X:=\check{\Sigma}^{\infty} X_{+}\ .$$
The suspension spectrum provides a functor \begin{equation}\label{fwefefewfewfewfewfewf234424324324}
\check{\Sigma}^{\infty}_{+}:\Top \to \check{\Sp}\ .
\end{equation}
We further use the notation
$$\check{\Sigma}^{\infty+k}_{+}X:=\check{\Sigma}^{\infty} \Sigma^{k}_{+}X\ .$$

\bigskip 

For $n\in \Z$ the $n$'th stable homotopy group of a spectrum $\check{E}$ is defined by 
$$\pi_{n}(\check{E}):=\colim_{k\ge -n} \pi_{k}(E_{n+k})\ .$$
A morphism $f:\check{E}\to \check{E}^{\prime}$ in $\check{\Sp}$ is a weak equivalence if the induced map between stable homotopy groups
$$\pi_{n}(\check{E})\to \pi_{n}(\check{E}^{\prime})$$
is an isomorphism for all $n\in \Z$. 
The suspension spectrum functor $\Sigma^{\infty}:\Top_{*}\to \check{\Sp}$ preserves
weak equivalences and therefore
descends to a functor \begin{equation}\label{kjklfjklewjflewfewfwef324}
\Spc_{*}\simeq \Nerve(\Top_{*})[W^{-1}]\to \Nerve(\check{\Sp})[W^{-1}]\ .\end{equation}
The category $\Nerve(\check{\Sp})[W^{-1}]$ is pointed and admits all small limits and colimits.
One now observes that the {suspension} endofunctor $${\Sigma}:  \Nerve(\check{\Sp})[W^{-1}]\to \Nerve(\check{\Sp})[W^{-1}]$$ (see \eqref{wjqdlkqwdjwlkdjwkldwqdwqdwqdqwdq} for a definition) is an auto-equivalence. Therefore $\Nerve(\check{\Sp})[W^{-1}]$ is  a stable $\infty$-category (Definition \ref{fjewfkljewfejelwfewfewfwf}). By the universal property of the $\infty$-category of spectra  $\Sp$ defined as the stabilization of $\Spc_{*}$ (Definition \ref{lekjflkewjlwfewf23i4p23i423o}){,} the suspension spectrum functor \eqref{kjklfjklewjflewfewfwef324} extends to a functor
\begin{equation}\label{kjkljdlkwjldwqou9872498273943241232r}\Sp=\Spc_{*}[{\Sigma}^{-1}] \to   \Nerve(\check{\Sp})[W^{-1}]\ .\end{equation}
This functor is actually an equivalence  and this fact justifies {our consideration of} $\check{\Sp}$ as a point-set model for $\Sp$. The shift $$\Sigma: \check{\Sp}\to \check{\Sp}\, \qquad    ((E_{n}),(\sigma_{n}))\mapsto ((E_{n+1}),(\sigma_{n+1}))$$ 
  in $\check{\Sp}$ is not an equivalence. But it preserves weak equivalences and descends to an autoequivalence of $\Nerve(\check{\Sp})[W^{-1}]$ which corresponds to the suspension   in $\Sp$  under the equivalence \eqref{kjkljdlkwjldwqou9872498273943241232r}.

For later use we introduce the functor
 between    $\infty$-categories  \begin{equation}\label{3224343hfweijfiowejfoewfjo}
\iota:\Nerve(\check{\Sp})\to \Nerve(\check{\Sp})[W^{-1}]\simeq \Sp
\end{equation}

\bigskip

The composition
$$\Sigma^{\infty}_{+}:\Spc\stackrel{\eqref{hdqhjwdqdwqddqdqdd134123}}{\simeq}\Nerve(\Top)[W^{-1}] \stackrel{\eqref{fwefefewfewfewfewfewf234424324324}}{\to} \Nerve(\check{\Sp})[W^{-1}] \simeq \Sp\ .$$
is the equivalent to the $\infty$-categorical version \eqref{qwjdkjqwdljqwkjdqlk89739423} of the  suspension spectrum functor.  
 
  }\hB\end{rem}

We use the symbol 
 \begin{equation}\label{may279}\check{\tr} (b):\check{\Sigma}^{\infty+k}_{+} B \to \check{\Sigma}^{\infty+k}_{+} W \end{equation} in order  to denote the associated   map between spectra  in $\check{\Sp}$ induced by \eqref{wqdqjwdqwjdwqjkdwjqkdkwqqwdwqdwqdw} by applying the point-set suspension spectrum functor \eqref{fwefefewfewfewfewfewf234424324324}. This 
{is  a precise point-set version of the map (\ref{sep0904}), where we add the argument $b$ in order to indicate the dependence on the transfer data}. 
We further write
\begin{equation}\label{fweffew3244rfwf23r22r23232r23}
\tr (b):=\Sigma^{-k}\iota(\check{\tr}_{k}(b)):\Sigma^{\infty}_{+}B \to \Sigma^{\infty}_{+}W 
\end{equation}
for the induced morphism   in $\Sp$, where
$\iota$ is the map \eqref{3224343hfweijfiowejfoewfjo},  {and} we use that the suspension becomes invertible in $\Sp$.

\begin{rem}{\rm If $M$ is a manifold, then we abbreviate
$\Sigma_{+}^{\infty} M_{top,\infty}$ to $\Sigma_{+}^{\infty}M$.
Similarly, for a topological space $X$ we write
$\Sigma_{+}^{\infty}X$ for $\Sigma_{+}^{\infty}X_{\infty}$.
}\hB
\end{rem}

 \bigskip
  
  The construction of the stable transfer map is compatible with stabilization.
Indeed, for
 $l\ge 0$, we have the equality $$\check{\tr} (S^{l}(b))=\Sigma^{l}\check{\tr} (b)\ .$$
 We thus get a natural equivalence $$\tr(b) \simeq \tr(S^{l}(b))$$
 which only depends on the stable transfer datum represented by $b$ up to contractible choice.

\begin{rem}{\rm 
The category $\check{\Sp}$ is enriched in topological spaces. 
 For every $k$ there is actually a space of transfer data of dimension $k$, and the map
   $b\mapsto \check{\tr}(b)$ is continuous. If $k$ is sufficiently large, then the space of transfer data is contractible and 
 (\ref{may279})   is well-defined up to contractible   choice. In particular,  
 \eqref{fweffew3244rfwf23r22r23232r23} is essentially unique independently of the choice of the stable transfer data. 
  {For the proof of Lemma \ref{sep2104}} we will {only need} the fact that any two choices of stable transfer data can be connected by a path; compare \cite[Sec. 3]{MR0377873}.
}
\hB
\end{rem}

We now consider a morphism $\phi:\pi^{\prime}\to \pi$ (i.e. a diagram of the shape  (\ref{jun105})) in $\Bundle$ and a choice of transfer data $b$  on $\pi$. Then there is a natural and functorial (with respect to compositions of morphisms in $\Bundle$) way to define transfer data  $b^{\prime}:=\phi^{*}b$ on $\pi^{\prime}$.
It has the same dimension as $b$.  The second component of the  fibrewise embedding of $b^{\prime}$ is given by the composition $$W^\prime\to W\to B\times \R^{k}\stackrel{\pr}{\to} \R^{k}\ ,$$ where the second map is the fibrewise embedding for $b$.  We get a canonical identification of the fibrewise normal bundles $N^{\prime}\stackrel{\cong}{\to} F^{*}N$, and the second component of the {open} embedding {of  $N^{\prime}$} for $b^{\prime}$ is given by
$$N^{\prime}\stackrel{\cong}{\to}{F^{*}N} \stackrel{\emb}{\to} B\times \R^{k}\stackrel{\pr}{\to}\R^{k}\ .$$
The first component of these maps is the natural map to $B^{\prime}$ in both cases.
We can now consider the bundle set
$$\Trans\in \Bundle(\Nerve(\Set))$$
which associates to a bundle $\pi:W\to B$ the set of stable  transfer data, 
and we let $$q:\Bundle_{trans}\to \Bundle$$ denote the associated Grothendieck construction.
As before, for $E\in \Sm(\bC)$  we write
$$\Funk(E)d:=\Funk(E)\circ \dom \circ q\ ,\quad \Funk(E){r}:=\Funk(E)\circ \ran \circ q$$
for the corresponding objects
in $\Bundle_{trans}(\bC)$.
\begin{lem}\label{kjdklqwdqwdqdwqdwqdqwdwqdqdqdqdq}For a spectrum $E\in \Sp$ 
we have a map
\begin{equation}\label{feb2702}\tr:\Funk(E)d\to\Funk(E)r\end{equation}
in $\Bundle_{trans}( \Sp)$.
\end{lem}
\proof
We use that the $\infty$-category $\Sp$ of spectra is tensored and cotensored over itself {(see Example \ref{dlkqjwdlkjqwljqwdqwdqwdqw})}. This structure extends the tensor and cotensor structure over $\Spc$, see  \eqref{jchjsachakscsacscaca}. In particular{,} the following diagram commutes
$$\xymatrix{\Spc^{op}\times \Sp\ar[rr]^{(X,E)\mapsto E^{X}}\ar[d]_{(X,E)\mapsto (\Sigma_{+}^{\infty}X,E)}&&\Sp\ar@{=}[d]\\
\Sp^{op}\times \Sp\ar[rr]^{(F,E)\mapsto \map(F,E)}&&\Sp}\ .$$
For a space $X$ and a spectrum $E$ we thus have a natural equivalence
$E^{X}\simeq \map(\Sigma^{\infty}_{+}X,E)$ of spectra. 
Consequently, the functor
$\Funk:\Sp\to \Sm(\Sp)$ introduced in Definition \ref{kjdlkqwwqddwqdqdq}
is equivalent to the functor
  given by
$$ M\mapsto \map(\Sigma^{\infty}_{+}M ,E)\ .$$

Above we have constructed functors
$$D,\: R: \Nerve(\Bundle)\to \Sp\ , \quad D (W\to B):= \Sigma^{\infty}_{+}W \ , \quad  R (W\to B):=\Sigma^{\infty}_{+}B \ .$$
Then \begin{equation}\label{wefewfijjk2uio23roiu3r2rr}
\Funk(E)d\simeq\map(-,E)\circ D\circ q\ ,\quad  \Funk(E)r\simeq \map(-,E)\circ R\circ q\ .
\end{equation}
For a bundle with a choice of transfer data $(\pi:W\to B),b)\in \ \Bundle_{trans}$ we have a
map in $\Sp$ \begin{equation}\label{fefefwefef252523532525}
\tr(b):R(q(\pi))\to D(q(\pi))
\end{equation}
 given by \eqref{fweffew3244rfwf23r22r23232r23}.
For a morphism $\phi:\pi^{\prime}\to \pi$ (i.e. a diagram \eqref{jun105})  and $b^{\prime}:=\phi^{*}b$
 the following diagram in $\check{\Sp} $ strictly commutes
\begin{equation}\label{sep2101}
\xymatrix{\check{\Sigma}^{\infty+k}_{+} W^{\prime} \ar[r]&\check{\Sigma}^{\infty+k}W \\
\check{\Sigma}^{\infty+k}_{+}B^{\prime} \ar[r]\ar[u]^{\check{\tr}(b^{\prime})}&\check{\Sigma}^{\infty+k}_{+} B \ar[u]^{\check{\tr}(b)}\ .
}\end{equation}
Consequently, 
 
\begin{equation}\label{sep2101nnn}
\xymatrix{ {\Sigma}^{\infty}_{+} W^{\prime} \ar[r]& {\Sigma}^{\infty}W \\
 {\Sigma}^{\infty}_{+}B^{\prime} \ar[r]\ar[u]^{\tr(b^{\prime})}& {\Sigma}^{\infty}_{+} B \ar[u]^{\tr(b)}
}\end{equation}
commutes in $\Sp$ in an essentially unique way. The diagrams \eqref{sep2101nnn} now  refine the collection of maps
\eqref{fefefwefef252523532525} for all $\pi$ to a transformation $$ \tilde \tr: R\circ q\to D\circ q$$
in $\Fun(\Bundle_{trans},\Sp)$ which induces the desired 
transformation $\tr$ by precomposition using the equivalences \eqref{wefewfijjk2uio23roiu3r2rr}. 
\hB


%

An inspection of the construction of $\tr$ shows that it is natural in the spectrum $E$. \begin{kor}\label{jdlqwdqwlkdjwqldqwdwqdqwdqwdqd}
We have constructed a map
$$\Sp\to \Map(\Funk(-)d,\Funk(-)r)\ .$$
\end{kor}
In particular, if 
 $f:E\to F$ is a map of spectra, then  
  we have a commutative diagram of $\Bundle_{trans}$-spectra
\begin{equation}\label{may285}
\xymatrix{\Funk(E)d\ar[d]^{\tr}\ar[r]^{\Funk(f)d}&\Funk(F)d\ar[d]^{\tr}\\
\Funk(E)r\ar[r]^{\Funk(f)r}&\Funk(F)r}\ .\end{equation}

\color{black}

\subsection{The left square in (\ref{may278}) and the  construction of $\hat \tr$}\label{jun145}

We   let $q:\Bundle_{geom,trans}\to \Bundle$ be the Grothendieck construction on the {bundle} set
$$\Riem\times \Trans\in \Bundle(\Nerve(\Set))\ .$$
The constructions given in Subsections \ref{jun143} and \ref{jun144}
provide maps \eqref{feb2702} and diagrams  {of the shape  \eqref{feb2701}} in $$\Bundle_{geom,trans}( \Sp)$$ via the obvious pull-backs along the functors {which forget} either the Riemannian or the transfer structures.

\bigskip

Recall 
the  construction of the map  (\ref{jun147}) 
$$\rat:\Funk(E)\to H(\Omega A)\ .$$
in $\Sm( \Sp )$. 
Because of (\ref{may285}) we have a commutative diagram
\begin{equation}\label{jul1714}\xymatrix{\Funk(E)d\ar[r]^{e_{\R}}\ar[d]^{\tr}&\Funk(E\R)d\ar[r]^{c}\ar[d]^{\tr}&\Funk(H(A))d\ar[d]^{\tr}\\
\Funk(E)r\ar[r]^{e_{\R}}&\Funk(E\R)r\ar[r]^{c}&\Funk(H(A))r}
\end{equation} in 
$\Bundle_{geom,trans}( \Sp )$.
Thus in order to complete the construction of the left square in (\ref{may278}) 
we must  show the following proposition.
\begin{prop}\label{jul1713}
For a chain complex $A\in \Ch_{{\R}}$ of real vector spaces
{there exists} a preferred commutative square
\begin{equation}\label{sep1109}\xymatrix{ H(\Omega A)d\ar[r]^{j}\ar[d]^{\int}&\Funk(H(A))d\ar[d]^{\tr}\\
H(\Omega A)r\ar[r]^{j}&\Funk(H(A))r}\end{equation}
in $\Bundle_{geom,trans}( \Sp )$.
\end{prop}
\proof
\begin{rem} {\rm
This proposition is non-trivial since it compares transfers in two different frameworks. On the one hand, the left vertical arrow in \eqref{sep1109} is defined via integration of differential forms. On the other hand, the right vertical arrow is defined using homotopy theory. In order to understand the transition between analysis and homotopy theory we decompose the square \eqref{sep1109} into a composite \eqref{may2911} of many smaller squares. Note that, {although we are actually} constructing an object in $\Bundle_{geom,trans}( \Sp )^{{\Nerve}(\Box)}$ with $\Box$ as in \eqref{ewfwfewfwj2lkr23iori32rio2ri322r2r2} (or even in 
$\Fun(\Nerve(\Ch_{{\R}}),\Bundle_{geom,trans}( \Sp )^{{\Nerve}(\Box)})$ if we take the dependence on $A$ into account) in order to make the argument readable we write down the construction evaluated at a proper submersion $\pi:W\to B$  with a geometry $g$ and transfer {data} $b$. In order to conclude naturality one must use e.g. \eqref{hckjjklejclejlejdldkjqwdwqdqwdqwd1}.
}\hB\end{rem}

\begin{rem}{\rm 
The assertion of Proposition \ref{jul1713} looks innocent. In homotopy it asserts that the Becker-Gottlieb transfer of a cohomology class  represented by a closed from can alternatively be expressed  by the integral of the product of this form with the Euler form of the bundle. This assertion is surely well-known. The main point of  the Proposition  is to find a natural homotopy between the two constructions.
The argument that we give here is unfortunately extremely long and technical. The whole Subsection \ref{jun148} is just added to provide some of the   details.   We would be happy to see a shorter  and more comprehensive argument.

}\hB
\end{rem}

We now explain the notation used in the diagram \eqref{may2911} below.
Let $p:V\to M$ be a vector bundle over a manifold $M$. We introduce the notions of   forms and smooth function objects   on $V$ 
with proper support over $M$.  
Let $A\in \Ch_{{\R}}$ be a chain complex of real vector spaces.

\begin{ddd}
A form $\omega\in \Omega A(V)$ has proper support over $M$ if  $$p_{|\supp(\omega)}:\supp(\omega)\to M$$ is proper (i.e. preimages of compact subsets are compact). We let \begin{equation}\label{wervhwekvhewvhewhiu32z984732846238462387feucj}
\Omega Ad_{c}(V\to M)\subseteq \Omega A(V)
\end{equation} denote the  
subcomplex of forms with proper support.
\end{ddd}
We refer to Subsubsection \ref{kkfkjewfjwlefo234uro2r23r32r} for more details.

\bigskip

 Let $\bC$ be a presentable $\infty$-category and $C\in \bC$. Recall that $M^{V}\in \Top_{*}$ denotes the Thom space of the vector bundle $V\to M$.
 
\begin{ddd}
We define   $$\Funk_{c }(C)(V\to M):=C^{M^{V}}\ .$$ 
\end{ddd}

More details are given in Subsubsection \ref{kljdkljlwqjdqjdwuo1iu2e12e12e12}.

\begin{rem}{\rm 
Since $M^{V}$ is a pointed topological space, but not a manifold in general we can not simply evaluate $\Funk(C)$ on this Thom space.
This is the reason that we go back to the $\infty$-categorical foundations \eqref{jchjsachakscsacscaca} in this definition.
%

Note that   
  $\Funk(C)(V)\simeq C^{V }$.
Hence
the natural inclusion $V_{top}\to M^{V}$ of the bundle into its Thom space  \eqref{ejg12e1g2hje21e2iu1e2i1oeu21oie21e21e21e}  induces
a morphism in $\bC$
\begin{equation}\label{ksjhvskvsdvsdvdsvdsv}
\Funk_{c }(C)(V\to M)\to \Funk(C)(V)
\end{equation}
which is the analog the inclusion \eqref{wervhwekvhewvhewhiu32z984732846238462387feucj}.
In this sense we can interpret $\Funk_{c}(C)(V\to M)$ as the evaluation of $\Funk(C) $ on $V$ with proper support 
over $M$. The precise statement is  Lemma {\ref{fowfopopipoifpowiofwf24545435}}.
}\hB \end{rem}
By Corollary \ref{fewklfjwelkfjlwefwecdce3rrrfwefewfewfewfewf} the de   Rham equivalence $H(\Omega A(V))\stackrel{j, \simeq}{\to} \Funk(H(A))(V)$  refines to a de Rham equivalence  with proper support 
$$j_{c}:H(\Omega Ad_{c}(V\to M))\stackrel{  \simeq}{\to} \Funk_{ c}(H(A))(V\to M)$$
which is compatible with \eqref{ksjhvskvsdvsdvdsvdsv} and  \eqref{wervhwekvhewvhewhiu32z984732846238462387feucj}.

 \bigskip
 
Let us now assume that the vector bundle has a metric and a connection    $\nabla^{V}$.
In   Subsection \ref{may2910} we will use the Mathai-Quillen formalism in order to define 
  a Thom form $$U(\nabla^{V})\in Z^{n}(\Omega_{c/M} (V,\Lambda))\ ,$$ see Equation \eqref{hfjkwehfejkwfhkewfoiufoiu}.  Here   $\Lambda\to V$ is the  orientation bundle  of $V\to M$ and $n:=\dim(V)$ is the dimension of $V$ as a vector bundle over $M$. 
 For a trivial vector bundle $M\times \R^{k}\to M$ we adopt the canonical trivialization of $\Lambda$. 
  
  The Euler form defined by \eqref{jul1603} is related with  the Thom form by \begin{equation}\label{fweflweflr4oekwf4lewkfewfewfewfewfewfew}
e(\nabla^{V}) =0_{V}^{*} U(\nabla^{V})\ ,
\end{equation}
 where $0_{V}:M\to V$ is the zero section. 

\bigskip 

In the diagram \eqref{may2911} below we shorten the notation and write
$\Funk_{c}(E)(V)$ instead of $\Funk_{c}(V\to M)$.

\bigskip
 
 
 {We decompose 
 \eqref{sep1109} as follows.} 
\begin{equation}\label{may2911}\xymatrix{&H(\Omega A(W))\ar@/^-2cm/[dddl]^{\cdots\wedge e(g)}\ar[r]^{j}\ar[d]^{\pr^{*}(\cdots)\wedge U(\nabla^{W\times \R^{k}})}&\Funk(H(A))(W)\ar[d]_{\simeq}^{susp}\ar@/^4cm/[ddddd]_{\tr(\pi,b)}\\
&H(\Omega Ad_{c}(W\times \R^{k})[k])\ar[r]^{j_{c}}\ar[d]^{\simeq}&\Funk_{c} (H(A[k])) (W\times \R^{k})\ar[d]^{\simeq}\\
 &H(\Omega Ad_{c}(T^{v}\pi\oplus N)[k])\ar[r]^{j_{c}} \ar[d]\ar[d]^{z^{*}}&\Funk_{c} (H(A[k])) (T^{v}\pi\oplus N)\ar[d]^{z^{*}}\\
H(\Omega A (W,\Lambda)[n])\ar[r]^{\pr^{*}(\cdots)\cup U(\nabla^{N})}\ar@/_1.3cm/[ddr]^{\int_{W/B}}&H(\Omega Ad_{c}(N)[k])\ar[r]^{j_{c}}\ar[d]^{\emb_{*}}&\Funk_{c} (H(A[k])) (N)\ar[d]^{\clps(\emb)^{*}}\\
&H(\Omega Ad_{c}(B\times \R^{k})[k])\ar[r]^{j_{c}}\ar[d]^{\int_{B\times \R^{k}/B}}&\Funk_{c} (H(A[k])) (B\times \R^{k})\ar[d]_{\simeq}^{desusp}\\
&H(\Omega A (B))\ar[r]^{j}&\Funk(H(A))(B)}\ .\end{equation}

The {normal} bundle $N\to W$  and the map $\emb:N\to B\times \R^{k}$ are parts of the transfer {data} $b$.  The morphism $\emb_{*}$ is defined in Corollary \ref{jekljdlkwjwkldjlqwkjwqldqwdwqdqwdqwdqwdqwd} 
and the corresponding square commmutes by Corollary \ref{fjhwjkf88z23430fwefwfewfewfewejkf8237458234r}.
 
 \bigskip
 
The map $z^{*}$ is induced
by the zero section of $T^{v}\pi$. It is {a} morphism of vector bundles.
  The square  involving $z^{*}$ is  a case of the commuting square \eqref{dhkjdhqkjwqwdiuoiuoiuoiu1}.
  
 \bigskip
  
  The suspension equivalence $susp$ and its inverse $desusp$ are defined in Definition \ref{fewklfwejkfweoiuo23u42398ruekjdnqwefewfwew4r}.
  
\bigskip

The diagram involves   involves the Thom forms for euclidean vector bundles with connection as mentioned above.  The choice of connections will be made precise during the proof of  Lemma \ref{jul1930}. 
The lower left square commutes already on the level of forms. 
{Here and in the text below the symbol $\Lambda$ stands for an appropriate orientation bundle; the precise object {will be clear} from the context in each case.}

\bigskip

The lower middle square commutes by Proposition \ref{ekjfwklejflfelwfjewfewfew42243242342342342342}.
 Commutativity of the upper  middle cell is provided by Proposition \ref{kldjhqwkdjqdwkljqwdiouqwodqwidqwdwqdwqdwqd}.

\bigskip

The commutativity of the {right} hexagon is the definition \eqref{feb2702} of $\tr(\pi, b)$, written in terms of function spectra with proper support in the Thom space picture introduced in Subsubsection \ref{kljdkljlwqjdqjdwuo1iu2e12e12e12}.

\bigskip

The outer square is the desired commutative square \eqref{sep1109}. We get it as the composite of the {smaller cells, all of which commute.

\bigskip

  However,} there  is one cell left   {which} does not obviously commute, namely
the upper left.
This will be discussed in  following lemma.

\begin{lem}\label{jul1930}
The upper left cell
\begin{equation}\label{feb2703}\xymatrix{H(\Omega A(W))\ar[rr]^(0.43){\pr^{*}(\cdots) \wedge U(\nabla^{W\times \R^{k}})}\ar[d]^{\cdots \wedge e(g)}&&H(\Omega A_{c}(W\times \R^{k})[k])\ar[r]^{\simeq}&H(\Omega A_{c}(T^{v}\pi\oplus N)[k])\ar[d]^{z^{*}}\\
H(\Omega A(W,\Lambda)[n])\ar[rrr]^{\pr^{*}(\cdots)\cup U(\nabla^{N})}&&&H(\Omega A_{c}(N)[k])
}\end{equation}
 in (\ref{may2911}) commutes.
\end{lem}
\proof
We have a natural isomorphism of   vector bundles \begin{equation}\label{fewfehwfkewfewfewfwfewr2342423rfee33564565}
W\times \R^{k}\cong T^{v}\pi\oplus N\ .
\end{equation}

The  bundle on the left-hand side is equipped with the trivial connection $\nabla^{W\times \R^{k}}$ and  metric. The right-hand side has the connection $\nabla^{T^{v}\pi}\oplus \nabla^{N}$ and the metric on $T^{v}\pi$ given by the Riemannian structure and {the}   metric on $N$ induced by the inclusion $N\to W\times \R^{k}$. Here $\nabla^{T^{v}\pi}$ is the connection induced by the Riemannian structure $g$, and $\nabla^{N}$ is the connection obtained by the projection of $\nabla^{W\times \R^{k}}$ to $N$. Using the isomorphism of vector bundles \eqref{fewfehwfkewfewfewfwfewr2342423rfee33564565}
we can consider the two Thom forms
$U(\nabla^{W\times \R^{k}})$ and $U(\nabla^{T^{v}\pi}\oplus \nabla^{N})$ as forms on $W\times \R^{k}$ with proper support over $W$. 
They are related  by a natural transgression form
$$\tilde U(\nabla^{W\times \R^{k}},\nabla^{T^{v}\pi}\oplus \nabla^{N})\in \Omega^{-1}_{c}(W\times \R^{k},\Lambda)$$ such that
$$d\tilde U(\nabla^{W\times \R^{k}},\nabla^{T^{v}\pi}\oplus \nabla^{N})=U(\nabla^{W\times \R^{k}})-U(\nabla^{T^{v}\pi}\oplus \nabla^{N})\ .$$
We further use that
$$U(\nabla^{T^{v}\pi}\oplus \nabla^{N})=U(\nabla^{T^{v}\pi})\wedge U( \nabla^{N})\ .$$
The pull-back  of the natural transgression form  by the inclusion $z:N\to T^{v}\pi\oplus N$ provides the chain homotopy for the {version of \eqref{feb2703}} obtained before application of {the Eilenberg-MacLane spectrum functor} $H$.
Here we use that
$e(g)=0^{*}_{T^{v}\pi} U(\nabla^{T^{v}\pi})$ by \eqref{fweflweflr4oekwf4lewkfewfewfewfewfewfew}.
The chain homotopy induces the required  homotopy after application of $H$, and this construction is natural in the bundle $\pi$.
\hB

%

 This finishes the construction of the diagram  
 (the functorial version of
 \eqref{may278})   \begin{equation}\label{flwejfwleif09uo3i23r323r3r32r}
\xymatrix{\Funk(E)d\ar[r]^{\rat}\ar[d]^{\tr}&H(\Omega A )d\ar[d]^{H(\int)}&H(\sigma \Omega A )d\ar[l]\ar[d]^{H(\int_{\sigma})}\\
 \Funk(E)r\ar[r]^{\rat}&H(\Omega A )r&H(\sigma \Omega A )r\ar[l]}
\end{equation}
  in $\Bundle_{geom,trans}( \Sp )$
 which induces the transfer map of bundle spectra
 $$\hat \tr :\Diff(E)d\to \Diff(E)r\ .$$
Its evaluation at a proper submersion $W\to B$ with Riemannian structure $g$ and transfer {data} $b$ provides, upon application of $\pi_{0}$, the differential transfer
$$\hat \tr(g,b) :\widehat E^{0}(W)\to \widehat E^{0}(B)$$  whose existence is asserted in Theorem \ref{may3102}. 
The diagram (\ref{may3103222})
commutes by construction.  In order to finish  the construction of $\hat \tr$ we must prove Lemma \ref{sep2104}. 

\bigskip

\proof{[of Lemma \ref{sep2104}]}
Given another choice of {stable transfer data} $b'$ on $\pi$, we can choose
 {stable transfer data} $\tilde b$ on
$\tilde \pi:=\pi\times \id:W\times [0,1]\to B\times [0,1]$ which interpolates between $b$ and $b^{\prime}$.
We equip {the  bundle $\tilde \pi$} with the Riemannian structure $\tilde g$ induced from $g$. We  let
$\pr_{W}:W\times [0,1]\to W$, $\pr_{B}:B\times [0,1]\to B$ be the projections and $i_{j}:B\to B\times [0,1]$, $j=0,1$ be the inclusions
at the endpoints of the interval. For the moment we let
$\hat \tr^{\prime}$ and $\tilde \tr$ denote the differential transfers for
$(\pi,g,b^{\prime})$ and $(\tilde \pi,\tilde g,\tilde b)$.
Using the naturality of the diagram (\ref{may278}) we get
$$\hat \tr=i_{0}^{*}\circ \tilde \tr\circ \pr^{*}\ , \quad \hat \tr^{\prime}=i_{1}^{*}\circ \tilde \tr\circ \pr^{*}\ .$$
We now apply the homotopy formula \eqref{feb2710} which for $x\in \widehat E^{0}(W)$ gives
$$\hat \tr^{\prime}(x)-\hat \tr(x)=a(\int_{[0,1]\times B/B} R(\tilde \tr(\pr_{W}^{*}x)))\ .$$
Since $e(\tilde g)=\pr^{*}_{W}e(g)$ we see that
$${\int_{[0,1]\times B/B}}R(\tilde \tr(\pr_{W}^{*}x))=\int_{[0,1]\times B/B}\pr_{B}^{*}\int_{[0,1]\times W/[0,1]\times B} R(x)\wedge e(g)=0$$
and therefore
$$\hat \tr(x)=\hat \tr^{\prime}(x)\ . $$
\hB

 An inspection of the proof of Proposition \ref{jul1713} together with the Corollaries \ref{lfjwlekfjeluo32ieifhqkf} and  \ref{jdlqwdqwlkdjwqldqwdwqdqwdqwdqd}
yields the  naturality of the transfer with respect to the differential data $(E,A,c)$. Let   $\Box\Box$ denote  the category of the shape
$$
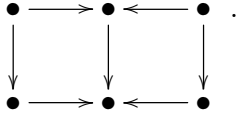
\xymatrix{\bullet\ar[d]\ar[r]&\bullet\ar[d]&\bullet\ar[l]\ar[d]\\\bullet\ar[r]&\bullet&\bullet\ar[l]}\ .$$
\begin{kor}
We have actually constructed a functor
$$\widehat{\Sp}\to \Bundle_{geom,trans}(\Sp)^{\Nerve(\Box\Box)}\ ,$$
which associates to the differential data $(E,A,c)$ the diagram \eqref{flwejfwleif09uo3i23r323r3r32r}.
In particular, for a morphism of differential data $(E,A,c)\to (E^{\prime},A^{\prime},c^{\prime})$ we get a commuting diagram of bundle spectra or (after evaluation at a bundle $W\to B$ with geometry), of  differential cohomology groups
$$\xymatrix{\Diff(E,A,c)d\ar[d]^{\hat \tr}\ar[r]&\Diff(E^{\prime},A^{\prime},c^{\prime})d\ar[d]^{\hat \tr}\\
\Diff(E,A,c)r\ar[r]&\Diff(E^{\prime},A^{\prime},c^{\prime})r}\ , \quad \xymatrix{\widehat{E}^{0}(W) \ar[d]^{\hat \tr}\ar[r]&\widehat{E}^{\prime,0}(W) \ar[d]^{\hat \tr}\\
\widehat{E}^{ 0}(B) \ar[r]& \widehat{E}^{\prime,0}(B) } \ .$$
 \end{kor}
This finishes the proof of Theorem  \ref{may3102}, 4, in particular of \eqref{aug1012}.

\subsection{Proof of  \eqref{jun0304}}\label{jul1710}

By Equation (\ref{hgr1}) we have {a natural isomorphism}
\begin{equation}\label{feb2750}E\R/\Z^{-1}(M)\cong \pi_{0}(\Sigma^{-1}\Funk(E\R/\Z)(M))\ ,\end{equation}
where the spectrum $\Sigma^{-1}\Funk(E\R/\Z)(M)$ fits into the homotopy pull-back (\ref{jul1711}).

{Now let} $W\to B$ be a proper submersion with Riemannian structure $g$ and transfer {data} $b$.
{Under the identification \eqref{feb2750}},
the Becker-Gottlieb transfer \begin{equation}\label{jul1715}
\tr^{*}:
E\R/\Z^{-1}(W)\to E\R/\Z^{-1}(B)\end{equation}
 corresponds   to the map induced on $\pi_{0}$ by the map
$$\tr: \Sigma^{-1}\Funk(E\R/\Z)(W)\to \Sigma^{-1}\Funk(E\R/\Z)(B),$$ which in turn is
induced by the map of pull-back {diagrams}
$$\xymatrix{\Funk(E)(W)\ar[r]\ar[d]^{\tr}&\Funk(E\R)(W)\ar[d]^{\tr}&0\ar[l]\ar[d]^{0}\\
\Funk(E)(B)\ar[r]&\Funk(E\R)(B)&0\ar[l]}\ .$$
Using Proposition \ref{jul1713} and the diagram (\ref{jul1714})
we see that the same map (\ref{jul1715}) is induced by the map of pull-back {diagrams}
\begin{equation}\label{jul1717}\xymatrix{\Funk(E)(W)\ar[r]^{\rat}\ar[d]^{\tr}&H(\Omega A (W))\ar[d]^{\int}&0\ar[l]\ar[d]^{0}\\
\Funk(E)(B)\ar[r]^{\rat}&H(\Omega A (B))&0\ar[l]} \ .\end{equation} By  
(\ref{eqq4111}) the  limit of  upper line in (\ref{jul1717}) is   $\Sigma^{-1}\Funk(E\R/\Z)(W)$, while {the limit of the lower is}  $\Sigma^{-1}\Funk(E\R/\Z)(B)$.

The map $0\to  \sigma \Omega A (W) $ induces a map of pull-back {diagrams}
$$\xymatrix{\Funk(E)(W)\ar[r]^{\rat}\ar@{=}[d]&H(\Omega A (W))\ar@{=}[d]&0\ar[l]\ar[d]^{0}\\
\Funk(E)(W)\ar[r]^{\rat}&H(\Omega A (W))&H(\sigma \Omega A (W))\ar[l]}$$ which in turn induces the map
$$\Sigma^{-1}\Funk(E\R/\Z)(W)\to \Diff(E)(W)\ .$$
We have a similar diagram for $B$.
In order to prove that \eqref{jun0304} commutes we must produce the commuting diagram
 $$\xymatrix{&\Funk(E)(W)\ar[rr]^(0.6){\rat}\ar[dd]&&H(\Omega A(W))\ar[dd]&&H(\sigma \Omega A(W))\ar[dd]^{\int_{W/B}\dots\wedge e(g)}\ar[ll]\\\Funk(E)(W)\ar[rr]^(0.6){\rat}\ar[dd]\ar@{=}[ur]&&
 H(\Omega A(W))\ar[dd]_(0.35){\int_{W/B}\dots\wedge e(g)}\ar@{=}[ur]&&0\ar@{=}[dd]\ar[ll]\ar[ur]&\\&\Funk(E)(B)\ar[rr]^(0.6){\rat} &
 &H(\Omega A(B))&&H(\sigma \Omega A(B))\ar[ll]\\\Funk(E)(B)\ar[rr]^(0.6){\rat}\ar@{=}[ur]&&
 H(\Omega A(B))\ar@{=}[ur]&&0\ar[ll]\ar[ur]
 }\ .
 $$
The left cube commutes by Proposition \eqref{jul1713} yielding the front and back faces. In order to get the right cube we apply 
$H$ to the strictly commutative diagram of chain complexes
 $$\xymatrix{&\Omega A(W)\ar[dd]&&\sigma \Omega A(W)\ar[dd]^{\int_{W/B}\dots\wedge e(g)}\ar[ll]\\
 \Omega A(W)\ar[dd]_{\int_{W/B}\dots\wedge e(g)}\ar@{=}[ur]&&0\ar@{=}[dd]\ar[ll]\ar[ur]&\\
 &\Omega A(B)&&\sigma \Omega A(B)\ar[ll]\\
 \Omega A(B)\ar@{=}[ur]&&0\ar[ll]\ar[ur]
 }\ .
 $$
 This produces the correct filler of the back face in view of the construction of \eqref{feb2701}.
\hB

\subsection{Functoriality of the transfer for iterated bundles}\label{fkwelfewfewfwef}

We consider two composable  proper submersions with closed fibres
$$W\xrightarrow{\pi} B\xrightarrow{\sigma} Z\ .$$ The composition
$\kappa:=\sigma\circ \pi$ is again a proper submersion with closed fibres. Let 
$g_{\pi}$ and $g_{\sigma}$ be Riemannian structures on $\pi$ and $\sigma$ (Definition \ref{may272}). Then we get a Riemannian structure $g_{\kappa}$ on the composition as follows. The decompositions
$$TW\cong T^{v}\pi\oplus \pi^{*}TB\ , \quad TB\cong T^{v}\sigma\oplus \sigma^{*}TZ$$ given by the horizontal distributions of $g_{\pi}$ and $g_{\sigma}$ induce
 a decomposition
$$T^{v}\kappa\cong T^{v}\pi\oplus \pi^{*}T^{v}\sigma\ ,$$ where $\pi^{*}T^{v}\sigma$ is canonically identified with a subbundle of $T^{h}\pi$. We define the vertical metric  of the Riemannian structure $g_{\kappa}$ by 
$$g^{T^{v}\kappa}:=g^{T^{v}\pi}\oplus \pi^{*}g^{T^{v}\sigma}\ .$$ Furthermore, we let
$$T^{h}\kappa:=\pi^{*}T^{h}\sigma \subset T^{h}\pi$$
be the horizontal subspace of $g_{\kappa}$.

We now have two connections on $T^{v}\kappa$, namely the connection
$\nabla^{T^{v}\kappa}$ induced by $g_{\kappa}$ and the direct sum connection
$\nabla^{\oplus}:=\nabla^{T^{v}\pi}\oplus \pi^{*}\nabla^{T^{v}\sigma}$.
In general they do not coincide. We define the transgression Euler form
\begin{equation}\label{uliapr1201}\eta:=\tilde e(\nabla^{T^{v}\kappa},\nabla^{\oplus})\in \Omega^{n-1}(W,\Lambda)\ ,\end{equation}
where $n:=\dim(W)-\dim(Z)$.
It satisfies
$$d\eta=e(g_{\kappa})-e(g_{\pi})\wedge \pi^{*}e(g_{\sigma})\ .$$

We now choose transfer data $b_{\sigma}$ and $b_{\pi}$ for $\sigma$ and $\pi$. Recall the definition of the point-set level transfer \eqref{may279} denoted by $\check{\tr}(-)$.
\begin{lem}\label{jul2012}
There is a canonical choice of transfer data $b_{\kappa}$ such that
$$\check{\tr}(b_{\sigma})\circ \check{\tr}(b_{\pi})=\check{\tr}(b_{\kappa})\ .$$
\end{lem}
\proof
 Let $W\hookrightarrow  B\times \R^{k_{\pi}}$ and $B\hookrightarrow Z\times \R^{k_{\sigma}}$ be the fibrewise embeddings for $b_{\pi}$ and $b_{\sigma}$.
They induce a fibrewise embedding
$W\hookrightarrow Z\times \R^{k_{\pi}+k_{\sigma}}$ by composition.
Let $\emb_{\pi}:N_{\pi}\hookrightarrow B\times \R^{k_{\pi}}$ and $\emb_{\sigma}:N_{\sigma}\hookrightarrow  Z\times \R^{k_{\sigma}}$ be the extensions to open embeddings of normal bundles.
Then we get an open embedding
$$\emb_{\kappa}:N_{\kappa}:=N_{\pi}\oplus \pi^{*}N_{\sigma}\hookrightarrow Z\times \R^{k_{\kappa}}$$
with $k_{\kappa}:=k_{\pi}+k_{\sigma}$.
With these choices the diagram of pointed topological spaces
\begin{equation}\label{jul2005}\xymatrix{\Sigma^{k_{\kappa}}_{+}Z\ar[d]^{\cong}\ar[r]^{\clps_{\kappa}}&W^{N_{\kappa}}\ar[rr]^{z_{\kappa}}&&\Sigma^{k_{\kappa}}_{+}W \\
\Sigma^{k_{\pi}}\Sigma^{k_{\sigma}}_{+}Z \ar[r]^{\clps_{\sigma}}&\Sigma^{k_{\pi}} B^{N_{\sigma}}\ar[r]^{z_{\sigma}}& \Sigma^{k_{\kappa}}_{+} B \ar[r]^{\clps_{\pi}}&W^{N_{\pi}}\ar[u]^{z_{\pi}}
}
\end{equation}
commutes (compare with \eqref{wqdqjwdqwjdwqjkdwjqkdkwqqwdwqdwqdw}). This implies the assertion.
\hB

We let
$\hat \tr_{\pi}$, $\hat \tr_{\sigma}$, and $\hat \tr_{\kappa}$ denote the differential Becker-Gottlieb transfers 
associated by the first part of  Theorem \ref{may3102} to the proper submersions $\pi$, $\sigma$, and $\kappa$ with the geometries $g_{\pi}$, $g_{\sigma}$ and $g_{\kappa}$ as above and transfer data 
$b_{\pi}$, $b_{\sigma}$ and  $b_{\kappa}$
as described in Lemma \ref{jul2012}.
The following proposition clarifies the functoriality of the differential Becker-Gottlieb transfer with respect to iterated fibre bundles.
{Recall the definition of the transgression Euler form $\eta$ in \eqref{uliapr1201}.}
\begin{prop}\label{jul2015}
For $x\in \widehat E^{0}(W)$ we have
\begin{equation}\label{jul1901}\hat \tr_{\kappa}(x)=\hat \tr_{\sigma}(\hat \tr_{\pi}(x))+a(\int_{W/Z} R(x)\wedge \eta)\ .\end{equation}
\end{prop}
\proof
The map
$\hat \tr_{\sigma}\circ \hat \tr_{\pi}$ is induced by the composition of  two diagrams  of the form
(\ref{may278})
\begin{equation}\xymatrix{\Funk(E)(W)\ar[d]^{\tr_{\pi}}\ar[r]^{\rat}&H(\Omega A (W))\ar[d]^{H(\int_{W/B}\cdots\wedge e(g_{\pi}))}&H(\sigma\Omega A (W))\ar[l]\ar@/_-0.3cm/[d]^{H(\int_{W/B}\cdots\wedge e(g_{\pi}))}\\\Funk(E)(B)\ar[d]^{\tr_{\sigma}}\ar[r]^{\rat}&H(\Omega A (B))\ar[d]^{H(\int_{B/Z}\cdots\wedge e(g_{\sigma}))}&H(\sigma \Omega A (B))\ar[l]\ar@/_-0.3cm/[d]^{H(\int_{B/Z}\cdots\wedge e(g_{\sigma}))}\\
\Funk(E)(Z)\ar[r]^{\rat}&H(\Omega A (Z))&H(\sigma \Omega A (Z))\ar[l]}\end{equation}
In general, given an element $\omega\in \Omega A^{-1}(M)$, the transformation
$$\dots+a(R(\dots)\wedge \omega):\widehat E^{0}(M)\to \widehat E^{0}(M)$$ is induced by the diagram
$$\xymatrix{\Funk(E)\ar@{=}[d]\ar[r]^{\rat}&H(\Omega A  (M))\ar@{=}[d]&H(\sigma \Omega A (M))\ar@/_-0.3cm/[d]^{\id_{\Omega A(M)}+\cdots \wedge d\omega}\ar[l]\ar@{:>}[dl]^{\omega}\\
\Funk(E)\ar[r]^{\rat}&H(\Omega A (M))&H(\sigma \Omega A (M))\ar[l]}\ ,$$
where the homotopy in the right square is induced from the chain homotopy between $\id_{\Omega A(M)}$ and $\id_{\Omega A(M)}+\dots\wedge d\omega$ given by $\omega$. 
We conclude that the right-hand side of Equation (\ref{jul1901}) is induced by 
\begin{equation}\xymatrix{\Funk(E)(W)\ar[d]^{\tr_{\pi}}\ar[r]^{\rat}&H(\Omega A  (W))\ar[d]^{H(\int_{W/B}\cdots\wedge e(g_{\pi}))}&H(\sigma\Omega A (W))\ar[l]\ar@/_-2cm/[dd]^{H(\int_{W/Z}\cdots\wedge e(g_{\kappa}))}\ar@{:>}[ddl]^{\eta} \\ \Funk(E)(B)\ar[d]^{\tr_{\sigma}}\ar[r]^{\rat}&H(\Omega A (B))\ar[d]^{H(\int_{B/Z}\cdots\wedge e(g_{\sigma}))}& \\
\Funk(E)(Z)\ar[r]^{\rat}&H(\Omega A (Z))&H(\sigma \Omega A (Z))\ar[l] }\end{equation}
It can be written as
\begin{equation}\label{jul1905}\xymatrix{\Funk(E)(W)\ar[d]^{\tr_{\pi}}\ar[r]^{\rat}&H(\Omega A (W))\ar[d]^{H(\int_{W/B}\cdots\wedge e(g_{\pi}))}\ar@/^3cm/@{->}[dd]^{H(\int_{W/Z}\cdots\wedge e(g_{\kappa})}&&H(\sigma\Omega A (W))\ar[ll]\ar@/_-2cm/[dd]^{H(\int_{W/Z}\cdots\wedge e(g_{\kappa}))}\\\Funk(E)(B)\ar[d]^{\tr_{\sigma}}\ar[r]^{\rat}&H(\Omega A (B))\ar[d]^{H(\int_{B/Z}\cdots\wedge e(g_{\sigma}))}&\ar@{:>}[l]^{\eta}& \\
\Funk(E)(Z)\ar[r]^{\rat}&H(\Omega A (Z))&&H(\sigma \Omega A (Z))\ar[ll] }\ , \end{equation}
where now the right square comes from a strictly commuting diagram of complexes. 
 We therefore must show that
the left part of the diagram (\ref{jul1905}) is equivalent to the square
$$\xymatrix{\Funk(E)(W)\ar[d]^{\tr_{\kappa}}\ar[r]^{\rat}&H(\Omega A_{\infty}(W))\ar[d]^{H(\int_{W/Z}\cdots\wedge e(g_{\kappa})}\\
\Funk(E)(Z)\ar[r]^{\rat}&H(\Omega A_{\infty}(Z))}$$
used in the construction of $\hat \tr_{\kappa}$.
Recall that we can expand this to
$$\xymatrix{\Funk(E)(W)\ar[d]^{\tr_{\kappa}}\ar[r]^{\rat}&\Funk(H(A))(W)\ar[d]^{\tr_{\kappa}}&\ar[l]^{j}H(\Omega A (W))\ar[d]^{H(\int_{W/Z}\cdots\wedge e(g_{\kappa}))}\\
\Funk(E)(Z)\ar[r]^{\rat}&\Funk(H(A))(Z)&H(\Omega A (Z))\ar[l]^{j}}\ .$$
Since $\tr_{\kappa}=\tr_{\sigma}\circ \tr_{\pi}$ is functorial in the spectrum argument of $\Funk$ the argument boils down to showing that the outer square of
\begin{equation}\label{jul1920}\xymatrix{&
H(\Omega A (W))\ar[d]_{H(\int_{W/B}\cdots\wedge e(g_{\pi}))}\ar@/_3cm/[dd]_{H(\int_{W/Z}\cdots\wedge e(g_{\kappa}))}\ar[r]&\Funk(H(A))(W)\ar[d]^{\tr_{\pi}}\ar@/^3cm/[dd]^{\tr_{\kappa}}&\\\ar@{:>}[r]^{\eta}&
H(\Omega A(B))\ar[d]_{H(\int_{B/Z}\cdots\wedge e(g_{\sigma}))}\ar[r]&\Funk(H(A))(B)\ar[d]^{\tr_{\sigma}}&\ar@{:>}[l]^-{\id}\\&
H(\Omega A_{\infty}(Z))\ar[r]&\Funk(H(A))(Z)&}\end{equation}
obtained by the composition of the four small squares of (\ref{jul1920})
is equivalent to the square constructed in Proposition \ref{jul1713} evaluated at the submersion
$\kappa:W\to Z$ with geometry $g_{\kappa}$. In detail the composition of the two inner squares
of  (\ref{jul1920}) expands to the composition of two diagrams of the form (\ref{may2911}), one for $\pi$ and another one for $\sigma$. One checks that the assertion follows from the   two  Lemmas \ref{jul2011}  and \ref{jul2013}.
\begin{lem}\label{jul2011}
The composition of the three homotopies
$$\xymatrix{&&\ar@{:>}[d]^{\eta}&&\\&\ar@{:>}[d]&&\ar@{:>}[d]&\\H(\Omega A (W))\ar@/^1.5cm/[rr]^{H(\int_{W/B}\dots\wedge e(g_{\pi}))}\ar@/^3cm/[rrrr]^{H(\int_{W/Z}\dots\wedge e(g_{\kappa}))}\ar[rr]_{f_{\pi}}&&H(\Omega A (B))\ar[rr]_{f_{\sigma}}\ar@/^1.5cm/[rr]^{H(\int_{B/Z}\dots e(g_{\sigma}))}&&H(\Omega A (Z))}$$  given by $\eta$ and two applications of Lemma \ref{jul1930} is homotopic to the homotopy 
$$\xymatrix{&\ar@{:>}[d]&\\H(\Omega A (W))\ar@/^1.5cm/[rr]^{H(\int_{W/Z}\dots\wedge e(g_{\kappa}))}\ar[rr]_{f_{\kappa}}&&H(\Omega A (Z))}$$
for the composition $\kappa$ (again by Lemma \ref{jul1930}), where $f_{\sigma}$, $f_{\pi}$, and $f_{\kappa}$ are the corresponding vertical compositions in the left column of (\ref{may2911}).
\end{lem}
\proof
This follows from the fact that all homotopies are obtained from transgressions and that the spaces of connections are contractible. The main point is   that the combination of transgressions $$\tilde e(\nabla_{0},\nabla_{1})+\tilde e(\nabla_{1},\nabla_{2})-\tilde e(\nabla_{0},\nabla_{2})$$ for three connections on a vector bundle is exact.
\hB

\begin{lem}\label{jul2013}
The composition of the two middle squares in (\ref{jul1920})
is equivalent to the corresponding middle square in (\ref{may2911}) for $\kappa$.
\end{lem}
\proof
We expand the small squares in (\ref{jul1920}) as
 a composition of two diagrams of the form (\ref{may2911}) (one on top of the other) and look at the sequence from up to down of small squares obtained from the middle parts.
The upper middle square of the second and the lower middle square of the first    are inverse to each other by Proposition \ref{ekjfwklejflfelwfjewfewfew42243242342342342342}  and cancel out. 
The remaining ladder defines a square
$$\xymatrix{H(\Omega A (W))\ar[r]^{j}\ar[d]&\Funk(H(A))(W)\ar[d]\\
H(\Omega A (Z))\ar[r]^{j}&\Funk(H(A))(Z)}\ .$$
Since (\ref{jul2005}) commutes  and using Remark \ref{hckjjklejclejlejdldkjqwdwqdqwdqwd} we conclude that this square is equivalent to the square for $\kappa$.
\hB 
This finishes the proof of Proposition \ref{jul2015}. \hB

\color{black}

%

\section{A transfer index conjecture}\label{sep2205}
\subsection{Introduction}

An index theorem relates two different constructions of fibre integration,
usually an analytic or geometric one and a homotopy theoretic version, see  (\ref{sep2401}). In the present section we discuss an index theorem for the fibre integration of locally constant sheaves of finitely generated projective $R$-modules, where $R$ the ring of integers in a number field. 

The index theorem which we are going to present here is not yet a proven theorem but the Conjecture \ref{may3103} about a differentials refinement of the  index theorems of  Dwyer-Weiss-Williams \cite{MR1982793} and its  characteristic class version  proven by Bismut-Lott   \cite{MR1303026}.   
{We will subsequently refer to it as the {\emph{transfer index conjecture}} (TIC).}
The main innovation of the present paper is the construction of all elements necessary to state
the  {transfer index conjecture} \ref{may3103}. We are far from having a proof which will probably be quite complicated. In order to support the validity of the TIC \ref{may3103} in Subsection \ref{sep0901} we discuss interesting consequences of the TIC which can be verified independently.

\subsection{The statement of the transfer index {conjecture}}\label{kljfelwkfwefwefwefewffewfef}


In this subsection give the statement of the transfer index {conjecture}, leaving the details of the construction of the analytic index to Subsection \ref{jun183}.

\bigskip

For a start  we formulate the Dwyer-Weiss-Williams Theorem \ref{dww} in the language of the present paper.
We consider a proper submersion $\pi:W\to B$. Let $\cV$ be a locally constant sheaf of finitely generated $R$-modules on $W$. It gives rise to a $K$-theory class
$[\cV]\in KR^{0}(W)$, see Definition \ref{may121311}. We have a Becker-Gottlieb transfer (\ref{jun182})
$$\tr^{*}:KR^{0}(W) \to  KR^{0}(B)\ ,$$
and we define
$$\ind^{top}(\cV):= \tr^{*} [\cV]\ .
$$
In order to recapture the naturality of this construction for pull-backs along cartesian squares
\begin{equation}\label{jun186}\xymatrix{W^{\prime}\ar[d]^{\pi^{\prime}}\ar[r]^{F}&W\ar[d]^{\pi}\\B^{\prime}\ar[r]^{f}&B}\end{equation} we will again use the language of $\Bundle_{geom}$-commutative monoids (see Subsection \ref{jun143} for more details about this notation).
We have the $\Bundle_{geom}$-commuative monoids
$\bloc(R) \: d$ and $KR^{0}r$ which associate to a proper submersion $\pi:W\to B$
the commutative monoid  $\bloc(R) (W)$ {(see \eqref{iwjiejfwofwefewfewflkewfjlfefwef} for a definition)} of isomorphism classes of   locally constant sheaves of finitely generated    $R$-modules on $W$, and
the commutative monoid $KR^{0}(B)$, respectively.
Then we can consider the topological index
as a map of $\Bundle_{geom}$-{commutative monoids}
 $$\ind^{top}:\bloc(R)\: d\to KR^{0}r\ .$$  
 
 \bigskip
 
 We can form the higher derived images $R^{i}\pi_{*}(\cV)\in \Sh_{\Mod(R)}(B)$ of the sheaf $\cV\in \Sh_{\Mod(R)}(W)$. By Lemma \ref{kjhdqkwdqwdqwdqwdqwdqwdwqd} they are again   locally constant sheaves of finitely generated  $R$-modules. {The analytic index is defined as  the $K$-theory class
\begin{equation}\label{may3116}\ind^{an}(\cV):=\sum_{i\ge 0} (-1)^{i} [R^{i}\pi_{*}(\cV) ]  \in KR^{0}(B) \ .\end{equation}}{Here we again use the topological cycle map given in Definition \ref{may121311}.}
In order to encode the naturality of the analytic index with respect to diagrams of the form (\ref{jun186})
 we   consider it as a map of $\Bundle_{geom}$-{commutative monoids}
$$\ind^{an}:\bloc(R)\: d\to KR^{0}r\ .$$
The index theorem of   
Dwyer-Weiss-Williams \cite{MR1982793} {(in particular Cor. 8.12 in the revised 2001 version)} now has as a consequence:
\begin{theorem}\label{dww}
We have the equality of the topological and analytical index 
$$\ind^{top}=\ind^{an}:\bloc(R) \:  d\to KR^{0}r\ .$$
\end{theorem}

The main purpose of the present section is to state a refinement of this theorem to differential algebraic $K$-theory. {We define the differential algebraic $K$-theory of the number ring $R$ using the data  $(KR,A,c)$ described  in Definition \ref{jhdkjdqwdqwdqwdqwdqwdwqd}.}

We fix an object $$(\pi: W\to B,g)\in \Bundle_{geom}$$ consisting of a proper submersion $\pi:W\to B$ with
a Riemannian structure $g\in \Riem(W\to B)$, see Definition \ref{may272}.
We furthermore consider  an isomorphism class of  locally constant sheaves of finitely generated    $R$-modules
  with geometry
$$(\cV,h^{\cV})\in \bloc_{geom}(R)(W)\ ,$$ see Definition \ref{defpre}. 
We have a cycle map $$\cycl: \bloc_{geom}(R)(W)\to \widehat{KR}^{0}(W)$$ (Theorem \ref{thecycleproj1}) and a transfer $$\hat \tr:\widehat{KR}^{0}(W)\to \widehat{KR}^{0}(B)$$(Theorem \ref{may3102}).
\begin{ddd}\label{jun181}
The topological index of $(\cV,h^{\cV})$ is the element
$$\hind^{top}(\cV,h^{\cV}):=\hat \tr(\cycl(\cV,h^{\cV}))\in \widehat{KR}^{0}(B)\ .$$
\end{ddd}
By construction, the topological index can be understood as a map of $\Bundle_{geom}$-{commutative monoids}
$$\hind^{top}: \bloc_{geom}(R)d\to \widehat{KR}^{0} r\ .$$
In {Subsection \ref{jun183}} we will construct an analytic index
$$\hind^{an}: \bloc_{geom}(R)d\to \widehat{KR}^{0} r\ .$$
The principal idea is that fibrewise Hodge theory provides geometries $h^{R^{i}\pi_{*}(\cV)}$ on the    locally constant sheaves of finitely generated    $R$-modules  $R^{i}\pi_{*}(\cV)$ so that we can define  (using the notation \eqref{may121207})
$$\hind^{an}(\cV,h^{\cV}):=\sum_{i\ge 0} (-1)^{i} \cycl \left(R^{i}\pi_{*}(\cV),h^{R^{i}\pi_{*}(\cV)}\right) +a(\cT)\ ,$$
where the correction $a(\cT)$ {involves the higher torsion form of Bismut-Lott.  It is added in} order to adjust the curvature of $\hind^{an}(\cV,h^{\cV})$ so that it coincides with that of $\hind^{top}(\cV,h^{\cV})$.

\bigskip

We have now all elements to formulate the TIC.
\begin{con}[{Transfer index conjecture}] \label{may3103}
{The two maps
$$\hind^{top}=\hind^{an}:\bloc_{geom}(R)d\to \widehat{KR}^{0} r$$
 of $\Bundle_{geom}$-{commutative monoids} are equal.}
\end{con}

This conjecture appears to be very natural. The equality holds true if we specialize to the curvature or the underlying algebraic $K$-theory classes. {Recall the structure maps $R$ and $I$ of differential cohomology, Definition \ref{jun216}.}
For the curvature we get essentially a reformulation of the theorem of Bismut-Lott \cite{MR1303026}:
\begin{theorem}\label{bl}
We have the equality
$$R\circ \hind^{top}=R\circ \hind^{an}:\bloc_{geom}(R)d\to Z^{0}(\Omega A) r\ .$$
\end{theorem}
We will verify this theorem at the end of Subsection \ref{jun183} after the construction of
the analytic index. On the level of underlying algebraic $K$-theory classes we get the  
Dwyer-Weiss-Williams Theorem \ref{dww}: 
\begin{theorem}
  $$I\circ \hind^{top}=I\circ \hind^{an}:\bloc_{geom}(R)d\to KR^{0}r\ .$$
  \end{theorem}
  \proof This immediately follows from
  $$I\circ \hind^{top}=\ind^{top}\circ \cF\ ,\quad I\circ \hind^{an}=\ind^{an}\circ \cF$$
  and Theorem \ref{dww}, where  
  $$\cF:\bloc_{geom}(R)\to \bloc(R)$$ is the map of symmetric monoidal stacks   which forgets the geometry.
  \hB
  
  \begin{rem}{\rm
The difference $\hat \delta:=\hind^{top}-\hind^{an}$ is therefore an additive natural transformation
$$\bloc_{geom}(R)(W)\ni (\cV,h^{\cV})\mapsto\hat \delta(\cV,h^{\cV})\in \frac{HA^{-1}(B)}{\im(KR^{-1}(B)\to   HA^{-1}(B))}\ .$$
It actually factors over a transformation
$$\delta:\bloc(R)\to  \frac{HA^{-1}(B)}{\im(KR^{-1}(B)\to   HA^{-1}(B))}\ ,$$ 
of $\Bundle$-{commutative monoids}.
The Conjecture \ref{may3103} can now be reformulated to say that
$$\delta=0\ .$$
One possible approach to a proof of Conjecture \ref{may3103} could be to show that such a transformation is necessarily zero. This approach works successfully in similar situations, e.g. for a
differential refinement of the index theorem of Atiyah-Singer in \cite{MR2664467}. The basic reason there was that the
rational cohomology of the classifying space of complex vector bundles is concentrated in even degrees.
But in the present case of algebraic $K$-theory of a number ring the situation is quite different since the  rational cohomology of the corresponding classifying space is generated as a ring by odd-degree classes.

}\hB\end{rem}

\subsection{The analytic index}\label{jun183}

 
In this subsection we provide the details of the construction of the analytic index.

\bigskip

{Let $R$ be a ring.}
For a manifold $W$ we let  $\Sh_{\Mod(R)}(W)$ denote the abelian categeory of sheaves of $R$-modules on the site of open subsets of $W$.    The proper submersion $\pi:W\to B$ induces a push-forward
$\pi_{*}:\Sh_{\Mod(R)}(W)\to \Sh_{\Mod(R)}(B)$. For a sheaf of $R$-modules $\cV\in \Sh_{\Mod(R)}(W)$ we have the  sequence of higher-derived images
$$R^{i}\pi_{*}(\cV)\in \Sh_{\Mod(R)}(B)\ ,\quad i\ge 0\ .$$

{
\begin{lem}\label{kjhdqkwdqwdqwdqwdqwdqwdwqd} If $\cV$ is a locally constant sheaf of finitely generated $R$-modules and $R$ is noetherian, then
for every $i\in \nat$, the sheaf $R^{i}\pi_{*}(\cV)$ is again a locally constant sheaf of finitely generated   $R$-modules.\end{lem}\proof}
We consider a point  $b\in B$ and let $W_{b}:=\pi^{-1}(b)$ denote the fibre of $\pi$ over $b$. For a sheaf $\cF\in \Sh_{\Mod(R)}(B)$ we let $\cF_{b}$ denote the stalk of $\cF$ at $b$.
Then we have a canonical isomorphism of $R$-modules \begin{equation}\label{hhjfewkfhewfooeewfwfwefewfe}
R^{i}\pi_{*}(\cV)_{b}\cong H^{*}(W_{b},\cV_{|W_{b}})\ .
\end{equation}
We now assume that $\cV$ is locally constant. {Since the bundle $W\to B$ is locally trivial and because of the homotopy invariance of the cohomology of a manifold with coefficients in a locally constant sheaf}
the sheaves $R^{i}\pi_{*}(\cV)$ are again locally constant sheaves of $R$-modules.
If $\cV$ is a locally constant sheaf of finitely generated $R$-modules, then
 so are  the higher derived images  $R^{i}\pi_{*}(\cV)$. Indeed, since $W_{b}$ is compact, the cohomology group \eqref{hhjfewkfhewfooeewfwfwefewfe}
can be calculated in the \v{C}ech picture using a finite good covering of $W_{b}$. It follows that
the cohomology is a subquotient of a finitely generated $R$-module. Then one uses the fact that $R$ is noetherian. 
\hB 

{From now on $R$ is again a number ring. It is noetherian so that Lemma \ref{kjhdqkwdqwdqwdqwdqwdqwdwqd} applies.} 
We  define the analytic index 
$$\ind^{an}:\bloc(R)\: d\to KR^{0}r$$ by the formula \eqref{may3116}.
 
 \bigskip
 
 Following the ideas of \cite{MR1724894} we now refine this construction to an analytic index in differential algebraic $K$-theory 
 $$\hind^{an}:\bloc_{geom}(R)d\to \widehat{KR}^{0}r\ .$$
   {For  an embedding $\sigma:R\to \C$   we let    
$H^{i}(W/B,V_{\sigma})\to B$ denote flat complex vector bundle} associated to the locally constant sheaf 
$R^{i}\pi_{*}(\cV)\otimes_{\sigma} \C$ of finite-dimensional complex vector spaces.
As the notation indicates, the fibre
$H^{i}(W/B,V_{\sigma})_{b}$  of this bundle at $b\in B$ can canonically be identified with
the cohomology $H^{i}(W_{b},\cV\otimes_{\sigma}\C)$ and hence with the $i$'th cohomology of the 
Rham complex $\Omega (W_{b},V_{\sigma|W_{b}})$ of $W_{b}$ twisted by the flat bundle
$V_{\sigma|W_{b}}$ associated to the local system $\cV\otimes_{\sigma} \C$. 
{Here we use the fact that $\C$ is a flat $R$-module in order to commute the operations of taking cohomology and
 tensoring with $\C$.}
 
  {The following constructions are similar as in  \cite[Sec.3(d)]{MR1303026}}. The Riemannian structure $g$ on $\pi$ and the geometry $h^{\cV}$ on $\cV$ supply  a Riemannian metric on $W_{b}$ and a metric on $V_{b}$. These structures induce an $L^{2}$-scalar product on the de Rham complex 
$\Omega (W_{b},V_{\sigma|W_{b}})$.
 We can now apply Hodge theory
  which provides an isomorphism of complex vector spaces
$$H^{i}(W/B,V_{\sigma})_{b}\cong \cH(W_{b},V_{\sigma|W_{b}})\ ,$$
where the right-hand side denotes the subspace $ \cH^{i}(W_{b},V_{\sigma|W_{b}})\subseteq \Omega^{i} (W_{b},V_{\sigma|W_{b}}) $ of harmonic $i$-forms on $W_{b}$ with coefficients in $V_{\sigma|W_{b}}$. The $L^{2}$-metric on $\cH^{i}(W_{b},V_{\sigma|W_{b}})$ thus induces a metric $h_{L^{2}}^{H^{i}(W/B,V_{\sigma})}(b)$ 
on $H^{i}(W/B,V_{\sigma})_{b}$. This metric depends smoothly on the base point and is therefore a hermitean metric on the bundle  $H^{i}(W/B,V_{\sigma})$.
The collection of metrics \begin{equation}\label{apr140113}h^{R^{i}\pi_{*}(\cV)}:=\left(h_{L^{2}}^{H^{i}(W/B,V_{\sigma})}\right)_{\sigma\in \Sigma}\end{equation}
 serves as a geometry (Definition \ref{defpre}) on $R^{i}\pi_{*}(\cV)$.
 Indeed, the condition (\ref{jun101}) for $h^{R^{i}\pi_{*}(\cV)}$ is easily implied by the corresponding condition for $h^{\cV}$.
 
As a first approximation {of the differential analytic index}, we define 
\begin{equation} \label{jun219}\hind_{0}^{an}(\cV,h^{\cV}):=\sum_{i}(-1)^{i}\cycl\left(R^{i}\pi_{*}(\cV),h^{R^{i}\pi_{*}(\cV)}\right) \\ \in \widehat{KR}^{0}(B) \ ,\end{equation}
{where $\cycl$  is as in    Definition  \ref{uliapr1210}. }
The subscript ''$0$'' indicates that this definition is not yet the final one. It has to be corrected  by the Bismut-Lott analytic torsion form $\cT(\pi,g,\cV,h^{\cV})\in \Omega A^{-1}(B)$ (see \eqref{jun104})
   in order to match the curvature forms. 
 
 \bigskip
 
  For each embedding $\sigma:R\to \C$ we have a Bismut-Lott analytic torsion form
 \cite[Def. 3.22]{MR1303026}
 $$\cT(T^{h}\pi,g^{T^{v}\pi},h^{V_{\sigma}})\in \Omega^{{ev}}(B)\ .$$
 By  \cite[Thm. 3.23]{MR1303026} it satisfies 
 \begin{equation}\label{jun103}d\cT(T^{h}\pi,g^{T^{v}\pi},h^{V_{\sigma}})=\int_{W/B} e(g)\wedge f(\nabla^{V_{\sigma}},h^{V_{\sigma}})-\sum_{i}(-1)^{i}f(\nabla^{H^{i}(W/B,V_{\sigma})},h_{L^{2}}^{H^{i}(W/B,V_{\sigma})})\ ,\end{equation}
 where the forms 
 $f(\nabla^{V_{\sigma}},h^{V_{\sigma}})\in \Omega(W)$ can be expressed in terms of {the  Kamber-Tondeur}  forms (\ref{jun102}) {as follows} (compare \cite[(3.102)]{MR1303026}):
 $$f(\nabla^{V},h^{V})=\sum_{k\ge 0} \frac{1}{k!} \omega_{2k+1}({V,\nabla^{V},h^{V}})\ . $$
 We consider the decomposition
 $$\cT(T^{h}\pi,g^{T^{v}\pi},h^{V_{\sigma}})=\sum_{k\ge 0} \cT_{2k}(T^{h}\pi,g^{T^{v}\pi},h^{V_{\sigma}})\ ,\quad \cT_{2k}(T^{h}\pi,g^{T^{v}\pi},h^{V_{\sigma}})\in \Omega^{2k}(B)$$
 into homogeneous components and
 define
 \begin{equation}\label{jun104}\cT(\pi,g,\cV,h^{\cV}):=\sum_{\sigma\in \Sigma}\sum_{k\ge 0} \frac{1}{k!}\: \cT_{2k}(T^{h}\pi,g^{T^{v}\pi},h^{V_{\sigma}}) \delta_{2k+1,\sigma}-\sum_{\sigma\in \Sigma}  \cT_{0}(T^{h}\pi,g^{T^{v}\pi},h^{V_{\sigma}})\kappa_{1} \in \Omega A^{-1}(B)\ ,\end{equation}
 {where the basis elements {$\delta_{2k+1,\sigma}\in \tilde A(R)^{-2k-1}$ were defined after Definition \ref{kjqwldqwdqwdwqdqwdqd}. We must add the term involving    $\kappa_{1}\in \tilde A(R)^{-1}$  given by \eqref{gdqjwhqjwduiiuiuiuiui} in order to ensure that  $\cT(\pi,g,\cV,h^{\cV})$ belongs to the subspace $\Omega A^{-1}(B)\subseteq \Omega \tilde A(R)^{-1}(B)$.}
{Using the forms introduced in Definition \ref{defpre1}  we can} now rewrite (\ref{jun103}) as follows:
 \begin{equation}\label{mar0101}d\cT(\pi,g,\cV,h^{\cV})=\int_{W/B} e(g)\wedge \omega({\cV,}h^{\cV})-\sum_{i}(-1)^{i}\omega({R^{i}\pi_{*}(\cV),} h^{R^{i}\pi_{*}(\cV)})\ .\end{equation}
  The final definition of the analytic index in differential algebraic $K$-theory will be the following.
\begin{ddd}\label{jun271}
If $W\to B$ is a proper submersion equipped with a geometry $g$ and $(\cV,h^{\cV})$ is a locally constant sheaf of finitely generated projective $R$-modules on $W$ {with geometry $h^{\cV}$}, then
we define
$$\hind^{an}(\cV,h^{\cV}):=\hind_{0}^{an}(\cV,h^{\cV})+a(\cT(\pi,g,\cV,h^{\cV}))\ .$$
\end{ddd}

\proof [of Theorem \ref{bl}]
We verify Theorem \ref{bl} by a curvature calculation.
 By  (\ref{m2159}) we have
$$R(\hind^{an}_{0}(\cV,h^{\cV}))=\sum_{i}(-1)^{i}\omega( h^{R^{i}\pi_{*}(\cV)} )\ .$$
From \eqref{mar0101} we therefore get  
\begin{equation}\label{may121310}R(\hind^{an}(\cV,h^{\cV}))=\int_{W/B}e(g)\wedge \omega(h^{\cV})\ .\end{equation}
Equation \eqref{may121310} implies Theorem \ref{bl} since
$$R(\hind^{top}(\cV,h^{\cV}))\stackrel{\eqref{may3103222}}{=} \int_{W/B}e(g)\wedge R(\cycl(\cV,h^{\cV}))\stackrel{\eqref{m2159}}{=}
\int_{W/B} e(g)\wedge   \omega(h^{\cV})\ .$$ \hB

The construction of the metric $h^{R^{i}\pi_{*}(\cV)}$
and the torsion form $\cT(\pi,g,\cV,h^{\cV})$ is compatible with pull-backs along morphisms in $\Bundle_{geom}$.  Hence our construction of the refined analytic index gives a map of $\Bundle_{geom}$-{commutative monoids}
$$\hind^{an}:\bloc_{geom}d\to \widehat{KR}^{0}r\ .$$

\subsection{Discussion of the transfer index conjecture}\label{sep0901}

\subsubsection{The relation with the work of Lott}\label{jul0481}

{The main result of the present paragraph is Proposition \ref{hkjh32kjhrk23r32r32r2r} which states that the Transfer Index Conjecture \ref{may3103}  implies a {version} of a conjecture of Lott \cite[Conj. 1]{MR1724894} in the case of  number rings. In \cite{MR1724894}  Lott introduces, for suitable rings $R$, secondary $K$-groups
of locally constant sheaves of finite-generated $R$-modules and constructs an analytic push-forward along proper submersions with closed fibres. He conjectures that this analytic push-forward is compatible with the Becker-Gottlieb transfer for the generalized cohomology theory represented by the fibre of a suitable regulator. 

In the present paragraph  we specialize to a number ring $R$. While Lott considers one embedding $R\hookrightarrow \C$ at a time here we will consider all embedding at {once}. We start {by} recalling Lott's definition of the secondary $K$-groups as a motivation for the definition of our version. We then state the corresponding version of Lott's Conjecture \ref{jkcnjaskcascaksacsacsacac}.}

\bigskip

Let $R$ be a number ring. 
In \cite{MR1724894} Lott constructed contravariant functors
$$\Mf^{op}\ni M\mapsto \bar K^{0}_{Lott,\sigma}(M)\in \Ab\ , \quad \Mf^{op}\ni M\mapsto \hat K^{0}_{Lott,\sigma}(M)\in \Ab$$
by cycles and relations. They depend on the choice of an embedding $\sigma:R\to \C$ which we added to the notation.
A cycle (\cite[Def. 2]{MR1724894}) for $\hat K^{0}_{Lott,\sigma}(M)$ is a triple $(\cV,h^{V_{\sigma}},\eta)$
of a locally constant sheaf {$\cV$} of finitely generated    $R$-modules on $M$, a metric $h^{V_{\sigma}}$ on the complex vector bundle $V_{\sigma}\to M$ associated to $\cV$ and $\sigma$, and a form $\eta\in \Omega^{ev}(M)/\im(d)$. The set of isomorphism classes of cycles becomes a monoid with respect to direct sum. The
group $\hat K^{0}_{Lott,\sigma}(M)$ is then defined by group completion and factorization 
 {modulo} the equivalence relation (\cite[Def. 3]{MR1724894}) generated by the following. If $(\cV_{i},h^{V_{i,\sigma}},\eta_{i})$, $i=0,1,2$, are three cycles and 
\begin{equation}\label{feb2810}0\to \cV_{0}\to \cV_{1}\to \cV_{2}\to 0\end{equation} is an exact sequence of sheaves, then
$$(\cV_{0},h^{V_{0,\sigma}},\eta_{0})+(\cV_{2},h^{V_{2,\sigma}},\eta_{2})\sim (\cV_{1},h^{V_{1,\sigma}},\eta_{1})$$ if
$$\eta_{0}+\eta_{2}-\eta_{1}=\cT_{\sigma}\in \Omega^{ev}(M)/\im(d)\ ,$$
were $\cT_{\sigma}$ is the analytic torsion form associated to the exact sequence \eqref{feb2810} equipped with metrics $h^{V_{i,\sigma}}$, see
\cite[A.2]{MR1724894} for details. {For every $k\in \nat$ its component  in degree $2k$ in particular satisfies the relation  \begin{equation}\label{fhwefhwekw298347234324234242}
d\cT_{\sigma,2k}=\omega_{2k+1}(\bV_{0,\sigma})+\omega_{2k+1}(\bV_{2,\sigma})-\omega_{2k+1}(\bV_{1,\sigma})\ ,
\end{equation}
  where $\bV_{i,\sigma}:=(V_{i,\sigma},\nabla^{V_{i,\sigma}},h^{V_{i,\sigma}})$.}
We let $[\cV,h^{V_{\sigma}},\eta]\in \hat K^{0}_{Lott,\sigma}(M)$ denote the class represented by the cycle $(\cV,h^{V_{\sigma}},\eta)$.

\bigskip

There are natural transformations
$$b:\hat K_{Lott,\sigma}(M)\to KR^{0}(M)\ ,\quad c:\hat K_{Lott,\sigma}(M)\to Z^{odd}(\Omega(M))$$
induced by 
$$(\cV,h^{V_{\sigma}},\eta)\mapsto \hat I(\cV)\ ,\quad (\cV,h^{V_{\sigma}},\eta)\mapsto \sum_{k\ge 0} \omega_{2k+1}(\bV_{\sigma}) -d\eta\ .$$
 See Definition \ref{may121311} for $\hat I$   and  \eqref{jun102}} for the definition of the forms $\omega_{2k+1}({\bV_\sigma})$. Furthermore, we have a
transformation
$$a:\Omega^{ev}(M)/\im(d)\to \hat K_{Lott,\sigma}(M)\ , \quad \eta\mapsto (0,0,-\eta)\ .$$
{Following \cite[Def. 5]{MR1724894} we define the}
functor $\bar K^{0}_{Lott,\sigma}$ by
$$\bar K^{0}_{Lott,\sigma}(M):=\ker\left(c: \hat K_{Lott,\sigma}(M)\to  Z^{odd}(\Omega(M))\right)\ ,$$
{which turns out to homotopy invariant}.

\begin{rem}{\rm
The two functors $\bar K^{0}_{Lott,\sigma}$ and $\hat K_{Lott,\sigma}$ bear some similarities with 
$\widehat{KR}^{0}_{flat}$ and $\widehat{KR}^{0}$.  As mentioned above, in order to set up a precise comparison we shall modify Lott's definition in order to get rid of the dependence {on}  the embedding $\sigma$. The idea is to consider all embeddings at once.
}\hB \end{rem}

 We define
the modified group $\widehat{KR}_{Lott}^{0}(M)$ again by cycles and relations as above, where now a cycle is
a triple $(\cV,h^{\cV},\eta)$ of a locally constant sheaf $\cV$ of finitely generated $R$-modules with geometry $h^{\cV}$ on $M$ and a form $\eta\in \Omega A^{-1}(M)$. This data is the same as a collection of cycles
$(\cV,h^{V_{\sigma}},\eta_{\sigma})_{\sigma\in \Sigma}$ with the fixed underlying  locally constant sheaf $\cV$ such that the relations  $\bar h^{V_{\sigma}}=h^{V_{\bar \sigma}}$, $\eta_{\sigma,2j+1}=(-1)^{j}\eta_{\bar \sigma,2j+1}$ and $\sum_{\sigma\in \Sigma} \eta_{\sigma,0}=0$ are satisfied. See Definition \ref{defpre} for the first and  Definition \ref{kjqwldqwdqwdwqdqwdqd} for the explanation of the latter two conditions.


 The equivalence relation for $\widehat{KR}_{Lott}^{0}(M)$ is generated as follows. If $(\cV_{i},h^{\cV_{i}},\eta_{i})$, $i=0,1,2$, are three cycles and 
$$0\to \cV_{0}\to \cV_{1}\to \cV_{2}\to 0$$ is an exact sequence of sheaves, then
$$(\cV_{0},h^{\cV_{0}},\eta_{0})+(\cV_{2},h^{\cV_{2}},\eta_{2})\sim (\cV_{1},h^{\cV_{1}},\eta_{1})$$ if
\begin{equation}\label{sep1601}\eta_{0}+\eta_{2}-\eta_{1}=\cT\in \Omega A^{odd}(M)/\im(d)\ .\end{equation}
The torsion form $\cT$ on the right-hand side of \eqref{sep1601} is given by   
\begin{equation}\label{may121003}\cT=\sum_{\sigma\in \Sigma}   \cT_{\sigma,2k} \: \bar \delta_{2k+1,\sigma} \ ,\end{equation}
where the components $\cT_{\sigma,2k}\in \Omega^{2k}(M)/\im(d)$ are determined by
$$\cT_\sigma=\sum_{k\ge 0} \cT_{\sigma,2k}\ ,$$
 and we set, for simplicity,  \begin{equation}\label{wkejfeljlwelkfjewlfjewlkfjewfkljewlfke902irp2ri2p3r32r32r32r23r}
\bar \delta_{2j+1,\sigma}:=\delta_{2j+1,\sigma}
\end{equation} for $j\ge 1$. 
{As a consequence of \eqref{fhwefhwekw298347234324234242} this version of the} torsion form satisfies the relation
\begin{equation}\label{mar0102}d\cT=\omega(\cV_{0},h^{\cV_{0}})+\omega(\cV_{2},h^{\cV_{2}})-\omega(\cV_{1},h^{\cV_{1}})\ .\end{equation}

As above there are natural transformations
$$b:\widehat{KR}^{0}_{Lott}(M)\to KR^{0}(M)\ ,\quad c:\widehat{KR}^{0}_{Lott}(M)\to Z^{0}(\Omega A(M))$$
induced by 
$$(\cV,h^{\cV},\eta)\mapsto \hat I(\cV)\ ,\quad (\cV,h^{\cV},\eta)\mapsto  \omega(\cV,h^{\cV}) -d\eta\ ,$$
{(see Definition \ref{jul1702} for $\hat I$ {and Definition \ref{defpre1} for $\omega(\cV,h^{\cV})$)}} and a transformation
$$a:\Omega A^{-1}(M)/\im(d)\to \widehat{KR}^{0}_{Lott}(M)\ , \quad \eta\mapsto (0,0,-\eta)\ .$$
The 
functor $M\mapsto \overline{ KR}^{0}_{Lott}(M)$ is defined as
$$\overline{KR}^{0}_{Lott}(M):=\ker\left(c: \widehat{ KR}^{0}_{Lott}(M)\to  Z^{0}(\Omega A(M))\right)\ .$$
{It is again  homotopy invariant.}

\bigskip

A cycle $(\cV,h^{\cV},\eta)$ gives rise to a class
\begin{equation}\label{jlwekjflkwefjklewfewfwf243243423423423432432}
\zz(\cV,h^{\cV},\eta):=\cycl(\cV,h^{\cV})+a(\eta)\in \widehat{KR}^{0}(M)
\end{equation}
such that
$$R(\zz(\cV,h^{\cV},\eta))=c([\cV,h^{\cV},\eta])\ ,\quad b(\zz(\cV,h^{\cV},\eta))=I([\cV,h^{\cV},\eta])\ .$$
The question is now whether this map $\zz$ factorizes over
$\widehat{KR}_{Lott}^{0}(M)$. To this end we must know that the following relation holds true in  $\widehat{KR}^{0}(M)$. 
\begin{theorem}[{Lott's relation}]\label{jun0301}
 If $(\cV_{i},h^{\cV_{i}})$,  $i=0,1,2$, are locally constant sheaves of finitely generated $R$-modules on a manifold $M$ with geometries and 
\begin{equation}\label{jul2021}0\to \cV_{0}\to \cV_{1}\to \cV_{2}\to 0\end{equation} is an exact sequence of sheaves, then we have the following relation in $\widehat{KR}^{0}(M)$:
\begin{equation}\label{jul2023}
 \cycl(\cV_{0},h^{\cV_{0}})+\cycl(\cV_{2},h^{\cV_{2}})-\cycl(\cV_{1},h^{\cV_{1}})=a(\cT)\ .\end{equation}
\end{theorem}

\proof
{A proof of this theorem for exact sequence of sheaves of finitely generated projective $R$-modules appears in  \cite{buta}.  In Lemma \ref{uliapr1501} we will show that one can drop the assumption of projectivity. \hB}

\begin{rem}{\rm
 The proof given in \cite{buta}   depends upon a generalization of the theory presented in the present paper. {Therefore, in Subsection \ref{aug0501}, we in addition}  verify some special cases and weaker statements which can easily be obtained in the framework of the present paper.  
 }\hB \end{rem}

{Lott's relation  Theorem \ref{jun0301} has the following consequence.
\begin{kor}\label{hdkjwhkjwqddwqdqwdwqdwqdq}
The map $\zz$ given by  \eqref{jlwekjflkwefjklewfewfwf243243423423423432432}  descends to a natural transformation
$$\phi:\widehat{KR}_{Lott}^{0}(M)\to \widehat{KR}^{0}(M)$$
such that the following diagram commutes
\begin{equation}\label{m2159hjkfhkjkjhhjkwef}\xymatrix{&& KR^{0}(M)\\ \widehat{KR}_{Lott}^{0}(M)\ar[urr]^{b}\ar[drr]_{c}\ar[r]^{\phi}&\widehat{KR}^{0}(M)\ar[ur]_{I}\ar[dr]^{R}\\
&&Z^{0}(\Omega A(M))} \ .
\end{equation}
\end{kor}
}

{\begin{rem}
{\rm Note that we de not expect that the map   
$$\phi:\widehat{KR}^{0}_{Lott}(M)\to \widehat{KR}^{0}(M)$$ given by Corollary \ref{hdkjwhkjwqddwqdqwdwqdwqdq} is an isomorphism since in general not all
classes in $KR^{0}(M)$ are represented by locally constant sheaves of finitely generated  $R$-modules on $M$ (consider e.g. the case where $M$ is simply connected!).
}\hB 
\end{rem}
}

{In view of the lower triangle in \eqref{m2159hjkfhkjkjhhjkwef} the map   $\phi$ maps the kernel of $c$ to the kernel of $R$. It hence restricts to  a transformation
\begin{equation}\label{hjhjfkwefewfwfwfwffef24}
\bar \phi:\overline{KR}^{0}_{Lott}(M)\to \widehat{KR}^{0}_{flat}(M)\xrightarrow{\cong\ Prop.\ \ref{hgr}} KR\R/\Z^{-1}(M)\ .
\end{equation}

\begin{rem}{\rm The construction of this map $\bar \phi$ solves a part of the conjecture \cite[Conj. 1]{MR1724894}. Namely it provides the map which in Lott's notation reads  $\overline{K}^{0}_{R}(M)\to \mathcal{F}^{0}_{R}(M)$.
}\hB \end{rem}
}

{One of the main results of \cite{MR1724894}  is the construction of the analytic index for $\widehat{KR}^{0}_{Lott}$ for proper submersions with a Riemannian structure.
In the following we present the  adaption of this  construction to our groups $\widehat{KR}_{Lott}^{0}$.
Let} $W\to B$ be a proper submersion with closed fibres and  with a  Riemannian structure $g$ as in Definition \ref{may272}.
Given a cycle
$(\cV,h^{\cV},\eta)$ for $\widehat{KR}_{Lott}^{0}(W)$ we set (compare \cite[Def. 7]{MR1724894})
\begin{eqnarray}\hind^{an}_{Lott}(\cV,h^{\cV},\eta)&:=&\sum_{i}(-1)^{i}[R^{i}\pi_{*}(\cV),h^{R^{i}\pi_{*}(\cV)},0] \label{klkjqkldjqwldkqjwlkjqwdjwqldjqwlkdkqwjdqwlkdjqwlkdwqdqwde3}\\
&&+[0,0,\int_{W/B}e(b)\wedge \eta-\cT(\pi,g,\cV,h^{\cV})]\nonumber\end{eqnarray}
{The same argument as for  \cite[Prop. 7]{MR1724894} shows that 
$$\hind^{an}_{Lott}:\widehat{KR}_{Lott}^{0}(W)\to \widehat{KR}^{0}_{Lott}(B)$$
is well-defined. {Moreover one checks  that $\hind^{an}_{Lott}$ restricts to a map
 $$\overline{\ind}_{Lott}^{an}:\overline{KR}^{0}_{Lott}(W)\to \overline{KR}^{0}_{Lott}(B)\ ,$$ and that this restriction is independent of the choice of the Riemannian structure. 
  }
}

\bigskip

{We can now formulate the remaining part of Lott's conjecture  \cite[Conj. 1]{MR1724894} (namely the commutativity of \cite[(85)]{MR1724894}) in a precise form.

\begin{con}[Lott] \label{jkcnjaskcascaksacsacsacac}
The following diagram commutes: \begin{equation}\label{wefjhejfhwekfewfewfewfewf23453434324}
\xymatrix{\overline{KR}^{0}_{Lott}(W)\ar[d]^{\overline{\ind}_{Lott}^{an}}\ar[r]^{\bar \phi}& KR\R/\Z^{-1}(W)\ar[d]_{\eqref{jun182}}^{\tr^{*}}\\
\overline{KR}^{0}_{Lott}(B)\ar[r]^{\bar \phi}&KR\R/\Z^{-1}(B)}\ .
\end{equation}
\end{con}
}

{
\begin{prop}\label{hkjh32kjhrk23r32r32r2r}
The Transfer Index Conjecture \ref{may3103} implies Lott's Conjecture \ref{jkcnjaskcascaksacsacsacac}.
\end{prop}
\proof
%
It immediately follows from the constructions of $\phi$ (in terms of $\zz$ given by \eqref{jlwekjflkwefjklewfewfwf243243423423423432432}),  $\hind^{an}_{Lott}$ (by \eqref{klkjqkldjqwldkqjwlkjqwdjwqldjqwlkdkqwjdqwlkdjqwlkdwqdqwde3}) and $\hind^{an}$ (Definition \ref{jun271}) 
that following diagram   commmutes
$$\xymatrix{\widehat{KR}^{0}_{Lott}(W) \ar[d]^{\hind^{an}_{Lott}}\ar[r]^{\phi}&\widehat{KR}^{0}(W)\ar[d]^{\hind^{an}}\\\widehat{KR}_{Lott}^{0}(B)\ar[r]^{\phi}&\widehat{KR}^{0}(B)}\ .$$
We now restrict to the flat part  and extend the diagram to
\begin{equation}\label{jun0307}\xymatrix{\overline{KR}^{0}_{Lott}(W)\ar@/^1cm/[rrr]^{\bar \phi}\ar[d]^{\overline{\ind}^{an}_{Lott}}\ar[r]^{\phi}&\widehat{KR}_{flat}^{0}(W)\ar[d]^{\hind^{an}}\ar@{=}[r]&\widehat{KR}_{flat}^{0}(W)\ar[d]^{\hind^{top}}\ar[r]^{\cong}&KR\R/\Z^{-1}(W)\ar[d]^{\tr^{*}}\\\overline{ KR}_{Lott}^{0}(B)\ar@/_1cm/[rrr]^{\bar \phi}\ar[r]^{\phi}&\widehat{KR}_{flat}^{0}(B)\ar@{=}[r]&\widehat{KR}^{0}_{flat}(B)\ar[r]^{\cong}&KR\R/\Z^{-1}(B)}\ .\end{equation}
The commutativity of the right square is  \eqref{jun0304}.
The Transfer Index Conjecture \ref{may3103}  implies that the middle square commutes.
The upper and lower horizontal compositions are  the definition \eqref{hjhjfkwefewfwfwfwffef24} of $\bar \phi$. 
Thus the outer square is exactly \eqref{wefjhejfhwekfewfewfewfewf23453434324}. \hB 
}

\subsubsection{The relation with the Cheeger-M\"uller theorem}\label{may121280}

{The goal of the present subsection is to support the  Transfer Index Conjecture through its compatibility with the Cheeger-M\"uller theorem {about the equality of Ray-Singer analytic torsion and Reidmeister torsion.}}

\bigskip

 Let $\pi:W\to B$ be a proper submersion with a Riemannian structure $g$ (Definition \ref{may272}).
We consider a locally constant sheaf    $\cV$  of finitely generated free (for simplicity) $R$-modules on $W$.
We further assume that we can choose the geometry $h^{\cV}$ (Definition \ref{defpre}) such that the metrics $h^{V_{\sigma}}$ are parallel for all places $\sigma\in \Sigma$ of $R$. {By  \eqref{cfr} and Proposition \ref{feb1704}, 4. the characterstic form of $(\cV,h^{\cV})$ is then given by  $$\omega(\cV,h^{\cV})=\kappa_{0}\dim(\cV)\ ,$$ where $\kappa_{0}$ is defined in   \eqref{gdqjwhqjwduiiuiuiuiui}.} 
{We cancel this curvature by substracting a multiple of the unit   class. To this end,}
using the projection $q:W\to *$, we  define \begin{equation}\label{fewjlkewfklelklkklwefwfewfw2355tgrg}
(\underline{R}_{W},h^{\underline{R}_{W}}) :=q^{*}(R,h^{R})\in  \bloc^{proj}(R)(W)\ , \quad \beins_{W}:=\cycl(\underline{R}_{W},h^{\underline{R}_{W}})\in \widehat{KR}^{0}(W)\ .
\end{equation} Then we have 
  $\beins_{W}:=q^{*}\beins\ $, where $\beins$ was defined in (\ref{jun201}).
We get the flat differential algebraic $K$-theory class
$$\cycl(\cV,h^{\cV})-\dim(\cV) \ \beins_{W} \in \widehat{KR}_{flat}^{0}(W)\cong KR\R/\Z^{-1}(W)\ .
$$
By   \eqref{jun0304} the restriction of the differential cohomology transfer $\hat \tr$ of the bundle $\pi$ to the flat subfunctor $\widehat{KR}^{0}_{flat}$ conincides with the Becker-Gottlieb transfer $\tr^{*}$ for the cohomology functor $KR\R/\Z^{-1}$. The latter is homotopy invariant.
 Hence
the element
$$\tr^{*}(\cycl(\cV,h^{\cV})-\dim(\cV)\  \beins_{W})\stackrel{ \eqref{jun0304}}{=}\hat \tr(\cycl(\cV,h^{\cV})-\dim(\cV)\ \beins_{W} )\in 
 \widehat{KR}_{flat}^{0}(B)$$
 does {not} depend on the Riemannian structure $g$.
 
 We now assume that $B$ is a point and that $W$ is connected. In this case the Becker-Gottlieb transfer can be calculated using a base point $s:*\to W$ as
\begin{equation}\label{regegggge34535566577}
\tr^{*}=\chi(W)\  s^{*}\ . 
\end{equation} 
{We use the Definition \ref{jun181} of the topological index in order to get the first equality of the following chain }   
 \begin{eqnarray}\label{jun189} \lefteqn{\hind^{top}(\cV,h^{\cV})-\hind^{top}(\underline{R}_{W},h^{\underline{R}_{W}})}\hspace{3cm}\nonumber\\&=&\hat \tr (  \cycl(\cV,h^{\cV})-\dim(\cV)  \beins_{W})  \nonumber \\&\stackrel{\eqref{regegggge34535566577}}{=}&\chi(W)\left(\cycl(s^{*}(\cV,h^{\cV})) -\dim(\cV)\ \beins\right)\label{fkjewfklewf9o2990r23u2309ufdoewfw}\\
 &\stackrel{Lemma \:\ref{jun1820}}{=}& \chi(W)  \left(\cycl(\det(s^{*}\cV),h^{\det(s^{*}\cV)})-\beins\right) \nonumber\ .\end{eqnarray}

We now consider the analytic index. For every $i\ge 0$ the cohomology group $H^{i}(W,\cV)$
of the compact  manifold $W$ with coefficients in the sheaf $\cV$ is a finitely generated $R$-module.
 Moreover,  we have a collection of metrics
$h^{H^{i}(W,\cV)}:=(h_{L^{2}}^{H^{i}(W,V_{\sigma})})_{\sigma\in \Sigma}$  induced by Hodge theory.
We set $h^{\cH(\cV)}:=(h^{H^{i}(W,\cV)})_{i=0,\dots,n}$.
We have by Definition \ref{jun271} {of the {analytic} index}
\begin{eqnarray*}\hind^{an}(\cV,h^{\cV})&=&\sum_{i\ge 0} (-1)^{i}\cycl\left(H^{i}(W,\cV),h^{H^{i}(W,\cV)}\right) \\
&& +a\left(\cT(\pi,g,\cV,h^{\cV})\right) 
\ .\end{eqnarray*}
{According to    \cite[Thm. 3.29]{MR1303026}, for $B$ a point,}   the   analytic torsion form is related with the logarithm of the classical Ray-Singer torsion  \cite{MR0295381}
$T^{an}(W,g^{TW}V_{\sigma},h^{V_{\sigma}})\in \R $ of the flat bundle $(V_{\sigma},\nabla^{V_{\sigma}},h^{V_{\sigma}})$ on the Riemannian manifold $(W,g^{TW})$ by
$$\cT(\pi,g,\cV,h^{\cV})=  T^{an}(W,g^{TW},\cV,h^{\cV}):= \sum_{\sigma}  T^{an}(W,g^{TW},V_{\sigma},h^{V_{\sigma}})\: \bar \delta_{1,\sigma}\in A^{-1}\ .$$ 
Observe that $\cT(\pi,g, \underline{R}_{W},h^{ \underline{R}_{W}})=0$ because of the usage of the projected basis vectors  $\bar \delta_{1,\sigma}$.
{We abbreviate
$$ D(\cV,h^{\cV}):=\sum_{i\ge 0} (-1)^{i}\cycl\left(H^{i}(W,\cV),h^{H^{i}(W,\cV)}\right)\ .$$} 
Then we have
\begin{eqnarray}
\lefteqn{\hind^{an}(\cV,h^{\cV})-\hind^{an}(\underline{R}_{W},h^{\underline{R}_{W}})}\hspace{3cm}\nonumber\\&=&D(\cV,h^{\cV})-\dim(\cV)\:D(\underline{R}_{W},h^{\underline{R}_{W}})\label{lkejdkljqlwkdqqwdddqwdqwdqwdqwdwqd}\\&&\hspace{-1cm}+a \left( T^{an}(W,g^{TW},\cV,h^{\cV})\right)- \dim(\cV)\ a\left( T^{an}(W,g^{TW},\underline{R}_{W},h^{\underline{R}_{W}})\right) \nonumber \\&=&D(\cV,h^{\cV})-\dim(\cV)\:D(\underline{R}_{W},h^{\underline{R}_{W}})+a \left( T^{an}(W,g^{TW},\cV,h^{\cV})\right)
\nonumber
\end{eqnarray}

 The Transfer Index Conjecture \ref{may3103} now  predicts that the   right-hand sides of 
 \eqref{fkjewfklewf9o2990r23u2309ufdoewfw} and \eqref{lkejdkljqlwkdqqwdddqwdqwdqwdqwdwqd} are equal, i.e. the equality
{
\begin{eqnarray}  \lefteqn{\chi(W)\left[\cycl\left(\det(s^{*}\cV),h^{\det(s^{*}\cV)}\right) -\beins\right]}&&\nonumber\\&=&D(\cV,h^{\cV})-\dim(\cV) \:D(\underline{R}_{W},h^{\underline{R}_{W}})  + a\left( T^{an}(W,g^{TW},\cV,h^{\cV})\right) \ .\label{jun1810} 
\end{eqnarray}}

\begin{prop}
The equality \eqref{jun1810} is implied by the  Cheeger-M\"uller theorem.  
\end{prop}
\proof {The flat bundle $(V_{\sigma},\nabla^{V_{\sigma}})$ of complex vector spaces
with the parallel  metric $h^{V_{\sigma}}$ gives rise to a Reidemeister torsion
$\tau_{RM}(W,g^{TW},\cV_{\sigma},h^{\cV_{\sigma}})\in \R^{+}$ defined as follows. Let $\cC(\cV)$ be the cochain complex of some finite smooth triangulation of $W$ with coefficients in $\cV$. It acquires a geometry $h^{\cC(\cV)}$ from
$h^{\cV}$. The de Rham isomorphism induces  isomorphisms
$H^{i}(\cC_{\sigma})\cong H^{i}(W,\cV_{\sigma})$ and therefore the isomorphisms $\kappa_{\sigma}$ in the following  chain of isomorphisms (see \eqref{ji3jfio2fj23987f9798f23f23f} for $\phi_{\sigma}(\cC(\cV))$)
$$\psi_{\sigma}:\det(\cC(\cV)_{\sigma})\stackrel{\phi_{\sigma}(\cC(\cV)),\cong}{\to}\det(H(\cC(\cV)_{\sigma}))\stackrel{\kappa_{\sigma}, \cong}{\to}
\det(H(W;\cV_{\sigma}))\ .$$
We let $h^{\det(\cC(\cV)_{\sigma})}$ and $h_{L^{2}}^{\det(H(W;\cV_{\sigma}))}$ denote the metrics induced on the determinants by the metrics $h^{\cC(\cV)_{\sigma}}$ and $h_{L^{2}}^{H^{i}(W;\cV_{\sigma})}$.
\begin{ddd} The Reidemeister torsion $\tau_{RM}(W,g^{TW},\cV_{\sigma},h^{\cV_{\sigma}})\in \R^{+}$ is uniquely defined such that the equality
$$\psi_{\sigma,*} h^{\det(\cC(\cV)_{\sigma})}= \tau_{RM}(W,g^{TW},\cV_{\sigma},h^{\cV_{\sigma}})^{2}\:h_{L^{2}}^{\det(H(W;\cV_{\sigma}))}$$ is satisfied.\end{ddd}
{The Reidemeister torsion does not depend on the choice of the triangulation. We refer to \cite[Sec. 1]{MR1189689} for more details.}
By comparison with 
  \eqref{may121210}  we get the equality
\begin{equation}\label{may121301}\tau_{RM}(W,g^{TW},\cV_{\sigma},h^{\cV_{\sigma}}) =\tau_{\sigma}(\cC(\cV),h^{C(\cV)},\kappa^{*}h^{\cH(\cV)}) \ .\end{equation}
In analogy with \eqref{hkhkfhkefhfhehfekjwhfkewfhewfwfwefwefwfwe} we define 
$$\ln \tau_{RM}(W,g^{TW},\cV,h^{\cV}):=\sum_{\sigma\in \Sigma } \ln \tau_{RM}(W,g^{TW},\cV_{\sigma},h^{\cV_{\sigma}}) \:\bar \delta_{1,\sigma}\in A^{-1}\ .$$
Again we have
$\tau_{RM}(W,g^{TW} ,\underline{R}_{W},h^{\underline{R}_{W}})=0$.
 We conclude that
\begin{eqnarray}
\lefteqn{\chi(W)\left[\cycl\left(\det(s^{*}\cV),h^{\det(s^{*}\cV)}\right) -\beins\right]
}&&\nonumber\\&\stackrel{Lemma\: \ref{jun1820}}{=}&
\sum_{i}(-1)^{i} \cycl( C^{i}(\cV),h^{C^{i}(\cV)})-\dim(\cV)\sum_{i}(-1)^{i} \cycl( C^{i}(\underline{R}_{W}),h^{C^{i}(\underline{R}_{W})})
\nonumber
\\&\stackrel{Prop. \:\ref{may121220} }{=}&D(\cV,h^{\cV})-\dim(\cV) D(\underline{R}_{W},h^{\underline{R}_{W}})\nonumber\\&&+
 a\left(  \ln  \tau(\cC(\cV),h^{\cC(\cV)},\kappa^{*}h^{\cH(\cV)})\right)  - \dim(\cV)\:a\left(  \ln  \tau(\cC(\underline{R}_{W}),h^{\cC(\underline{R}_{W}},\kappa^{*}h^{\cH(\underline{R}_{W})} )\right)    \nonumber\\
&\stackrel{\eqref{may121301}}{=}&D(\cV,h^{\cV})-\dim(\cV) D(\underline{R}_{W},h^{\underline{R}_{W}})\nonumber\\ 
&&+a\left(  \ln  \tau_{RM}(W,g^{TW},\cV ,h^{\cV})\right)  -\dim(\cV)\:a\left(  \ln  \tau_{RM}(W,g^{TW},\underline{R}_{W},h^{\underline{R}_{W}})   \right) \nonumber\\&=&D(\cV,h^{\cV})-\dim(\cV) D(\underline{R}_{W},h^{\underline{R}_{W}})+a\left(  \ln  \tau_{RM}(W,g^{TW},\cV ,h^{\cV})\right)
\label{may121230}
\end{eqnarray}
}

{
The Cheeger-M\"uller theorem   \cite{MR528965}, \cite{MR498252} (see \cite{MR1189689}, \cite{MR1185803} for generalizations), applied to the present situation, states 
\begin{theorem}\label{sep1301}
For every $\sigma\in \Sigma$ we have the equality
$$T^{an}(W,g^{TW},\cV_{\sigma},h^{\cV_{\sigma}})=\ln \tau_{RM}(W,g^{TW},\cV_{\sigma},h^{\cV_{\sigma}})\ .$$
\end{theorem}
Plugging this equality into \eqref{may121230} we get exactly the right-hand side of  \eqref{jun1810}.  \hB

\begin{rem}{\rm 
The Cheeger-M\"uller theorem  therefore  implies a non-trivial consequence of   the transfer index conjecture.
Vice versa, observe that the Transfer Index Conjecture implies an $\R/\Z$-version of the Cheeger-M\"uller theorem.
More precisely,  it implies that 
$$ a\left(  \ln  \tau_{RM}(W,g^{TW},\cV ,h^{\cV})\right)  = a\left( T^{an}(W,g^{TW},\cV,h^{\cV})\right) $$
 The TIC thus  only implies an identity between torsion invariants after application of $a$, i.e. modulo the image of the regulator
$K_{1}(R)\to A^{-1}$.

}\hB \end{rem}}


\subsubsection{$S^{1}$-bundles}\label{jul0110}

The analytic  torsion form for $U(1)$-bundles equipped with a flat line bundle has been calculated in \cite[Prop. 4.13]{MR1303026}. In this paragraph we discuss the consequences of this calculation for the Transfer Index Conjecture \ref{may3103}.
In particular we see that the TIC implies  Proposition \ref{jun281} which can be verified independently by  methods from completely different fields. We consider this again as an example supporting the TIC.

\bigskip

We consider an $U(1)$-principal bundle $p:P\to B$ over a simply connected base  manifold $B$ with non-trivial first Chern class $c_{1}(p)\in H^{2}(B;\Z)$. The second cohomology group $H^{2}(B;\Z)$ of the base space $B$ is torsion-free, and we assume that the Chern class $c_{1}(p)$ is a primitive element. This implies that the total space $P$ of the bundle is simply connected, too.
For a prime $r\ge 3$ we have an embedding of groups $$\Z/r\Z\hookrightarrow U(1)\ ,\quad [n]\mapsto e^{\frac{2\pi i n}{r}}\ .$$ The group $\Z/r\Z$ therefore acts on $P$ via restriction of the $U(1)$-action. The quotient
$$\pi:W:=P/(\Z/r\Z)\to B$$ again has the structure of an $U(1)$-principal bundle. Its first Chern class satisfies
$c_{1}(\pi)=rc_{1}(p)$. 

\bigskip

We consider the ring $$R:=\Z[\xi]/(1+\xi+\dots+\xi^{r-1})\ .$$ 
This is the ring of the integers in the cyclotomic number field $\Q(\xi)$ obtained from    $\Q$   by  adjoining an $r$'th root of unity, {see Example \ref{fkjwelfwefewfwf}}. We get a representation
$$\Z/r\Z\to GL(1,R)\ , \quad [n]\mapsto \xi^{n}\ .$$  This representation determines a local system
$\cV$ of free $R$-modules of rank one on $W$.

Each choice of a nontrivial complex $r$'th root of unity $\sigma(\xi)\in \C$ defines a place $\sigma\in \Sigma$ of $\Q(\xi)$. Since we assume that $r$ is prime and
$r\ge 3$, all these places are complex. Therefore $$r_{\R}=0\ , \quad r_{\C}=\frac{r-1}{2}$$ {by Example  \ref{fkjwelfwefewfwf1}.}
For each place $\sigma$ the image of
$\Z/r\Z\to GL(1,R)\stackrel{\sigma}{\to} GL(1,\C)$ is contained in the subgroup $U(1)\subset GL(1,\C)$ so that
  the bundle $V_{\sigma}\to W$ associated to $\sigma$ and $\cV$  has a canonical flat hermitean metric $h^{V_{\sigma}}$. The collection of these metrics forms a geometry $h^{\cV}$ on $\cV$.
  
  \bigskip
  
In the following we calculate the {analytic} and topological index of $(\cV,h^{\cV})$ explicitly. 
For the topological side it turns out to be simpler to consider a relative situation. 
We therefore consider the class
\begin{equation} \label{jun218}\cycl(\cV,h^{\cV})-\beins_{W}\in \widehat{KR}^{0}(W)\ ,\end{equation}
 where $\beins_{W}$ is defined in \eqref{fewjlkewfklelklkklwefwfewfw2355tgrg}.
This  difference is flat. 

We equip the proper submersion $\pi:W\to B$ with the fibrewise Riemannian metric $g^{T^{v}\pi}$ which is $U(1)$-invariant and normalized such that the fibres have unit volume. We furthermore choose a principal bundle connection $\nabla^{W}$ on $W$. It defines and is defined by a $U(1)$-invariant horizontal distribution $T^{h}\pi$.
The pair $g:=(g^{T^{v}\pi},T^{h}\pi)$ is a Riemannian structure on the proper submersion $\pi$. These choices together fix the transfer
$$\hat \tr:\widehat{KR}^{0}(W) \to \widehat{KR}^{0}(B)$$ in differential algebraic $K$-theory by Theorem \ref{may3102}.

Since the class (\ref{jun218}) is flat we can interpret it in $KR\R/\Z^{-1}(W)$, and  by  \eqref{jun0304}
we can express the differential cohomology transfer of flat classes in terms of the usual Becker-Gottlieb transfer:
  $$\hat \tr\left(\cycl(\cV,h^{\cV})-\beins_{W}\right) =\tr^{*}\left(\cycl(\cV,h^{\cV})-\beins_{W}\right)\ .$$ 
The fundamental vector field generating the $U(1)$-action is a vertical vector field on $W$ without zeros.  This implies (compare e.g. \cite{MR2248973}) that the Becker-Gottlieb transfer for $\pi$ vanishes. Using \eqref{fewjlkewfklelklkklwefwfewfw2355tgrg}
it follows that
\begin{equation}\label{jun2110}0=\hat \tr\left(\cycl(\cV,h^{\cV})-\beins_{W}\right)=\hind^{top}( \cV,h^{\cV} )-\hind^{top}( \underline{R}_{W},h^{\underline{R}_{W}} )\ .\end{equation}

We now consider the analytic indices.
The holonomy of $\cV$ along the fibres $W_{b}\cong S^{1}$, $b\in B$, is the non-trivial element $\xi$.
We get an exact sequence
$$0\to H^{0}(W_{b},\cV_{|W_{b}})\to R\xrightarrow{1-\xi} R\to H^{1}(W_{b},\cV_{|W_{b}})\to 0\ .$$
Since $R$ is an integral domain we see that $H^{0}(W_{b},\cV_{|W_{b}})=0$. 
It follows that $H^{1}(W_{b},\cV_{|W_{b}})$ is torsion. 
{Recall the decomposition \eqref{apr123013} of a sheaf of finitely generated $R$-modules into a torsion and a projective part.} We conclude that
$$\Proj(R^{i}\pi_{*}(\cV))=0\ , \quad i\ge 0\ .$$
We define the torsion $R$-module $L$ by
$$0\to R\stackrel{1-\xi}{\to} R\to L\to 0\ .$$
Since $B$ is simply connected we have
$$R^{1}\pi_{*}(\cV){\cong} \Tors(R^{1}\pi_{*}(\cV))\cong \underline{L}_{B}\ ,$$
{where $\underline{L}_{B}$ denotes the locally constant sheaf of $R$-modules generated by $L$.}
Let $q:B\to *$ be the projection and $\hat \cZ(L)\in \widehat{KR}^{0}(*)$, {see \eqref{may121207}}. {Combining \eqref{may121260} with Definition \ref{jun271}}  we can write
$$\hind^{an}(\cV,h^{\cV})=a(\cT(\pi,g,\cV,h^{\cV}))+q^{*} 
\hat \cZ(L)\ .$$
Next we calculate $  \hind^{an}( \underline{R}_{W},h^{ \underline{R}_{W}})$.   We have
$$R^{i}\pi_{*}(\underline{R}_{W})\cong \left\{\begin{array}{cc} \underline{R}_{B}&i=0,1\\0&i\ge 2
\end{array}\right .\ .$$ 
Since we have normalized the volume of the fibres to one it is  easy to check that
$h_{L^{2}}^{H^{i}(W/B,R_{\sigma})}$ is the standard metric on $(\underline{R}_{B})_{\sigma}$
for $i=0,1$. This implies
$$\cycl\left( R^{i}\pi_{*}(\underline{R}_{W}),h^{R^{i}\pi_{*}(\underline{R}_{W})}\right)=\beins_{B},\quad i=0,1\ .$$
The contributions of $i=0$ and $i=1$ to $\hind_{0}(1_{W})$ thus cancel out in (\ref{jun219}),  and
we get $$\hind^{an}(\beins_{W})=a(\cT(\pi,g,\underline{R}_{W},h^{\underline{R}_{W}}))\ .$$ 
These torsion forms have explicitly been calculated by  Bismut-Lott \cite[Prop. 4.13]{MR1303026} and 
Lott \cite{MR1297676}. In the case of non-trivial holonomy we get
 \begin{eqnarray*}
\cT(\pi,g,\cV,h^{\cV})&=&\sum_{\sigma\in \Sigma} \left(\sum_{j even} (-1)^{j/2} \frac{1}{(2\pi)^{j}}\: \frac{(2j+1)!}{2^{2j}(j!)^{2}} \:\Ree(\Li_{j+1}(\sigma(\xi))) \: c_{1}^{j}(\nabla^{W}) \:\: \bar \delta_{2j+1,\sigma}\right.  \\&&\left. +\sum_{j odd} (-1)^{(j-1)/2} \frac{1}{(2\pi)^{j}}\: \frac{(2j+1)!}{2^{2j}(j!)^{2}} \:\Imm(\Li_{j+1}(\sigma(\xi))) \:c_{1}^{j}(\nabla^{W})\:\:\bar \delta_{2j+1,\sigma}\right)
\end{eqnarray*}
(see {\eqref{dh12j1l2jlkjdljlj21ldjl12jd21d} and} \eqref{wkejfeljlwelkfjewlfjewlkfjewfkljewlfke902irp2ri2p3r32r32r32r23r} for $\bar \delta_{2j+1,\sigma}$).
 Similarly, in the case of trivial holonomy we get
$$  
\cT(\pi,g,{\underline{R}}_{W},h^{{\underline{R}}_{W}})=\sum_{\sigma\in \Sigma} \sum_{j even} (-1)^{j/2} \frac{1}{(2\pi)^{j}}\: \frac{(2j+1)!}{2^{2j}(j!)^{2}} \:\Li_{j+1}(1) \: c_{1}^{j}(\nabla^{W}) \:\: \bar \delta_{2j+1,\sigma}\ \ .$$
The Transfer Index Conjecture \ref{may3103} together with the vanishing (\ref{jun2110}) of the topological index now implies that
$$a(\cT(\pi,g,\cV,h^{\cV}))+{q^{*} 
\hat \cZ(L)}=a(\cT(\pi,g,\underline{R}_{W},h^{\underline{R}_{W}}))\ .$$
It has been checked in Subsection \ref{may121280} that {this, as a consequence of the Cheeger-M\"uller theorem},   holds true after restriction to a point. 
Therefore {it is now interesting to consider}   the contribution of higher-degree forms. {Then the contribution of $q^{*}\hat \cZ(L)$ drops out.}     
The transfer index conjecture predicts that
$$ a\left(\sum_{\sigma\in \Sigma}   \frac{1}{(2\pi)^{j}}\: \frac{(2j+1)!}{2^{2j}(j!)^{2}} \:\Imm(\Li_{j+1}(\sigma(\xi))) \:c_{1}^{j}(\nabla^{W})\:\:\bar \delta_{2j+1,\sigma}\right)=0$$
for odd $j$, and
$$a\left(
\sum_{\sigma\in \Sigma}    \frac{1}{(2\pi)^{j}}\: \frac{(2j+1)!}{2^{2j}(j!)^{2}} \:\left(\Ree(\Li_{j+1}(\sigma(\xi)))-\Li_{j+1}(1)\right) \: c_{1}^{j}(\nabla^{W})\:\: \bar\delta_{2j+1,\sigma}\right)=0$$
for even $j\ge 1$.
Recall that we have an exact sequence \eqref{exats111}
$$KR^{-1}(B)\stackrel{c}{\to} H(A)^{-1}(B) \stackrel{a}{\to} \widehat{KR}^{0}(B)\stackrel{R}{\to} Z^{0}(\Omega A
(B))\ .$$
We now use the fact that $\C\P^{\infty}:=\colim_{n} \C\P^{n}$ classifies $U(1)$-principal bundles. 
If we apply the above reasoning to  the bundles $S^{2n+1}\to \C\P^{n}$ for all $n$ in place of $P\to B$, then we see  
that the Transfer Index Conjecture \ref{may3103} implies:
\begin{con}\label{sep1401}
There exists an element
$x\in KR^{-1}(\C\P^{\infty})$ such that
$$c(x)=\sum_{j\ge 1} u_{j} \ c_{1}^{j}\in H(A)^{-1}(\C\P^{\infty})$$
with \begin{equation}\label{hdkqhwdkjhkhqwdjqwdwqdwqdwqwqdwqdw}
u_{j} := \frac{1}{(2\pi)^{j}} \frac{(2j+1)!}{2^{2j}(j!)^{2}}\left\{\begin{array}{cc}
\sum_{\sigma\in \Sigma }\Imm(\Li_{j+1}(\sigma(\xi))) \ \bar \delta_{2j+1,\sigma}& j \: odd\\[0.3cm]
\sum_{\sigma\in \Sigma }\big(\Ree(\Li_{j+1}(\sigma(\xi))) -  \Li_{j+1}(1)\big)\ \bar \delta_{2j+1,\sigma}&j\ge 1\: even
\end{array}\right\} \in A_{2j+1}
\end{equation}

and the generator $c_{1}\in H^{2}(\C\P^{\infty};\Z)$.
\end{con}

{
\begin{rem}{\rm
While the ranks of the algebraic $K$-groups of a space like $\C\P^{\infty}$ can easily be calculated using the theorem of Borel 
the determination of the rational or integral structure is a very difficult matter. In this respect it is very interesting to predict the existence of elements in $KR^{-1}(\C\P^{\infty})$ whose regulator is explicitly prescribed. In the following we discuss a consequence for  the image of the regulator  $K_{*}(R)\to A_{*}$ itself.

}\hB
\end{rem}
}

{Using the Atiyah-Hirzebruch spectral sequence} we can calculate
$$KR\Q^{*}(\C\P^{\infty})\cong KR\Q^{*}[[c_{1}]]\ .$$
Therefore the image $x_\Q\in KR\Q^{*}(\C\P^{\infty})$ of the class $x\in  KR^{-1}(\C\P^{\infty})$
  can be written in the form
$$x_{\Q}=\sum_{j\ge 1}  \frac{(2j+1)!}{2^{2j}(j!)^{2}} \ x_{2j+1} \ c_{1}^{j}$$
with $x_{2j+1}\in K_{2i+1}(R)\otimes \Q$ satisfying
$$\frac{1}{(2\pi)^{j}}\frac{(2j+1)!}{2^{2j}(j!)^{2}}\sum_{\sigma\in \Sigma }\langle x_{2j+1},\omega_{2j+1}(\sigma)\rangle \ \bar \delta_{2j+1,\sigma}=u_{j}\in A_{2j+1}
\ ,$$
where {the cohomology classes} $\omega_{2j+1}(\sigma)$ are given by \eqref{kqedlkjdlqkwjlwqkjqo32jle}.

\bigskip

{Finally}, let us express the sum on the right-hand side in terms of the Borel regulator. 
{We use the target $X_{2j+1}(R)$ defined in \eqref{jul0230} usually considered in arithmetic geometry. The isomorphism with $A_{2j+1}$ is fixed by Proposition \ref{jul0231}.} We define the element
$$\Li^{R}_{j+1}\in X_{2j+1}(R)$$  such that
$$\Li^{R}_{j+1}(\sigma):=[\Li_{j+1}(\sigma(\xi))]\in \C/\R(j+1)\cong \R(j)$$
for all $\sigma\in \Sigma$.
We define
$$\Li^\Z_{j+1}\in X_{2j+1}(\Z)$$ by the same formula with $\xi=1$.
The inclusion $i:\Z\to R$ induces a map $$i_{*}:X_{2j+1}(\Z)\to X_{2j+1}(R)\ .$$
By Proposition \ref{jul0231} and {comparison with \eqref{hdkqhwdkjhkhqwdjqwdwqdwqdwqwqdwqdw}}
we have in $A_{2j+1}${,}
$$\psi(\Li^{R}_{j+1}-i_{*}\Li^{\Z}_{j+1}) = (-1)^{j} j! \   2^{2j} \ u_{j}\ .$$
Observe that the prefactors are rational. Recall the definition
(\ref{jul0232}) of the normalized Borel regulator map $r_{Bor}$.
Then 
  Conjecture   \ref{sep1401} implies {the assertion of the following proposition, {for which we provide} an independent proof.} 
\begin{prop}\label{jun281}
For every $j\ge 1$ there exists an element $y_{2j+1}\in K_{2j+1}(R)\otimes \Q$ such that
$$r_{Bor}(y_{2j+1})=\Li^{R}_{j+1}-i_{*}\Li^{\Z}_{j+1}\ .$$
 \end{prop}
 \proof
 This fact has indeed been proven in arithmetic geometry. In fact, one can realize both terms, 
$ \Li^{R}_{j+1}$ and $i_{*}\Li^{\Z}_{j+1}$, separately as Borel regulators.

For the second,
 applying the main result of Borel \cite{MR0506168} in the case of the field $\Q$,
 we conclude that there exists an element
 $z_{2j+1}\in K_{2j+1}(\Z)\otimes \Q$ such that
 $r_{Bor}(z_{2j+1})=[\Li^{\Z}_{j+1}]\in X_{2j+1}(\Z)$. 
 We then have  
 $$i_{*}\Li^{\Z}_{j+1}=r_{Bor}(i_{*}z_{2j+1})$$
 by the naturality of the Borel regulator.
 
 \bigskip
 
 The existence of an element $w_{2j+1}\in K_{2j+1}(R)\otimes \Q$ with $r^{Bor}(w_{2j+1})=\Li^{R}_{j+1}$ is ensured by a Theorem of Beilinson, see \cite[Thm. 5.2.1]{MR2002643}, \cite{MR862627}, \cite{MR944995}, \cite{MR1014822}\footnote{{We thank} G. Kings for explaining this result.}.
 In these papers an element in motivic cohomology with  prescribed Beilinson regulator is constructed.
 Now one can identify this particular motivic cohomology group with $K_{2j+1}(R)\otimes \Q$ so that the Beilinson regulator is half of the Borel regulator \cite{MR1869655}.
 Then the combination $y_{2j+1}:=w_{2j+1}-i_{*} z_{2j+1}$ does the job. \hB
 
{ 
 \begin{rem}{\rm The elements $y_{2j+1}\in K_{2j+1}(R)\otimes \Q$ belong to the very few examples of non-trivial elements in rational  higher $K$-theory of number rings with known regulators. 
  }
 \end{rem}
}

 \begin{rem}{\rm 
As an alternative to the arithmetic geometric construction of the elements $w_{2k+1}$ and $y_{2j+1}$ above  one could also obtain the existence of these elements by an argument based on Fact \ref{jul1621} and Theorem \ref{jul1640}. The details will be written up elsewhere.\footnote{{This is a result from the diploma thesis of Aron Strack.}}
}\hB\end{rem}

 The validity of the  consequence \ref{jun281} further supports the Transfer Index Conjecture \ref{may3103}.
   
 \subsubsection{Higher analytic and  Igusa-Klein torsion}
 
 {
In this subsection we collect further  consequences of the Transfer Index Conjecture \ref{may3103}.
  They involve the topological version of the higher torsion defined by Igusa and Klein \cite{MR1945530} and its comparison with the analytic torsion due to Bismut and Goette \cite{MR1867006}, \cite{MR2674876}.

These consequences come out of a comparison between the analytical and topological index for geometric locally constant sheaves of finitely generated projective $R$-modules whose underlying topological $K$-theory classes are in the image of the unit \eqref{jul3101}. This property allows an explicit calculation, as we shall see below.

The interesting feature of the discussion is that the TIC implies the existence of algebraic $K$-theory classes with regulator given by
 topological expressions. We interpret the fact that one can show, in some cases, the existence of such classes by independent means as a support for the validity of the TIC.
}

\bigskip

 The unit of the algebraic $K$-theory spectrum is a map of spectra  
 \begin{equation}\label{jul3101}\epsilon:S\to KR\ .\end{equation} 
 
 {
 \begin{rem}{\rm 
 The algebraic $K$-theory spectrum of a commutative ring $R$ is actually a commutative ring spectrum, i.e. an object
of $\CAlg(\Sp)$. The unit is part of this structure. We refer to \cite{2013arXiv1305.4550G}, \cite{Bunke:2013kx} for a construction of $KR$ as a commutative ring spectrum  using ring completion. 

In order to define the unit using the approach to algebraic $K$-theory via Definition \ref{nov0101} one could proceed as follows. Let $\mathbf{Fin}$ be the symmetric monoidal category of finite sets with the monoidal structure given by  the disjoint union. By the Baratt-Priddy theorem the sphere spectrum can be presented as the $K$-theory of $\mathbf{\mathbf{Fin}}$:
$$S\simeq  \spp(\Omega B(\Nerve(i\mathbf{Fin})))\ .$$
For a ring $R$ we can consider the symmetric monoidal functor
$$\mathbf{Fin}\to \bP(R)\ , \quad X\mapsto R[X]\ ,$$
where $R[X]$ denotes the $R$-module freely generated by the set $X$.
This functor induces the unit by application of the composition $\spp \circ \Omega B\circ \Nerve\circ i$.
   }\hB\end{rem}
 }
 
  We define the differential sphere spectrum $\Diff(S)$ using the canonical {differential data} introduced in Definition \ref{candef} and obtain the differential stable cohomotopy theory $\hat{S}^{0}$ by Definition \ref{difdef1}. Since the sphere spectrum is rationally even, by \cite{MR2608479}  the axioms stated in Definition \ref{jun216}  determine the functor $\hat S^{0}$  up to unique isomorphism.  We can extend the unit
 $\epsilon$ to a map of canonical {differential} data (unique up the action of $K_{1}(R)\otimes \Q$) and thus obtain a map
   $$\hat\epsilon:\Diff(S)\to \Diff(KR)$$
   of differential function spectra.  In the following we denote the underlying class map $I$ (see Definition \ref{jun216}, 2.) for the differential cohomology theory $\hat S^{0}$ by $I_{S}$ in order to distinguish it from the corresponding transformation for $\widehat{KR}^{0}$ still denoted by $I$. Since the homotopy groups of the sphere spectrum are finite in non-zero degrees (by a theorem of Serre){,}
the $\Z$-graded group $\pi_{*}(S)\otimes \R\cong \R$ is trivial in non-zero degrees.
{It follows from the exact sequence \eqref{exats111} that}  $I_{S}$ induces an isomorphism of smooth groups
 \begin{equation}\label{sep1405}I_{S}:\hat S^{0}\stackrel{\cong}{\to}S^{0}\ .\end{equation} 
 
 \bigskip

 We consider a proper submersion $\pi:W\to B$  with closed fibres and a choice of a Riemannian structure (Definition \ref{may272}). We  assume, for simplicity, that $B$ is connected.
 It follows from the naturality of the differential transfer in the data Theorem \ref{may3102}, 4. 
 that {on $\hat S^{0}$  we have} \begin{equation}\label{r4zgrrr872r83r23kjfo23r32lkfnfkjhwgrfjzwffewfe}
\hat \tr\circ \hat\epsilon=\hat \epsilon\circ I_{S}^{-1}\circ   \tr^{*}\circ I_{S}\ .
\end{equation} {For a manifold $M$ let $1_{M}\in S^{0}(M)$ be the unit and note that $\beins_{M}=\hat \epsilon(I_{S}^{-1}(1_{M}))$}, see \eqref{fewjlkewfklelklkklwefwfewfw2355tgrg} for $\beins_{M}$.
  It follows from  the relation $\tr^{*}(1_{W})=\chi (F) 1_{B}$ that we have \begin{equation}\label{fwefkjflkewfjewflkj3493495835435435}
\hat \tr(\beins_{W})=\hat \tr (\hat\epsilon( I^{{-1}}_{S}(1_{W})))\stackrel{\eqref{r4zgrrr872r83r23kjfo23r32lkfnfkjhwgrfjzwffewfe}}{=}\hat \epsilon(I^{{-1}}_{S}(\tr^{*}(1_{W})))= \hat\epsilon(I^{{-1}}_{S}( \chi(F)1_{B}))=\chi(F)\beins_{B}\ ,
\end{equation}
 where $\chi(F)$ is the Euler characteristic of the fibre of
$\pi:W\to B$.

\bigskip

 We now consider a geometric locally constant sheaf $(\cV,h^{\cV})$ of {finitely generated} $R$-modules such that
$[\cV]={\dim(\cV)}  1_{W}\in K R^{0}(W)$. By \eqref{exats111}  there exists a form $\eta\in \Omega A^{-1}(W)$ such that
\begin{equation}\label{mar0103}\cycl(\cV,h^{\cV})+a(\eta)={\dim(\cV)}\beins_{W}\end{equation}
holds true in $\widehat{KR}^{0}(W)$.
In this case, using \eqref{may3103222} and \eqref{fwefkjflkewfjewflkj3493495835435435}, we can calculate the topological index explicitly:
\begin{equation}\label{mar0104}\hind^{top}(\cV,h^{\cV})=\chi(F){\dim(\cV)} \beins_{B} -a\left(\int_{W/B} \eta\wedge e(g)\right)\ .\end{equation}

In order to calculate the analytic index $\ind^{an}(\cV,h^{\cV})$ we need the notion of an unipotent sheaf.
\begin{ddd} A locally constant sheaf $\cV$ of finitely generated
$R$-modules  
 is called unipotent if it has
 a decreasing filtration 
$$0=F^{r}\cV\subseteq F^{r-1}\cV\subseteq\dots\subseteq F^{1}\cV\subseteq F^{0}\cV=\cV$$
by locally constant subsheaves such that the quotients $F^{i-1}\cV/F^{i}\cV$  are trivialized, i.e. identified with $\underline{R^{n_{i}}}_{M}$ for all $i\ge 1$.\end{ddd}
If $\cV$ is an  unipotent {locally constant sheaf of finitely generated $R$-modules} on a manifold $M$, then  the trivializations  provide canonical flat geometries $h^{F^{i-1}\cV/F^{i}\cV}$ on the subquotients.  Let $h^{\cV}$ be any geometry on $\cV$. It induces for all $i$ geometries $h^{F^{i}\cV}$
on the subsheaves $F^{i}\cV$ by restriction. 
{Let $-\cT_{i}\in \Omega A^{-1}(M)/\im(d)$ be the torsion form of the sequence (starting in degree $0$)
$$0\to F^{i+1}\cV\to F^{i}\cV\to F^{i}\cV/F^{i+1}\cV\to 0$$
with the geometries fixed above.}
Using   Lott's relation   Theorem \ref{jun0301} 
\begin{equation}\label{sep1410}\cycl(F^{i+1}\cV,h^{F^{i+1}\cV})+\cycl(F^{i}\cV/F^{i+1}\cV,h^{F^{i}\cV/F^{i+1}\cV})-\cycl(F^{i}\cV,h^{F^{i}\cV})=-a(\cT_{i} )\end{equation}
and the {equalities} $$\cycl(F^{i}\cV/F^{i+1}\cV,h^{F^{i}\cV/F^{i+1}\cV})=\dim(F^{i}\cV/F^{i+1}\cV) \beins_{M}$$
in order to {obtain the equality}
$$\cycl(\cV,h^{\cV})=\dim(\cV) \beins_{M}+a(\cT(F^{*},h^{\cV}))\ , \quad \cT(F^{*},h^{\cV}):=\sum_{i} \cT_{i}  \ .$$

We now come back to our original situation and assume that the geometry $h^{\cV}$ on the sheaf $\cV$ on $W$ is flat and  that the sheaves
$R^{i}\pi_{*}(\cV)$ are unipotent with filtrations $F_{i}^{*}$. 
{Using the geometries \eqref{apr140113}
we} define the form
$$\tau(\pi,g,\cV,h^{\cV}):=\cT(\pi,g,\cV,h^{\cV})+ \sum_{i}(-1)^{i} a(\cT(F_{i}^{*},h_{L^{2}}^{R^{i}\pi_{*}(\cV)}))  \in \Omega A^{-1}(B)\ .$$
One checks using \eqref{mar0101} and \eqref{mar0102} that $\tau(\pi,g,\cV,h^{\cV})$
 is closed. By the usual homotopy argument its cohomology class $$\tau(\pi,\cV)\in H(A)^{-1}(B)$$  does not depend on the Riemannian geometry $g$ of $\pi$  and the flat geometry $h^{\cV}$. 
 \begin{ddd}
 Let $\pi:W\to B$ be a proper submersion with closed fibres.
Assume that  $\cV$ is a locally constant sheaf of finitely generated $R$-modules on $W$  which admits a flat geometry $h^{\cV}$ and whose  cohomology sheaves $R^{i}\pi_{*}(\cV) $ are unipotent. Then we define  the  higher analytic torsion class   by { 
 $$\tau(\pi,\cV):=[\tau(\pi,g,\cV,h^{\cV})]\in H(A)^{-1}(M)\ .$$ }   \end{ddd} A similar definition has been given in \cite[Def. 2.8]{MR2674876}.
Note that here we use the normalization fixed by Bismut-Lott.  
We have by the above construction
  \begin{equation}\label{mar0105}\hind^{an}(\cV,h^{\cV})=\chi(F){\dim(\cV)}  \beins_{B} +a(\tau(\pi,\cV) )\ .\end{equation}

   Assuming the Transfer Index Conjecture \ref{may3103} {by comparison with \eqref{mar0104}}  we get the following consequence:
  \begin{equation}\label{apr140413}a\left(\int_{W/B}\eta\wedge e(g))+\tau(\pi,\cV)\right)=0\ .\end{equation}
  In particular, if we take the canonical geometry $h^{\underline{R_{W}}}$ on
  $\cV=\underline{R_{W}}$, then in \eqref{mar0103} we can take $\eta=0$ and get
 \begin{equation}\label{apr140313}a(\tau(\pi,\underline{R}_{W}))=0
 \ .\end{equation}
\begin{ddd} We call a proper submersion  $\pi:W\to B$ unipotent if the cohomology sheaves $R^{i}\pi_{*}(\underline{R_{W}})$ are unipotent for all $i\ge 0$.\end{ddd}
  Note that any proper submersion over a simply connected base manifold $B$ is unipotent.
  {From \eqref{apr140213} and  \eqref{apr140313} we conclude}
 \begin{kor}\label{kjdlqwdqdwddwdwdqd}
 We assume that $\pi:W\to B$ is a unipotent proper submersion with closed fibres. {If the 
  Transfer Index Conjecture \ref{may3103}  holds true, then}
   there exists an element $x\in KR^{-1}(B)$ such that
 $c(x)=\tau(\pi,\underline{R}_{W})$ in $H(A)^{-1}(B)$.  \end{kor}

\begin{rem}{\rm 
On the one hand this   is an integrality statement for the
 higher analytic torsion class.
 On the other hand, if the torsion is known,  Corollary \ref{kjdlqwdqdwddwdwdqd} predicts again the existence of certain algebraic $K$-theory classes of $R$ with known regulator. In order to proceed along these lines we must specialize to cases where the torsion can be calculated.
 }\hB \end{rem}

{We now return to the situation where  $\cV$ is a locally constant sheaf of finitely generated projective $R$-modules on $W$. In addition we   assume}
   that the sheaves $ R^{i}\pi_{*}(\cV) $ are trivializable for all $i\in \nat$.
This is the case e.g. if the base manifold $B$ is simply connected and the cohomology groups of the fibres of $\pi$ with coefficients in $\cV$ are free. 
We further assume that $h^{\cV}$ is flat so that we can take $\eta=0$.

If $(V,\nabla^{V})$ is a flat complex vector bundle with trivializable cohomology bundles along $\pi$, then   Igusa-Klein torsion $\tau^{Igusa}_{IK}(\pi, V)$  is defined  \cite[Def. 4.4]{MR2674876}, \cite{MR2365656}, \cite{MR1945530}, and (\ref{jul1608}). The superscript indicates the Igusa normalization. Using (\ref{jul1605}) we define for all $j\ge 0$ the renormalized version
$$\tau_{IK}(\pi, V)_{2j}:= N_{Igusa}(2j+1)\tau^{Igusa}_{IK}(\pi, V)_{2j}\ .$$
To our sheaf $\cV$ of $R$-modules we can now associate the class
$$\tau_{IK}(\pi,\cV):=\sum_{\sigma\in \Sigma } \sum_{j=1}^{\infty}\tau_{IK}(\pi,V_{\sigma})_{2j} \ {\bar \delta_{2j+1,\sigma}}\in  H(A)^{-1}(B)\ .$$
We further define the  Igusa normalized   characteristic class of a real vector bundle $U\to M$
$${}^{0}J^{Igusa}(U):=  \sum_{j=1}^{\infty} (-1)^{j}\zeta(2j+1)\  [\ch(U\otimes \C)]_{4j} \ .$$
Its Bimut-Lott normalization is then 
$${}^{0}J(U):=  \sum_{j=1}^{\infty}    N_{Igusa}(4j+1)\  {}^{0}J^{Igusa}(U)_{4j} \ .$$
We furthermore define
$${}^{0}J_{R}(U):=  \sum_{j=1}^{\infty} {}^{0}J(U)_{4j}\   \kappa_{4j+1}\ \in H(A)^{-1}(M)\ .$$
The following theorem has been anounced by Goette \cite[Thm. 5.5]{MR2674876} (note that he uses the Chern normalization).
\begin{theorem}[Goette]\label{jul1640}
$$\tau(\pi,\cV)=\tau_{IK}(\pi,\cV)+\tr^{*}\ {}^{0}J_{R}(T^{v}\pi)\in H(A)^{-1}(M)\ .$$
\end{theorem}
  
If the fibre dimension is odd, then the transfer $\tr^{*}$   vanishes.   {Since $\eta=0$ from} {\eqref{apr140213} and} \eqref{apr140413} we {now} conclude:
 
\begin{kor}\label{jul1633}
Assume that the fibre dimension of $\pi:W\to B$ is odd, {that $\cV$ satisfies $[\cV] = \dim(\cV)1_{W}$,}
that the  geometry $h^{\cV}$ is flat, and that
the sheaves $R^{i}\pi_{*}(\cV)$ are trivializable for all $i\in \nat$. {If the Transfer Index Conjecture  {\ref{may3103}} holds true, 
then  there exists an element
$x\in KR^{-1}(B)$ such that $c(x)=\tau_{IK}(\pi,\cV)$. }
\end{kor}

{\begin{rem}{\rm 
At the moment we can verify {the existence of  an element
$x\in KR^{-1}(B)$ such that $c(x)=\tau_{IK}(\pi,\cV)$ }  in certain cases {by independent means}, see {Proposition} \ref{sep0507}.
But we can show in general {that such an element indeed exists rationally. This again supports the validity of the TIC.} }\hB\end{rem}}
\begin{prop}\label{jul1621} {Assume that the fibre dimension of $\pi:W\to B$ is odd, {that $\cV$ satisfies $[\cV] = \dim(\cV)1_{W}$,}
that the  geometry $h^{\cV}$ is flat, and that
the sheaves $R^{i}\pi_{*}(\cV)$ are trivializable for all $i\in \nat$.   Then} {there exists an element
 $x\in KR^{-1}(B)\otimes \Q$ such that $c(x)=\tau_{IK}(\pi,\cV)$. }
\end{prop}
 \proof
We first recall some details of the construction of the Igusa-Klein torsion. The unit of the algebraic $K$-theory spectrum $\epsilon:S\to KR$ extends to a fibre sequence
$$\dots \Omega KR\to Wh^{R}(*)\to S\xrightarrow{\epsilon} KR\to\Sigma Wh^{R}(*)\dots\ ,$$
defining the Whitehead spectrum $Wh^{R}(*)$.  Since the homotopy groups $\pi_{i}(S)$ for $i\ge 1$ are torsion the map of infinite loop spaces
$$\Omega^{\infty+1} KR\to \Omega^{\infty} Wh^{R}(*)$$
is a rational equivalence. In particular we have an isomorphism of real cohomology groups
\begin{equation}\label{jul1607}H^{*}(\Omega^{\infty} Wh^{R}(*);\R) \stackrel{{\cong}}{\rightarrow} H^{*}(\Omega^{\infty+1} KR;\R) \ .\end{equation}
The spectrum cohomology classes (\ref{kqedlkjdlqkwjlwqkjqo32jle})
$$\omega_{2k+1}(\sigma)\in H^{2k+1}(KR;\R)\ ,\quad \sigma\in \Sigma$$
 induce classes
$\tau_{2k}(\sigma)\in H^{2k}(\Omega^{\infty+1}KR;\R)$. Under the isomorphism (\ref{jul1607}) they correspond to the universal
torsion classes
$$\tau_{IK,2k}(\sigma)\in H^{2k}(\Omega^{\infty} Wh^{R}(*);\R)\ .$$
The main ingredient of the construction of the Igusa-Klein torsion of a bundle
$\pi:W\to B$ equipped with a locally constant sheaf $\cV$ of finitely generated projective $R$-modules is the construction of a well-defined homotopy class of maps
$$m^{\cV}_{IK}:B\to \Omega^{\infty}Wh^{R}(*)$$
using the theory of framed Morse functions.
Then by definition 
\begin{equation}\label{jul1608} \tau_{IK}(\pi,V_{\sigma})_{2k}=(m^{\cV}_{IK})^{*}\tau_{IK,2k}(\sigma)\ .\end{equation}
We now consider the diagram {of spaces} 
$$\xymatrix{&\Omega^{\infty+1}KR\ar[d]&\\B\ar@{-->}[ur]^{\check{f}}\ar@{.>}[dr]^{\chi}\ar@/_{2cm}/@{.>}[drr]^{\chi_{\Q}}\ar[r]^-{m^{\cV}_{IK}}&\Omega^{\infty}Wh^{R}(*)\ar[d]&\\
&\Omega^{\infty} S\ar[r]^{\epsilon_{{H\Q}}}&\Omega^{\infty} H\Q}\ .$$
The class $[\chi_{\Q}]\in H^{0}(B;\Q)\cong \Q$ is the Euler characteristic
of the fibre which vanishes by assumption. Since $\epsilon_{{H\Q}}$ is a rational
equivalence we conclude that the class
$[\chi]\in S^{0}(B)$ vanishes rationally. We choose $0\not=l\in \nat$ such that
$l [\chi]=0$. Then there exists a lift $\check{f}:B\to \Omega^{\infty+1} KR$ of the  map $l \ m^{\cV}_{IK}$.
We consider this map as a class
$[\check f]\in KR^{-1}(B)$.
Then we have
$l\  \tau_{IK}(\pi,\cV)= c([\check f])$. {
Hence $$x:=[\check f]\otimes l^{-1}\in KR^{-1}(B)\otimes \Q$$ does the job.}
\hB

We now consider the special case of bundles over the base $B=S^{2j}$.
Then we have $KR^{-1}(S^{2j})\cong K_{2j+1}(R)$.
In this case {Corollary} \ref{jul1633} is equivalent to:
\begin{kor}\label{jul1612}
Let $\pi:W\to S^{2j}$ be a proper submersion with odd-dimensional fibres and $(\cV,h^\cV)$ be a locally constant sheaf of finitely generated projective $R$-modules with flat geometry $h^{\cV}$. {We assume that the  cohomology of the fibres of $\pi$ with coefficients in $\cV$ is free.}  {If the Transfer Index Conjecture {\ref{may3103}} holds true, then} there exists an element $x\in K_{2j+1}(R)$
such that $$\omega_{2j+1}(\sigma)(x)=  \langle \tau_{IK}(\pi,V_\sigma),[S^{2j}]\rangle$$ for all $\sigma\in \Sigma$.
\end{kor}

Hatcher's construction provides examples of bundles $\pi:W\to S^{4k}$ with non-trivial Igusa-Klein torsion. 
Let
$$a_{k}:=\denom(\frac{B_{k}}{4k})\in \nat\ ,$$
where $B_{k}$ is the $k$'th Bernoulii number. By a Theorem of Adams \cite{MR0159336}
this number is the order of {the image of} the $J$-homomorphism
$$\pi_{4k-1}(O)\to \pi_{4k-1}(S)\ .$$
Let $$\epsilon_{k}:=\left\{\begin{array}{cc} 1&\mbox{$k$ is odd}\\
\frac{1}{2} &\mbox{$k$ is even}\end{array} \right\}\ �.$$
Then we have the following theorem shown in \cite[Cor. 1.2.4]{2010arXiv1011.4653G}:
\begin{theorem}\label{jul1610}
 For sufficiently large odd $n$  there exists a bundle
$\pi:W\to S^{4k}$ with fibre $S^{n}$ such that
$$\langle \tau^{Igusa}_{IK}(\pi,\C)_{4k},[S^{4k}]\rangle =a_{k}\epsilon_{k} \zeta(2k+1)$$
\end{theorem}
{For those examples we can verify the consequence of the Transfer Index Conjecture formulated in Corollary \ref{jul1612} by independent means.}
\begin{prop}\label{sep0507}
{ If the bundle $\pi:W\to S^{4k}$ obtained by Hatcher's construction  and $\cV=\underline{R}_{W}$,
then  there exists an element $x\in K_{4k+1}(R)$
such that $$\omega_{4k+1}(\sigma)(x)=  \langle \tau_{IK}(\pi,V_\sigma),[S^{4k}]\rangle$$ for all $\sigma\in \Sigma$.}
  \end{prop}
\proof
Homotopy theoretically, Hatcher's example can be understood by the following diagram which we take from \cite[Sec. 8]{MR2739777} with slight modifications.
\begin{equation}\xymatrix{&&&\Omega^{\infty+1} KR\ar[d]\\
&SG/SO\ar[d]\ar[r]&\Omega^{\infty+1} Wh^{Diff}(*)\ar[r]\ar[d]&\Omega^{\infty}Wh^{R}(*)\ar[d]\\
S^{4k}\ar@{.>}[ur]^{\tilde f}\ar[r]^{f}\ar@/^2cm/@{.>}[urrr]^{\hat{f}}\ar@/^3cm/@{.>}[uurrr]^{\check{f}}&BSO\ar[r]^{J}\ar[d]&\Omega^{\infty}S\ar@{=}[r]\ar[d]&\Omega^{\infty}S\ar[d]\\
& BSG\ar[r]&\Omega^{\infty} A(*)\ar[r]^{\lambda^{R}}&\Omega^{\infty}KR}
\end{equation} Here $SG=SL_{1}(S)\subseteq GL_{1}(S)$ is the connected component of the identity of  the group of units of the sphere spectrum, $A(*)$ is Waldhausen's $K$-theory of spaces \cite{MR802796} evaluated at the point, and
$Wh^{Diff}(*)$ is the spectrum appearing in Waldhausen's splitting \cite{MR802796}, \cite{MR921487}:
\begin{equation}\label{jul1609}A(*)\cong S\vee Wh^{Diff}(*)\ .\end{equation}
The map $\lambda^{R}:A(*)\to KR$ is the linearization, and
the map $\Omega^{\infty+1} Wh^{Diff}(*)\to \Omega^{\infty}Wh^{R}(*)$ is induced from the lower right commutative square.

 The map $J:BSO\to \Omega^{\infty} S$ induces the $J$-homomorphism \cite{MR0159336}. 
We choose $f\in \pi_{4k}(BSO)$ such that it represents an element the kernel 
of the $J$-homomorphism.
Because of the splitting (\ref{jul1609}) the induced map $S^{4k}\to BSG$ is zero homotopic and we can find a lift $\tilde f:S^{4k}\to SG/SO$. Since $J\circ f$ is zero homotopic we can further find a factorization of $\hat f:S^{4k}\to \Omega^{\infty} Wh^{R}(*)$ through a map
$\check{f}:S^{4k}\to \Omega^{\infty+1} KR$. This map represents a class
$$[\check{f}]\in K_{4k+1}(R)\ .$$
  
If we build Hatcher's example from the map $f$, then we get a bundle
$W\to S^{4k}$ such that
$$m^{\underline{R}_{W}}_{IK}\sim \hat f\ .$$
This implies that
$$\tau_{IK}(\pi,\underline{R}_{W})=\sum_{\sigma\in \Sigma}\hat{f}^{*} \tau_{4k}(\sigma) \:\bar \delta_{4k+1,\sigma}$$
and therefore
$$\langle \tau_{IK}(\pi,\underline{R}_{W}),[S^{4k}]\rangle= \sum_{\sigma\in \Sigma}\omega_{4k+1}(\sigma)([\check{f}]) \:\bar \delta_{4k+1,\sigma}\ .$$
\hB


%

\subsection{Discussion of Lott's relation}\label{aug0501}

In this subsection we show Theorem  \ref{jun0301}. For the case of exact sequences of locally constant sheaves of projective $R$-modules we refer to \cite{buta}. It remains to remove the assumption ``projective''. This is done in  Lemma
 \ref{uliapr1501}. Since the proof of Lott's relation given  \cite{buta} is quite involved at the end of this subsection
we give independent arguments in some special cases.

\bigskip

{
We first remove the projectivity assumption in Lott's relation and thus finish the proof of Theorem \ref{jun0301}.}

{
\begin{lem}\label{uliapr1501}
If Lott's relation \eqref{jul2023} holds true for exact sequence of sheaves of finitely generated projective
$R$-modules, then it also holds true  for exact sequences of sheaves of finitely generated  $R$-modules.
\end{lem}
\proof
We consider an exact sequence
$$\cV: 0\to \cV_{0}\to \cV_{1}\to \cV_{2}\to 0$$
of finitely generated $R$-modules. It induces a diagram of sheaves (see \eqref{apr123013} for notation)
$$\xymatrix{
&&0\ar[d]&0\ar[d]&0\ar[d]&\\&0\ar[r]&\Tors(\cV_{0})\ar[r]\ar[d]&\Tors(\cV_{1})\ar[r]^{!}\ar[d]&\Tors(\cV_{2})\ar[r]\ar[d]&0\\\cV:&0\ar[r]&\cV_{0}\ar[r]\ar[d]^{!!}&\cV_{1}\ar[r]\ar[d]&\cV_{2}\ar[r]\ar[d]&0\\\Proj(\cV):&0\ar[r]&\Proj(\cV_{0})\ar[r]\ar[d]&\Proj(\cV_{1})\ar[r]\ar[d]&\Proj(\cV_{2})\ar[r]\ar[d]&0\\
&&0&0&0&}$$
with exact raws and columns. Indeed, the arrow {named} $!$ is surjective since the arrow {named} $!!$ locally splits.}

{
We now assume that we are given geometries $h^{\cV_{i}}$ for $i=0,1,2$.
By Definition  \ref{uliapr1210} we have the relation
$$\sum_{i=0}^{2} (-1)^{i}\cycl(\cV_{i},h^{\cV_{i}})=\sum_{i=0}^{2} (-1)^{i}\cycl(\Proj(\cV_{i}),h^{\cV_{i}})
+\sum_{i=0}^{2} (-1)^{i} \hat \cZ(\Tors(\cV_{i}))\ .$$
Since the torsion forms for the sequences with geometries $(\cV,h^{\cV})$ and $(\Proj(\cV),h^{\cV})$
coincide it remains to show that
\begin{equation}\label{uliapr1510}\sum_{i=0}^{2} (-1)^{i} \hat \cZ(\Tors(\cV_{i}))=0\ .\end{equation}
To this end we consider the diagram  \eqref{vwvwnenwvewvnewvnwvnw23849823u43242343242}
$$\xymatrix{
&&0\ar[d]&0\ar[d]&0\ar[d]&\\\cB:&0\ar[r]&E(\Tors(\cV_{0}))\ar[r]\ar[d]&E(\Tors(\cV_{1}))\ar[r]\ar[d]&\cB_{2}\ar[r]\ar[d]&0\\\cA:&0\ar[r]&F(\Tors(\cV_{0})) \ar[r]\ar[d]&F(\Tors(\cV_{1}))\ar[r]\ar[d]&\cA_{2}\ar[r]\ar[d]&0\\&0\ar[r]&\Tors(\cV_{0})\ar[r]\ar[d]&\Tors(\cV_{1})\ar[r]\ar[d]&\Tors(\cV_{2})\ar[r]\ar[d]&0\\
&&0&0&0&}\ .$$
In particular, the right vertical sequence is a projective resolution of $\Tors(\cV_{2})$.
%
}

{
We now choose geometries  $h^{\cA_{i}}$, $i=0,1,2$ for the sheaves in the complex $\cA$. They induce geometries $h^{\cB_{i}}$ on the sheaves in the complex $\cB$ by restriction.
Then we have by Lemma \ref{may121201}
$$\sum_{i=0}^{2} (-1)^{i} \hat \cZ(\Tors(\cV_{i})) =\sum_{i=0}^{2} (-1)^{i} \cycl(\cA_{i},h^{\cA_{i}})-\sum_{i=0}^{2} (-1)^{i} \cycl(\cB_{i},h^{\cB_{i}})\ .$$
By \cite[Thm. 5.25]{buta} we can assume Lott's relation  \eqref{jul2023} for the sequences of projective $R$-modules $\cA$ and $\cB$. Since their torsion forms coincide (since their complexifications for every place $\sigma$  are isometric), we conclude  the equality \eqref{uliapr1510}. \hB
}

\bigskip

{We now discuss some cases of Lott's relation  \eqref{jul2023} which can be verified independently of \cite{buta}.}
We consider  three locally constant sheaves of finitely generated projective $R$-modules on a manifold $M$ with geometries $(\cV_{i},h^{\cV_{i}})$,  $i=0,1,2$,  which fit into an exact sequence
\begin{equation}\label{jul2021nnn}  \bV:0\to \cV_{0}\to \cV_{1}\to \cV_{2}\to 0\ .\end{equation} 
Then we define
\begin{equation}\label{aug0502}
 \delta(\bV):= \cycl(\cV_{0},h^{\cV_{0}})+\cycl(\cV_{2},h^{\cV_{2}})-\cycl(\cV_{1},h^{\cV_{1}})-a(\cT)\in \widehat{KR}^{0}(M)\ ,\end{equation}
 where $\cT$ is the analytic torsion form defined in (\ref{may121003}).
{Lott's relation asserts  the equality
  $\delta(\bV)=0$}. 
   We start with some general properties of $\delta(\bV)$.

\begin{lem}\label{aug1601}
We have
$\delta(\bV)\in \im(H(A)^{-1}(M)\xrightarrow{a} \widehat{KR}^{0}(M))$.
In particular, {the class $\delta(\bV)$}  only depends on the sequence $\bV$ and not on the choice of geometries.
\end{lem}
\proof
We have the relation $[\cV_{0}]+[\cV_{2}]=[\cV_{1}]$ in $KR^{0}(M)$ and therefore
$I(\delta(\bV))=0$. Furthermore, by \cite[Eq. (292)]{MR1724894} we have the equality \eqref{mar0102}
$$\omega(\cV_{0},h^{\cV_{0}})+\omega(\cV_{2},h^{\cV_{2}})-\omega(\cV_{1},h^{\cV_{1}})=d\cT$$ and therefore
$R(\delta(\bV))=0$. This implies that
$\delta(\bV)\in  \im(H(A)^{-1}(M)\to \widehat{KR}^{0}(M))$.
Since this image is a homotopy invariant functor we conclude that $\delta(\bV)$ does not depend on the geometries by a standard deformation argument.
\hB 

\begin{lem}\label{aug0503}
We have
$\delta(\bV)=0$ if the sequence $\bV$ splits.
\end{lem}
\proof
We use the isomorphism $\cV_{1}\cong \cV_{0}\oplus \cV_{2}$, the additivity of the cycle map, and that the torsion form
$\cT$ vanishes for a metrically split sequence in order to conclude that
$$\delta(\bV)=\cycl(\cV_{0},h^{\cV_{0}})+\cycl(\cV_{2},h^{\cV_{2}})-\cycl(\cV_{1},h^{\cV_{0}}\oplus h^{\cV_{2}})=0\ .$$ 
\hB
Here are further cases in which we know that $\delta(\bV)=0$ independently of \cite{buta}.
 \begin{lem}\label{aug0504}
 Assume that $R$ is the ring of integers in a totally real number field (resp. arbitrary number field ). Then we have
 $\delta(\bV)=0$ if $\bV$ is pulled back from a $CW$-complex of rational cohomological dimension $\le 3$  ($\le 1$).
 \end{lem}
\proof
We discuss the totally real case. Let $X$ be a $CW$-complex of rational cohomological dimension  $\le 3$ with a sequence $\bV^{\prime}$ and a map $f:M\to X$ such that $\bV\cong f^{*}\bV^{\prime}$. By approximation of a finite skeleton of $X$
we can assume that $X$ is a manifold with $H^{l}(X;\Q)=0$ for $4\le l\le \dim(M)$, and  that $f$ is smooth. Then we have
$\delta(\bV)=\delta(f^{*}\bV^{\prime})=f^{*}\delta(\bV^{\prime})$.
We claim that $f^{*}\delta(\bV^{\prime})=0$.
Note that by Theorem \ref{borel} we have $K_i(R)\otimes \R=0$ for $i=2,3,4$. It follows that 
\begin{equation}\label{aug0403}
  f^{*}(H(A)^{-1}(X))\cong f^{*}(H^0(X;\R)\otimes K_1(R))\ .
\end{equation}
 By (\ref{aug0403}) it suffices to show that $\delta(\bV^{\prime})$ is zero after restriction to a point. But after restriction to a point the sequence $\bV^\prime$ splits and we know that the conjecture holds true in this case by Lemma \ref{aug0503}. \hB

\section{Technicalities}\label{sep2210}

\subsection{Categories with weak equivalences and $\infty$-categories} \label{jul0501}

In this subsection we review the basic $\infty$-categorical constructions which we use throughout the paper.

\bigskip

In the present paper we work with the notion of $\infty$-categories as developed in detail in \cite{MR2522659}. Thus an $\infty$-category $\bC$ is a simplicial set $\bC\in \sSet$ which satisfies an inner horn filling condition. The set of zero simplices $\bC[0]$ in an $\infty$-category $\bC$  is interpreted as  the set of its objects, and the set of one-simplices $\bC[1]$ is considered as the set of one-morphisms.

An $\infty$-category $\bC$ gives rise to a homotopy category $\Ho(\bC)$ which is an ordinary category.  The homotopy category has the same set of objects as $\bC$, but its morphisms are equivalence classes of one-morphisms in $\bC$, see \eqref{jjkqdhwkqdkhjwqdqwd}. 
A morphism in $\bC$ is an equivalence if and only if it induces an isomorphism in $\Ho(\bC)$.

 A functor $\bC\to \bC^{\prime}$ between $\infty$-categories $\bC$ and $\bC^{\prime}$ is a map of simplicial sets, {and the $\infty$-category $\Fun(\bC,\bC^{\prime})$ of functors from $\bC$ to $\bC^{\prime}$ is the simplicial mapping space $\Map(\bC,\bC^{\prime})$, i.e. the  simplicial set whose set of $n$-simplices is the set of maps $\Delta^n\times\bC\to\bC^{\prime}$ in $\sSet$.}

If $\bC$ is an ordinary category, then its nerve $\Nerve(\bC)\in \sSet$ is an example of an $\infty$-category. If $\bD$ is a second ordinary category, then we have the relation
\begin{equation}\label{lkdjlkqwjdldwqdqwd}\Nerve(\Fun(\bC,\bD))\simeq \Fun(\Nerve(\bC),\Nerve(\bD))\ .\end{equation}
In the present paper, {most} $\infty$-categories arise from localizations. If $\bC$ is an $\infty$-category and $W$
 is a collection of arrows of $\bC$, then the localization
$\bC\to \bC[W^{-1}]$ is characterized by the following universal property:  For all $\infty$-categories
$\bD$ the induced map 
$$\Fun(\bC[W^{-1}],\bD)\to \Fun_{W^{-1}}(\bC,\bD)$$   is an equivalence, where
$ \Fun_{W^{-1}}(\bC,\bD)\subseteq  \Fun(\bC,\bD)$ is the subcategory of those functors which map arrows in $W$ to equivalences in $\bD$.
If necessary we use the notation $\iota$ for the functor $\bC\to \bC[W^{-1}]$. Sometimes we write the application of this functor to objects as $X\mapsto X_{\infty}$. 

\begin{ex}{\rm

In the present paper we mainly consider localizations which arise from
a category with weak equivalences,  i.e. a pair $(\bC,W)$ consisting of an ordinary category $\bC$ and a  collection of arrows $W$ of $\bC$. We consider the set of arrows as a subset $W\subseteq \Nerve(\bC)[1]$ of the one-simplices of the nerve of $\bC$.
 The $\infty$-category associated to the category with weak equivalences $(\bC,W)$ is the localization
$
\Nerve(\bC)[W^{-1}]
$
of $\Nerve(\bC)$.


If $\bC$ is an ordinary  category with a  collection of weak equivalences    $W$, then 
we have an equivalence between ordinary categories
$$\Ho(\Nerve(\bC)[W^{-1}])\simeq \bC[W^{-1}]\ ,$$
 although these may only be categories in a larger universe.
{Indeed, enlarging the universe if necessary, if $\bD$ is an ordinary category then, by the universal property of $\bC[W^{-1}]$, we see that
 \begin{eqnarray*}
\Nerve(\Fun(\bC[W^{-1}],\bD))&\simeq&\Nerve(\Fun_{W^{-1}}(\bC,\bD))\\&\stackrel{\eqref{lkdjlkqwjdldwqdqwd}}{\simeq}& \Fun_{W^{-1}}(\Nerve(\bC),\Nerve(\bD))\\&\simeq&\Fun(\Nerve(\bC)[W^{-1}],\Nerve(\bD))\\&\simeq& 
\Fun(\Nerve(\Ho(\Nerve(\bC)[W^{-1}])),\Nerve(\bD))\\&\simeq&
\Nerve(\Fun(\Ho(\Nerve(\bC)[W^{-1}]),\bD))\ .
 \end{eqnarray*}}
 
 }\hB \end{ex}
   



\begin{ex}{\rm In the category $\Top$ of topological spaces a continuous map $X\to Y$ {is} defined to be    a weak equivalence, if it induces a bijection between the sets of connected components  $\pi_{0}(X)\to \pi_{0}(Y)$ and isomorphisms of higher homotopy groups $\pi_{k}(X,x)\to \pi_{k}(Y,f(x))$ for every base point $x\in X$ and $k\ge 1$. 
Using this choice of weak equivalences we can then form the $\infty$-category $\Nerve(\Top)[W^{-1}]$.}\hB
\end{ex}

\begin{ex}\label{fjeklfjlewfewfewfewfewfewfefewfew3453453567}{\rm  Let $\sSet$ be the category of simplicial sets. We have a geometric realization functor $|-|:\Set\to \Top$. A morphism $f:X\to Y$ in $\sSet$ is defined to be a weak equivalence, if its geometric realization $|f|:|X|\to |Y|$  is a weak equivalence in $\Top$.
 \begin{ddd}\label{lkdjlkqwjwqlkdjwqdwqdwqd} The $\infty$-category of spaces is defined by $$\Spc:=\Nerve(\sSet)[W^{-1}]\ .$$\end{ddd}
   }\hB
\end{ex}

\begin{ex}\label{kdjqwdkljqlwdjlqwjdlwqdwqdwqdqwdwqdwqddqwdqd}{\rm 

In this example we present a    calculation of a colimit in $\Spc$.  The well-known result  will be used at various places in the present paper.

Let $c:\Set\to \sSet$ be the functor which maps a set $Z$ to the constant simplicial set $c(Z)\in \sSet$. Furthermore we let    $\iota:\Nerve(\sSet)\to \Spc$ be the localization map. 
To a set  $Z\in \Set$ we can now associate the space $Z_{\infty}:=\iota(c(Z)) \in \Spc$.
Applying this construction to the evaluations $X[n]\in \Set$ of a simplicial set $X\in \sSet$ for all $[n]\in \Delta$ we 
  form the simplicial space $X_{\infty}\in \Fun(\Nerve(\Delta)^{op},\Spc)$. We are interested in the colimit $\colim_{\Nerve(\Delta)^{op}} X_{\infty}$ in $\Spc$. It is calculated by the
   natural equivalence \begin{equation}\label{dweidweiduewoiduoiweudoiu239087e9237e93}
\iota(X)\simeq \colim_{\Nerve(\Delta)^{op}} X_{\infty}\ .
\end{equation}
Classically one would say that $X$ is weakly equivalent to
$ \hocolim_{\Delta^{op}} c(X)$. 
For an argument note that  homotopy colimit of the bisimplicial set $c(X)$   can be represented by its diagonal which isomorphic to $X$. Compare e.g. with \cite[Remark 2.1]{MR1982793}

 }\hB
\end{ex}

Let $\bC$ and $\bC^{\prime}$ be $\infty$-categories with distinguished sets of arrows $W$ and $W^{\prime}$.
If $f:\bC\to \bC^{\prime}$ is a functor  
  which carries the arrows in $W$ to the arrows in $W'$, then we obtain a functor
\[
f:\bC[W^{-1}]\longrightarrow \bC^{\prime}[W'^{-1}]
\]
between the localizations.

\begin{ex}\label{jhdqwkjdhqkdhqwkjdhwdwqdqwdwqd}{\rm The functor $\pi_{0}:\sSet\to \Set$ induces a functor between $\infty$-categories
$\pi_{0}:\Nerve(\sSet)\to \Nerve(\Set)$. Since by the definition of weak equivalences in $\Nerve(\sSet)$ (see Example \ref{fjeklfjlewfewfewfewfewfewfefewfew3453453567}) it sends weak equivalences to isomorphisms it induces a functor
$$\pi_{0}:\Spc=\Nerve(\sSet)[W^{-1}]\to \Nerve(\Set)\ .$$
We thus   have a well-defined set of connected component of a space.
A similar remark applies to higher homotopy groups.

\bigskip

One can use homotopy groups in order to detect equivalences in $\Spc$.
Indeed, a morphism $X\to Y$ in $\Spc$ is   an equivalence if and only if it induces a bijection
$\pi_{0}(X)\to \pi_{0}(Y)$ and isomorphisms $\pi_{k}(X,x)\to \pi_{k}(Y,f(x))$ for all $k\ge 1$ and $x\in X$.
}\hB
\end{ex}

\begin{ex}\label{lkejfewlkfewfwffewfwe34}
{\rm By the definition of the weak equivalences in  $\sSet$, the geometric realization functor $|-|:\sSet\to \Top$ preserves weak equivalences.
It therefore induces a morphism between $\infty$-categories \begin{equation}\label{udhiquwdhwqdqhwrwerr23rwrd}
\Spc=\Nerve(\sSet)[W^{-1}]\stackrel{ }{\to} \Nerve(\Top)[W^{-1}]\end{equation}
This morphism is actually an equivalence between $\infty$-categories  \begin{equation}\label{udhiquwdhwqdqhwd}\Spc\simeq \Nerve(\Top)[W^{-1}]\ .\end{equation}
}\hB
\end{ex}

\bigskip

\begin{ddd}\label{fkjwelfkwefwefwefewfewe3455}
An $\infty$-category $\bC$ is called pointed if it admits a zero object $*$, i.e. an object which is initial and final at the same time. \end{ddd}
If $\bC$ is pointed and admits finite {co}limits, then it has a {suspension} endofunctor \begin{equation}\label{wjqdlkqwdjwlkdjwkldwqdwqdwqdqwdq}
\Sigma:\bC\to \bC\end{equation}
which maps an object $X\in \bC$ to the {co}limit of the canonically associated diagram
$${*\leftarrow X\rightarrow *}\ .$$ The following definition can be found e.g. in \cite[1.1.1.9]{highalg}.
\begin{ddd}\label{fjewfkljewfejelwfewfewfwf}
An $\infty$-category is called stable, if it is pointed, admits finite limits, and if its  {suspension} endofunctor is an equivalence.
\end{ddd}
The homotopy category of a stable $\infty$-category has a natural structure of a triangulated category.

\bigskip

If $\bC$ is pointed and admits finite {co}limits
then we can define its stabilization. This is an $\infty$-category $\bC[{\Sigma}^{-1}]$ with a functor
\begin{equation}\label{fewfwfewkjfefkjewfhkewjhfkewfhkewhfkjefewfewfewfwef}
\bC\to \bC[{\Sigma}^{-1}] 
\end{equation}
which is initial among functors from $\bC$ to stable $\infty$-categories.
{If $\bC$ is presentable, then a stabilization is known to exist.}

\begin{ex}\label{dkqwjdlkqwdjldjldjqwdoi2u1oe21e12e}{\rm  
We can form the category of pointed spaces $\Spc_{*}$. It is pointed and admits all small limits and colimits. 
 \begin{ddd}\label{lekjflkewjlwfewf23i4p23i423o}
The category of  spectra is defined as the stabilization of the $\infty$-category of pointed spaces
$$\Sp:=\Spc_{*}[{\Sigma}^{-1}]\ .$$\end{ddd}
The triangulated homotopy category $\Ho(\Sp)$ is the  stable homotopy category.

\bigskip

The functor \eqref{fewfwfewkjfefkjewfhkewjhfkewfhkewhfkjefewfewfewfwef} in this case is denoted by $\Sigma^{\infty}:\Spc_{*}\to \Sp$ and called the suspension spectrum functor.  It is 
the left adjoint of an adjunction
$$\Sigma^{\infty}:\Spc_{*}\leftrightarrows \Sp:\Omega^{\infty}\ .$$
We furthermore have the following version of the suspension spectrum functor
\begin{equation}\label{qwjdkjqwdljqwkjdqlk89739423} \Sigma^{\infty}_{+}:\Spc\stackrel{X\mapsto X_{+}}{\to} \Spc_{*}\stackrel{\Sigma^{\infty}}{\to} \Sp\ .\end{equation}
\hB
In the construction of the Becker-Gottlieb transfer we use an alternative construction of the category of spectra. Its relation with the definition above is explained in some detail in Remark \ref{dqwjdklqwdwqwqdqdwq}.
}\end{ex}

\begin{ex}\label{lkjklfjqwlfjqlwkjfwqfwqfq}{\rm 
We let $\Ch$ denote the category of chain complexes of $\Z$-modules. We index our objects homologically
$$ \dots \xrightarrow{d}  A_{i+1}\xrightarrow{d} A_{i}\xrightarrow{d} A_{i-1}\xrightarrow{d}\dots\ .$$
and adopt the convention that $A^{i}:=A_{-i}$ in order to consider cochain complexes as chain complexes. By definition, a weak equivalence in $\Ch$ is a quasi-isomorphism, i.e. a morphism $X\to Y$ between chain complexes which induces an isomorphism in homology groups $H_{k}(X)\to H_{k}(Y)$ for all $k\in \Z$. The resulting $\infty$-category
$\Nerve(\Ch)[W^{-1}]$ is  stable. Its triangulated homotopy category is the derived category of the abelian category of abelian groups.
}\hB
\end{ex}

\bigskip

\begin{ex}\label{fjfwlefjlwefjlewfewfewffewfewfef}{\rm

We must be able to calculate certain limits and colimits  in $\Ch[W^{-1}]$. For a chain complex $A\in \Nerve(\Ch)$ let $A_{\infty}$ be its image in $\Nerve(\Ch)[W^{-1}]$.
 If $$A[-]:=([n]\mapsto A[n]) \in \Fun(\Nerve(\Delta),\Nerve(\Ch))$$  is a cosimplicial chain complex, then we apply this notation objectwise and get the cosimplicial object
$$A[-]_{\infty}:=([n]\mapsto A[n]_{\infty})\in \Fun(\Nerve(\Delta),\Nerve(\Ch)[W^{-1}])$$ of   $\Nerve(\Ch)[W^{-1}]$.
We want to calculate the limit $\lim_{\Nerve(\Delta)}A[-]_{\infty}$ in $\Nerve(\Ch)[W^{-1}]$.

\bigskip

By the Dold-Kan correspondence a cosimplicial group $G[-]\in \Fun(\Delta^{op},\Ab)$ gives rise to a chain complex  $DK(G[-])$ such that $DK(G[-])^{n}=G[n]$.  Its differential is given by the usual formula which we will not write out here.
Similarly, the cosimplicial chain   complex $A[-]$ gives rise to a double complex which we also denote by
$DK(A[-])$. In bidegree $(n,m)$ it is given by $DK(A[-])^{m,n}=A[m]^{n}$. 
We can form  the total complex   $\Tot(  A[-])  \in \Nerve(\Ch)$ of this double complex using products so that $$\Tot(A[-])^{j}:=\prod_{k+l=j} A[k]^{l}\ .$$ 

The limit in question is calculated by 
  the equivalence \begin{equation}\label{dqwhqwdiwjhdwqjdwqdijwqdoijwqdoiwqdu987987983243424}
\lim_{  \Nerve(\Delta)} A[-]_{\infty}\simeq \Tot(A[-])_{\infty}
\end{equation}
in $\Nerve(\Ch)[W^{-1}]$. We refer to \cite[Problem 4.23]{skript} for an argument.

\bigskip

A similar statement holds true for colimits of simiplicial chain complexes. We   consider simplicial chain complex $B[-]\in \Fun(\Nerve(\Delta)^{op},\Nerve(\Ch))$ and want to calculate
the colimit $\colim_{  \Nerve(\Delta)^{op}} B[-]_{\infty}$ in $\Nerve(\Ch)[W^{-1}]$.

\bigskip

The version of the Dold-Kan correspondence for simplicial chain complexes associates to a simplicial  abelian group $G$ the chain complex $DK(G)$ with $DK(G)^{n}:=G[-n]$. Similarly, for a simplicial chain complex $B[-]$  we get a double complex $DK(B[-])$. Then we can form the total complex $\mathrm{tot}(B[-])$ of the associated double complex using direct sums, i.e.
$$\mathrm{tot}(B[-])^{i}:=\bigoplus_{k+l=i} B[-k]^{l}\ .$$ The colimit in question is calculated by 
 the equivalence
\begin{equation}\label{dqwhqwdiwjhdwqjdwqdijwqdoijwqdoiwqdu9879879832434241}\colim_{  \Nerve(\Delta)^{op}} B[-]_{\infty}\simeq \mathrm{tot}(B[-])_{\infty}\ .\end{equation}
The proof is simpler, see e.g. \cite[Problem 4.24]{skript}

}
\hB\end{ex}

\begin{ex}{\rm 
We can consider the one-category of small categories $\Cat$. After enlarging the universe it can again be considered as a small category. The morphisms in $\Cat$ are functors. The collection of weak equivalences in $\Cat$ consists of functors which extend to  equivalences of categories. Then we define the $\infty$-category $\Nerve(\Cat)[W^{-1}]$ of categories. We have a localization map
$$\Nerve(\Cat)\to \Nerve(\Cat)[W^{-1}]\ , \quad \bC\mapsto \bC_{\infty}\ .$$
The nerve functor $\Nerve:\Cat\to \sSet$ preserves weak equivalences and therefore descends to a functor \begin{equation}\label{dkqwjkldjwqldjklwqjd8897193123123123}
\Nerve:\Nerve(\Cat)[W^{-1}]\to \Nerve(\sSet)[W^{-1}]=\Spc\ .\end{equation}
In this way we can use categories in order to present spaces.

For example, let $G$ be a group which we consider as a category $\mathbf{G}$ with one object.
Then the space \begin{equation}\label{dqwkdwjkwqhdqwkdwqdwd789}BG:=\Nerve(\mathbf{G}_{\infty})\end{equation}
is the classifying space of $G$.
A more general version of this construction is a basic ingredient in the construction of the $K$-theory spectrum in Subsection \ref{jun27333}.
}\hB
\end{ex}


\bigskip

For two objects $x,y$ in an $\infty$-category $\bC$ we write 
\begin{equation}\label{mar0701}
\Map(x,y)\in \Spc\end{equation}
for the mapping space between {$x$ and $y$. In terms of the mapping space we can describe the set of  morphisms between two objects $x,y$ in   the homotopy category $\Ho(\bC)$ as follows:
\begin{equation}\label{jjkqdhwkqdkhjwqdqwd}
\Hom_{\Ho(\bC)}(x,y)\cong\pi_{0}(\Map(x,y))\ .\end{equation}

If $\bC$ is stable, then we can define the mapping spectrum $$\map(x,y)\in \Sp\ .$$
The mapping spectrum and the mapping space are related by a natural equivalence
$$\Omega^{\infty}\map(x,y)\simeq \Map(x,y)$$
(where we implicitly forget the base point on the left-hand side).
One can consider the mapping space and mapping spectrum as functors \begin{equation}\label{kjdkqwjdlqkwdjlwqdqwdqwd}
\Map:\bC^{op}\times \bC\to \Spc\ , \quad \map:\bC^{op}\times \bC\to \Sp\ .\end{equation}
}

\begin{ex}\label{jcsdkjcdksjhcdkschkdscer}{\rm 
For a chain complex $C\in \Ch$ let $C_{\infty}\in \Nerve(\Ch)[W^{-1}]$ denote its image under the localization $\Nerve(\Ch)\to \Nerve(\Ch)[W^{-1}]$. 
We consider the abelian group $\Z$ as a chain complex in degree zero and obtain $\Z_{\infty}\in \Nerve(\Ch)[W^{-1}]$.

\begin{ddd}\label{jdjwqldjwqldjwqljdwqopiop12e}
We define the Eilenberg-MacLane functor
$$H:\Nerve(\Ch)[W^{-1}]\to \Sp\ , \quad C\mapsto \map(\Z_{\infty},C)\ .$$
\end{ddd}
In order to make the definition of $H$ as a functor precise one uses \eqref{kjdkqwjdlqkwdjlwqdqwdqwd}.
It immediately follows from this definition that $H$ preserves all small limits.

For a chain complex $A\in \Ch$ and $k\in \Z$ we have a canonical isomomorphism
$$\pi_{k}(H(A))\cong H_{k}(A)$$ of abelian groups. 
In order to show this we must be able to calculate the mapping spectra in $\Nerve(\Ch)[W^{-1}]$. 
At this point it is convenient to use a model for $\Nerve(\Ch)[W^{-1}]$ based on $dg$-nerves
$$\Nerve(\Ch)[W^{-1}]\simeq \Nerve^{dg}(\Ch_{dg}^{free})\ .$$
Here $\Ch_{dg}^{free}$ is the dg-category of chain complexes of free abelian groups.
Then for $k\in \Z$ and $X,Y\in \Ch_{dg}^{free}$ we have isomorphisms of abelian groups
$$\pi_{k}(\map(X,Y))\cong H_{k}(\Hom_{dg}(X,Y))\ .$$
Since $\Z$ is a free abelian group, if $A^{\prime}\to A$ is a quasi-isomorphism with $A^{\prime}\in \Ch^{free}_{dg}$, then we get
$$\pi_{k}(\map(\Z_{\infty},A))\cong H_{k}(\Hom_{dg}(\Z,A^{\prime}))\cong H_{k}(\Hom_{dg}(\Z,A))\cong H_{k}(A)\ .$$

 }\hB
\end{ex}

\begin{ex}\label{dlkqjwdlkjqwljqwdqwdqwdqw} {\rm A particularly nice class of $\infty$-categories is given by  the presentable $\infty$-categories, see \cite[Def. 5.5.0.1]{MR2522659}. 
A presentable $\infty$-category admits all small limits and colimits. 
One can characterize the $\infty$-category of spaces $\Spc$ as the universal presentable $\infty$-category generated by the one-point category $*$. If $X\in \Spc$, then we  have an equivalence \begin{equation}\label{cweckjhcjewkc879234324}
X\simeq \colim_{X}*\ .
\end{equation}

A presentable $\infty$-category  it is tensored and cotensored over the category of spaces. This means that we have functors
$$\bC\times \Spc\to \bC\ ,\quad (C,X)\to C\otimes X$$
 and \begin{equation}\label{dkjkjwqlwjwljdwldjwldjwqdq546456234}\Spc^{op}\times \bC\to \bC\ , \quad (X,C)\mapsto C^{X}\ .\end{equation}
 For a fixed space $X$ these functors fit into the adjunction
 $$-\otimes X:\bC\leftrightarrows \bC:(-)^{X}\ .$$
 In particular, $(-)^{X}$, as a right-adjoint, preserves limits.  
 
   For a given object $C\in \bC$ we have the adjunction
 $$C\otimes (-):\Spc\leftrightarrows \bC:\Map(C,(-)^{X})\ .$$
 This shows that
 $C\otimes (-)$ preserves colimits.  For $X\in \Spc$ and using $C\otimes *\simeq C$ we have by \eqref{cweckjhcjewkc879234324}  
 $$C\otimes X\simeq \colim_{X}C\ .$$

 In fact, {the}  functor $C\otimes (-)$ is 
 the left Kan-extension of the functor $*\to \bC$ given by $C$ along $*\to \Spc$.
 
Similarly,  for a given object $D\in \bC$ we have the adjunction
 $$\Map(-,D):\bC\leftrightarrows\Spc^{op}:D^{(-)}$$
 showing that  the contravariant functor $D^{(-)}$ maps colimits in $\Spc$ to limits in $\bC$. This functor
 is the right Kan extension of the functor $*^{op}\to \bC$ given by $D$ along $*^{op}\to \Spc^{op}$. 
 Using \eqref{cweckjhcjewkc879234324} for $C\in \bC$ and $X\in \Spc$ we get the equivalence
 \begin{equation}\label{chcjkwhkecheckehkeecewcewcewc}
C^{X}\simeq \lim_{X} C\ .
\end{equation}

  If the presentable $\infty$-category $\bC$ is stable (Definition \ref{fjewfkljewfejelwfewfewfwf}),   then the tensor and cotensor structures
  over $\Spc$ extend to tensor and cotensor structures over $\Sp$ along the suspension spectrum functor $\Sigma^{\infty}_{+}:\Spc\to \Sp$, see \eqref{qwjdkjqwdljqwkjdqlk89739423}. This structure fits into similar adjunctions and has the corresponding limit- and colimit preserving properties.
 
Most of the $\infty$-categories considered in the present paper are presentable. Here is a list:
\begin{enumerate}
\item $\Spc$, $\Sp$, $\Nerve(\Set)$, $\Nerve(\sSet)$, \item $\Nerve(\Ch)$, $\Nerve(\Ch)[W^{-1}]$, \item $\Nerve(\Cat)$, $\Nerve(\Cat)[W^{-1}]$, \item $\CAlg(\bC)$ for a presentable symmetric monoidal category (see Subsection \ref{mar0601}), so in particular \item $\CommMon(\Set)$, $\CommMon(\Spc)$, \item and also
$\CommGroups(\Set)$, $\CommGroups(\Spc)$. \end{enumerate}

}\hB
\end{ex}

\subsection{Commutative algebras and monoids} \label{mar0601}
 
 \newcommand{\Comm}{{\mathbf{Comm}}}
 \newcommand{\Ass}{{\mathcal{A}ss}}
 \newcommand{\cO}{{\mathcal{O}}}
 \newcommand{\Alg}{\bf{Alg}}
 
 In this subsection we review symmetric monoidal categories, commutative algebras and monoids, and the group completion.

 \bigskip

  A symmetric monoidal   $\infty$-category  is defined as a cocartesian fibration
$\bC^{\otimes}\to \Comm^{\otimes}$ over the commutative $\infty$-operad  {$\Comm^{\otimes}$}, see \cite[Def. 2.0.0.7]{highalg}. 
A symmetric monoidal   structure on an $\infty$-category $\bC$ is then an embedding of $\bC$ into a {cocartesian fibration $\bC^\otimes\to\Comm^\otimes$ such that we have a pull-back}
$$\xymatrix{\bC\ar[d]\ar[r]
&\bC^{\otimes}\ar[d]\\
\mbox{}*\ar[r]
&\Comm^{\otimes}}
\ .$$

The $\infty$-category $\CAlg(\bC)$ of commutative {algebras} in $\bC^\otimes$ is  defined as a full subcategory of the $\infty$-category $\Fun_{\Comm^\otimes}(\Comm^{\otimes},\bC^{\otimes})$ of sections of $\bC^\otimes\to\Comm^\otimes$ (see  \cite[Def. 2.4.2.1]{highalg} for details). 
\begin{rem}{\rm A more consistent notation would be $\CAlg_{\Comm^{\otimes}}(\bC^{\otimes})$ in order to indicate the dependence of this notion on the symmetric monoidal structure but we prefer the simpler symbolics. But one must be careful if there is more than one possible structure as e.g. in Example \ref{flkwejweklfjelkewfewfw4324234234}.} \hB\end{rem}
The forgetful functor
$$\CAlg(\bC)\to  \bC$$ which forgets the commutative algebra structure is then given by pullback along
$$*\to   \Comm^{\otimes}\ .$$
{If $\bC$ is an $\infty$-category with finite products, then $\bC$ naturally determines a symmetric monoidal $\infty$-category $\bC^\times\to\Comm^\otimes$ in which the multiplication is the cartesian product.
In this case, we may write
\[
\CommMon(\bC):=\CAlg(\bC) \]
and refer to the commutative algebra objects as commutative monoids.}

\begin{ex}\label{flkwejweklfjelkewfewfw4324234234}{\rm 
We consider the symmetric monoidal category of abelian groups $\Ab^{\otimes}$ with the usual tensor product and the corresponding symmetric monoidal $\infty$-category $\Nerve(\Ab^{\otimes})$.
Then $\CAlg(\Nerve(\Ab^{\otimes}))$ is equivalent to the $\infty$-category $\Nerve(\Comm\Rings)$  of  commutative rings.

In contrast, if we would consider the symmetric monoidal structure $\oplus$ on $\Ab$, then the forgetful functor $\CAlg(\Nerve(\Ab^{\oplus}))\to \Nerve(\Ab)$ is an equivalence. So every abelian group is naturally a commutative algebra with respect to this symmetric monoidal structure. This simply  reflects the fact that for $A\in \Ab$
$$A\oplus A\to A\ , \quad (a,b)\mapsto a+b$$
is a homomorphism of abelian groups.

}\hB
\end{ex}

\begin{ex}\label{kdjqwkldjqwdljlwqwqd1341323213213123}{\rm
The tensor product  of chain complexes induces a symmetric monoidal structure on $\Nerve(\Ch)$. 
Commutative algebras in $\Nerve(\Ch)$ are commutative differential graded algebras.

Since the tensor product in general does not preserve weak equivalences it does not directly descend to $\Nerve(\Ch)[W^{-1}]$. In order to resolve this problem we consider the  subcategory  of chain complexes $\Ch^{free}\subset \Ch$ of free abelian groups.
Then the inclusion $\Ch^{free}\to \Ch$ induces an equivalence of $\infty$-categories
\begin{equation}\label{kdqjwkdjqwlkdjlqwdqwdqwdqw}\Nerve(\Ch^{free})[W^{-1}]\stackrel{\simeq}{\to} \Nerve(\Ch)[W^{-1}]\ .\end{equation} The tensor product on $\Ch^{free}$ directly descends to $\Nerve(\Ch^{free})[W^{-1}]$, and we define the symmetric monoidal structure on $\Nerve(\Ch)[W^{-1}]$ by transport over the equivalence
\eqref{kdqjwkdjqwlkdjlqwdqwdqwdqw}.

Note that {the image} $\Z_{\infty}$ {of $\Z$ under the localization $\Nerve(\Ch)\to \Nerve(\Ch)[W^{-1}]$} is the tensor unit of the symmetric monoidal structure on $\Nerve(\Ch)[W^{-1}]$. 
Consequently we get a commutative algebra
$$H\Z:=\map(\Z_{\infty},\Z_{\infty})\in \CAlg(\Sp)\ .$$
By construction, the Eilenberg-MacLane functor (Definition \ref{jdjwqldjwqldjwqljdwqopiop12e})
has a factorization
$$\Nerve(\Ch)[W^{-1}]\stackrel{H_{\Z}}{\to} \Mod(H\Z)\to \Sp\ ,$$
where the second map takes the underlying spectrum of a $H\Z$-module. 
It is known that $H_{\Z}$ is an equivalence of $\infty$-categories, see \cite{MR2306038},  \cite[Proposition 7.1.2.7]{highalg}.

\bigskip

In Subsection \ref{aug1006} we will consider a rational version of this. Let $\Ch_{\Q}$ be the category of chain complexes of rational vector spaces. In this case the tensor product $\otimes_{\Q}$ preserves weak equivalences and therefore directly descents from $\Nerve(\Ch_{\Q})$ to a symmetric monoidal structure on $\Nerve(\Ch_{\Q})[W^{-1}]$.
We
 consider the commutative ring spectrum $$H\Q:=\map(\Q_{\infty},\Q_{\infty})\in  \CAlg(\Sp)\ .$$ Then we  have an equivalence
$$H_{\Q}:\Nerve(\Ch_{\Q})[W^{-1}]\stackrel{\simeq}{\to} \Mod(H\Q)\ , \quad C\mapsto \map(\Q_{\infty},C)\ .$$
If we equip the domain and the target with their natural symmetric monoidal structures $\otimes_{\Q}$ and $\wedge_{H\Q}$, then $H_{\Q}$ refines to an equivalence between symmetric monoidal categories.

The same applies to the analog $H_{\R}$ for chain complexes of real vector spaces.
}
\hB
\end{ex}

\begin{ex}{\rm 
We consider the $\infty$-category $\Nerve(\Cat)[W^{-1}]$ with its cartesian symmetric monoidal structure. Symmetric monoidal categories present objects in $\CAlg(\Nerve(\Cat)[W^{-1}])$.
To make this precise we first observe that we have symmetric monoidal functors
$$\Nerve(\Cat)\to \tilde\Nerve(\Cat_{2})\to \Nerve(\Cat)[W^{-1}])\ ,$$
where $\Cat_{2}$ is the two-category of categories, functors and natural isomorphisms,  and $\tilde \Nerve$ is an appropriate nerve functor.
We have induced functors between the categories of commutative algebras
$$\CAlg(\Nerve(\Cat))\to \CAlg(\tilde\Nerve(\Cat_{2}))\to \CAlg(\Nerve(\Cat)[W^{-1}])\ .$$
A commutative monoid in $\Set$ like $\nat$ can be considered as a symmetric monoidal category with only identity morphisms and hence as an object of $\CAlg(\Nerve(\Cat))$. More generally,   objects of $\CAlg(\Nerve(\Cat))$ are   symmetric monoidal categories whose associator and symmetry constraints are identities. A general symmetric monoidal one-category is naturally an object of $\CAlg(\tilde\Nerve(\Cat_{2}))$. All these examples provide objects in $\CAlg(\Nerve(\Cat)[W^{-1}])$.

The nerve functor \eqref{dkqwjkldjwqldjklwqjd8897193123123123} preserves the cartesian structure and therefore induces a functor
$$\CAlg(\Nerve(\Cat)[W^{-1}])\to \CAlg(\Spc)=\CommMon(\Spc)\ .$$
Using this functor we can present commutative monoids in spaces by symmetric monoidal categories.

If $A$ is an abelian group, then we get a commutative monoid $A_{\infty}\in \CommMon(\Spc)$.
We can also view $A$ is a symmetric category $A^{\delta}$ with only identity morphisms\footnote{We use the superscript $\delta$ in order to avoid a conflict with the notation $A_{\infty}$ used in the next paragraph.}. Then we have the equivalence
$$\Nerve(A^{\delta}_{\infty})\simeq A_{\infty}\ .$$
On the other hand we can consider $A$ as a symmetric monoidal category $A$ with one object. This symmetric monoidal category belongs to $\CAlg(\Nerve(\Cat))$ and induces 
an object $A_{\infty}\in \CAlg(\Nerve(\Cat)[W^{-1}])$.
The additional  symmetric monoidal structure thus refines the classifying space \eqref{dqwkdwjkwqhdqwkdwqdwd789} of $A$ to
 a commutative monoid (actually a group, see below) in spaces $BA\simeq \Nerve(A_{\infty})\in \CommMon(\Spc)$.
}\hB
\end{ex}

\bigskip

In order to discuss the group completion we
  specialize to the category $\bC=\Spc$ with the cartesian symmetric monoidal structure. We have a symmetric monoidal functor $$\pi_{0}:\Spc\to \Nerve(\Set)\ .$$ If $X\in \CommMon(\Spc)$ is a monoid in the $\infty$-category of simplicial sets, then
$\pi_{0}(X)$ is a monoid in the ordinary sense. 
\begin{ddd}\label{djqwldjlqwdqwd324234234} The $\infty$-category \begin{equation}\label{qwegdqjhwgjqwhdgwqjgdjwqdqwdqwqd34}
\CommGroups(\Spc)\subset \CommMon(\Spc) 
\end{equation}of commutative groups in $\Spc$ is defined as the full subcategory of $\CommMon(\Spc)$ of those commutative monoids $X$ which have  the property that $\pi_{0}(X)$ is a group.\end{ddd}

\begin{lem} The inclusion  of commutative groups into commutative monoids is part of an adjunction
\begin{equation}\label{khfkqhwfkhqwkfhqwkfhkhwqfqw34124124}\Omega B:\CommMon(\Spc)\leftrightarrows\CommGroups(\Spc):inclusion\ .\end{equation}
\end{lem}
\proof  {For a proof  of this fact we could refer to  the arguments leading to \cite[Cor. 4.4]{2013arXiv1305.4550G}. But for educational reasons we will give an argument\footnote{The first author thanks Thomas Nikolaus for the kind explanation.} which points out the essential feature of the problem.}

The $\infty$-category of commutative monoids in $\Spc$ is  presentable. 
We will characterize the subcategory of commutative groups as an accessible localizations.
{To this end we must characterize commutative groups as commutative monoids which are local with respect to a set of morphisms. In the present case this set will contain one element.}
We have an adjunction $$\Free:\Spc\leftrightarrows\CommMon(\Spc):forget\ .$$
We let $F(a,b):=\Free(\{a,b\}_{\infty})$ be the free commutative monoid on the space with two points $a,b$. Using the universal property of the free commutative  monoid we can  consider the shearing map
$$s:F(a,b)\to F(a,b)\ , \quad (a,b)\mapsto (a,ab)\ .$$
We claim that
a commutative monoid $M\in \CommMon(\Spc)$ is a commutative
group if an only if $$\Map(s,M):\Map(F(a,b),M)\to \Map(F(a,b),M)$$ is an equivalence of spaces, where the mapping spaces are taken   in $\CommMon(\Spc)$.
Assuming the claim we can conclude that
$\CommGroups(\Spc)$ is the full subcategory of objects which are local with respect to the shearing map. Using the presentability of $\CommMon(\Spc)$ we conclude that 
 the inclusion \eqref{qwegdqjhwgjqwhdgwqjgdjwqdqwdqwqd34}
is the right-adjoint of a localization.

We now show the claim. Using the universal property of the free commutative monoid we have an equivalence $\Map(F(a,b),M)\simeq M\times M$ of spaces under which   $\Map(s,M)$ goes to  the usual shearing map $M\times M\to M\times M$ for the commutative monoid $M$. If we apply $\pi_{0}$, then we get the shearing map on the level of commutative monoids $\pi_{0}(M)$ in sets. It is an isomorphism if and only if $\pi_{0}(M)$ is a group.
So if $\Map(s,M)$ is an equivalence of spaces, then $\pi_{0}(M)$ is a group. 
Now assume that $\pi_{0}(M)$ is a group. In order to show that
the shearing map $s:M\times M\to M\times M$ is an equivalence we must show that it induces an isomorphism in higher homotopy groups (see Example \ref{jhdqwkjdhqkdhqwkjdhwdwqdqwdwqd}). We fix a section $\phi:\pi_{0}(M)_{\infty}\to M$ of the canonical map  $[-]:M\to \pi_{0}(M)_{\infty}$. Then we consider the composition
\begin{equation}\label{dhwkjdhkwqdqwdqwdqwdqwd}M\times M\stackrel{s}{\to} M\times M\xrightarrow{(x,y)\mapsto (x,\phi([x]^{-1})y} M\times M\ .\end{equation}
Now  we have a decomposition of abelian groups $$\pi_{k}(M\times M,(x_{0},y_{0}))\cong \pi_{k}(M,x_{0})\times \pi_{k}(M,y_{0})\ .$$
Under this decomposition
the morphism \eqref{dhwkjdhkwqdqwdqwdqwdqwd} has the structure
$$\left(\begin{array}{cc}\id&0\\a&\kappa\end{array}\right):\left(\begin{array}{c}\pi_{k}(M,x_{0})\\\pi_{k}(M,y_{0})\end{array}\right)\to \left(\begin{array}{c}\pi_{k}(M,x_{0})\\\pi_{k}(M,\phi([x_{0}]^{-1})x_{0}y_{0})\end{array}\right) \ , $$ where $\kappa$ is an isomorphism.
 Consequently \eqref{dhwkjdhkwqdqwdqwdqwdqwd} is an equivalence. A similar discussion applies to the composition of the two morphisms in  \eqref{dhwkjdhkwqdqwdqwdqwdqwd}  in the other order.
 \hB

\begin{ddd}
The functor $\Omega B$ is called the group completion functor
\end{ddd}

{
We will frequently use the fact that the $\infty$-category $\CommGroups(\Spc)$ is a model for the
$\infty$-category of connective spectra. More precisely, the  $\infty$-loop space functor is part of an adjoint equivalence of $\infty$-categories
\begin{equation}\label{feb1703} \spp  :\CommGroups(\Spc)\leftrightarrows \Sp_{\ge 0}:\Omega^{\infty}\ .
\end{equation}}
 
{
\begin{ex} {\rm 
We have two functors $$\Phi,\Psi:\Nerve(\Ab)\to \Sp\ .$$
The first is the composition
$$\Phi:\Nerve(\Ab)\simeq\CommGroups(\Set)\stackrel{(-)_{\infty}}{\to} \CommGroups(\Spc)\stackrel{\spp}{\to} \Sp\ ,$$
and the second is given by
$$\Psi:\Nerve(\Ab)\to \Nerve(\Ch)\to \Nerve(\Ch)[W^{-1}]\stackrel{H}{\to}\Sp\ .$$
\begin{lem}\label{fjelfjwkfjekjfwlipo31}
There is a natural equivalence $\Phi\simeq \Psi$.
\end{lem}
\proof
This is \cite[Remark 2.13]{buta}. \hB

} 
\end{ex}
}

\subsection{Smooth objects}\label{smoothobjects}

In this subsection we introduce the notion of  smooth objects of an $\infty$-category, the descent condition and the basic construction of smooth objects as smooth function objects.

\bigskip

 {In the present paper a smooth manifold is a smooth manifold with corners in the sense of  \cite[3.3.2-3.3.6]{MR2191484}.
We require that the underlying topological space of a smooth manifold is metrizable and second countable.
The model for a corner of codimension $n $ in a manifold of dimension $n+m$ is $[0,\infty)^{n}\times \R^{m}$. 
 We require that the transition functions of an atlas preserve the product structure so that we have a product structure near the boundary of a smooth manifold. 
 Moreover, a smooth map between manifolds is required to be compatible with the product structure    \cite[ 3.3.6]{MR2191484}.
 This convention allows us to define the notion of a product structure for connections and metrics on bundles over a manifold, and it ensures that the pull-back of such structures along maps between manifolds preserves such product structures.
We let $\Mf$ denote the category of smooth manifolds and smooth maps.
 The category $\Mf$ contains intervals $[0,1]$, rays $[0,\infty)$, and the simplices $\Delta^{n}$ for $n\in \nat$ (see \cite[3.3.8]{MR2191484}). Furthermore, it is closed under forming products.}

\bigskip

Let $\bC$ be some $\infty$-category. In the present paper we consider the category of presheaves on the category of smooth manifolds $\Mf$ with values in $\bC$. We will call such presheaves smooth objects in $\bC$.
\begin{ddd}
The $\infty$-category of smooth objects in $\bC$ is defined by 
$$\Sm(\bC):=\Fun(\Nerve(\Mf^{op}),\bC)\ .$$
\end{ddd}

\begin{ex}\label{dqjlwdjqwdwqdwqdwqdqwdwq3432423423}{\rm
Given an object $C\in \bC$ we can consider the constant smooth object $$\underline{C}\in \Sm(\bC)\ , \quad M\mapsto  C 
\ , \quad  (f:M\to M^{\prime})\mapsto \id_{C}   \ .$$
We get a constant smooth object functor
$$\underline{(-)}:\bC\to\Sm(\bC)\ .$$
}\hB
\end{ex}

\begin{ex}{\rm
If $\bC$ is a one-category, then in view of \eqref{lkdjlkqwjdldwqdqwd} smooth objects in $\Nerve(\bC)$ are nothing else than contravariant functors from $\Mf$ to $\bC$, i.e. we have an equivalence of $\infty$-categories
\begin{equation}\label{fhewfhwefhewfiuizuii3u3244234234ewddewd}
\Nerve(\Fun(\Mf^{op},\bC))\simeq \Sm(\Nerve(\bC))\ .
\end{equation}
For example, smooth abelian groups  are classically called presheaves of abelian groups. But note that these presheaves are not just defined on the open subsets of a given manifold, but on all smooth manifolds.
}\end{ex}

\begin{ex}\label{fklwjklewfwfwf45453453453wfwfwf}{\rm 
  A  typical example is the presheaf  of rings $C^{\infty}$ which associates to a smooth manifold the ring of real-valued smooth functions $C^{\infty}(M)$. This presheaf is represented by the manifold $\R$ with its usual algebraic structures, i.e. $C^{\infty}(M)\cong \Hom_{\Mf}(M,\R)$.
In the notation introduced above we have $C^{\infty} \in \Sm(\CAlg(\Nerve(\Ab)))$. This object satisfies descent in the sense of Definition \ref{jun215111} below.
}\hB
\end{ex}

\begin{ex}\label{djkjwlkdjlwqdqwdqwdwqdwdwqd}{\rm 
Let $R$ be a ring and consider the one-category $\Mod(R)$ of $R$-modules. For an $R$-module $V\in \Mod(R)$ we can consider the smooth $R$-module $\cV$ on $\Mf$ which associates to a manifold $M$ the $R$-module of continuous functions $M\to V$, where we consider $V$ as a discrete space. This   
 is first of all an object of $\Sm(\Nerve(\Mod(R)))$.
In addition it satisfies descent   the sense of Definition \ref{jun215111} and is homotopy invariant in the sense of Definition \ref{efjwflewjfl239874297329cjsajkas}.   
\bigskip

The symmetric monoidal stack $\loc(R)$ of sheaves of $R$-modules {introduced in Definition \ref{fjewjflwefjlwfoiiouoiuoifwefwefwef}} is an object of $\Sm(\CAlg(\Nerve(\Cat)[W^{-1}]))$. Being a stack is again an additional descent property.

\bigskip

For example, if the $R$-module  $V$ above is finitely generated, then
$\cV_{|M}\in \loc(R)(M)$ for every manifold $M\in \Mf$.

\bigskip

Finally, we can form the smooth monoid $\bloc(R):=\pi_{0}(\loc(R))$  of isomorphism classes in $\loc(R)$. It is an object of $\Sm(\CommMon(\Set))$.
}\hB 
\end{ex}

\begin{ex}{\rm If $\bC$ is a presentable $\infty$-category (Example \ref{dlkqjwdlkjqwljqwdqwdqwdqw}), then we have a very natural method to construct smooth objects in $\bC$.
We consider the composition
$$\Nerve(\Mf)\stackrel{M\mapsto M_{top}}{\to} \Nerve(\Top)\stackrel{(-)_{\infty}}{\to} \Nerve(\Top)[W^{-1}]\stackrel{\eqref{udhiquwdhwqdqhwd}}{\simeq} \Spc\ , \quad  M\mapsto M_{top,\infty}\ .$$ If we precompose this composition with
the cotensor structure \eqref{dkjkjwqlwjwljdwldjwldjwqdq546456234}, then we obtain the functor
\begin{equation}\label{72938djkekdj83e2e2}\Nerve(\Mf^{op})\times \bC\to \bC\ . \end{equation}
\begin{ddd}
We define the smooth function object functor
$$\Funk:\bC\to \Sm(\bC)$$
as the adjoint of \eqref{72938djkekdj83e2e2}.
\end{ddd}
Thus for $C\in \bC$ and $M\in \Mf$ we have \begin{equation}\label{lfkjlekfjewlfjewlfjfelwjfljlekwf9038409832423424}\Funk(C)(M)\simeq C^{M_{top,\infty}}\ .\end{equation}
Note that we often abbreviate the notation and write $C^{M}$ for this object.
}\hB\end{ex}


 \begin{ex}{\rm 
The smooth object  $\cV\in \Sm(\Nerve(\Mod(R)))$ considered in Example \ref{djkjwlkdjlwqdqwdqwdwqdwdwqd}
can be written in the form $\cV\simeq \Funk(V)$.
}\hB\end{ex} 


 \begin{ex}{\rm 
If $\bC=\Spc$ and $X\in \Spc$, then $\Funk(X)(M)$ has the homotopy type of the space of continuous maps from $M$ to $X$. So for example $\Funk(X)(S^{1}) $ has the homotopy type of the free loop space of the space $X$.

Similarly, for a group $G$ we can identify $\pi_{0}(\Sm(BG)(M))$ (see \eqref{dqwkdwjkwqhdqwkdwqdwd789} for $BG$) with the set of isomorphism classes of 
$G$-principal bundles over $M$.
}\hB\end{ex} 


 \begin{ex}{\rm 
In the case $\bC=\Sp$ and for $E\in \Sp$ the spectrum  $\Funk(E)(M)$ has the homotopy type of the function spectrum from $M$ to $E$. In particular, for every $k\in \Z$ we have a natural  
isomorphism of abelian groups
$$\pi_{k}(\Funk(E)(M))\cong E^{-k}(M)\ .$$
}\hB 
\end{ex}

\begin{rem}\label{ijlkdjqwldqwdqwdwqdwqdqwd}{\rm
We assume that $\bC$ is a presentable and   symmetric monoidal $\infty$-category.  Then the smooth function object functor
$$\Funk: \bC\to \Sm(\bC)$$
is lax symmetric monoidal as we will explain in the following.
Let $\Fun^{lim}(\Spc^{op},\bC)$ be the $\infty$-category of limit preserving functors. We consider the functor
$$\bC\to\Fun^{lim}(\Spc^{op},\bC)\ , \quad C\mapsto  C^X \ .$$ Since $\Spc$ is freely generated under colimits by the one-point space (see Example \ref{dlkqjwdlkjqwljqwdqwdqwdqw}), this functor is actually  an equivalence, hence symmetric monoidal.
 We can now write $\Funk$ as the composition of lax symmetric monoidal functors
 $$\bC\simeq \Fun^{lim}(\Spc^{op},\bC)\to \Fun (\Spc^{op},\bC)\to \Sm(\bC)\ ,$$
 where the second is the natural inclusion and the last is the restriction along $\Nerve(\Mf)\to \Spc$ given by  $M\mapsto M_{top,\infty}$.

\bigskip

We consider objects $C,D\in \bC$ and spaces $X,Y\in \Spc$.
The lax symmetric monoidal structure on the power functor $C\mapsto C^{X}$ provides  a natural transformation
\begin{equation}\label{kfjkfjklwejfwejfewopipoiopiopieqweqeqweeqeeqwe}
C^{X}\otimes D^{Y}\to (C\otimes D)^{X\times Y}\ .
\end{equation}
This map is defined as the composition$$
 C^{X}\otimes D^{Y} \to  C^{X\times Y}\otimes D^{Y}
\to C^{X\times Y}\otimes D^{X\times Y}
\to (C\otimes D)^{X\times Y}\ ,
$$
where the first two maps are induced by the projections $X\times Y\to X$ and $X\times Y\to Y$, and the last comes from the aforementioned lax symmetric monoidal structure.
}\hB
%
\end{rem}

{The category $\Mf$ has a Grothendieck pre-topology given by covering families of manifolds.
Here a covering family of a manifold $M$ is an at most countable family of open subsets $(U_{\alpha})_{\alpha\in A}$ such that $\bigcup_{\alpha\in A} U_{\alpha}=M$.  We obtain an induced Grothendieck topology on the $\infty$-category $\Nerve(\Mf)$, see \cite[Def. 6.2.2.1 and Rem. 6.2.2.3]{MR2522659}. Using this Grothendieck topology we can define the full subcategory of sheaves $$\Sm^{desc}(\bC)\subseteq \Sm(\bC)\ .$$ We first recall the standard definition of sheaves which works   for arbitrary sites
and then explain that in our special case of manifolds it is equivalent to Definition \ref{jun215}. The latter will be used in praxis.

\bigskip

{
We assume that $\bC$ is presentable.  Then we can extend  smooth objects in $\bC$ to functors from smooth spaces to $\bC$,  i.e. we have a map \begin{equation}\label{kjkl3jkl32j2e3r32r3r332r23r32r}
\Sm(\bC)\simeq \Fun^{lim}(\Sm(\Spc)^{op},\bC)\to \Fun(\Sm(\Spc)^{op},\bC)\ ,\quad F\mapsto \tilde F
\end{equation}
where $\Fun^{lim}$ denotes the full subcategory of limit-preserving functors. The first equivalence is determined by the relation 
\begin{equation}\label{cewcjwelckjwlkcjlwekjclkwejcewcecewcecc}
F(M)\simeq \tilde F(Y(M)) 
\end{equation} for every manifold $M$, where \begin{equation}\label{fwefewfjewflkjiou3ro23r23r32r32r32r2r3}
 Y(M):=y(M)_{\infty}\in \Sm(\Spc) 
\end{equation} 
with $y(M)\in \Sm(\Nerve(\Set))$  given by  usual one-categorical Yoneda embedding, i.e. 
$$y(M)(N):=\Hom_{\Mf}(N,M) \ ,$$ {and for a set $X$ the symbol $X_{\infty}$ denotes the space represented by the set $X$ considered as a constant simplicial set, see Example \ref{kdjqwdkljqlwdjlqwjdlwqdwqdwqdqwdwqdwqddqwdqd}.} Note  that $\Sm(\Spc)$ is generated under colimits by representables.   
The second map in \eqref{kjkl3jkl32j2e3r32r3r332r23r32r} is the inclusion of limit preserving functors into all functors.

\begin{rem}{\rm  In order to avoid mistakes below note that $y(\emptyset) \not\cong \underline{\emptyset}$  since $y(\emptyset)(\emptyset)\cong *$ while $\underline{\emptyset}(\emptyset)=\emptyset$. }\hB
\end{rem}

If $\cU:=(U_{\alpha})_{\alpha\in A}$ is a covering family of $M$, then we can define the \v{C}ech nerve
$\cU_{M}\in \Sm(\Spc)$ as follows. We first define the smooth simplicial set $\cU^{\bullet}_{M}\in \Sm(\Nerve(\sSet))$
such that for $n\in \nat$ \begin{equation}\label{gergergkjerkgneoi4r354354353453er}
\cU^{n}_{M}:=\underbrace{\left(\coprod_{\alpha\in A}y(U_{\alpha})\right)\times_{y(M)} \dots \times_{y(M)} \left(\coprod_{\alpha\in A}y(U_{\alpha})\right)}_{n+1\:\: factors}\ .
\end{equation}
 Using the localization $\iota:\sSet\to  \Spc$ we then set \begin{equation}\label{edghjegdgjqhwdwqiudziuwzdiuziuqwdqwdqwdqd3}
\cU_{M}:=\iota(\cU^{\bullet}_{M})\ .
\end{equation}
 This smooth space comes with a natural augmentation
$\cU_{M}\to Y(M)$.
}

%
\begin{ddd} \label{jun215111}
A smooth object $F\in \Sm(\bC)$ of $\bC$ satisfies descent  (or is called a sheaf) if the  map 
\begin{equation}\label{frefkjefkerjferlkfrefreferfferferf}
\tilde F(Y(M))\longrightarrow   \tilde F(\cU_{M})
\end{equation}
induced by the augmentation
is an equivalence in $\bC$ for all $M$ and coverings $\cU$ of $M$. \end{ddd}

The site $\Mf$ has following special properties:
\begin{enumerate}
\item If $(U_{\alpha})_{\alpha\in A}$ is a covering family of $M$, then
$U:=\coprod_{\alpha} U_{\alpha}$ exists in $\Mf$.
\item If $(U_{\alpha})_{\alpha\in A}$ is an at most countable family of manifolds, then it is a covering family of $M:=\coprod_{\alpha\in A} U_{\alpha}$.
\item Coproducts in $\Mf$ are disjoint.
\end{enumerate}
The following Lemma is related to \cite[Thm. 10.3]{MR2034012}.
\begin{lem}
These properties imply that   Definition \ref{jun215111} is equivalent to   Definition \ref{jun215}.\end{lem}
\proof
Assume that $F\in \Sm(\bC)$
  satisfies descent according to Definition \ref{jun215111}. We then verify the conditions of
  Definition \ref{jun215}.  We first show that $F$ maps disjoint unions to products. 

\bigskip

  Note that the empty manifold $\emptyset\in \Mf$ admits the empty covering family $\cU$. In this case the associated \v{C}ech nerve 
  $\cU_{\emptyset} $ is equivalent to the   empty smooth space, i.e. an initial object in smooth spaces. Since $\tilde F$ preserves limits we conclude that
  $\tilde F(\cU_{\emptyset})\simeq *$, i.e. $\tilde F(\cU_{\emptyset})$ is a final object in $\Spc$, a   contractible space. We finally get  \begin{equation}\label{g3g34jnn43kgj434389gjkkjregergerg345}
  \tilde F(Y(\emptyset))\stackrel{\eqref{frefkjefkerjferlkfrefreferfferferf}}{\simeq} \tilde F(\cU_{\emptyset})\simeq *\ .
\end{equation}

  Let $(U_{\alpha})_{\alpha\in A}$ be an at most countable family of manifolds. Then $
 \cU:=(U_{\alpha})_{\alpha\in A}$ is a covering family of $M:=\coprod_{\alpha\in A} U_{\alpha}$.
 Using that coproducts in $\Mf$ are disjoint we  calculate
$$U_{\alpha_{0}}\times_{M} \dots \times_{M} U_{\alpha_{n}}\cong \left\{\begin{array}{cc} U_{\alpha_{0}}&\alpha_{0}=\dots=\alpha_{n}\\\emptyset&else\end{array}\right. \ .$$
Using the fact that the Yoneda embedding $y$ preserves limits (so in particular fibre products over $M$), by specializing \eqref{gergergkjerkgneoi4r354354353453er}  we get  
$$\cU_{M}^{n}\cong \coprod_{\alpha\in A} y(U_{\alpha})\sqcup \coprod_{A^{n+1}\setminus \diag(A)} y(\emptyset)\ .$$
 We use \eqref{dweidweiduewoiduoiweudoiu239087e9237e93} in order to calculate the \v{C}ech nerve $\cU_{M}$ defined by \eqref{edghjegdgjqhwdwqiudziuwzdiuziuqwdqwdqwdqd3}. We get \begin{equation}\label{gergeklgrklekrjlrege345t3}
\cU_{M}\simeq \colim_{[n]\in \Nerve(\Delta)^{op}}\cU_{M,\infty}^{n} \ .
\end{equation}
 
 Since $\tilde F$ preserves limits (and colimits in $\Sm(\Spc)$ are limits in $\Sm(\Spc)^{op}$) and  using \eqref{fwefewfjewflkjiou3ro23r23r32r32r32r2r3}  we get   \begin{equation}\label{vervnener89798745928knjr3emfwffewfhewjkffewf}
 \tilde F(\cU_{M,\infty}^{n})\simeq \prod_{\alpha\in A} \tilde F(Y(U_{\alpha}))\times \prod_{A^{n}\setminus \diag(A)} \tilde F(Y(\emptyset))\stackrel{\eqref{g3g34jnn43kgj434389gjkkjregergerg345}}{\simeq}  \prod_{\alpha\in A} \tilde F(Y(U_{\alpha}))\ . 
\end{equation}
 Therefore $\tilde F(\cU_{M,\infty}^{\bullet})$ is equivalent to the constant simplicial object
with value  $\prod_{\alpha\in A} \tilde F(Y(U_{\alpha}))$. For a constant simplicial object $\underline{C}$
with value $C$ we have by \eqref{chcjkwhkecheckehkeecewcewcewc} 
 \begin{equation}\label{iou23ioru2o3r32ru3o2iru32r32r32r32r2r}
\lim_{\Nerve(\Delta)^{op}} \underline{C}\simeq C^{ \Nerve(\Delta)^{op} }\simeq C
\end{equation}
since $ \Nerve(\Delta)^{op} \simeq *$.
  This gives finally the desired equivalence 
  \begin{eqnarray*}
 F(M)&\stackrel{\eqref{fwefewfjewflkjiou3ro23r23r32r32r32r2r3}}{\simeq} &\tilde F(Y(M))\\&\stackrel{\eqref{frefkjefkerjferlkfrefreferfferferf}}{\simeq}& \tilde F(\cU_{M})\\&\stackrel{\eqref{gergeklgrklekrjlrege345t3}, \eqref{vervnener89798745928knjr3emfwffewfhewjkffewf}}{\simeq}& \lim_{[n]\in \Nerve(\Delta^{op})}   \prod_{\alpha\in A} \tilde F(Y(U_{\alpha}))\\&\stackrel{\eqref{iou23ioru2o3r32ru3o2iru32r32r32r32r2r}}{\simeq} &\prod_{\alpha\in A} \tilde F(Y(U_{\alpha}))\\&\stackrel{\eqref{fwefewfjewflkjiou3ro23r23r32r32r32r2r3}}{\simeq} & \prod_{\alpha\in A} F( U_{\alpha} )\ .\end{eqnarray*}
 
 
 For the second condition of  Definition \ref{jun215} we consider a covering $\cU=(U_{\alpha})_{\alpha\in A}$ of a manifold $M$. We define the simplicial manifold $U_{M}^{\bullet}$ for the map
 $U:=\coprod_{\alpha\in A} U_{\alpha}\to M$ by \eqref{f4uhfjkhhfk3hf3fmnj334f487fz3487fz34ff34f3f}. Note that \begin{equation}\label{htrkhrtklhrjh09905363454543543543534}
U_{M}^{n}\cong \coprod_{(\alpha_{0},\dots,\alpha_{n})\in A^{n+1}}  U_{\alpha_{0}}\times_{M}\dots\times_{M} U_{\alpha_{n}}\ .
\end{equation}
 By \eqref{dweidweiduewoiduoiweudoiu239087e9237e93} we have an equivalence
 $ \cU_{M} \simeq \colim_{\Nerve(\Delta)^{op}}  \cU_{M,\infty}^{\bullet}$. Hence 
   \begin{equation}\label{efewfjkhwekfjhefkewfho2u349283u432424324243242}
\tilde F(\cU_{M})\simeq \tilde F(\colim_{\Nerve(\Delta)^{op}} \cU_{M,\infty}^{\bullet}))\stackrel{!}{\simeq}
 \lim_{\Nerve(\Delta)}\tilde F(  \cU_{M,\infty}^{\bullet})\ .
\end{equation} Since the Yoneda embedding $y$ preserves limits (used at the equivalence marked by $!!$) and $F$ and $\tilde F$ map disjoint unions to products (used at the equivalence marked by $!!!$) we further get   the equivalences
 \begin{eqnarray} \tilde F(  \cU_{M,\infty}^{n}) &\stackrel{!!}{\simeq}& \tilde F( \coprod_{(\alpha_{0},\dots,\alpha_{n})\in A^{n+1}} {Y}( U_{\alpha_{0}}\times_{M}\dots\times_{M} U_{\alpha_{n}} ))\label{lidjqwlkdjqwdwqdqwdwqdwqdqwd} \\
 &\stackrel{!!!}{\simeq} &\prod_{(\alpha_{0},\dots,\alpha_{n})\in A^{n+1}} \tilde F({Y}( U_{\alpha_{0}}\times_{M}\dots\times_{M} U_{\alpha_{n}} ))\nonumber \\ &\stackrel{\eqref{cewcjwelckjwlkcjlwekjclkwejcewcecewcecc}}{\simeq} &\prod_{(\alpha_{0},\dots,\alpha_{n})\in A^{n+1}} F(  U_{\alpha_{0}}\times_{M}\dots\times_{M} U_{\alpha_{n}}  ) \nonumber\\&\stackrel{!!!}{\simeq}&  F(  \coprod_{(\alpha_{0},\dots,\alpha_{n})\in A^{n+1}}  U_{\alpha_{0}}\times_{M}\dots\times_{M} U_{\alpha_{n}}) \nonumber \\&
 \stackrel{\eqref{htrkhrtklhrjh09905363454543543543534}}{\simeq} &F(U_{M}^{n})\ .\nonumber\end{eqnarray}
  Combining these equivalences we get the equivalence
$$F(M) \stackrel{\eqref{cewcjwelckjwlkcjlwekjclkwejcewcecewcecc}}{\simeq}   \tilde F(Y(M)) \stackrel{\eqref{frefkjefkerjferlkfrefreferfferferf}}{\simeq}   \tilde F(\cU_{M})\stackrel{\eqref{efewfjkhwekfjhefkewfho2u349283u432424324243242}, \eqref{lidjqwlkdjqwdwqdqwdwqdwqdqwd}}{ \simeq}   \lim_{\Nerve(\Delta)} F(U_{M}^{\bullet}) $$
as required.

\bigskip

We now consider the opposite direction 
 and assume that $F$ satisfies the conditions of Definition \ref{jun215}.
Then using the calculations above we  get the condition \eqref{frefkjefkerjferlkfrefreferfferferf}:
\begin{eqnarray*}\tilde F(\cU_{M})\simeq \tilde F(\colim_{\Nerve(\Delta)^{op}} \cU_{M,\infty}^{\bullet})
&\stackrel{!}{\simeq} & \lim_{ \Nerve(\Delta)}  \tilde F(\cU_{M,\infty}^{\bullet})\\&\stackrel{\eqref{lidjqwlkdjqwdwqdqwdwqdwqdqwd}}{\simeq} &  \lim_{ \Nerve(\Delta)}  F(U_{M}^{\bullet})\\&\stackrel{Def. \ref{jun215}, 2.}{\simeq} &F(M)\\&\stackrel{\eqref{cewcjwelckjwlkcjlwekjclkwejcewcecewcecc}}{\simeq}& \tilde F(Y(M)) \ .\end{eqnarray*}
\hB
}

\begin{ddd}
{For a presentable $\infty$-category $\bC$ we   define   $$\Sm^{desc}(\bC)\subseteq \Sm(\bC)$$ as} the full subcategory of smooth objects which satisfy descent. 
\end{ddd}
This inclusion is the right adjoint of an adjunction
\begin{equation}\label{sep0508}L:\Sm(\bC)\leftrightarrows \Sm^{desc}(\bC):i\ ,\end{equation}
where $L$ is the sheafification functor \cite[Def. 6.2.2.6, Lem. 6.2.2.7, Not. 6.2.2.8]{MR2522659}.
An object $E\in \Sm(\bC)$ satisfies descent if and only if the natural morphism
$E\to L(E)$ given by the unit of the adjunction \eqref{sep0508} is an equivalence.

\begin{rem}\label{dqwkldjqlkdjqwldjlwqdjwqlkdjwqldj}{\rm 
For every $n\in \nat$ we consider the poset $\mathcal{P}_{n}$ of open neighbourhoods of $0\in \R^{n}$.
For $E\in \Sm(\bC)$ we define the stalk $$E(\R^{n},0):=\colim_{U\in \mathcal{P}_{n}} E(U)\ .$$
Sheaves on the site $\Mf$ {with values in $\bC\in \{\Spc,\Sp,\Nerve(\Ch)[W^{-1}]\}$ or $\bC=\Nerve(\cC)$ for some complete and cocomplete category $\cC$} 
have the following very convenient properties:
\begin{enumerate}
\item Local equivalences can be checked on stalks. More precisely, consider a morphism $E\to F$  in $\Sm(\bC)$. Then $L(E)\to L(F)$ is an equivalence if and only if the induced morphism
$E(\R^{n},0)\to F(\R^{n},0)$ between the stalks is an equivalence for all $n\in \nat$.
\item Sheafification preserves stalks, i.e. for $E\in \Sm(\bC)$ the natural map $E\to L(E)$ induces equivalences
\begin{equation}\label{fwefjwefliwefjlewfjiofueoeu9u23424326t45z5444647586786}
E(\R^{n},0)\to L(E)(\R^{n},0)
\end{equation}
for all $n\in \nat$.
\end{enumerate}
{Note that the second property follows from the first.
Indeed, the unit $E\to L(E)$ is a local equivalence since $L\simeq L\circ L$. 
Evaluating on $(\R^{n},0)$ and using the ''only if''  part of 1. we see that  \eqref{fwefjwefliwefjlewfjiofueoeu9u23424326t45z5444647586786}
is an equivalence.
}

{Let us explain the reason for 1.
in more detail. One could start with presheaves  $\Sm(\bC)$ and declare a morphism to be a weak equivalence if  it induces an equivalence on stalks. Using the ideas in \cite{MR2034012} one can show that the corresponding localization of  $\Sm(\bC)$ is equivalent to the full subcategory $\Sm^{hyper-desc}(\bC)$ of presheaves which satisfy descent for all hypercoverings. 
Now for manifolds it is known that hypercover descent and \v{C}ech descent are equivalent conditions.  Indeed, since the covering dimension of an $n$-dimensional manifold $M$ is {bounded by} $\le n$
one knows that the $\infty$-topos  of $\Spc$-valued sheaves on the small site of open subsets of $M$
has homotopical dimensional dimension $\le n$. In particular, every $\Spc$-valued sheaf on this small site is equivalent to the limit of its Postnikov tower.  By a comparison of the big site $\Mf$ and the small sites of all connected manifolds in $\Mf$ one shows that every object in $\Sm^{desc}(\Spc)$ 
is equivalent to the limit of its Postnikov tower. This implies hypercompleteness of $\Sm^{desc}(\Spc)$.\footnote{See the discussion in MathOverflow \url{http://mathoverflow.net/questions/130999/is-the-site-of-smooth-manifolds-hypercomplete}}

Note that these results are  shown for space-valued presheaves. {In order to extend them to other targets one must generalize the  notion of Postnikov towers. If $\bC$ is the $\infty$-category of spectra $\Sp$ or of chain complexes $\Nerve(\Ch)[W^{-1}]$, then this is possible. Details will be given in the appendix of \cite{bst}.}}

\bigskip

Evaluation at $0$ provides a morphism \begin{equation}\label{lijlkerjlrekjfrlefjljreferferferfreferfrefrfref}
E(\R^{n},0)\to E(*)\ .
\end{equation}
If $n\ge 1$, then it is not an equivalence in general.  But see Example \ref{1234}.

\bigskip 

Consider the sheaf of  rings $C^{\infty}$ of smooth real-valued functions,  Example \ref{fklwjklewfwfwf45453453453wfwfwf}. Then $C^{\infty}(\R^{n},0)$ is the ring of  germs of such functions at $0\in \R^{{n}}$, while $C^{\infty}(*)\simeq \R$. }\hB
\end{rem}

\begin{ex}\label{fwlkfejewlkfewfewfewf}{\rm  For a   set $X\not\cong \{*\}$ the constant smooth object
$\underline{X}\in \Sm(\Nerve(\Set))$ (see Example \ref{dqjlwdjqwdwqdwqdwqdqwdwq3432423423}) does not satisfy descent. Indeed, a smooth object with descent must send the empty manifold to the final object of the target, i.e. the  one point set in this case. But $\underline{X}(\emptyset)=X\not\not\cong\{*\}$.
 
 \bigskip
 
On the other hand, $\Funk(X)(M)$ is the set of functions from $M_{top}$ to the discrete space $X$. One easily sees that $\Funk(X)$ satisfies descent. 
 This fact generalizes, see Lemma \ref{kdqqwdwqjldwqdqdqwdq343224342342242}. Note that ${\Funk}(X)(\emptyset)\cong X^{\emptyset}\cong \{*\}$. 
}\hB

\end{ex}

\begin{ex}\label{kldjqkldjqwlkdjqwdqwdqwdqwdwqdqwdqwdq}{\rm As an application of Remark \ref{dqwkldjqlkdjqwldjlwqdjwqlkdjwqldj} we note the following fact. 
Let $$ \Nerve(\Ch)\to \Nerve(\Ch)[W^{-1}] \ , \quad X\mapsto X_{\infty}$$ be {the} localization. Then for $X\in \Sm( \Nerve(\Ch))$ the natural morphism \begin{equation}\label{fwefhewfkehfkewhfkwfh8923u984234234324}
L( X_{\infty})\to L(  L(X)_{\infty})
\end{equation}
is an equivalence.
By Remark \ref{dqwkldjqlkdjqwldjlwqdjwqlkdjwqldj} we can check this on stalks. We must show that
$$ X_{\infty}(\R^{n},0)\to  L(X)_{\infty}(\R^{n},0)$$
is an equivalence for every $n\in \nat$. Now the stalk is defined via a filtered colimit. Since filtered colimits of chain complexes preserve
quasi-isomorphisms they are preserved by the localization.
Consequently it suffices to show that
$X (\R^{n},0)\to  L(X) (\R^{n},0)$ is an isomorphism, but this is clear.
}\hB

\end{ex}

{Let $\bC$ be a presentable $\infty$-category.}
\begin{lem}\label{kdqqwdwqjldwqdqdqwdq343224342342242}  {For every  
 $C\in \bC$ the corresponding smooth function object}  $\Funk(C)\in \Sm(\bC)$ satisfies descent.
\end{lem}
\proof{
We verify the two conditions of Definition \ref{jun215}. For the first observe that
the map
$M\to M_{top,\infty}$ preserves coproducts, and that the cotensor $\Spc\ni X\mapsto C^{X}\in \bC$ sends colimits to limits.
Hence
$\Funk(C)$ maps countable coproducts to products.
}

\bigskip

{
We now consider the second condition in  Definition \ref{jun215}. Let $(U_{\alpha})_{\alpha\in A}$ be a covering family of a manifold $M$ and form  the map
 $U:=\coprod_{\alpha\in A}U_{\alpha}\to M$.} Then we have the simplicial manifold $U_{M}^{\bullet}\in \Fun(\Delta^{op},\Mf)$ {(see \eqref{f4uhfjkhhfk3hf3fmnj334f487fz3487fz34ff34f3f} for a formula for $U_{M}^{n}$)} and therefore the simplicial space $U_{M,top,\infty}^{\bullet}\in \Fun(\Nerve(\Delta^{op}),\Spc)$. We consider  the augmentation $U_{M}^{\bullet}\to M$ as a map between simplicial manifolds,  where we consider $M$ as a constant simplicial manifold. Then we will use the crucial fact   that
$$\colim_{\Nerve(\Delta^{op})} U_{M,top,\infty}^{\bullet} \to M_{top,\infty}$$ is an equivalence in $\Spc$. Classically this  is known as the fact that the underlying topological space of a manifold is weakly equivalent to the homotopy colimit of the \v{C}ech nerve of the covering family, see e.g. \cite[Prop 4.1]{MR0232393}.
Since the cotensor maps colimits in $\Spc$ into limits in $\bC$ we get the equivalence \begin{equation}\label{kdhdkhqkwhdqkjdhjhd821e912e}
C^{M_{top,\infty}}\to \lim_{\Nerve(\Delta)} C^{U_{M,top,\infty}^{\bullet}}\ .\end{equation}
This is precisely the {desired second condition.} \hB

\begin{kor}\label{jfljlejlwf232343242342342342}
The smooth function object refines to a functor
$$\Funk:\bC\to \Sm^{desc}(\bC)\ .$$
\end{kor}

\bigskip

\begin{ex}\label{fwlkfejewlkfewfewfewfwww}{\rm

Some mathematical objects can be considered as objects of different categories. The fact that a smooth object satisfies descent may depend on the category in which it is considered. 
Here is an example.
 For a chain complex  $A\in \Nerve(\Ch)$  we  use the notation  $ A_{\infty}\in\Nerve(\Ch)[W^{-1}]$ for its image under localization. 
We apply the same notation for the object wise localization of smooth chain complexes.  
We consider $\Z$ as a chain complex concentrated in degree zero and form the smooth function object  
$\Funk(\Z)\in \Sm^{desc}( \Nerve(\Ch))$. Then 
$ \Funk(\Z)_{\infty}$ does not satisfies descent. Indeed,
$$H^{1}(L( \Funk(\Z)_{\infty})(S^{1}))\cong H^{1}(S^{1};\Z)\cong \Z\ , \quad H^{1}( \Funk(\Z)_{\infty} (S^{1}))\cong 0\ .$$
This example also shows that
$\Funk(\Z)_{\infty}\not\simeq \Funk(\Z_{\infty})$, since the latter  satisfies descent.
}\hB\end{ex}

\subsection{Homotopy invariance}\label{kjdkljkdljlqwdwqdqwdwqddqwdqdwqdwqd}
In this subsection we discuss the homotopy invariance of sheaves and presheaves and its consequences.
Let $\bC$ be some $\infty$-category. We have the unit interval $I:=[0,1]\in \Mf$.

\begin{ddd}\label{efjwflewjfl239874297329cjsajkas}
A smooth object $E\in \Sm(\bC)$ is called homotopy invariant if the projection $M\times I\to M$ induces an equivalence {$E(M )\to E(M\times I)$}  in $\bC$ for all manifolds $M$.
\end{ddd}
We let $$\Sm^{h}(\bC){\subseteq} \Sm(\bC)\ , \quad \Sm^{desc,h}(\bC)\subseteq \Sm^{desc}(\bC)$$ denote the full subcategories of homotopy invariant presheaves and sheaves.

\begin{ex}{\rm 
A constant smooth object $\underline{C}$ (see Example \ref{dqjlwdjqwdwqdwqdwqdqwdwq3432423423}) is homotopy invariant.

\bigskip

We consider the $\infty$-category $\Nerve(\Comm\Rings)$ of commutative rings (Example \ref{flkwejweklfjelkewfewfw4324234234}). 
The structure sheaf $C^{\infty}\in \Sm^{desc}(\Nerve(\Comm\Rings))$ of $\Mf$ associates to every manifold $M$ the commutative ring of smooth real-valued functions $C^{\infty}(M)$.
This sheaf is obviously  not homotopy invariant.

\bigskip

On the other hand, the sheaf $\cV\in \Sm^{desc}(\Nerve(\Mod(R)))$ of $R$-modules  considered in Example \ref{djkjwlkdjlwqdqwdqwdwqdwdwqd} is homotopy invariant.
}\hB
\end{ex}

\begin{ex}\label{1234}
{\rm If $E\in \Sm(\bC)$ is homotopy invariant, then the natural morphism
 \eqref{lijlkerjlrekjfrlefjljreferferferfreferfrefrfref}  
is an equivalence. Indeed, we can restrict the index category for the colimit defining the stalk to the poset of balls $B\subset \R^{n}$ in $\R^{n}$ centered at $0$. Since a ball is contractible we conclude that the map $*\stackrel{0}{\to} B$ induces an equivalence $E(B)\stackrel{\simeq}{\to} E(*)$.
}\hB
\end{ex}

Let $\bC$ be a presentable $\infty$-category. By Corollary \ref{jfljlejlwf232343242342342342}, for $C\in \bC$ we have the smooth function object $\Funk(C)\in \Sm^{desc}(\bC)$.

\begin{lem}\label{fjelwjfw82uro23ur3r22r23r23r2rr}
The object $\Funk(C)$ is homotopy invariant.
\end{lem}
\proof
The projection $(M\times I)_{top}\to M_{top}$ is a homotopy equivalence.
Consequently, $(M\times I)_{top,\infty}\to M_{top,\infty}$ is an equivalence in $\Spc$.
It follows that $C^{M_{top,\infty}}\to C^{(M\times I)_{top,\infty}}$ is an equivalence in $\bC$ for every manifold $M$. In view of \eqref{lfkjlekfjewlfjewlfjfelwjfljlekwf9038409832423424} this is exactly the condition for the  homotopy invariance  of $\Funk(C)$ according to Definition \ref{efjwflewjfl239874297329cjsajkas}. \hB 
  \begin{kor}
 The smooth function object construction refines to a functor
 $$\bC\to \Sm^{desc,h}(\bC)\ .$$
\end{kor}

\begin{lem}\label{oiwioewoifwoifjwfoj245239852345}
If $\bC$ is presentable, then we have adjunctions
$$\cH^{pre}:\Sm(\bC)\leftrightarrows \Sm^{h}(\bC):inclusion\ , \quad \cH:\Sm^{desc}(\bC)\leftrightarrows \Sm^{desc,h}(\bC):inclusion\ .$$
The left-adjoints are called homotopification functors.
\end{lem}
\proof
We reformulate the homotopy invariance condition as a locality condition.
We use that the (co)tensor structure of $\bC$ over $\Spc$ extends to a (co)tensor structure over $\Sm(\Spc)$ 
\begin{equation}\label{kkjdqwkdjwqldjwdq9894324234}\Sm(\bC)\times \Sm(\Spc)\to \Sm(\bC)\ , \quad  \Sm(\Spc)^{op}\times \Sm(\bC)\to \Sm(\bC)\ .\end{equation} 
For a manifold $M$ we let  $Y(M)\in \Sm(\Spc)$  be as in \eqref{fwefewfjewflkjiou3ro23r23r32r32r32r2r3}.
Our first reformulation is as follows.
A presheaf $E$ is homotopy invariant if and only if the morphism
$E\to E^{Y(I)}$ induced  by the projection of $Y(I)\to *$ is an equivalence.
Indeed, if we evaluate this map on a manifold $M$ and use the equivalence
$E^{Y(N)}(M)\simeq E(M\times N)$, then we get exactly the defining condition for homotopy invariance.

For every object $F\in \Sm(\bC)$ we consider the morphism 
$$p_{F}:F\otimes Y(I)\to F$$
induced by the projection of $Y(I)\to *$ and the tensor structure \eqref{kkjdqwkdjwqldjwdq9894324234}. Then in our   second reformulation 
an object $E\in \Sm(\bC)$ is homotopy invariant, if and only if
$$\Map(p_{F},E):\Map(F,E)\to \Map(F\otimes Y(I),E)$$ is an equivalence for every $F$.
To this end we simply rewrite $$\Map(F\otimes Y(I),E)\simeq \Map(F, E^{Y(I)})$$ and use the first reformulation.

Therefore $E\in \Sm(\bC)$ is homotopy invariant if and only if it is local with respect to the morphisms
$p_{F}$ for all $F\in \Sm(\bC)$. 
It suffices to require this locality  for $F$ running through a set of generators of $\Sm(\bC)$. This implies that the inclusion of homotopy invariant presheaves into all presheaves is the right-adjoint of an accessible localization.  This finishes the argument for presheaves.

\bigskip

For sheaves we argue similarly. \hB

An explicit model for the homotopification functor $\cH^{pre}$ has been discussed in \cite[Lemma 7.5]{Bunke:2013uq}. Let $\Delta^{\bullet}\in \Fun(\Delta,\Mf)$ be the cosimplicial manifold given by the collection of standard simplices. We then define the functor \begin{equation}\label{dqwhdkjqwhdwqdqdhwqkd89ue12e12ee1e21e1e12e12}
\bs:\Sm(\bC)\to \Sm(\bC)\ , \quad E\mapsto \colim_{\Nerve(\Delta)^{op}}E^{Y(\Delta^{\bullet})}\ .
\end{equation}
We have a natural augmentation $  \Delta^{\bullet}\to *$ which induces a morphism
$\id\to  \bs$. By Lemma \cite[Lemma 7.5 (2)]{Bunke:2013uq} we have an equivalence
\begin{equation}\label{ufeufioewufouaweofuoewfuoewifuewofui}\cH^{pre}\simeq \bs\ .\end{equation}

In general, for $E\in \Sm^{desc}(\bC)$ the object $\cH^{pre}(E)\in \Sm^{h}(\bC)$ does not satisfies descent (but see \cite[Proposition 7.6]{Bunke:2013uq} which says that it does for finite coverings if $\bC$ is stable). Before we provide the model for $\cH$ we show some general facts about homotopy invariant presheaves and sheaves on the category of manifolds. They are essentially due to Dugger, see e.g. \cite[Prop. 8.3]{MR1870515} and consequences of the fact that manifolds are locally contractible.

 We have an adjunction
\begin{equation}\label{kljwdklwqjdjqqlo32ue}\underline{(-)}:\bC\leftrightarrows \Sm(\bC):\ev_{*}\end{equation}
where $\ev_{*}$ denotes the evaluation at the one-point manifold $*$ and $\underline{(-)}$ is the constant presheaf functor
is as in Example \ref{dqjlwdjqwdwqdwqdwqdqwdwq3432423423}.

\begin{ass}\label{wijfwijeoifjeowifjw}
{We assume that $\bC$ is presentable and such that Remark \ref{dqwkldjqlkdjqwldjlwqdjwqlkdjwqldj} applies, i.e. that we can detect \v{C}ech local equivalences in $\Sm(\bC)$  on stalks.}
\end{ass}

\begin{lem}\label{lfehwefh984984234324}
{Assume that $\bC$ satisfies \ref{wijfwijeoifjeowifjw}.}
\begin{enumerate}
\item For every $C\in \bC$ the natural morphism
$L(\underline{C})\to \Funk(C)$ is an equivalence.
\item For every $E\in \Sm^{h}(\bC)$ the natural map
$L(\underline{\ev_{*}(E)})\to L(E)$ is an equivalence.
\item The sheafification  of a homotopy invariant presheaf is  still homotopy invariant.
\end{enumerate}
\end{lem}
\proof
We have an equivalence of functors from $\bC$ to $\bC$
 $$\id_{\bC}\stackrel{\simeq}{\to} \ev_{*}\circ \Funk\ .$$ From the adjunction \eqref{ufeufioewufouaweofuoewfuoewifuewofui} we then get
  a morphism between functors $$\underline{(-)}\to \Funk(-)$$
  from $\bC$ to $\Sm(\bC)$.
  Since its target takes values in sheaves it induces the  morphism
$$L(\underline{(-)})\to \Funk(-)$$ in the statement 1.
In order to show that it is an equivalence we insert an object $C\in \bC$ and show that the induced morphism on stalks 
$$L(\underline{C})(\R^{n},0)\to  \Funk(C)(\R^{n},0)$$
is an equivalence for every $n\in \nat$. To this end we consider the commutative diagram
$$\xymatrix{\underline{C}(\R^{n},0)\ar[d]_{\simeq }^{\eqref{fwefjwefliwefjlewfjiofueoeu9u23424326t45z5444647586786}}\ar[r]^{\simeq}&C\\L(\underline{C})(\R^{n},0)\ar[r]& \Funk(C)(\R^{n},0)\ar[u]^{\eqref{lijlkerjlrekjfrlefjljreferferferfreferfrefrfref} }_{\simeq}}\ .$$
The upper equivalence uses that the colimit defining $\underline{C}(\R^{n},0)$, by definition of the constant presheaf, is a colimit over the constant diagram with value $C$. The right vertical map is an equivalence by Example \ref{1234} and Lemma \ref{fjelwjfw82uro23ur3r22r23r23r2rr}.

\bigskip

The map in statement 2. of the Lemma is induced by the counit 
$$\underline{\ev_{*}(E)}\to E$$
of the adjunction \eqref{kljwdklwqjdjqqlo32ue}.
In view of Remark \ref{dqwkldjqlkdjqwldjlwqdjwqlkdjwqldj}, 1. we must check that for every $n\in \nat$ the induced map
$$\underline{\ev_{*}(E)}(\R^{n},0)\to E(\R^{n},0)$$
is an equivalence.  Since $\underline{\ev_{*}(E)} $ and $E$  are homotopy invariant,  by Example \ref{1234} this map is equivalent to the canonical map
$\underline{\ev_{*}(E)}(*)\to E(*)$, i.e. the identity on $E(*)$.


\bigskip

For 3. we consider a homotopy invariant presheaf $E$. Then by 1. and 2. we have equivalences
$$L(E)\simeq L(\underline{\ev_{*}(E)})\simeq \Funk(\ev_{*}(E))\ .$$
Since the smooth function object is homotopy invariant by Lemma \ref{fjelwjfw82uro23ur3r22r23r23r2rr} we conclude that $L(E)$ is homotopy invariant, too. 
\hB

We can now give a formula for the homotopification on the level of sheaves.

 \begin{lem}\label{eklwfwefefewfewfewf} {Assume that $\bC$ satisfies \ref{wijfwijeoifjeowifjw}.}
 We have an equivalence
 \begin{equation}\label{jlkjelwjfweklfwpfoiopi32pr23r32r2r}\cH\simeq  L\circ \cH^{pre}\ .\end{equation}
 \end{lem}
 \proof For $E\in \Sm^{desc}(\bC)$ we know  by Lemma \ref{lfehwefh984984234324}, 3. that $L(\cH^{pre}(E))$ is homotopy invariant.
 We have therefore a functor $L\circ \cH^{pre}:\Sm^{desc}(\bC)\to \Sm^{desc,h}(\bC)$. In order to see that this functor represents the homotopification  we observe that for every  $F\in \Sm^{desc,h}(\bC)$ we have an equivalence of mapping spaces
 $$\Map(L(\cH^{pre}(E)),F)\simeq \Map(\cH^{pre}(E),F)\simeq \Map(E,F)$$
 induced by the map $E\to L(\cH^{pre}(E))$. 
%
%
%
%
%
 \hB

  We consider a morphism $E\to F$  in $\Sm^{desc}(\bC)$ with a homotopy invariant target $F$.
Recall the Definition \ref{dlqwkjqlwkdjwqldjwqdou231ildejodidqw} of the property that this morphisms presents $F$ as the homotopification of $E$. The following Lemma provides a usable criterion for this property.

\begin{lem}\label{dkhkqlwdjqlwdwqdwdwqdqwdw} {Assume that $\bC$ satisfies \ref{wijfwijeoifjeowifjw}.} If the induced morphism
 $\bs(E)\to \bs(F)$ is an equivalence, then  $E\to F$  presents $F$ as the homotopification of $E$.
\end{lem}
\proof
Since $F$ is homotopy invariant the natural map $F\to \bs(F)$ is an equivalence. In particular $\bs(E)$ satisfies descent and is homotopy invariant. This acounts for the first two equivalences in the upper line of the following diagram:
$$\xymatrix{E\ar[r]\ar[d]&\bs(E)\ar[d]^{\simeq}\ar[r]^{\simeq}&\cH(\bs(E))\ar[d]^{\simeq}\ar[r]^{\simeq}& \cH(L(\bs(E)))\ar[r]_{\eqref{jlkjelwjfweklfwpfoiopi32pr23r32r2r}, \eqref{ufeufioewufouaweofuoewfuoewifuewofui}}^{\simeq}\ar[d]^{\simeq}&\cH(E)\ar[d]^{\simeq} \\
F\ar[r]^{\simeq}&\bs(F)\ar[r]^{\simeq}&\cH(\bs(F))\ar[r]^{\simeq}&\cH(L(\bs(F)))\ar[r]^{\simeq}&\cH(F)}$$
In view of Definition  \ref{dlqwkjqlwkdjwqldjwqdou231ildejodidqw} the outer square implies the assertion. 
\hB

In the following Lemma we consider the case $\bC=\Spc$. More generally, the argument applies if $\Sm(\bC)$ is an $\infty$-topos. Note that $\cH$ preserves colimits since it is a left-adjoint. It does not preserve general limits. But we have the following special result.
\begin{lem}\label{ehfwkejhfewkfwdwqwqdqdqd} {Assume that $\bC$ satisfies \ref{wijfwijeoifjeowifjw} and}
 consider a cartesian square  {in} $\Sm^{desc}(\bC)$
$$\xymatrix{P\ar[r]\ar[d]&Y\ar[d]\\X\ar[r]&Z}\ .$$
If $Z$ is homotopy invariant, then the induced square
$$\xymatrix{\cH(P)\ar[r]\ar[d]&\cH(Y)\ar[d]\\\cH(X)\ar[r]&Z}$$
is cartesian.
\end{lem}
\proof Since $Z$ is homotopy invariant we have equivalences $Z\to Z^{Y(\Delta^{n})}$ for all $n\in \nat$.
We calculate $\cH^{pre}$ using the  model \eqref{ufeufioewufouaweofuoewfuoewifuewofui}.
 Using that the diagonal $$\Delta^{op}\to \Delta^{op}\times \Delta^{op}$$ is cofinal, and that colimits in an $\infty$-topos are universal, we get
\begin{eqnarray*}
\cH^{pre}(P)&\simeq& \colim^{pre}_{\Nerve(\Delta^{op})} (X\times_{Z}Y)^{Y(\Delta^{\bullet})}\\
&\simeq&  \colim^{pre}_{\Nerve(\Delta^{op})} (X^{Y(\Delta^{\bullet})}\times_{Z^{Y(\Delta^{\bullet}})}Y^{Y(\Delta^{\bullet})})\\
&\simeq&  \colim^{pre}_{\Nerve(\Delta^{op})} (X^{Y(\Delta^{\bullet})}\times_{Z}Y^{Y(\Delta^{\bullet})})\\
&\simeq& \colim^{pre}_{\Nerve(\Delta^{op} )\times \Nerve(\Delta^{op})} (X^{Y(\Delta^{\bullet})}\times_{Z }Y^{Y(\Delta^{\bullet})}) \\ &\simeq&
 \colim^{pre}_{\Nerve(\Delta^{op})} X^{Y(\Delta^{\bullet})}
 \times_{  Z } \colim^{pre}_{\Nerve(\Delta^{op})} Y^{Y(\Delta^{\bullet})} \\&\simeq&
 \cH^{pre}(X)\times_{Z}\cH^{pre}(Y)\ .
\end{eqnarray*}
Here the superscript $pre$ indicates that we take the colimit in $\Sm(\bC)$.
{S}ince $L$ commutes with finite limits {we have the equivalence}
$$(L\circ \cH^{pre})(P)\simeq  (L\circ \cH^{pre})(X)\times_{Z} (L\circ \cH^{pre})(Y) \ .$$
We now use \eqref{jlkjelwjfweklfwpfoiopi32pr23r32r2r} in order to conclude that 
$$\cH(P)\simeq  \cH(X)\times_{Z} \cH(Y) \ .$$
 \hB

{\begin{rem}{\rm  The proof of   \cite[Thm. 4.2.3]{Asok:2015uq} gives a hint how one could weaken the assumptions on $Z$ in Lemma \ref{ehfwkejhfewkfwdwqwqdqdqd}.}\hB
\end{rem}
}

 We finally show Lemma \ref{djqwjkqwdqwdqdqd213}. We recall its statement

\begin{lem}\label{djqwjkqwdqwdqdqd2131}
 {Assume that $\bC$ satisfies \ref{wijfwijeoifjeowifjw}.} The functors  $$\Funk:\bC\to \Sm^{desc,h}(\bC)\ ,\quad \ev_{*}:\Sm^{desc,h}(\bC)\to \bC$$   are inverse to each other equivalences of $\infty$-categories. In particular, for $F\in \Sm^{desc,h}(\bC)$ 
 there is a natural equivalence \begin{equation}\label{dkqwkdqwdjdljdlkwd9879kwjefkjefewjfhewkf}
\Funk(F(*))\stackrel{\simeq}{\to} F\ . \end{equation}
  \end{lem}
\proof
If $E\in \Sm^{desc,h}(\bC)$, then combining the equivalences in Lemma \ref{lfehwefh984984234324}, 1. and 2., and the equivalence $E\stackrel{\simeq}{\to} L(E)$ we get the desired natural equivalence \begin{equation}\label{jhffhskdhfkshfkwhekhfkwfwefwfwefwefweff}
\Funk(\ev_{*}(E))\stackrel{\simeq}{\leftarrow} L(\underline{\ev_{*}(E)})\stackrel{\simeq}{\to} L(E)\stackrel{\simeq}{\leftarrow} E\ .
\end{equation}
  
 We furthermore have the obvious natural equivalence
 $$\ev_{*}(\Funk(C))\simeq C$$
 for every $C\in \bC$.
 In a first step we can conclude that the functors $\ev_{*}$ and $\Funk$ induce inverse to each other equivalences between the homotopy categories
 $\Ho(\Sm^{desc,h}(\bC))$ and $\Ho(\bC)$. It is now a general fact about $\infty$-categories that this implies that these functors also induce
inverse to each other equivalences between $\infty$-categories. \hB

\subsection{The de Rham complex}\label{jhdkjqwdqwdwqdwqdqd}

In this subsection we define the de Rham complex $\Omega A$ and study its descent and homotopy invariance properties.

\bigskip

We consider the category $\Sm^{desc}(\Nerve(\Ch_{\R}))$  of sheaves of  chain complexes of real vector spaces. 
\begin{ex}{\rm A chain complex $A\in \Ch_{\R}$ of real vector spaces gives rise to a locally constant sheaf of chain complexes $L(\underline{A})\in \Sm^{desc}(\Nerve(\Ch_{\R}))$. Another example of a sheaf of chain complexes  is the de Rham complex $\Omega\in \Sm^{desc}(\Nerve(\Ch_{\R}))$ which associates to a smooth manifold $M$ the chain complex $\Omega(M)$ of smooth real differential forms.  }\hB \end{ex}

\bigskip

The tensor product of chain complexes of real vector spaces induces a symmetric monoidal structure on $\Nerve(\Ch_{\R})$ and therefore symmetric monoidal structures $\otimes_{\R}^{pre}$ on $ \Sm (\Nerve(\Ch_{\R}))$ and  $\otimes_{\R}$ on $\Sm^{desc}(\Nerve(\Ch_{\R}))$. The latter two are related by
\begin{equation}\label{vsdvhsdkvjhsdkvhksdvhsdkvhdskviowrwerewrwerwer}
X\otimes_{\R} Y\cong L(X\otimes_{\R}^{pre} Y)
\end{equation} for $X,Y\in \Sm^{desc}(\Nerve(\Ch_{\R}))$.  
 The tensor product $\otimes_{\R}$ preserves quasi-isomorphisms and therefore descents directly to
 a symmetric monoidal structure on $\Nerve(\Ch_{\R})[W^{-1}]$, {we have the equivalence
 \begin{equation} \label{jfewkljflwjlewooioi}(X\otimes_{\R}Y)_{\infty}\simeq X_{\infty}\otimes_{\R}Y_{\infty}\ .\end{equation}} We  get an induced symmetric monoidal structure 
 $\otimes_{\R}^{pre}$ on $\Sm( \Nerve(\Ch_{\R})[W^{-1}])$.
 
 \bigskip
\begin{ddd}\label{deRham}
{Let $A\in \Ch_{\R}$ be a chain complex of real vector spaces.}
The de Rham complex $\Omega A$ with coefficients in $A$ is defined as  the tensor product in $\Sm^{desc}(\Nerve(\Ch_{\R}))$
\begin{equation}\label{regreerreg34543536777754z45}
\Omega A:=\Omega\otimes_{\R} L(\underline{A})\ .
\end{equation}
\end{ddd} 

In most applications we will forget the $\R$-vector space structure and use the same symbol $\Omega A$ for the corresponding object of $ \Sm^{desc}(\Nerve(\Ch))$. 

\begin{rem}\label{jkdhwqjkdhkwqdwqhdqwdkjhwqdqwdwqdwqdwd}{\rm In this remark we motivate the definition of $\Omega A$ using the tensor product of sheaves. In fact, if
 $A$ is degree wise finite dimensional or if the manifold $M$ is compact, then we simply have  an isomorphism $$\Omega A(M)\cong \Omega(M)\otimes_{\R}A\ .$$  But in general, the map $M\mapsto \Omega(M)\otimes_{\R}A$ does not define a sheaf and we have a proper inclusion $\Omega(M)\otimes A\to \Omega A(M)$. The sheafification in the definition of the tensor product between sheaves \eqref{vsdvhsdkvjhsdkvhksdvhsdkvhdskviowrwerewrwerwer} acts as a sort of completion. But note that we do not use any kind of topology on $A$. 
 
 \bigskip
 
 If $A$ is degree wise finite-dimensional, then we can identify
 $\Omega A^{i}(M)$ with the smooth sections of the vector bundle $$\bigoplus_{p+q=i} \Lambda^{p}T^{*}M \otimes_{\R} A^{q}\ .$$

\bigskip

Here is a typical example. Assume that $A=\R^{\Z}$  and let $\chi\in C_{c}^{\infty}(\R)$. For $i\in \Z$ let $T_{i}:\R\to \R$ be the translation $x\mapsto T_{i}(x):=x+i$ and $\delta_{i}\in \R^{\Z}$ be the characteristic function of the subset $\{i\}\in \Z$.
Then we can interpret the infinite sum
$\sum_{i\in \Z} T_{i}^{*}\chi \otimes \delta_{i}$
naturally as an {element in} 
$\Omega A^{0}(\R)$. But note that it does not belong to the subspace $(\Omega(\R)\otimes_{\R} A)^{0}$.

}\hB\end{rem}

\bigskip

 We now consider the image  \begin{equation}\label{mar0801}\Omega A_{\infty}\in \Sm(\Nerve(\Ch )[W^{-1}])\end{equation}  of $\Omega A$ under the localization $\Sm(\Nerve(\Ch))\to  \Sm(\Nerve(\Ch)[W^{-1}])$.  
    We  further consider the sheaf of subcomplexes $$\sigma \Omega A\subseteq \Omega A\ , \quad \sigma\Omega A\in  \Sm^{desc}(\Nerve(\Ch))$$ of cochains of non-negative total degree, see \eqref{gerggrklgrjeglregjelrgjerlre}.  Its localization is denoted by \begin{equation}\label{mar0801111}\sigma \Omega A_{\infty}\in  \Sm(\Nerve(\Ch)[W^{-1}])\ .\end{equation}

    \begin{rem}{\rm 
In this construction we must be careful with the order. We first form the subcomplex   of an object  $C\in \Nerve(\Ch)$ and then apply the localization $\Nerve(\Ch)\to \Nerve(\Ch)_{\infty}$. A more appropriate notation would
be $(\sigma (C))_{\infty}$ but we drop the brackets in order  to simplify the notation.
On $ \Nerve(\Ch)[W^{-1}] $ the operation $\sigma$ is not well-defined, see Example \ref{kewjfwklefjlwekfewfwefefefwefwef}.
}\hB \end{rem}

  \begin{ex}\label{kewjfwklefjlwekfewfwefefefwefwef}{\rm 
 Assume that $A=\R[k]$, the cochain complex $\R$ with $\R$ in degree $-k$.
 Then $\Omega A\cong \Omega[k]$ is the shifted de Rham  complex where $\ell$-forms are in degree $\ell-k$.
 We then have 
 $$\sigma \Omega A\cong (\Omega^{\ge k})[k]\ .$$
 Note that
 $$H^{0}(\sigma \Omega A(M))\cong Z^{k}(\Omega(M))\ ,$$
where $Z^{k}(\Omega(M))$ is the space of closed $k$-forms. In degree zero cohomology the inclusion
$\sigma \Omega A\to \Omega A$ induces the map
$$[-]:Z^{k}(\Omega(M))\to H^{k}(\Omega(M))$$
which sends a closed $k$-form to the cohomology class it represents.

\bigskip

Let us fix a connected  manifold $M$ of positive dimension without boundary and consider the complex
 $\Omega_{-\infty}[k](M)$  of distributional forms shifted by $k$ for some $k\in\nat$ such that   $\dim(M)\ge k\ge 1$.
 Then the inclusion
$$\Omega[k](M)\to \Omega_{-\infty}[k](M)$$ is a quasi-isomorphism,
but
$$\sigma \Omega[k](M)\to \sigma\Omega_{-\infty}[k](M)$$
is not. In fact, the inclusion
$$Z^{k}(\Omega(M))\to Z^{k}(\Omega_{-\infty}(M))$$
is proper.}\hB
 
 \end{ex}

 \begin{rem}{\rm 
 Note that  in \eqref{mar0801} and \eqref{mar0801111} we drop the superscript $desc$ on purpose. In general,   the image of an object of $\Sm^{desc}(\Nerve(\Ch))$ in $\Sm(\Nerve(\Ch)[W^{-1}])$ does not satisfy descent. See Example \ref{fwlkfejewlkfewfewfewfwww}.
 But by Lemma \ref{feb1601} below  the localizations of the de Rham complex and its truncation  still satisfy descent. 
 } \hB
 \end{rem}

%
 
\begin{lem}\label{feb1601}
If $A$ is a chain complex of real vector spaces, then $\Omega A_{\infty}$ and $\sigma \Omega A_{\infty}$ satisfy descent.
\end{lem}
\proof
This is  \cite[Lemma   A.4]{buta}. For completeness we  give a detailed version of  the proof.\footnote{This corrects the proof given in an earlier preprint version of the present paper which only works if $A$ is bounded below.} 
We use the characterization  of sheaves given in Definition \ref{jun215} and  write out the details  in the case of $\Omega A_{\infty}$. The argument for $\sigma \Omega A_{\infty}$ is the same.

First of all, since $\Omega A$ sends disjoint unions to products, so does $\Omega A_{\infty}$. We now verify the second condition. Let $ \cU:=(U_{\alpha})_{\alpha\in A}$ be a covering family of a   manifold $M$ and $U^{\bullet}_{M}$ be the associated simplicial manifold \eqref{f4uhfjkhhfk3hf3fmnj334f487fz3487fz34ff34f3f}. Without loss of generality we can assume that $M$ is connected so that $\dim(M)$ is defined. We must show that the natural map
\begin{equation}\label{eq11} \Omega A_{\infty}(M)\to\lim_{\Nerve(\Delta)} \Omega A_{\infty}(U^{\bullet}_{M})\end{equation}
 is an equivalence.  Now  $\Omega A(U^{\bullet}_{M})$ is a cosimplicial chain complex. It gives rise to a double complex $\check{C} (\cU,\Omega A):=DK(\Omega A(U^{\bullet}_{M})) $ (see Example \ref{fjfwlefjlwefjlewfewfewffewfewfef} for notation), the \v{C}ech complex associated to the covering family $\cU$. By \eqref{dqwhqwdiwjhdwqjdwqdijwqdoijwqdoiwqdu987987983243424} the limit in (\ref{eq11}) is realized as the associated total complex $\Tot (
\check{C} (\cU,\Omega A))$ (defined using products), i.e. we have an equivalence  
\begin{equation}\label{ewfewfewf234423455346634f1}
\lim_{\Nerve(\Delta)} \Omega A_{\infty}(U^{\bullet}_{M})\simeq  \Tot (
\check{C} (\cU,\Omega A))_{\infty}\ .\end{equation}

If $E$ is a  sheaf of $C^{\infty}$-modules, then using a partition of unity we can define a chain homotopy inverse of  the natural map  $E(M) \to \check{C}(\cU,E)$. Note that  this construction involves  restrictions and multiplication by smooth functions. Hence
more generally, if $E$ is a complex of $C^{\infty}$-modules (i.e. its homogeneous components are sheaves of $C^{\infty}$-modules and the differentials are morphisms of $C^{\infty}$-modules), then we can use  the same construction   in order to define a  homotopy inverse of $E(M)\to \Tot( \check{C}(\cU,E))$. 

The differential of $\Omega A$ is the sum of the de Rham differential and the differential of $A$. Because of the presence of the de Rham differential it is not a complex of $C^{\infty}$-modules. In order to proceed we  introduce
  decreasing filtration of $\Omega A$ by 
 $$F^{p}\Omega A^{n}:=\bigoplus_{q\ge p} \Omega^{q}\otimes A^{n-q}\ .$$ Then 
 $\mathrm{Gr}^{p} (\Omega A)$  only involves the differential of $A$ and is therefore a complex of $C^{\infty}$-modules.
We conclude that  $$\mathrm{Gr}^{p}(\Omega A(M))\to \Tot(\check{C} (\cU,\mathrm{Gr}^{p}( \Omega A)))\cong \mathrm{Gr}^{p} (\Tot(\check{C} (\cU,  \Omega A)))$$ admits a chain homotopy inverse for all $p\in \nat$.
Since the filtration is finite (its length is bounded by the dimension of $M$)  we conclude that
$\Omega A(M)\to  \Tot(\check{C}(\cU,  \Omega A))$ is a chain homotopy equivalence.
This implies \begin{equation}\label{ewfewfewf234423455346634f}
 \Omega A(M)_{\infty}\simeq  \Tot\,(\check{C} (  \cU,\Omega A))_{\infty}\ . 
\end{equation}
The combination of \eqref{ewfewfewf234423455346634f} and \eqref{ewfewfewf234423455346634f1} shows that \eqref{eq11} is an equivalence as desired. \hB




 {
\begin{lem}\label{djqwkdjqwlwsqdqw}
The sheaf $\Omega A_{\infty}$ is homotopy invariant.
\end{lem}
\proof
We first argue that $\Omega_{\infty}$ is homotopy invariant.
 We have a natural map
 $\underline{\R}\to \Omega$. It induces
 $\underline{\R}_{\infty}\to \Omega_{\infty}$ by localization.
 Since the target is a sheaf by Lemma \ref{feb1601}, we   get the morphism
 $L(\underline{\R}_{\infty})\to \Omega_{\infty}$.
This morphism  is an equivalence. Indeed, using  Remark \ref{dqwkldjqlkdjqwldjlwqdjwqlkdjwqldj} 
this can be checked on stalks where the argument boils down to the usual Poincar\'e Lemma.
By Lemma \ref{lfehwefh984984234324}, 1. and Lemma \ref{fjelwjfw82uro23ur3r22r23r23r2rr} the sheaf 
$L(\underline{\R}_{\infty})$ is homotopy invariant. Hence $ \Omega_{\infty}$ is homotopy invariant.
Note that, in contrast, $\Omega$ is not homotopy invariant.

\bigskip

Again by Lemma \ref{lfehwefh984984234324}, 1. and Lemma \ref{fjelwjfw82uro23ur3r22r23r23r2rr}, for a chain complex $A$ of real vector spaces, the sheaf $L(\underline{A})$ is homotopy invariant.
It follows that the presheaf
$ \Omega_{\infty}\otimes_{\R}^{pre}L(\underline{A})_{\infty}  $
is homotopy invariant. By Lemma \eqref{lfehwefh984984234324} 3. the sheaf
$L( \Omega_{\infty}\otimes_{\R}^{pre}L(\underline{A})_{\infty})$ is homotopy invariant.
We now have the chain of equivalences
\begin{eqnarray*}
\Omega A_{\infty}&\stackrel{\eqref{regreerreg34543536777754z45}}{\simeq}&(\Omega \otimes_{\R}L(\underline{A}))_{\infty}\\&\stackrel{Lemma \, \ref{feb1601}}{\simeq} &L((\Omega \otimes_{\R}L(\underline{A}))_{\infty})\\
\\&\stackrel{\eqref{vsdvhsdkvjhsdkvhksdvhsdkvhdskviowrwerewrwerwer}}{\simeq}&L(L(\Omega \otimes^{pre}_{\R}L(\underline{A}))_{\infty})\\
&\stackrel{ \eqref{fwefhewfkehfkewhfkwfh8923u984234234324}}{\simeq}&
L(( \Omega \otimes^{pre}_{\R}L(\underline{A}))_{\infty})\\&{\stackrel{\eqref{jfewkljflwjlewooioi}}{\simeq}}&
L(  \Omega_{\infty} \otimes^{pre}_{\R}L(\underline{A} )_{\infty})\ .
\end{eqnarray*}
We conclude that $\Omega A_{\infty}$ is homotopy invariant. \hB

\begin{rem}{\rm Note that $\sigma \Omega A_{\infty}$ is not homotopy invariant, in general.}\hB
\end{rem}

\begin{rem}{\rm 
In this remark we explain why we use such an abstract argument to prove Lemma \ref{djqwkdjqwlwsqdqw}.
In fact, if $A$ is degree wise finite-dimensional, then we could employ the usual homotopy argument in order to show that
$$\Omega A(M)\to\Omega A(M\times I)$$ 
is a quasi-isomorphism for every manifold $M$. But without that finiteness condition on $A$ it requires some work   to extend   this argument  directly. Note that the homotopy formula involves integrals which would be well-defined on elements of $\Omega(M\times I)\otimes_{\R} A$. We would have to extend this formula
   to sections of the in general larger space $\Omega A(M\times I )$, see also Remark \ref{jkdhwqjkdhkwqdwqhdqwdkjhwqdqwdwqdwqdwd}. This extension could be constructed using similar methods as in Subsubsection \ref{kwjfweklfjewlfjelwfewopipi3ir23r20pr32r32r32r}.
 }\hB 

\end{rem}

}

Recall Definition \ref{dlqwkjqlwkdjwqldjwqdou231ildejodidqw}.
{
\begin{lem}\label{djqlwdjqwldjwqldwqdwqdwqdwqdq}
The map $\sigma \Omega A_{\infty}\to \Omega A_{\infty}$ presents $ \Omega A_{\infty}$ as the homotopification of $\sigma \Omega A_{\infty}$.
\end{lem}
\proof 
By Lemma  \ref{djqwkdjqwlwsqdqw} the sheaf  $\Omega A_{\infty}$ is homotopy invariant. It therefore   suffices to show that
the natural morphism \begin{equation}\label{fwefkjhewkjew89279834u3kejwfwei8ru32423r2}
\cH(\sigma \Omega A_{\infty})\to \cH(\Omega A_{\infty})
\end{equation} is an equivalence.  This follows from \cite[Lemma 7.15]{Bunke:2013uq}. For completeness we give the argument. 

For $k\in \Z$ let $\sigma^{\ge -k} \Omega A$ be the truncation to degrees $\ge -k$ defined similary as $\sigma=\sigma^{\ge 0}$.  Formally we can write  $\sigma^{\ge -k} \Omega A:= (\sigma \Omega A[-k])[k]$. Hence we can apply Lemma \ref{feb1601} and conclude that $\sigma^{\ge -k} \Omega A_{\infty}$ is a sheaf.

We have an isomorphism in $\Sm^{desc} (\Nerve(\Ch))$
$$\Omega A\cong \colim_{k\in \nat} \sigma^{\ge -k} \Omega A\ .$$
Since the localization preserves filtered colimits we get an equivalence \begin{equation}\label{vdsvsdvksdvdklvdklvlsdvsdvdv}
\Omega A_{\infty}\simeq \colim_{k\in \Nerve(\nat)} \sigma^{\ge -k} \Omega A_{\infty} 
\end{equation} in $\Sm^{desc} (\Nerve(\Ch)[W^{-1}])$.

For every
$k\in \nat$ we have a short exact sequence of sheaves of chain complexes
$$0\to \sigma^{\ge -k} \Omega A \to \sigma^{\ge -k-1} \Omega A \to \Omega A^{-k-1} \to 0\ ,$$
where $\Omega A^{-k-1}$ is considered as a sheaf of complexes  with only one non-trivial component in degree $-k-1$.
The localization sends short exact sequences of chain complexes in $\Ch$ to cofibre  sequences in $\Nerve(\Ch)[W^{-1}]$. Therefore we get a cofibre sequence
\begin{equation}\label{dqwdlkjqwldkqwdqwdwqdqwdqw}
\sigma^{\ge -k} \Omega A_{\infty}\to \sigma^{\ge -k-1}\Omega A_{\infty}\to   \Omega A^{-k-1}_{\infty} \end{equation}
in the stable $\infty$-category $\Sm(\Nerve(\Ch)[W^{-1}])$.
Since, as observed above,   the first two entries of \eqref{dqwdlkjqwldkqwdqwdwqdqwdqw} are sheaves  the third, $\Omega A^{-k-1}_{\infty}$, is a sheaf, too.
Using the fact, that in a stable $\infty$-category like  $\Sm^{desc} (\Nerve(\Ch)[W^{-1}])$    fibre sequences and cofibre sequences are the same and fibre sequences in sheaves are detected objectwise, we  conclude that that
\eqref{dqwdlkjqwldkqwdqwdwqdqwdqw} is a cofibre sequence in $\Sm^{desc}(\Nerve(\Ch)[W^{-1}])$.
As a left-adjoint the homotopification $\cH$ preserves cofibre sequences. We get a cofibre sequence
 $$\cH(\sigma^{\ge -k} \Omega_{\infty})\to \cH(\sigma^{\ge -k-1}\Omega A_{\infty}) \to \cH( \Omega A_{\infty}^{-k-1})\ .$$
 We claim that the first map    is an equivalence. To this end
 it suffices to show that $ \cH( \Omega A_{\infty}^{-k-1})\simeq 0$.
 Here we could apply   \cite[Lemma 7.13]{Bunke:2013uq}. 
 In detail, because of \eqref{ufeufioewufouaweofuoewfuoewifuewofui} and \eqref{jlkjelwjfweklfwpfoiopi32pr23r32r2r} it suffices to show that $\bs(\Omega A_{\infty}^{-k-1})\simeq 0$.
 Explicitly, by \eqref{dqwhdkjqwhdwqdqdhwqkd89ue12e12ee1e21e1e12e12} the evaluation of this presheaf on a manifold $M$ is given by
 $$\colim_{\Nerve(\Delta)^{op}} \Omega A_{\infty}^{-k-1}(M\times \Delta^{\bullet})\ .$$
 Using \eqref{dqwhqwdiwjhdwqjdwqdijwqdoijwqdoiwqdu9879879832434241} this is equivalent to the localization of 
 $$ \mathrm{tot} (\Omega A^{-k-1}(M\times \Delta^{\bullet}))  \ .$$
 Using partitions of unity one    easily checks that this complex is acyclic.
 
 We now have shown that the    natural embeddings induce a chain of equivalences
\begin{equation}\label{fwhefkjwehk32423424ewfewfewfweewfff}
\cH(\sigma^{\ge 0}\Omega A_{\infty})\stackrel{\simeq}{\to} \cH(\sigma^{\ge -1}\Omega A_{\infty}) \stackrel{\simeq}{\to} \cH(\sigma^{\ge -2}\Omega A_{\infty})  \stackrel{\simeq}{\to} \dots\ .
\end{equation}   As a left-adjoint $\cH$ commutes with colimits. We conclude that
 $$\cH(\sigma \Omega A_{\infty})\stackrel{\eqref{fwhefkjwehk32423424ewfewfewfweewfff}}{\simeq} \colim_{\Nerve(\nat)}  \cH(\sigma^{\ge -k}\Omega A_{\infty})\simeq
 \cH( \colim_{\Nerve(\nat)}  \sigma^{\ge -k}\Omega A_{\infty})\stackrel{\eqref{vdsvsdvksdvdklvdklvlsdvsdvdv}}{\simeq} \cH(\Omega A_{\infty})\ ,$$ i.e.
 the morphism \eqref{fwefkjhewkjew89279834u3kejwfwei8ru32423r2} 
 is an equivalence as required. 
 \hB

}


\begin{rem}{\rm  In this  remark we discuss the symmetric monoidal aspect of the construction of $\Omega A$ and the de Rham equivalence.     

   \bigskip
   
We first observe that   
the sheafification $L:\Sm(\Nerve(\Ch_{\R}))\to \Sm^{desc}(\Nerve(\Ch_{\R}))$ is symmetric monoidal.
Indeed, for $C,D\in \Sm(\Nerve(\Ch_{\R}))$ the natural morphism
\begin{equation}\label{gergjerlgregregregreogjpo34ip34t4t3t}L(C\otimes^{pre}_{\R}D)\to L(C)\otimes_{\R}L(D)
\end{equation}
is an isomorphism. By Remark \ref{dqwkldjqlkdjqwldjlwqdjwqlkdjwqldj} this can   checked on stalks, where the argument is straightforward.

\bigskip

 Note that the wedge product
  turns the de Rham complex $\Omega$  into a sheaf of commutative differential graded algebras over $\R$. Consequently,
the formula \eqref{regreerreg34543536777754z45}      defines a lax symmetric monoidal functor
\begin{equation}\label{wjkegdjkhkqwdqdqwdqwdqwdqwdqd}
\Nerve(\Ch_{\R})\to \Sm^{desc}(\Nerve(\Ch_{\R}))\ .
\end{equation}
In particular, 
the monoidal structure map is given by the composition
\begin{eqnarray*}\Omega A\otimes_{\R} \Omega B &\stackrel{\eqref{regreerreg34543536777754z45}}{\cong} & \Omega \otimes_{\R} L(\underline{A})\otimes_{\R} \Omega \otimes_{\R} L(\underline{B})\\&\cong& \Omega\otimes_{\R} \Omega\otimes_{\R} L(\underline{A})\otimes_{\R}L(\underline{B})\\&\to&  \Omega\otimes_{\R} L(\underline{A \otimes_{\R}B})\\&\stackrel{\eqref{regreerreg34543536777754z45}}{\cong}&   \Omega(A\otimes_{\R}B)\ ,
\end{eqnarray*}
where the second map uses symmetries and in the third map we use the product $\wedge:\Omega\otimes_{\R}\Omega\to \Omega$ together with the isomorphism $ L(\underline{A \otimes_{\R}B})\stackrel{ }{\cong} L(\underline{A})\otimes_{\R}L(\underline{B})
$   implied by the isomorphism \eqref{gergjerlgregregregreogjpo34ip34t4t3t}.

\bigskip

We have two functors
$$H(\Omega-)\ , \Funk(H(-)):\Nerve(\Ch_{\R})\to \Sm^{desc}(\Sp)$$
given by
$$A\mapsto H(\Omega A)\ , \quad A\mapsto \Funk(H(A))\ .$$
 The first is a composition of lax symmetric monoidal functors
$$\hspace{-0.5cm}\Nerve(\Ch_{\R})\stackrel{\eqref{wjkegdjkhkqwdqdqwdqwdqwdqwdqd}}{\to} \Sm^{desc}(\Nerve(\Ch_{\R}))  \stackrel{X\mapsto L(X_{\infty})}{\to}
 \Sm^{desc}(\Nerve(\Ch_{\R})[W^{-1}])\stackrel{H_{\R}}{\simeq} \Sm^{desc}(\Mod(H\R))\to \Sm^{desc}(\Sp)\ ,$$
 where the last map forgets the $H\R$-module structure.
 The second can be written as a composition of lax symmetric monoidal functors
 $$\Nerve(\Ch_{\R})\stackrel{H_{\R}}{\simeq} \Mod(H\R)\to \Sp\stackrel{\Funk}{\to} \Sm^{desc}(\Sp)\ ,$$
 where we use Remark \ref{ijlkdjqwldqwdqwdwqdwqdqwd} for the last map.

\bigskip

By Corollary \ref{lqejhdlwdjqwodu321oidwqd} we have a natural de Rham equivalence
$$j:H(\Omega-)\stackrel{\simeq}{\to} \Funk(H(-))$$
between the underlying functors.
\begin{prop}\label{flkijweflwefewfewfewfwe23442434}
The de Rham equivalence refines to an equivalence  between lax symmetric monoidal functors. In particular, for chain complexes $A,A^{\prime}\in \Ch_{\R}$ we have a commuting diagram
 \begin{equation}\label{fewlkjfwlkejlwefewfewfwefewfewfwe}
\xymatrix{H(\Omega A) \wedge H(\Omega A^{\prime})\ar[r]^{j\wedge j}\ar[d]&\Funk(HA)\wedge \Funk(HA^{\prime})\ar[d]\\
 H(\Omega (A\otimes_{\R} A^{\prime}))\ar[r]^{j}&\Funk(H(A\otimes_{\R}A))}\ .
\end{equation}

\end{prop}
\proof
We first  claim that the equivalence $$ \Funk(\ev_{*}(-))\stackrel{\eqref{dkqwkdqwdjdljdlkwd9879kwjefkjefewjfhewkf} }{\to} \id$$ of  functors on $\Sm^{desc,h}(\Sp)$ is lax symmetric monoidal.  In order see this observe that all transformations   in the composition \eqref{jhffhskdhfkshfkwhekhfkwfwefwfwefwefweff} are  symmetric monoidal transformations between lax symmetric monoidal functors. We apply the equivalence to $H(\Omega -)$ and note that \begin{equation}\label{xqwijisoiqwosoiqwiswq90809123132}
 \Funk(\ev_{*}(H(\Omega -)))\simeq \Funk(H(-)) \end{equation}
is an equivalence of lax symmetric monoidal functors.
Hence $j$ given as the composition
$$H(\Omega(-))\stackrel{\eqref{dkqwkdqwdjdljdlkwd9879kwjefkjefewjfhewkf}}{\simeq}
\Funk(\ev_{*}(H(\Omega -)))\stackrel{\eqref{xqwijisoiqwosoiqwiswq90809123132}}{\simeq}\Funk(H(-))$$
is an equivalence of lax symmetric monoidal functors.
\hB


}\hB \end{rem}

\subsection{Function spectra with proper support}\label{jun148}

In this subsection, for a map $W\to B$, we consider the evaluation of smooth objects on $W$   with proper support over $B$. We extend the de Rham equivalence to this context.  In the present paper the contents of this subsection is only  used in the proof   of Proposition \ref{jul1713}. But  it lays the foundations for the construction of integration maps in differential cohomology theories in general  as sketched in \cite{skript}.

\subsubsection{Evaluation of smooth objects with proper support}\label{kkfkjewfjwlefo234uro2r23r32r}

We consider a map $\pi:W\to B$ between manifolds. A subset $K\subseteq W$ is called proper over $B$ if
the restriction $\pi_{|K}:K\to B$ is proper, i.e. it has property that for every compact subset $L\subseteq B$ the subset
$\pi^{-1}(L)\cap K$ is compact.  We define a subcategory $$\Subm\subset \Mor(\Mf)$$ as follows:
\begin{enumerate}\item $\Subm$ has the same objects as $\Mor(\Mf)$. \item The morphisms in $\Subm$     are  certain morphisms in $\Mor(\Mf)$, namely commutative diagrams
\begin{equation}\label{jun105111}\xymatrix{W^{\prime}\ar[d]^{\pi^{\prime}}\ar[r]^{F}&W\ar[d]^{\pi}\\B^{\prime}\ar[r]^{f}&B}\end{equation}  such that $F$ is relatively proper. This condition means that for every subset $K\subseteq W$ which is proper over $B$ the subset $F^{-1}(K)\subseteq W^{\prime}$ is proper over $B^{\prime}$.
 \end{enumerate}
 
 \begin{ex}\label{ewhkjfhkfhewkfhhefkhefhfkwhfkhfoiuoiu}{\rm 

\begin{enumerate}
\item  If the square  \eqref{jun105111} is cartesian, then $F$ is relatively proper. Consequently, 
 the category $\Subm$ contains $\Bundle$ introduced in Subsection \ref{jun143} as a subcategory.
 
\item  
 Our main example of an object $(W\to B)\in \Subm$ which does not belong to $\Bundle$  is a real   vector bundle. 
 It does not belong {to} $\Bundle$ since the fibres are not compact.
\item  
If $W\to B$ is an object of $\Bundle$ (and not the identity), then the fibrewise diagonal
$(W\to B)\to (W\times_{B}W\to B)$ is a morphism in $\Subm$ which does not belong to $\Bundle$ since
the corresponding square is not cartesian.
\item If 
  $U\to B$ and $V\to B$ are vector bundles and $z:U\to U\oplus V$ is induced by the zero section of $V$, then
$(\id_{B},z):(U\to B)\to (U\oplus V\to B)$ is a morphism in $\Subm$. 
\item The morphism 
$(\R\to *)\to (*\to*)$ in $\Mor(\Mf)$ does not belong to $\Mor(\Subm)$ since it is not relatively proper.
\end{enumerate}
}
\end{ex}

Following our previous conventions, for an $\infty$-category $\bC$ we define the $\infty$-category
$$\Subm(\bC):=\Fun(\Nerve(\Subm)^{op},\bC)$$ of $\Subm$-objects in $\bC$.
We extend the  constructions \eqref{may281} from $\Bundle$ to $\Subm$.  
In particular we consider the functors
$$(-)d\ , (-)r:\Sm(\bC)\to \Subm(\bC)\ , \quad X\mapsto Xd\ ,Xr\ .$$
\begin{ex}{\rm For $X\in \Sm(\bC)$ and $(W\to B)\in \Subm$ we have  $$Xd(W\to B)\simeq X(W)\ , \quad Xr(W\to B)\simeq X(B)\ .$$ }\hB\end{ex}
Our next goal is  the definition of an object $Xd_{c}\in \Subm(\bC)$   such that
$Xd_{c}(W\to B)$ is the evaluation of $Xd$ on $W$ with proper support over $B$.

\begin{rem}\label{dhqwkjdhqwkdhqwdkqwhwqkddwqdwqdwdqwddqwdd}{\rm In this remark we explain the construction of $Xd_{c}(W\to B)$ morally.
If  $K\subseteq W$   proper over $B$,  then  $W\setminus K\subseteq W$  is an open subset.
Under appropriate conditions on $\bC$ we can form the fibre
$$Xd_{K}(W\to B):=\Fib(X(W)\to X(W\setminus K))\ .$$
{This should be interpreted as sections of $X$ on $W$ which are supported on $K$. We then define the sections of $X$ on $W$ with proper support over $B$ by}
$$Xd_{c}(W\to B):=\colim_{K} Xd_{K}(W\to B)\ ,$$
where the colimit is taken over the subsets $K\subseteq W$ which are proper over $B$.
The main point of the following discussion is to make this idea precise.
 }\hB

\end{rem}

\bigskip

First of all we must assume that $\bC$ allows the definition of fibres. To this end
we   assume that $\bC$ is a pointed     $\infty$-category (Definition \ref{fkjwelfkwefwefwefewfewe3455})   which admits finite limits. 
\begin{ex}{\rm 
Examples of pointed $\infty$-categories with finite limits  are stable $\infty$-categories like $\Sp$ or $\Nerve(\Ch)[W^{-1}]$, but also  pointed spaces $\Spc_{*}$ or chain complexes $\Nerve(\Ch)$.
}\hB
\end{ex}
Under this assumption on $\bC$
we can then define a functor
$$\Fib:\Mor(\bC)\to  \bC $$
which maps a morphism $f:C\to D$ to  its fibre.
  
\begin{rem}{\rm 
Here are more details on this construction. 
In a first step, using the universal property an initial object $0$,  we can define a functor
$$\Mor(\bC)\to \bC^{\Nerve(\bullet\to\bullet\leftarrow \bullet)}$$ which sends 
a morphism $C\to D$ to the pull-back diagram $C\to D\leftarrow 0$. 
We define $\Fib$ as the composition of this functor with $\lim_{\Nerve(\bullet\to\bullet\leftarrow \bullet)}$.
}\hB
\end{rem}

\begin{ex}\label{kljfewlkfjweklewfouoi23r23r32r23rwefwefwef}{\rm
 In this example we explain and compare the fibres of a morphism  between of chain complexes and its localization.
Thus let  $\phi:C\to D$ is a morphism in $\Nerve(\Ch)$. Then its fibre is given by  the kernel
$\Fib(\phi):=\ker(\phi) $ of $\phi$. Explicitly, this is the chain complex whose group in degree $i\in \Z$ is
$\ker(\phi:C^{i}\to D^{i})$, and whose differential is induced from the differential of $C$ by restriction.

We now consider the fibre of the localization $\phi_{\infty}:C_{\infty}\to D_{\infty}$ of $\phi$. It can be represented by the mapping cone of $\phi$, i.e. we have an equivalence
$$\Fib(\phi_{\infty})\simeq \Sigma^{-1}\Cone(\phi)_{\infty}\ .$$ Here $\Sigma^{-1}\Cone(\phi)$ is the complex
with underlying $\Z$-graded group
$C\oplus D[-1]$. A degree $i$-cochain will be written as a pair
$(\omega,\alpha)$ with $\omega\in C^{i}$ and $\alpha\in D^{i-1}$.
The differential of $\Sigma^{-1}\Cone(\phi)$ is given by $d(\omega,\alpha):=(d\omega,\phi(\omega)-d\alpha)$.
Observe that we have a natural map of chain complexes
$$\ker(\phi)\to \Sigma^{-1}\Cone(\phi)\ , \quad \omega\mapsto (\omega,0)\ .$$
It induces the natural map
$$\Fib(\phi)_{\infty}\to \Fib(\phi_{\infty})\ ,$$ which we will abstractly 
construct  in Lemma \ref{jdhdjdkdkwqhdqwkjhdwqkhdwqdwqdwqdwqdwqdwqdwq} below.
In general this map  
 is not a quasi-isomorphism.
}\hB\end{ex}

For an object $(W\to B)\in \Subm$ we consider the poset \begin{equation}\label{thjkhtjkkjer8guo3ut3ot45}P(W\to B):= \{K\subseteq W\:|\:\mbox{$K$ is proper over $B$}\}\ .
\end{equation}
Note that $P(W\to B)$ is filtered since for $K,L\in P(W\to B)$ we have
$K\cup L\in P(W\to B)$.  By the  definition of a morphism \eqref{jun105111}   in $\Subm$  we get a map of posets
$$F^{-1}:P(W\to B)\to P(W^{\prime}\to B^{\prime})\ , \quad K\mapsto F^{-1}(K)\ .$$
In particular we have defined a functor
$$P:\Subm^{op}\to \Cat\ .$$

We let $\widetilde{\Subm}$ be the Grothendieck construction of this functor. \begin{rem}{\rm  An object of  
$\widetilde{\Subm}$ is a pair $(W\to B,K)$ such that $K\in P(W\to B)$. A morphism $(W^{\prime}\to B^{\prime},K^{\prime})\to  (W\to B,K)$ is a morphism \eqref{jun105111}  in $\Subm$ such that $  F^{-1}(K)\subseteq K^{\prime}$. }\hB \end{rem}

We have a forgetful functor
$$q:\widetilde{\Subm}\to  \Subm\ .$$
We can now consider the following constructions:
\begin{enumerate}
\item We have a functor
$$\widetilde{\Subm}\to \Mor(\Mf)\ , \quad (W\to B,K)\to (W\setminus K\to W)\ . $$  
\item
Pull-back along this functor provides
$$\Sm(\bC)\to \Fun( \Nerve(\widetilde{\Subm})^{op},\Mor(\bC))\ , \quad X\mapsto [(W\to B,K)\mapsto (X(W)\to X(W\setminus K))]\ .$$
\item If we compose this with $\Fib$, then we get  functor 
\begin{equation}\label{ferfjrefirefjio9irore23r32r23r32r32r31}\Sm(\bC)\to  \widetilde{\Subm}(\bC)\ , \quad X\mapsto \tilde X\ .\end{equation}
We have \begin{equation}\label{ferfjrefirefjio9irore23r32r23r32r32r3}
 \tilde X(W\to B,K)\simeq Xd_{K}(W\to B):= \Fib(X(W)\to X(W\setminus K))\ .\end{equation}

 \end{enumerate}

\begin{ddd}\label{dijqwidjqwldjqwlqwdqwdqwdwqdwqdqwdqw} 
We define the functor
$$Xd_{c}:\Sm(\bC)\to \Subm(\bC)$$
as the composition
$$\Sm(\bC)\stackrel{ \eqref{ferfjrefirefjio9irore23r32r23r32r32r31}}{\to}    \widetilde{\Subm}(\bC)\stackrel{q_{*}}{\to} 
 \Subm (\bC)\ , $$ where $q_{*}$ is the  
  left Kan extension  along the functor $q$.\end{ddd}


   \begin{rem}{\rm Note that $\tilde X(W\to B,K)$ is equivalent to the object $Xd_{K}$ considered in  Remark \ref{dhqwkjdhqwkdhqwdkqwhwqkddwqdwqdwdqwddqwdd}. 
 By the pointwise formula for the left Kan extension we get the equivalence \begin{equation}\label{fwefwejfewlkfewfoifuewoiewfjewifjewoifjeioj32io3u4284u2394u2398423r3kjwehfwkejfnwelkfwe}
Xd_{c}(W\to B)\simeq \colim_{K \in P(W\to B) }Xd_{K}(W\to B)  \end{equation}  as required in Remark \ref{dhqwkjdhqwkdhqwdkqwhwqkddwqdwqdwdqwddqwdd}.
}\hB\end{rem}

\begin{ex}\label{fjkhefjkwefkjewfewfewfewfwefwef}{ \rm
For $n\in \Z$ and a spectrum $E$ the cohomology $$H^{n}_{c/B}(W;E):=\pi_{-n}(\Funk(E)d_{c}(W\to B))$$ is the $E$-cohomology of $W$ with proper support over $B$. We refer to Example \ref{fjeljweljewlfewiofofewufewfewfef}  for a further discussion.}\hB
\end{ex}

%

\subsubsection{Application to the de Rham complex}

In this subsection we apply the evaluation with proper support to the de Rham complex. We analyse the compatibility with localization.

\bigskip

\begin{rem} {\rm We consider   a map $\pi:W\to B$ between manifolds.
For an $A$-valued form    $\omega\in \Omega A(W)$  on $W$ we have a well-defined notion of support $\supp(\omega)\subseteq W$.
A point $w\in W$ belongs to the necessarily closed subset $\supp(\omega)$ if and only if $\omega$ does not vanish in every  neighbourhood of $w$. The subcomplex
$\Omega Ad_{c}(W\to B)\subset \Omega A(W)$ consists of all forms $\omega\in \Omega A(W)$ such that $\pi_{|\supp(\omega)}:\supp(\omega)\to B$ is proper. 

In contrast, using the cone model for the fibre as in Example \ref{kljfewlkfjweklewfouoi23r23r32r23rwefwefwef},  one could think of an element of $\Omega A_{\infty}d_{c}(W\to B)$ as being represented by a form  $\omega\in \Omega A(W)$ together with a form $\alpha\in  \Omega A(W\setminus K)$ for some $K\subseteq W$ which is  proper over $B$ and such that $d\alpha=\omega_{|W\setminus K}$.
In particular, $\omega$ does not have necessarily proper support over $B$. It has proper support ''up to a specified homotopy''.}\hB
\end{rem}

\begin{lem}\label{jdhdjdkdkwqhdqwkjhdwqkhdwqdwqdwqdwqdwqdwqdwq}
We have a natural morphism 
\begin{equation}\label{hdjkjhkwqwqdqwdqwdwqd}
\Fib(C\to D)_{\infty}\to \Fib(C_{\infty}\to D_{\infty}) 
\end{equation} in $\Nerve(\Ch)[W^{-1}]$
\end{lem}
\proof
The filler of the cartesian diagram
$$\xymatrix{\Fib(C\to D)\ar[r]\ar[d]&0\ar[d]\\
C\ar[r]&D}$$
induces a filler
of the diagram
 $$\xymatrix{\Fib(C\to D)_{\infty}\ar[r]\ar[d]&0_{\infty}\ar[d]\\
C_{\infty}\ar[r]&D_{\infty}}\ .$$
We get the required morphism  \eqref{hdjkjhkwqwqdqwdqwdwqd} from the universal property of $\Fib(C_{\infty}\to D_{\infty})$ as left upper corner of the cartesian diagram $$\xymatrix{\Fib(C_{\infty}\to D_{\infty})\ar[r]\ar[d]&0_{\infty}\ar[d]\\
C_{\infty}\ar[r]&D_{\infty}}\ .$$   An inspection of this proof shows that  \eqref{hdjkjhkwqwqdqwdqwdwqd}   actually refines to a transformation between functors
$$\Mor(\Nerve(\Ch))\to \Mor(\Nerve(\Ch)[W^{-1}])\ .$$
\hB

We apply this to the de Rham complex. In view of \eqref{ferfjrefirefjio9irore23r32r23r32r32r3} we get a morphism $$(\widetilde{\Omega A})_{\infty}\to \widetilde{\Omega A_{\infty}}$$
between functors from $\Nerve(\widetilde{\Subm})^{op}$ to $\Nerve(\Ch)[W^{-1}]$.
By Definition \ref{dijqwidjqwldjqwlqwdqwdqwdwqdwqdqwdqw}, applying $q_{*}$ to this morphism, we get a morphism
\begin{equation}\label{fjhefjwehfkwefewfew3244234324244324}
(\Omega A d_{c})_{\infty}\to   \Omega A_{\infty}d_{c}\ .
\end{equation}
For the de Rham complex the conditions of proper support and proper support up  to a specified homotopy are equivalent.
More precisely we have the following Lemma.  \begin{lem}
The morphism \eqref{fjhefjwehfkwefewfew3244234324244324} is an equivalence.
\end{lem}
\proof
The assertion is a consequence of the fact that $\Omega A$ consists of level wise
fine sheaves. We fix a map $W\to B$ between manifolds and must show that
$$(\Omega A d_{c})_{\infty}(W\to B)\to  \Omega A_{\infty}d_{c}(W\to B)$$
is an equivalence in $\Nerve(\Ch)[W^{-1}]$. We thus must show that for every $ n\in \nat$   the induced map in cohomology
$$H^{n}(\Omega A d_{c}(W\to B))\to  H^{n}(\Omega A_{\infty}d_{c}(W\to B))$$
  is an isomorphism.
  Now note that cohomology commutes with filtered colimits and the 
  poset  $P(W\to B)$ defined in \eqref{thjkhtjkkjer8guo3ut3ot45} is filtered. Using \eqref{fwefwejfewlkfewfoifuewoiewfjewifjewoifjeioj32io3u4284u2394u2398423r3kjwehfwkejfnwelkfwe} we therefore must show that \begin{equation}\label{fewfwfewjkfwefkljeflejflfjoiu23ou2323424232rfeewfwfwefwe}
\colim_{K\in P(W\to B)}H^{n}(\Omega A d_{K}(W\to B))\to \colim_{K\in P(W\to B)} H^{n}(\Omega A_{\infty}d_{K}(W\to B))
\end{equation}
is an isomorphism for every $ n\in \nat$.

\bigskip
 
We first show that \eqref{fewfwfewjkfwefkljeflejflfjoiu23ou2323424232rfeewfwfwefwe}   is injective. For the explicit calculation we use the identification 
$$\Omega A_{\infty}d_{K}(W\to B)\simeq \Sigma^{-1}\Cone(\Omega A(W)\to \Omega A(W\setminus K))_{\infty}$$
(see Example \ref{kljfewlkfjweklewfouoi23r23r32r23rwefwefwef}).
 Fix $K\in P(W\to B)$ and consider 
$[\omega]\in H^{n}(\Omega A d_{K}(W\to B))$ represented by $\omega\in \Omega A^{n}(W)$ such that $\supp(\omega)\subseteq K$.  If $[\omega]$ maps to zero under  \eqref{fewfwfewjkfwefkljeflejflfjoiu23ou2323424232rfeewfwfwefwe}, then there exists a subset $K^{\prime}\in P(W\to B)$ containing $K$   and a pair of forms $(\alpha,\beta)$ such that $\alpha\in \Omega A^{n-1}(W )$ and $\beta\in \Omega A^{n-2}(W\setminus K^{\prime} )$  and
$\omega=d\alpha$, $\alpha_{|W\setminus K^{\prime}}=d\beta$. We can find a subset $K^{\prime\prime}\in P(W\to B) $  which  contains $K^{\prime}$ in its interior. Using a partition of unity we can extend the restriction of  $\beta$ to $ W\setminus K^{\prime\prime}$   to  a form $\tilde \beta\in \Omega A^{n-2}(W)$.
 Then $\supp(\alpha-d\tilde \beta)\subseteq K^{\prime\prime}$ and $d  (\alpha-d\tilde \beta)=\omega$. Hence the image of 
$[\omega]$ in $H^{n}(\Omega A d_{K^{\prime\prime}}(W\to B))$ and therefore in $\colim_{K\in P(W\to B)}H^{n}(\Omega A d_{K}(W\to B))$ vanishes. 

\bigskip 
Next we show that \eqref{fewfwfewjkfwefkljeflejflfjoiu23ou2323424232rfeewfwfwefwe}   is surjective.
Let $(\alpha,\beta)\in  \Sigma^{-1}\Cone(\Omega A(W)\to \Omega A(W\setminus K))^{n}$ represent a 
class in $\colim_{K\in P(W\to B)} H^{n}(\Omega A_{\infty}d_{K}(W\to B))$. 
Then $\alpha\in \Omega A^{n} (W)$ and $\beta\in \Omega^{n-1}(W\setminus K)$ satisfy
$d\alpha=0$ and $\alpha_{|W\setminus K}=d\beta$. 
We can find  a subset $K^{\prime}\in P(W\to B)$ which contains $K$ in its interior.
Using a partition of unity we can choose an  extension $\tilde \beta\in \Omega A^{n-1}(W)$ of  the restriction of $\beta$ to $K^{\prime}$.
Then $\alpha-d\tilde \beta$ is supported in $K^{\prime}$. The class
$[\alpha-d\tilde \beta]\in H^{n}(\Omega A d_{K^{\prime }}(W\to B))$ represents the required preimage.
Indeed, we have the relation $(\alpha-d\tilde \beta,0)-(\alpha,\beta)=-d(\tilde \beta,0)$
in $\Sigma^{-1}\Cone(\Omega A(W)\to \Omega A(W\setminus K^{\prime}))^{n}$. \hB

Recall the de Rham equivalence $j:H(\Omega A)\stackrel{\simeq}{\to} \Funk(H(A))$, see \eqref{hgjhgdjhgqwd873iu2kjd31}.  The functoriality of $X\mapsto Xd_{c}$ (Definition \ref{dijqwidjqwldjqwlqwdqwdqwdwqdwqdqwdqw})   implies:
\begin{kor}\label{fewklfjwelkfjlwefwecdce3rrrfwefewfewfewfewf}
We have a  de Rham equivalence with proper support
$$j_{c}: H(\Omega Ad_{c})\stackrel{\eqref{fjhefjwehfkwefewfew3244234324244324}}{\simeq} H( \Omega A)d_{c} \to \Funk(H(A))d_{c}\ .$$
For a morphism \eqref{jun105111} in $\Subm$ the following diagram commutes
\begin{equation}\label{gergkedtjtjt3k4tjl3t34t053405} \xymatrix{H( \Omega Ad_{c})(W \to B )\ar[r]^{F^{*}}\ar[d]^{j_{c}}& H( \Omega Ad_{c})(W^{\prime}\to B^{\prime}) \ar[d]^{j_{c}}\\
\Funk(H(A))d_{c}(W \to B )\ar[r]^{(f,F)^{*}}&\Funk(H(A))d_{c}(W^{\prime}\to B^{\prime}) }\ .
\end{equation}
\end{kor}

\subsubsection{Extension by zero}

In this subsection we will use the sheaf property and therefore assume that $\bC$ is pointed and presentable.
In particular it then has finite limits so that the previous theory applies.

We consider a morphism \eqref{jun105111} in $\Subm$ such that $F$ is an open embedding and $f$ is proper. We will call such a map a relatively proper embedding. Under this assumption, if $K^{\prime}\subseteq W^{\prime}$
is proper over $B^{\prime}$, then $F(K^{\prime}
)\subseteq W$ is proper over $B$. Indeed, if $L\subseteq B$ is compact, then
$$F(K^{\prime})\cap \pi^{-1}(L)=F(K^{\prime}\cap F^{-1} \pi^{-1}(L))=F(K^{\prime}\cap \pi^{\prime,-1}f^{-1}(L))$$
is compact.

\begin{ex}{\rm The map $\emb$ occuring in Definition \ref{jul2001}, 3. is an example of  a relatively proper embedding. In the proof of 
Proposition \ref{jul1713} we need the covariant functoriality $\emb_{*}$ and its compatibility with
the de Rham equivalence.
}\hB
\end{ex}

For a smooth object $X\in \Sm^{desc}(\bC)$ we want to define \begin{equation}\label{}
F_{*}:Xd_{c}(W^{\prime}\to B)\to Xd_{c}(W \to B)\ .
\end{equation}
 To this end we first show excision.
\begin{lem} Assume that $X\in \Sm^{desc}(\bC)$  and $(f,F)$  is a relatively proper embedding. If $K^{\prime}\subseteq W^{\prime}$ is proper over
$B^{\prime}$, then the natural map \begin{equation}\label{hdjhkjhwqkdhwqdo9udowqdjwqldkjqdlqwjdlqwdqwdqwd}
F_{K^{\prime}}^{*}:Xd_{K}(W\to B)\to Xd_{K^{\prime}}(W^{\prime}\to B)
\end{equation}
is an equivalence, where $K:=F(K^{\prime})$.
\end{lem}
\proof
We consider the covering family $\cU:=(W^{\prime},W\setminus K)$ of $W$ and let $U^{\bullet}_{M}$ be the corresponding simplicial manifold. The sheaf condition for $X$ now provides the cartesian  square
$$\xymatrix{X(W)\ar[r]\ar[d]&X(W\setminus K)\ar[d]\\X(W^{\prime})\ar[r]&X(W^{\prime}\setminus K^{\prime})}\ .$$ 
It induces an equivalence of the fibres of the horizontal maps which is exactly \eqref{hdjhkjhwqkdhwqdo9udowqdjwqldkjqdlqwjdlqwdqwdqwd}. 
\hB

If $K^{\prime}\subseteq K_{1}^{\prime}$ and
$K_{1}^{\prime}$ is still proper over $B$, then the corresponding diagram
$$\xymatrix{Xd_{K}(W\to B)\ar[d]\ar[r]_{\simeq}^{F_{K^{\prime}}^{*}}&Xd_{K^{\prime}}(W^{\prime}\to B)\ar[d]\\
Xd_{K_{1}}(W\to B)\ar[r]_{\simeq}^{F_{K^{\prime}_{1}}^{*}} &Xd_{K_{1}^{\prime}}(W^{\prime}\to B)}$$ commutes. 
 We can thus define the map $F_{*}$ as the composition 
\begin{eqnarray*}
 Xd_{c}(W^{\prime}\to B)&\simeq &\colim_{K^{\prime}\in P(W^{\prime}\to B^{\prime})} Xd_{K^{\prime}}(W^{\prime}\to B)\\
 &\simeq &\colim_{K^{\prime}\in P(W^{\prime}\to B^{\prime})} Xd_{F(K^{\prime})}(W \to B)\\
 &\to&\colim_{K\in P(W\to B) } Xd_{ K}(W \to B)\\
& \simeq& Xd_{c}(W\to B)
\end{eqnarray*}

An inspection of this construction shows that it is natural in $X\in   \Sm^{desc}(\bC)$.
In other words:

\begin{kor}\label{jekljdlkwjwkldjlqwkjwqldqwdwqdqwdqwdqwdqwd} For a relatively proper embedding $(f,F)$   we have constructed  a morphism
$$F_{*}:\Sm^{desc}(\bC)\to \Mor(\bC)\ .$$
\end{kor}

\begin{ex}{\rm 
For $\bC=\Nerve(\Ch)$ the functor $F_{*}$ is simply given by extension by zero. This applies e.g. to
$$F_{*}:\Omega Ad_{c}(W^{\prime}\to B^{\prime})\to \Omega A d_{c}(W\to B)\ .$$
The point of the construction above is to extend the idea of extension by zero to the $\infty$-categorical setting.
}\hB\end{ex}

\begin{rem}\label{hckjjklejclejlejdldkjqwdwqdqwdqwd}{\rm 
An inspection of the construction of $F_{*}$ further shows that
if $(f^{\prime},F^{\prime}):(W\to B)\to (W^{\prime\prime}\to B^{\prime\prime})$ is a second  relatively proper embedding, then we have a natural equivalence
\begin{equation}\label{kcndncnlklkqklejwqlkejlqw}
F_{*}^{\prime}\circ F_{*}\simeq (F^{\prime}\circ F)_{*}\ .
\end{equation}
Furthermore, given a commuting diagram in $\Subm$
$$\xymatrix{(W_{11} \to B_{11} )\ar[d]^{(g_{ 1},G_{ 1})}\ar[r]^{(f_{1},F_{1})}&(W_{01}\to B_{01})\ar[d]^{(g,G)}\\
(W_{10} \to B_{10})\ar[r]^{(f,F)}&(W\to B)}$$
where $(f,F)$ and $(f_{1},F_{1})$ are relatively proper embeddings, 
we have a natural equivalence
\begin{equation}\label{hckjjklejclejlejdldkjqwdwqdqwdqwd1}
F_{1,*}\circ (g_{1},G_{1})^{*}\simeq (g ,G)^{*}\circ F_{*}\ .
\end{equation}
} \hB
\end{rem}

We apply Corollary \ref{jekljdlkwjwkldjlqwkjwqldqwdwqdqwdqwdqwdqwd}     to the de Rham equivalence $j\in \Mor(\Sm^{desc}(\Sp))$.
Since $F_{*}$ obviously commutes with \eqref{fjhefjwehfkwefewfew3244234324244324} we get:
 \begin{kor}\label{ihdwqdqwdlqwjdwqldjqwldjwlkdjwdoiuoiudoiqwudqwidwqdwqdqwdwqd} Let $F$ be a relatively proper embedding.
Then we have a commuting square 
$$\xymatrix{H(\Omega Ad_{c})(W^{\prime}\to B^{\prime})\ar[r]^{j_{c}}\ar[d]^{F_{*}}&\Funk(H(A))d_{c}(W^{\prime}\to B^{\prime})\ar[d]^{F_{*}}\\
H(\Omega Ad_{c})(W \to B)\ar[r]^{j_{c}}&\Funk(H(A))d_{c}(W\to B)}\ .$$
\end{kor}

 \subsubsection{The exterior product}

We now consider the exterior product. 

\begin{rem}{\rm 
The category of spectra $\Sp$ is symmetric monoidal with respect to $\wedge$.  
By Remark \ref{ijlkdjqwldqwdqwdwqdwqdqwd}, for   topological spaces $X,Y$ and spectra $E,F$  we get a natural map \begin{equation}\label{fkwfkjewfewki2oiuio3h32fhiu2f23f23f32f}
E^{X}\wedge F^{Y}\to (E\wedge F)^{X\times  Y}
\end{equation}
by specializing \eqref{kfjkfjklwejfwejfewopipoiopiopieqweqeqweeqeeqwe}.
  }\hB \end{rem}

We consider spectra $E,F\in \Sp$,
 a manifold $B$ and a map
$W\to C$ between manifolds. Then we define a map
$$\Funk(E)(B)\wedge \Funk(F)d_{c}(W\to C)\to \Funk(E\wedge F)d_{c}(B\times W\to B\times C)\ .$$
It is induced by
$$\Funk(E)(B)\wedge \Funk(F)d_{K}(W\to C)\to \Funk(E\wedge F)d_{B\times K}(B\times W\to B\times C)$$
which in turn is the induced map on the fibres of the horizontal maps in the square
$$\xymatrix{E^{B}\wedge F^{W}\ar[r]\ar[d]^{\eqref{fkwfkjewfewki2oiuio3h32fhiu2f23f23f32f}}&E^{B}\wedge F^{W\setminus K}\ar[d]^{\eqref{fkwfkjewfewki2oiuio3h32fhiu2f23f23f32f}}\\
(E\wedge F)^{B\times W}\ar[r]&(E\wedge F)^{B\times (W\setminus K)} }\ .$$
This product is natural in the data $E,F$ and $B, W\to C$, but we refrain from writing out the corresponding diagrams.

Similarly, for chain complexes $A,A^{\prime}\in \Ch_{\R}$
the wedge product of forms induces an exterior product
$$\Omega A(B)\otimes_{\R} \Omega A^{\prime}d_{c}(W\to C)\to \Omega (A\otimes_{\R} A^{\prime})d_{c}(B\times W\to B\times C)\ .$$
It is again natural in the  data.

\begin{lem}\label{fkjewkfjfklewjjewflejflewfjejfewlfewfwefewfewfewfe}
The de Rham equivalence is compatible with
the exterior products. In particular, we have a commuting diagram
$$\xymatrix{H(\Omega A)(B)\wedge H( \Omega A^{\prime}d_{c})(W\to C)\ar[d]\ar[r]^{j\wedge j_{c}}&
\Funk(H(A))(B)\wedge \Funk(H(A^{\prime}))d_{c}(W\to C)\ar[d]\\
H(\Omega(A\otimes A^{\prime})d_{c})(B\times W\to B\times C)\ar[r]^{j_{c}}&\Funk(H(A\otimes_{\R} A^{\prime}))d_{c}(B\times W\to B\times C)}\ .$$
\end{lem}
\proof
As a first step we observe that
$$\xymatrix{H(\Omega A)(B)\wedge H( \Omega A^{\prime} d_{c}) (W\to C)\ar[d]\ar[r]^{\id\wedge \eqref{fjhefjwehfkwefewfew3244234324244324}} &H(\Omega A{)}(B)\wedge H( \Omega A^{\prime} )d_{c}(W\to C) 
  \ar[d]\\
H(\Omega(A\otimes A^{\prime})d_{c})(B\times W\to B\times C)\ar[r]^{\eqref{fjhefjwehfkwefewfew3244234324244324}} &H(\Omega(A\otimes A^{\prime}))d_{c}(B\times W\to B\times C)\ }$$
commutes. It remains to show that
\begin{equation}\label{gerergiljkigrejreoigrgugoirerioureoiure345}
\xymatrix{H(\Omega A)(B)\wedge H( \Omega A^{\prime} )d_{c}(W\to C)\ar[d]\ar[r]^{j\wedge j_{c}}&
\Funk(H(A))(B)\wedge \Funk(H(A^{\prime}))d_{c}(W\to C)\ar[d]\\
H(\Omega(A\otimes A^{\prime}) )d_{c}(B\times W\to B\times C)\ar[r]^{j_{c}}&\Funk(H(A\otimes_{\R} A^{\prime}))d_{c}(B\times W\to B\times C)}\ .
\end{equation}
commutes. If $K\subseteq W$ is proper over $C$, then as a consequence of  \eqref{fewlkjfwlkejlwefewfewfwefewfewfwe} and the fact that the smash product with
a spectrum preserves fibre sequences we get 
$$\xymatrix{H(\Omega A)(B)\wedge H( \Omega A^{\prime} )d_{K}(W\to C)\ar[d]\ar[r]^{j\wedge j_{c}}&
\Funk(H(A))(B)\wedge \Funk(H(A^{\prime}))d_{K}(W\to C)\ar[d]\\
H(\Omega(A\otimes A^{\prime}) )d_{B\times K}(B\times W\to B\times C)\ar[r]^{j_{c}}&\Funk(H(A\otimes_{\R} A^{\prime}))d_{B\times K}(B\times W\to B\times C)}\ .$$ These diagrams are compatible with the inclusions of subsets $K\in P(W\to C)$ which are proper over $C$.
Taking the colimit over $K\in P(W\to C)$ we get \eqref{gerergiljkigrejreoigrgugoirerioureoiure345}. \hB

  \subsubsection{Integration}\label{kwjfweklfjewlfjelwfewopipi3ir23r20pr32r32r32r}

In order to discuss integration we must consider oriented submersions.
If $\pi:W\to B$ is a submersion, then we have a vertical bundle $T^{v}\pi:=\ker(d\pi)$. An orientation of $\pi$ is an orientation of the real vector bundle $T^{v}\pi$.
We let $\Submor$ be the category of submersions with smooth fibres and with a given orientation of $\pi$.
A  morphim in $\Submor$ is a diagram \eqref{jun105111} which is in addition cartesian and such that the map
$dF_{|T^{v}\pi^{\prime}}:T^{v}\pi^{\prime}
\to T^{v}\pi$ preserves orientations. 
We have a forgetful functor
$$\Submor\to \Subm\ .$$
We can thus define the functors
$$(-)d\ ,(-)r \ ,(-)d_{c}:\Sm(\bC)\to \Submor(\bC)$$ by restriction along this forgetful functor. 

\bigskip

Let $\pi:W\to B$ be an oriented  submersion with fibre dimension $n$. 
  If $\omega\in \Omega A(W)$ is such that $  \supp(\omega) $ is proper over $B$, then we can   define an integral \begin{equation}\label{fwehfkejfekfewfewfewfewfwfewfewfwefefw}
\int_{W/B}\omega\in \Omega A(B)[-n]\ .
\end{equation}
Consider a relatively compact  open subset   $U\subseteq B$. Then   $\supp(\omega)\cap \pi^{-1}(U)\subseteq W$ is relatively compact in $W$. 
  As noted in Remark \ref{jkdhwqjkdhkwqdwqhdqwdkjhwqdqwdwqdwqdwd} we have
  $$\omega_{|\pi^{-1}(U)}\in \Omega(\pi^{-1}(U)  )\otimes_{\R} A\ .$$
In particular we can write  $\omega$ as a finite linear combination \begin{equation}\label{fwefewefwfw2344rr34z4657868563r}
 \omega=\sum_{i\in I} \omega_{i}\otimes a_{i} 
\end{equation}
for $\omega_{i}\in \Omega(\pi^{-1}(U) )$ and $a_{i}\in A$. 
Then we can define the integral by
$$(\int_{W/B}\omega)_{|U}:=\sum_{i\in I}(\int_{W/B}\omega_{i })\otimes a_{i}\in \Omega A(U)[-n]\ .$$
The form $(\int_{W/B}\omega)_{|U}$ is well-defined  independently of the presentation \eqref{fwefewefwfw2344rr34z4657868563r}. Moreover, these local forms patch together and define
\eqref{fwehfkejfekfewfewfewfewfwfewfewfwefefw}. It follows from Stokes' theorem that 
$$d\int_{W/B}\omega=\int_{W/B}d\omega\ .$$
For every $(W\to B)\in \Submor$ we thus have constructed a morphism of complexes
$$\int_{W/B}:\Omega Ad_{c}(W\to B  )\to \Omega Ar(W\to B)[-n]\ .$$
The integration is compatible with morphisms in $\Submor$. 
We have therefore defined a transformation
$$\int_{rel}:\Omega A d_{c}\to  \Omega Ar [-n]\ .    $$
We add the subscript $rel$ in order to distinguish it from the integration \eqref{hdjwdkqjwdhkjhkhdkjhkqwdoiuoiqdqwdqwdqdd}
which is similar but involves the Euler form. Applying the localization $\Nerve(\Ch)\to \Nerve(\Ch)[W^{-1}]$
 we get an induced morphism
$$\int_{rel,\infty }:(\Omega A d_{c})_{\infty}\to  (\Omega Ar)_{\infty} [-n]\ .$$
 We finally define
 $$\int^{\prime}_{rel} :\Omega A_{\infty} d_{c} \to  (\Omega Ar)_{\infty} [-n]$$
 such that
 $$\xymatrix{(\Omega A d_{c})_{\infty}\ar[d]^{\int_{rel,\infty }}\ar[r]^{\eqref{fjhefjwehfkwefewfew3244234324244324}}_{\simeq}&\Omega A_{\infty} d_{c}\ar[d]^{\int_{rel}^{\prime}}\\
 (\Omega Ar)_{\infty} [-n])\ar@{=}[r]&(\Omega Ar)_{\infty}[-n]}$$
 commutes.
 
 \newcommand{\Cofib}{\mathrm{Cofib}}

\subsubsection{The Thom space approach}\label{kljdkljlwqjdqjdwuo1iu2e12e12e12}

We now consider an alternative of the evaluation with proper support which is defined for vector bundles and
homotopy invariant  smooth function objects, i.e. objects of the form $\Funk(C)\in \Sm(\bC)$. We   assume that the {pointed} $\infty$-category $\bC$ is presentable.

\bigskip

\newcommand{\RBundle}{\mathbf{VectBdl}} 

We consider the  subcategory of $$\RBundle\subset \Subm$$ consisting of vector bundles
and vector bundle morphisms. For $X\in \Sm(\bC)$ we use the notation $Xd$, $Xr$ and $Xd_{c}$ also for the
restriction of these functors along this inclusion.

\bigskip

The Thom space construction (see Remark \ref{djqwdkjqwdkdjqwldjqwdqwdqwdwqd}) provides a functor $$ \RBundle\to \Top_{*}\ , \quad (W\to B)\mapsto B^{W}\ ,$$
where $\Top_{*}$ denotes the category of pointed topological spaces. As before,
  we omit the subscript $top$ so that
  $B^{W}$ abbreviates the more  precise  notation
$B_{top}^{W_{top}}$. 

We further simplify the notation and write $C^{X}$ instead of $\Fib(C^{X}\to C^{*})$ for $C\in \bC$ and $X\in \Top_{*}$, where $*\to X$ is the base point.

\begin{ddd}\label{flkjwelfkewfjewlfjewfewfewfwefewfewf}
We define a functor $$\Funk_{c}:\bC\to \RBundle(\bC)\ , \quad \Funk_{c}(C)(W\to B):=C^{ B^{W} }\ .$$
\end{ddd}

\begin{ex}{\rm For $n\in \nat$ we consider the bundle $\R^{n}\to *$ and $\bC=\Sp$. The Thom space  of the bundle is given by 
$*^{\R^{n}}\cong S^{n}_{*}$.   Consequently, for a spectrum $E$  we get an equivalence  $$\Funk_{c}(E)(\R^{n}\to *)\simeq E^{S^{n}_{*}}\simeq \Sigma^{-n}E\ .$$
In particular we have the isomorphism
\begin{equation}\label{fjkefhwekjfhkewfewfewkfhewfw}
\pi_{-n}(\Funk_{c}(E)(\R^{n}\to *))\cong  \pi_{-n}(\Sigma^{-n}E)\cong \pi_{0}(E)\ .
\end{equation}
}\hB
\end{ex}

 \begin{lem}\label{fowfopopipoifpowiofwf24545435}
We have an equivalence \begin{equation}\label{rfhjf287489223rfefef4}
\Funk (-)d_{c}\stackrel{\simeq}{\to} \Funk_{c}(-)
\end{equation}
of functors 
$\bC\to \RBundle(\bC)$.
\end{lem}
\proof
For every subset $K\subseteq W$ which is proper over $B$ we have a natural map
\begin{equation}\label{ghjegdhjghjqwgdjqwgwquidwqiudwqdqwdqwdqwd1}B^{W}\to W/(W\setminus K)\end{equation}  between pointed topological spaces. This induces a morphism
\begin{equation}\label{ghjegdhjghjqwgdjqwgwquidwqiudwqdqwdqwdqwd}
C^{W/(W\setminus K)}\to C^{B^{W}}
\end{equation} in $\bC$.
For a second subset $ K^{\prime}\subseteq W$ such that $K^{\prime}$ proper over $B$ and contains $K$ the diagram in $\Top_{*}$
$$\xymatrix{B^{W}\ar[dr]\ar[r]&W/(W\setminus K^{\prime})\ar[d]\\
&W/(W\setminus K)}$$
commutes. Hence the diagram
 $$\xymatrix{ C^{W/(W\setminus K )}\ar[r]\ar[d]&C^{B^{W}}\\C^{W/(W\setminus K^{\prime})}\ar[ur]&}
 $$ in $\bC$
commutes. By taking the colimit over all $K\in P(W\to B)$ we get the desired morphism
\begin{equation}\label{fkejrfklerflerfoi34pr34435435}
\Funk(C)d_{c}(W\to B)\to \Funk_{c}(C)(W\to B)\ .
\end{equation} 
By an inspection of the construction we see that it comes as the specialization  at $C$ of a morphism between functors \eqref{rfhjf287489223rfefef4}.
If $K$ is a bundle of discs, then \eqref{ghjegdhjghjqwgdjqwgwquidwqiudwqdqwdqwdqwd1} is a homotopy equivalence between pointed topological spaces. Consequently,  
\eqref{ghjegdhjghjqwgdjqwgwquidwqiudwqdqwdqwdqwd} is an equivalence.
Since disc bundles form a cofinal subset of $P(W\to B)$ we conclude that 
\eqref{fkejrfklerflerfoi34pr34435435} is an equivalence. 
 \hB

\begin{kor}\label{dhkjdhqkjwqwdiuoiuoiuoiu}
We have a de Rham equivalence between $\RBundle$-objects   \begin{equation}\label{ferfriufreiof89789r89r78947r4fkjfewfewfewfe}
H(\Omega Ad_{c})\stackrel{\eqref{fjhefjwehfkwefewfew3244234324244324}}{\simeq} H(\Omega A)d_{c}\stackrel{j_{c}}{\simeq }\Funk(H(A))d_{c}\stackrel{\eqref{rfhjf287489223rfefef4}}{\simeq}  \Funk_{c}(H(A))
\end{equation}
which we also denote by $j_{c}$.
It is natural in $A$ and for a morphism \eqref{jun105111} in $\RBundle$ the diagram
\begin{equation}\label{dhkjdhqkjwqwdiuoiuoiuoiu1}
\xymatrix{H(\Omega Ad_{c})(W\to B)\ar[r]^{j_{c}}\ar[d]^{(f,F)^{*}}&\Funk_{c}(H(A))(W\to B)\ar[d]^{(f,F)^{*}}\\
H(\Omega Ad_{c})(W^{\prime}\to B^{\prime})\ar[r]^{j_{c}} &\Funk_{c}(H(A))(W^{\prime}\to B^{\prime}) 
}\end{equation}
commutes. 
  \end{kor}

\begin{rem}{\rm 
The map \eqref{fkwfkjewfewki2oiuio3h32fhiu2f23f23f32f} naturally induces a map 
 \begin{equation}\label{fkwfkjewfewki2oiuio3h32fhiu2f23f23f32f1}
E^{X}\wedge F^{Y}\to (E\wedge F)^{X\wedge  Y}
\end{equation}
for pointed spaces $X,Y$.}\hB
\end{rem}

For a vector  bundle $W\to C$ and a manifold $B$
we have a homeomorphism of pointed topological spaces \begin{equation}\label{ewfwewnmbzui87z3289r9328r32}
(B\times C)^{B\times W}\simeq B_{+}\wedge C^{W}\ .
\end{equation} 
This gives
$$ E^{B_{+}}\wedge F^{C^{W}}\stackrel{\eqref{fkwfkjewfewki2oiuio3h32fhiu2f23f23f32f1}}{\to} 
  (E\wedge F)^{B_{+}\wedge C^{W}} \stackrel{\eqref{ewfwewnmbzui87z3289r9328r32}}{\simeq} 
     (E\wedge F)^{(B\times C)^{B\times W}}\ ,$$
i.e. \begin{equation}\label{kjkefjekwlopi2iro2r23r783465847345}
\Funk(E)(B)\wedge \Funk_{c}(F)(W\to C)\to \Funk_{c}(E\wedge F)(B\times W\to C\times W)\ .
\end{equation}

\begin{lem}\label{kjdkdjqwkljdlkwjdqwdjwqjldqwdqwdwq254345}
The   de Rham equivalence \eqref{ferfriufreiof89789r89r78947r4fkjfewfewfewfe} is compatible with the exterior product. More precisely we have a commuting diagram
$$\xymatrix{H(\Omega A(B))\wedge H( \Omega A^{\prime}d_{c})(W\to C)\ar[d]\ar[r]^{j\wedge j_{c}}&
\Funk(H(A))(B)\wedge \Funk_{c}(H(A^{\prime})) (W\to C)\ar[d]^{\eqref{kjkefjekwlopi2iro2r23r783465847345}}\\
H(\Omega(A\otimes A^{\prime})d_{c})(B\times W\to B\times C)\ar[r]^{j_{c}}&\Funk_{c}(H(A\otimes A^{\prime})) (B\times W\to B\times C)}\ .$$
\end{lem}
\proof
This follows from Lemma \ref{fkjewkfjfklewjjewflejflewfjejfewlfewfwefewfewfewfe} and the almost immediate fact that \eqref{rfhjf287489223rfefef4} is compatible with the exterior product. 
\hB

\begin{ex}\label{fjeljweljewlfewiofofewufewfewfef}{\rm 
We assume that $E$ is a ring spectrum with multiplication $mult:E\wedge E\to E$.
Let $\pi:W\to B$  be an $n$-dimensional vector bundle.
An $E$-orientation of this bundle is determined by a Thom class
$U\in \pi_{-n}(\Funk_{c}(E)(W\to B))$, i.e. a properly supported cohomology class such that
$$(i_{b}^{*},I_{b})^{*}U\in  \pi_{-n}(\Funk_{c}(E)(\R^{n}\to*))\stackrel{\eqref{fjkefhwekjfhkewfewfewkfhewfw}}{\simeq} \pi_{0}(E)$$ is a unit of the ring $\pi_{0}(E)$ for every point $b\in B$. Here  the morphism
$(i_{b},I_{b}):(\R^{n}\to *)\to (W\to B)$ in $\RBundle$ is a trivialization of the fibre over $b$.
Multiplication by the Thom class induces   the first map in the composition 
\begin{equation}\label{dlkwejdelkdewdewdwedwedededd}\Funk(E)(B)\stackrel{\times U, \eqref{kjkefjekwlopi2iro2r23r783465847345}}{\to} \Funk_{c}(E\wedge E)(B\times W\to B\times B)\stackrel{(d,D)^{*}, mult}{\to}  \Funk_{c}(E)( W\to B )\ ,
\end{equation}
where the morphism $(d,D):(W\to B)\to (B\times W\to B\times B)$ in $\RBundle$ is the diagonal given by
$D(w)=(\pi(w),w)$ and $d(b):=(b,b)$.
The map in homotopy induced by \eqref{dlkwejdelkdewdewdwedwedededd} is the Thom isomorphism, which in other notation could be written as
$$Thom^{U}:H^{*}(B;E)\stackrel{\cong}{\to} H^{*+n}_{c/B}(W;E)\ .$$

}\hB
\end{ex}

\subsubsection{The collapse map}

We consider a relatively proper embedding $(f,F): (W^{\prime}\to B^{\prime})\to (W\to B)$.  Here we assume that the domain and target are vector bundles, but the map is just a map in $\Subm$.

\bigskip

\begin{ex}{\rm
The map $\emb$ occuring in Definition \ref{jul2001}, 3. is an example of  a  morphism between vector bundles  considered as objects of $\Subm$  which is a relatively proper embedding. It is not a vector bundle morphism.

}\hB
\end{ex}

\bigskip

For  a relatively proper embedding between vector bundles we define the collapse map $$\clps(F):B^{W}\to (B^{\prime})^{W^{\prime}}\ .$$
The map $F$ induces an embedding
$$W^{\prime}\cong F(W^{\prime})\subseteq W\subset B^{W}\ .$$
The collapse is defined in the image of this map as the    inverse.
It further sends the complement of the image to the base point of $(B^{\prime})^{W^{\prime}}$. 
Pull-back along the collapse yields
$$\clps(F)^{*}:\Funk_{c}(C)(W^{\prime}\to B^{\prime})\to \Funk_{c}(C)(W\to B)\ .$$

\begin{lem}\label{fjhkjkewf8978923423434324324324}
The following diagram commutes
$$\xymatrix{\Funk_{c}(C)(W^{\prime}\to B^{\prime})\ \ar[r]^{\clps(f,F)^{*}}& \Funk_{c}(C)(W\to B) \\\Funk(C)d_{c}  (W^{\prime}\to B^{\prime})\ar[u]^{\eqref{rfhjf287489223rfefef4}}_{\simeq}\ar[r]^{F_{*}}&\Funk(C)d_{c}(W\to B)\ar[u]^{\eqref{rfhjf287489223rfefef4}}_{\simeq}
}
$$
\end{lem}
\proof
We simplify the notation and write $\clps:=\clps(F)$.
We consider the diagram
$$\xymatrix{B^{W}\ar[d]\ar[r]^{\clps}&(B^{\prime})^{W^{\prime}}\ar[d]\\
(W/W\setminus F(K^{\prime}))&W^{\prime}/(W^{\prime}\setminus K^{\prime})\ar[l]^{\cong}}\ ,$$
where $K^{\prime}\subset W^{\prime}
$ is proper over $B^{\prime}$.
Note that the lower map is a homeomorphism between pointed topological spaces.
We now use that $C^{(-)}$ sends colimits to limits in order to conclude that
$$\Fib(C^{W}\to C^{W\setminus K})\simeq C^{\Cofib(W\setminus K\to W)}\simeq C^{W/(W\setminus K)}\ .$$
Hence we get a commuting diagram
$$\xymatrix{ \Funk_{c}(C)(W\to B)  &\Funk_{c}(C)(W^{\prime}\to B^{\prime}) \ar[l]^{\clps^{*}}\\ \Funk(C)d_{F(K^{\prime})}(W\to B)\ar[u]\ar[r]^{\simeq}&\Funk(C)d_{ K^{\prime} }(W^{\prime}\to B^{\prime})\ar[u]
 }\ ,$$
These diagram are natural with respect to inclusions in $P(W^{\prime}\to B^{\prime})$.
We first take the colimit over $K^{\prime}\in P(W^{\prime}\to B^{\prime})$ and then 
replace the left lower corner by $\Funk(C)d_{c}(W\to B)$. We then  get a diagram
 $$\xymatrix{ \Funk_{c}(C)(W\to B)  &\Funk_{c}(C)(W^{\prime}\to B^{\prime}) \ar[l]^{\clps^{*}}\\  \Funk(C)d_{c}(W\to B)\ar[u]\ar[r]^{F_{*} }&\Funk(C)d_{ c}(W^{\prime}\to B^{\prime})\ar[u]
 }\ .$$
We observe that the resulting lower horizontal map is $F_{*}$. 
\hB 

\begin{rem}{\rm The analog of Remark \ref{hckjjklejclejlejdldkjqwdwqdqwdqwd} applies.} \hB
\end{rem}

If we combine Lemma \ref{fjhkjkewf8978923423434324324324} with Corollary \ref{ihdwqdqwdlqwjdwqldjqwldjwlkdjwdoiuoiudoiqwudqwidwqdwqdqwdwqd}
we conclude:

\begin{kor}\label{fjhwjkf88z23430fwefwfewfewfewejkf8237458234r}
The following diagram commutes:
 \begin{equation}\label{dkjlkdqwdiuqdoqwd9878947873232hkdjds} \xymatrix{
H(\Omega Ad_{c})(W^{\prime}\to B^{\prime})\ar[r]^{j_{c}}\ar[d]^{F_{*}}&\Funk_{c}(H(A))(W^{\prime}\to B^{\prime})\ar[d]^{\clps(F)^{*}}\\
H(\Omega Ad_{c})(W \to B )\ar[r]^{j_{c}}&\Funk_{c}(H(A))(W \to B )}
\end{equation}
\end{kor}

\subsubsection{The suspension equivalence}

For $k\in \nat$ we consider 
  vector bundle $ \R^{k}\to *$. Then we have a homeomorphism of pointed spaces $*^{\R^{k}}\simeq S^{k}_{*}$. 
  For a spectrum $E\in \Sp$ this homeomorphism induces an equivalence
$$\Funk_{c}(E)(\R^{k}\to *)\simeq  E^{S^{k}_{*}}   \simeq  \Sigma^{-k}E$$
in $\Sp$.
Setting $E:=\Sigma^{k}S$, the $k$-fold shift of the sphere spectrum {$S$},  then we get  the  equivalence  
$$susp_{0}:S\simeq  \Funk_{c}(\Sigma^{k}S)(\R^{k}\to *)$$ in $\Sp$.

In the following definition we use the identification   $E\wedge \Sigma^{k}S\simeq \Sigma^{k}E$.

 \begin{ddd}\label{fewklfwejkfweoiuo23u42398ruekjdnqwefewfwew4r}
 We define the suspension equivalence 
 $$susp:\Funk ( E)(B)\simeq \Funk_{c}( \Sigma^{k}E)(B\times \R^{k}\to B)$$
 as the composition
 \begin{eqnarray*}\Funk ( E)(B)&\simeq& \Funk ( E)(B)\wedge \Funk(S)(*)\\&\stackrel{\id\wedge susp_{0}}{\to} & \Funk ( E)(B)\wedge 
  \Funk_{c}(\Sigma^{k}S)(\R^{k}\to *)\\& \stackrel{\eqref{kjkefjekwlopi2iro2r23r783465847345}}{\to} &
  \Funk_{c}( \Sigma^{k}E)(B\times \R^{k}\to B)\ .\end{eqnarray*}
   \end{ddd}
   The suspension is 
    an equivalence 
 between functors
 $$\Sp\times \Nerve(\Mf^{op})\to \Sp\ .$$
 Its inverse will also be called desuspension and denoted by $desusp$.
 
 \bigskip

%
%
%
%
%

We fix a form $U\in \Omega^{k}_{c}(\R^{k})$ such that \begin{equation}\label{jfklewjflkwefjklewffewfe7498234234242}
\int_{\R^{k}} U=1
\end{equation}
We consider $U$ as a map
$$U :S\simeq\Funk(S)(*)\to H(\Omega[k]d_{c})(\R^{k}\to *)\ .$$

\begin{prop}\label{kldjhqwkdjqdwkljqwdiouqwodqwidqwdwqdwqdwqd}
We have a canonical commuting square
\begin{equation}\label{wefjkhewfjke89724332}
\xymatrix{H(\Omega A )(B)\ar[r]^{j}\ar[d]^{\times U}&\Funk(H(A))(B)\ar[d]^{susp}\\
H(\Omega A[k]d_{c})(B\times \R^{k}\to B)\ar[r]^{j_{c}}&\Funk_{c}(H(A[k]))(B\times \R^{k}\to B)}
\end{equation} \end{prop}\proof
 We expand the diagram as follows \begin{equation}\label{fewkjhjfkfefjewflkewfjewfewf4343534535}
\xymatrix{H(\Omega A )(B)\ar[r]^{j} \ar[d]&\Funk(H(A))(B)\ar[d] \\
H(\Omega A )(B)\wedge  \Funk(S)(*)\ar[r]^{j\wedge \id}\ar[d]^{\id\wedge  j_{c}(U) }\ar@/_4cm/[dd]_{\id\wedge U}&\Funk(H(A))(B)\wedge \Funk(S)(*)\ar[d]^{\id\wedge susp_{0}}\\
H(\Omega A )(B)\wedge H(\Funk_{c}(H\R[k] ))(\R^{k}\to *)    &\Funk(H(A))(B)\wedge \Funk_{c}(\Sigma^{k}S)(\R^{k}\to *)\ar[d]^{\id\wedge \epsilon} \ar[l]^{j^{-1}\wedge\epsilon}  \\ 
H(\Omega A )(B)\wedge H( \Omega [k]d_{c})(\R^{k}\to *)  \ar[u]^{\id\wedge j_{c}} \ar[r]^{j\wedge j_{c}}\ar[d]&  \Funk(H(A))(B)\wedge \Funk_{c}(H\R[k])(\R^{k}\to *)\ar[d]^{\eqref{kjkefjekwlopi2iro2r23r783465847345}}     \\
H(\Omega A [k]d_{c} )(B\times \R^{k}\to B)\ar[r]^{j_{c}}& \Funk_{c}( H(A[k])) (B\times \R^{k}\to B)
}
\end{equation}
The lower square commutes by Lemma \ref{kjdkdjqwkljdlkwjdqwdjwqjldqwdqwdwq254345}.
The first and third cell (counting from the top) obviously commute. 
In order to get commutativity of the second cell it remains to  provide the equivalence
$$j_{c}  (U)\simeq   \epsilon\circ susp_{0} :S\simeq   \Funk(S)(*) \to \Funk_{c}( H\R[k])(\R^{k}\to *) \simeq H\R\ ,$$ {where $\epsilon$ is induced by the shift of the unit $S\to H\R$ of the ring spectrum $H\R$.}
To this end we observe that the map
$$\pi_{0}:\Map(S,H\R)\to \pi_{0}(H\R)_{\infty}\cong \R^{\delta}_{\infty}$$ is an  equivalence  of spaces, where $\R^{\delta}$  is $\R$ with the discrete topology. 
The condition \eqref{jfklewjflkwefjklewffewfe7498234234242} ensures that
both maps belong to  the same component. Hence they are equivalent in an essentially unique way. \hB

\begin{prop}\label{ekjfwklejflfelwfjewfewfew42243242342342342342}
We have a natural commuting square
\begin{equation}\label{djkhdhqwkjhqwkjdkjdwqdwqdwqdwqdw}
\xymatrix{ H(\Omega A[k]d_{c})(B\times \R^{k}\to B)\ar[d]^{\int^{\prime}_{rel}}\ar[r]^{j_{c}}&\Funk_{c}(H(A[k]))(B\times \R^{k}\ar[d]^{desusp}\to B)\\H(\Omega A )(B)\ar[r]^{j} &\Funk(H(A))(B) 
}\end{equation}
 which is inverse to \eqref{wefjkhewfjke89724332}.
\end{prop}
\proof
We consider the following diagram
{\scriptsize
$$\hspace{-2cm}\xymatrix{\ar@/_2cm/[dddd]^(.3){\id}H(\Omega A)(B)\ar[dr]^{\id\wedge U}\ar[rr]^{j}\ar[dd]^{\times U}&&\Funk(H(A))\ar[dd]^(.3){susp}\ar[dr]^{\epsilon\circ susp_{0}}\ar@/^9cm/[dddd]_{\id}\\
&H(\Omega A)(B)\wedge H(\Omega[k]d_{c})(\R^{k})\ar[rr]^(0.3){j\wedge j_{c}}\ar[dl]\ar[dd]^(.3){\id\wedge \int}&&\Funk(H(A))(B)\wedge \Funk_{c}(H\R[k])(\R^{k})\ar[dl] \\H(\Omega Ad_{c})(B\times \R^{k})\ar[rr]^(.6){j_{c}}\ar[dd]^{{\int^{\prime}_{rel}}}&&\Funk_{c}(H(A[k]))(B\times \R^{k})\ar[dd]^(.3){desusp}&\\
&H(\Omega A)(B)\wedge H\R\ar[rr]^(.3){j\wedge \id}\ar[dl]&&\Funk(H(A))(B)\wedge H\R \ar[uu]^{\id\wedge (susp_{0}\wedge H\R)}\\  H(\Omega A)(B)\ar[rr]^{j}&& \Funk(H(A))(B)  \ar[ur]^{\epsilon}   & }\ .$$}
The square in question \eqref{djkhdhqwkjhqwkjdkjdwqdwqdwqdwqdw} appears in the lower middle.
All other cells canonically commute by the preceeding discussion.  For example the cells on the left commute on the form level. The cells on the right commute by the definition of $desusp$ as the inverse of $susp$.
The remaining middle squares appear in \eqref{fewkjhjfkfefjewflkewfjewfewf4343534535}.
We therefore define the square 
\eqref{djkhdhqwkjhqwkjdkjdwqdwqdwqdwqdw} such that this big diagram commutes.
The outer square expresses the assertion that \eqref{djkhdhqwkjhqwkjdkjdwqdwqdwqdwqdw}  is inverse to 
\eqref{wefjkhewfjke89724332}. \hB

An inspection of the proofs shows:

 \begin{kor}
 The squares \eqref{wefjkhewfjke89724332} and \eqref{djkhdhqwkjhqwkjdkjdwqdwqdwqdwqdw} depend  functorially on $B$ and $A$. More precisely we have constructed them as   functors
 $$\Nerve(\Ch_{\R})\times \Nerve(\Mf)^{op}\to \Sp^{\Nerve(\Box)}\ ,$$
 where $\Box$ is as in \eqref{ewfwfewfwj2lkr23iori32rio2ri322r2r2}. 
Moreover, they are inverse to each other in this functorial sense. 
 \end{kor}

\subsection{Thom and Euler forms}\label{may2910}

{
In this subsection, which is a review of \cite[Ch. 1.6]{MR1215720}, we recall the constructions of the Euler and Thom forms. 
 }
 
 \bigskip

If $\pi:V\to M$ is an $n$-dimensional real  vector bundle, then we can form the orientation bundle $\Lambda\to M$. It is a one-dimensional real vector bundle with structure group reduced to the subgroup $\{\pm 1\}\subset GL(1,\R)$. If $\Fr(V)\to  M$ denotes the $GL(n,R)$-principal bundle of frames  of $V$, then we can write the orientation bundle as an associated vector bundle
\begin{equation}\label{jul1601}\Lambda:=\Fr(V)\times_{GL(n,R),\sign(\det)}\R\ ,\end{equation}
where $\sign(\det):GL(n,\R)\to \{\pm 1\}$
is the ``sign of the determinant'' representation. 
The main reason for introducing the orientation bundle is that we have an integration map  {
\begin{equation}\label{jul1602}\int_{V/M}:\Omega_{c/M}(V,\pi^{*}\Lambda)[n]\to \Omega(M) \end{equation} without any condition of orientability. Here $\Omega_{c/M}(V,\pi^{*}\Lambda)[n]$ denotes the de Rham complex of differential forms on $V$ with coefficients in $\Lambda$ which are properly supported over $M$, and which is shifted such that $n$-forms are in degree $0$.}

\begin{rem}{\rm 
In the notation introduced in Definition \ref{dijqwidjqwldjqwlqwdqwdqwdwqdwqdqwdqw} the complex  $\Omega_{c/M}(V)$ would have been denoted by $\Omega d_{c}(V\to M)$, but we find this notation inappropriate at this point.}\hB
\end{rem}

The bundle $\Lambda$ has a natural flat connection, and by $\ulambda$ we denote 
 the associated locally constant sheaf of parallel sections. The Thom class of $V$ is {a}  cohomology class
$$U(V)\in H_{c/M}^{\dim(V)}(V,\pi^{*}\ulambda)\ ,$$ where the subscript $c/M$ indicates  cohomology with proper support over $M$ as explained in Example \ref{fjkhefjkwefkjewfewfewfewfwefwef}, and we have an isomorphism
$$H_{c/M}^{\dim(V)}(V,\pi^{*}\ulambda)\cong H^{\dim(V)}(\Omega_{c/M}(V,\pi^{*}\Lambda))\ .$$
The Thom class $U(V)$ is
uniquely determined by the condition that $\int_{V/M} U(V)=1$,
where  the integration in cohomology is induced by the  map (\ref{jul1602}) {on} form level.

\begin{rem}{\rm
This can be compared to Example \ref{fjeljweljewlfewiofofewufewfewfef}. Since the present definition of a Thom class
involves the twist $\pi^{*}\underline{\Lambda}$ the Thom class  considered here induces a twisted $H\R$-orientation of $V$. It provides an orientation in the sense of Example \ref{fjeljweljewlfewiofofewufewfewfef} if the twist is trivialized.
 }\hB
\end{rem}

The restriction of the Thom class to the zero section
 $$e(V):=0_{V}^{*}U(V)\in H^{\dim(V)}(M,\ulambda)$$   is,  by definition, the Euler class of $V$.

\bigskip

\newcommand{\so}{\mathbf{so}}
\newcommand{\vol}{\mathrm{vol}}
We need a refinement of the Thom class on the level of differential forms for real  vector bundles with metric and metric connection $(V,h^{V},\nabla^{V})$.  
 We let $$MQ(\nabla^{V})\in \Omega^{\dim(V)}(V,\pi^{*}\Lambda)$$ be the Mathai-Quillen form  {(see
\cite[Ch. 1.6, (1.37)]{MR1215720} but note the different normalization).} {For the  sake of  completeness  we briefly recall its construction.
We consider the bigraded algebra $A$ with constituents $A^{i,j}:=\Omega^{i}(V,\pi^{*}\Lambda^{j}V)$.
We have the tautological section $x\in A^{0,1}$. Using the metric
on $V$ we define  the {section} $|x|^{2}\in A^{0,0}$.
We further have $\nabla^{\pi^{*}V}x\in A^{1,1}$. The curvature
$R^{\nabla^{V}}$ can be considered as an element in $\Omega^{2}(M,\so(V))$.
Using the metric again we identify $\so(V)\cong \Lambda^{2}V$. With this identification we get
$F:=\pi^{*} R^{\nabla^{V}}\in A^{2,2}$.
The metric on $V$ induces a canonical element $\vol\in C^{\infty}(M,\Lambda^{\dim(V)} V^{*}\otimes \Lambda)$. The Berezin  integral $T:A\to \Omega^{*}(V,\Lambda)$ is given by the composition of projection from $A$ to $A^{*,\dim(V)}$ and the pairing with $\vol$.  In terms of these constructs
the Mathai-Quillen form is given by 
$$MQ(\nabla^{V}):=\frac{1}{(2\pi)^{\dim(V)/2}}  T \exp\left(-\frac{|x|^{2}}{2} +i\nabla^{\pi^{*}V} x +F\right)\ .$$
}

{The Mathai-Quillen form} does not have proper support over $M$, but it decays rapidly at infinity. We modify this form in order to obtain a form with proper support over $M$.
We define the diffeomorphism $f:B_{1}(V)\to V$ of the unit ball bundle $B_{1}(V)\subset V$  with $V$  by 
$$f(v):=\frac{v}{(1-|v|^{2})^{1/2}}\ , \quad v\in B_{1}(V)\ .$$   
The form
$f^{*}MQ(\nabla^{V})$ is first defined on $B_{1}(V)$ but can be extended by zero smoothly to all of $V$.
We define\begin{equation}\label{hfjkwehfejkwfhkewfoiufoiu}
 U(\nabla^{V}):=f^{*}MQ(\nabla^{V})\in \Omega_{c/M}^{\dim(V)}(V,\pi^{*}\Lambda)\ .
\end{equation}The Euler form $e(\nabla^{V})\in \Omega^{\dim(V)}(M,\Lambda)$ of the real euclidean vector bundle $(V,h^{V})$ with connection $\nabla^{V}$ is defined  in terms of the curvature $R^{\nabla^{V}}\in \Omega^{2}(W,\so(V))$  by
\begin{equation}\label{jul1603}e(\nabla^{V}):= {\frac{1}{(2\pi)^{\dim(V)/2}} T\exp(F)}\ ,\end{equation}
see \cite[Sec. 1.6]{MR1215720} for details.
In the following we list some obvious properties of the Thom form.
\begin{enumerate}
\item
 We have $$[U(\nabla^{V})]=U(V)\in H_{c/M}^{\dim(V)}(V,{\pi^{*}}\ulambda)\ .$$
 \item For a {pull-back diagram}
$$\xymatrix{V^{\prime}\ar[r]^{F}\ar[d]^{{\pi^{\prime}}}&V\ar[d]^{{\pi}}\\M^{\prime}\ar[r]^{{f}}&M}$$
we have an isomorphism $\Lambda^{\prime}\cong {f}^{*}\Lambda$ and (under this isomorphism)
$$F^{*}U(\nabla^{V})=U(\nabla^{V^{\prime}})\in \Omega_{c/M^{\prime}}^{\dim(V)}(V^{\prime},{\pi^{\prime*}}\Lambda^{\prime})\ .$$
Here we assume that $V^\prime$ has the induced metric and metric connection $\nabla^{V^{\prime}}=f^{*}\nabla^{V}$.
\item We have 
$$0_{V}^{*}(U(\nabla^V))=e(\nabla^{V})\in \Omega^{\dim(V)}(M,\Lambda)\ ,$$
where $e(\nabla^{V})$ denotes the { Euler form \eqref{jul1603}}.
\item For  two bundles with metric and connection $(V,h^{V},\nabla^{V})$ and $(V^{\prime},h^{V^{\prime}},\nabla^{V^{\prime}})$ we have
$$0_{V}^{*}U(\nabla^{V}\oplus \nabla^{V^{\prime}})=\pi^{\prime,*}e(\nabla^{V})\wedge U(\nabla^{V^{\prime}})\ .$$
\end{enumerate}

\subsection{The normalized Borel regulator map}\label{jun274333}

The goal of the present subsection is to understand the relation  between the Kamber-Tondeur forms introduced in (\ref{jun102})  and the normalized Borel regulator map.  
This becomes important if one wants to connect to developments in arithmetic geometry. The results of this subsection are only used in Subsubsection \ref{jul0110}.
 
 \bigskip

  We  follow \cite{MR1869655} for the description of the normalized Borel regulator map.
For $n\ge 1$ we consider $G_{n}:=SL(n,\C)$ as a semisimple real Lie group. Its maximal compact subgroup is $K_{n}=SU(n)$, and we let $X_{n}:=G_{n}/K_{n}$ denote the associated symmetric space.
The complexification of $G_{n}$ can be identified with
$G_{n,\C}:=SL(n,\C)\times SL(n,\C)$ with the maximal compact subgroup $G_{u,n}:=SU(n)\times SU(n)$. The embedding $G_{n}\to G_{n,\C}$ is given by 
$$SL(n,\C)\ni A\mapsto (A,\bar A)\in SL(n,\C)\times SL(n,\C)\ .$$ 
The space $X_{u,n}:=G_{u,n}/K_{n}$ is the compact dual of $X_{n}$. 
We have
$X_{u,n}=SU(n)\times SU(n)/SU(n)$ with the action $(A,B)U=(AU,B\bar U)$ and therefore
a diffeomorphism $X_{u,n}\stackrel{\sim}{\to} SU(n)$ given by $(A,B)\mapsto AB^{t}$. 
The induced left action of $SU(n)\times SU(n)$ on $SU(n)$ is given by $(A,B)X=AXB^{t}$.

The rational cohomology of the $h$-space $SU(n)$ is the exterior algebra
$$H^{*}(SU(n);\Q)\cong \Lambda(\beta(n)_{3}^{\Q},\dots,\beta(n)_{2n-1}^{\Q})$$ generated by primitive
elements $\beta(n)_{2j-1}^{\Q}\in H^{2j-1}(SU(n);\Q)$, where $x\in  H^{*}(SU(n);\Q)$ is called primitive if
$\mu^{*}(x)=x\times 1+1\times x$ for the multiplication map $\mu:SU(n)\times SU(n)\to SU(n)$.

We consider $\C$ as a real vector space with an action of $\Z/2\Z$ by complex conjugation and form the invariant subspace $\R(j):=(2\pi i)^{j}\R\subset \C$. We  define normalized elements
$$\beta(n)_{2j-1}:= c(n)_{2j-1}\beta(n)^{\Q}_{2j-1}\in H^{2j-1}(SU(n);\R(j))\ ,\quad  c(n)_{2j-1}\in \C$$ such that
$$\beta(n)_{2j-1}(\pi_{2j-1}(SU(n))=\Z(j):=(2\pi i)^{j}\Z\ .$$ This fixes these elements up to sign. We choose the signs such that for $3\le j\le m\le n$ we have the compatibility
$(\beta(n)_{2j-1})_{|U(m)}=\beta(m)_{2j-1}$. An explicit choice of the {classes} $\beta(n)_{2j-1}$ will be described in (\ref{jul0210}).
{In order to distingish the cohomology of a Lie group as a space from group cohomology we will denote the group-cohomology by $H^{*}_{gr}$. We write $H^{*}_{cont}$ for the continuous group cohomology.}  There is a canonical morphism
\begin{equation}\label{grgegergregrgre34t534334657657}
\phi_{n}: H^{2j-1}(SU(n);\R(j))\to H_{gr}^{2j-1}(SL(n,\C);\R(j-1))
\end{equation}
defined as the composition of the following maps:
\begin{enumerate}
\item The de Rham isomorphism and Hodge theory identify the cohomology of a compact Lie group with  its  biinvariant differential forms. Hence we have 
\begin{equation}\label{jul0206}H^{{2j-1}}(SU(n);\R(j))=(2\pi i)^{j}H^{{2j-1}}(SU(n);\R) \cong (2\pi i)^{j}\Omega^{2j-1}(SU(n))^{SU(n)\times SU(n)}\ ,\end{equation}
{viewed as subspaces in the corresponding complexifications.}
\item  We have an isomorphism
\begin{equation}\label{jul0207}(2\pi i)^{j} \Omega^{2j-1}(SU(n))^{SU(n)\times SU(n)}=(2\pi i)^{j} \Omega^{2j-1}(X_{u,n})^{G_{u,n}}\cong (2\pi i)^{j} H^{2j-1}(\frg_{n,u},K_{n};\R)\ ,\end{equation}
where the last group is relative Lie algebra cohomology. We use the general convention that small letters denote the Lie algebras of the groups denoted by the corresponding capital symbols.
\item We have Cartan decompositions of the Lie algebras $\frg_{u,n},\frg_{n}\subset \frg_{n,\C}$
$$ \frg_{u,n}= \frk_{n}\oplus  \frp_{u,n}\ , \quad \frg_{n}= \frk_{n}\oplus  \frp_{n}\ ,$$
and the relation $ \frp_{u,n}=i  \frp_{n}$. Note that
\begin{equation}\label{jul0208}H^{2j-1}( \frg_{u,n},K_{n};\R)\cong (\Lambda^{2j-1} \frp_{u,n})^{K_{n}} ,\quad  H^{2j-1}( \frg_{n},K_{n},\R)\cong (\Lambda^{2j-1}  \frp_{n})^{K_{n}}\ .\end{equation}
By complex linear extension and restriction of forms
we now have an isomorphism
\begin{equation}\label{jul12}(2\pi i)^{j} H^{2j-1}( \frg_{u,n},K_{n};\R)\cong  (2\pi i)^{j-1} H^{2j-1}( \frg_{n},K_{n};\R)\ .  \end{equation}  
\item 
We have the  van Est isomorphism 
\begin{equation}\label{jul0203} H^{2j-1}( \frg_{n},K_{n};\R)\stackrel{\sim}{\to} H_{cont}^{2j-1}(G_{n};\R)= H_{cont}^{2j-1}(SL(n,\C);\R)\ .\end{equation}
\item  We have a morphism 
\begin{equation}\label{jul0204}(2\pi i)^{j-1}H_{cont}^{2j-1}(SL(n,\C);\R)\to (2\pi i)^{j-1}H_{gr}^{2j-1}(SL(n,\C);\R)= H_{gr}^{2j-1}(SL(n,\C);\R(j-1))\ .\end{equation}
\end{enumerate}

\begin{ddd}For $j\ge 2$
a class $c^{Bor}_{2j-1}\in H^{2j-1}_{gr}(GL(\C);\R(j-1))$ is called a Borel element 
if its restriction to $SL(\C,n)$ coincides with
$\phi_{n}(\beta(n)_{2j-1})$ for some (hence all) $n\ge j$.
\end{ddd}
This definition fixes the Borel elements up to decomposeables.



%

%

%

For every complex embedding $\sigma\in \Sigma$ we get
a map
$GL(R)\to GL(\C)$ which induces a map
\begin{equation}\label{jul0201} H^{2j-1}(GL(\C);\R(j-1))\to H^{2j-1}(GL(R);\R(j-1)) \to \Hom(K_{2j-1}(R);\R(j-1))\ ,\end{equation}
where the last map is induced by
the Hurewich map, the fact that  Quillen's $+$-constrution produces a homology isomorphism
$$K_{*}(R)\xrightarrow{Hurewich} H_{*}(BGL(R)^{+};\Z)\cong H_{*}(BGL(R);\Z)\ ,$$
and the evaluation of cohomology against homology. 
We let
$$\tilde c^{Bor}_{2j-1}(\sigma):K_{2j-1}(R)\to \R(j-1)$$ denote the image under (\ref{jul0201}) of a Borel element
$c^{Bor}_{2j-1}$. This homomorphism $\tilde c^{Bor}_{2j-1}(\sigma)$ does not depend on the choice of the Borel element since decomposeable classes evaluate trivially against classes in the image of the Hurewich map.

  We have a $\Z/2\Z$-invariant decomposition
$$\C\cong \R(j)\oplus \R(j-1)$$ and therefore an isomorphism of real vector spaces $\C/\R(j)\cong  \R(j-1)$. 
We define the real vector space {of $\Z/2\Z$-invariants}
\begin{equation}\label{jul0230}X_{2j-1}(R):=((\C/\R(j)^{\Sigma})^{\Z/2\Z}\end{equation} where the group $\Z/2\Z$ acts by complex conjugation on both,  { $\C/\R(j)$ and $\Sigma$}.
Then we define the Borel regulator map
\begin{equation}\label{jul0232}r_{Bor} :K_{2j-1}(R)\to X_{2j-1}(R)\ ,\quad r_{Bor}(x)(\sigma):=\tilde c^{Bor}_{2j-1}(\sigma)(x)\ .\end{equation}
 Theorem \ref{borel} of Borel \cite{MR0387496} {and its refinement Proposition \ref{jun261}} can be reformulated to say that the Borel regulator map
$r_{Bor}$ induces an isomorphism 
$$K_{2j-1}(R)\otimes \R\xrightarrow{\cong} X_{2j-1}(R)$$ for $j\ge 2$.

 
 Recall the Definition \ref{kjqwldqwdqwdwqdqwdqd} of the graded vector space $A$.
 We have the following commutative diagram
 \begin{equation}\label{jul0220}\xymatrix{&&X_{2j-1}(R)\ar[dd]^{\psi}\\K_{2j-1}(R)\ar[urr]^{r_{Bor}}\ar[drr]^{c}&&\\&&A_{2j-1}}
\end{equation} with the isomorphism of real vector spaces $\psi$.
We will need the explicit form of the isomorphism
$\psi$ in the discussion \ref{jul0110}.   
We have  
$$\psi(f)= \frac{1}{(2\pi i)^{j-1}} \sum_{\sigma,\sigma^{\prime}\in \Sigma }  \psi_{\sigma,\sigma^{\prime}} f(\sigma^{\prime}) \delta_{2j-1,\sigma}$$
for some real matrix $(\psi_{\sigma,\sigma^{\prime}})_{\sigma,\sigma^{\prime}\in \Sigma }$.
In order to fix the matrix uniquely we will assume the relations  \begin{equation}\label{ffehewjkfhkewf89u2or23r332r32r}
\psi_{\sigma,\sigma}=\psi_{\bar \sigma,\bar\sigma}\ , \quad \psi_{\sigma,\bar\sigma}=0
\end{equation} for all complex places $\sigma$.

\begin{prop}\label{jul0231}
{We assume that $j\ge 2$. Then
we} have $$\psi_{\sigma,\sigma^{\prime}}= \frac{(-1)^{j-1}(2j-1)!}{(j-1)!  } \delta_{\sigma,\sigma^{\prime}}\ .$$
\end{prop}
\proof
By \eqref{kdqwkdjlwqkdjwqkdjq987}  we have 
$$\psi(r_{Bor}(x))=\sum_{\sigma\in \Sigma } \omega_{2j-1}(\sigma)(x) \delta_{2j-1,\sigma}
= \frac{1}{(2\pi i)^{j-1}}\sum_{\sigma,\sigma^{\prime}\in \Sigma }   \psi_{\sigma,\sigma^{\prime}} \tilde c^{Bor}_{2j-1}(\sigma^{\prime})(x) \delta_{2j-1,\sigma}\ ,$$
and 
hence the relation (see \eqref{dkqwdjqwkldlwqkwdwqdooiuo123213213} for $\omega_{2j-1}(\sigma)$)
$$ \omega_{2j-1}(\sigma) =  \frac{1}{(2\pi i)^{j-1}}\sum_{\sigma,\sigma^{\prime}\in \Sigma }   \psi_{\sigma,\sigma^{\prime}} \ \tilde c^{Bor}_{2j-1}(\sigma^{\prime})\ ,\quad \forall \sigma\in \Sigma$$ together with \eqref{ffehewjkfhkewf89u2or23r332r32r}
determines the matrix $(\psi_{\sigma,\sigma^{\prime}})_{\sigma,\sigma^{\prime}\in \Sigma }$. 

In order to compare $\omega_{2j-1}(\sigma)$ and $\tilde c^{Bor}_{2j-1}(\sigma)$ we give a description of the latter which is similar to the former. 

\bigskip

We fix $n\ge j$.
The standard hermitean metric $h_{0}$ on $\C^{n}$ induces a { normalization of the volume form}. We let $\Met_{1}(n)$ be the space of hermitean metrics on $\C^{n}$ {inducing the normalized volume form}. The group
$SL(n;\C)$ acts on $\C^{n}$ and therefore on $\Met_{1}(n)$. The stabilizer of $h_{0}$ is
$SU(n)$ so that we get an identification $X_{n}\cong \Met_{1}(n)$. The trivial vector bundle
$$V:=\Met_{1}(n)\times \C^{n}\to \Met_{1}(n)$$ has a trivial connection $\nabla^{V}$ and a tautological metric
$h_{taut}^{V}$. Both structures are    $SL(n,\C)$-invariant. Hence the associated  Kamber-Tondeur form \eqref{jun102} is an invariant form  $$\omega_{2j-1}(h^{V}_{taut})\in \Omega^{2j-1}(\Met_{1}(n))^{SL(n,\C)}\ .$$ We have an isomorphism
\begin{equation}\label{jul0205}\Omega^{2j-1}(\Met_{1}(n))^{SL(n,\C)}\cong \Omega^{2j-1}(X_{n})^{G_{n}}\cong (\Lambda^{2j-1}\frp_{n})^{K_{n}}\ .\end{equation}
We let
$\tilde \omega_{2j-1}\in H_{gr}^{2j-1}(SL(n,\C);\R)$ be the image of the form
$\omega_{2j-1}(h^{V}_{taut})$ under the composition of (\ref{jul0205}) with (\ref{jul0203}) and (\ref{jul0204}). {It is clear that the following Lemma implies Proposition \ref{jul0231}.} Recall the description of $\phi_{n}$,  \eqref{grgegergregrgre34t534334657657}. 
\begin{lem}\label{mar1001}
In $H_{gr}^{2j-1}(SL(n,\C);\R)$ we have
$$\tilde \omega_{2j-1}= \frac{(-1)^{j-1}(j-1)!}{(2j-1)!(2\pi i)^{j-1}} \phi_{n}(\beta(n)_{2j-1})\ .$$
\end{lem}
 \proof
 Let $\hat \beta_{2j-1}\in (\Lambda^{2j-1}\frp_{n})^{K_{n}}$ be the image of $\beta(n)_{2j-1}$ under
 (\ref{jul0206}), (\ref{jul0207}), (\ref{jul12}),  and (\ref{jul0208}).
 It suffices to compare $\hat \beta_{2j-1}$ with the image $\hat \omega_{2j-1}$ of $\omega_{2j-1}(h^{V}_{taut})$ under
 (\ref{jul0205}). 
 The Cartan decomposition $\frg_{n}=\frk_{n}\oplus \frp_{n}$ is the decomposition of complex
 trace-free $n\times n$-matrices into the anti-hermitean  and hermitean parts.
 We thus identify 
 $$\frp_{n}\cong \{A\in \Mat(n,\C)\:|\:A=A^{*}\ ,\tr A=0\}\ .$$
 We use the composition of the exponential map and the projection 
 $\frp_{n}\to G_{n}\to X_{n}$ as a chart. In these coordinates the metric $h_{taut}^{V}$ is a symmetric matrix valued function $$h_{taut}^{V}(A)=\exp(A)\exp(A)^{*}=1+A+A^{*}+\dots=1+2A+\dots\ .$$ We have
 $\nabla^{V}=d$ and $\nabla^{V,*}-\nabla^{V}=(h_{taut}^{V})^{-1} dh^{V}_{taut}$. At the origin $A=0$ we get $$
 (\nabla^{V,*}-\nabla^{V})(0)(X)= 2X\ ,\quad  X\in \frp_{n}\ .$$ 
{We insert this into  
  \eqref{jun102} and get}
 $$\hat \omega_{2j-1}(X_{1},\dots, X_{2j-1})=\frac{1}{(2\pi i)^{j-1}}  \sum_{s\in \Sigma^{2j-1}} \sign(s) \Tr(X_{s(1)}\dots X_{s(2j-1)})\ .$$
 We now calculate $\hat \beta_{2j-1}$ explicitly.
 We first give a cohomological description of a primitive element in $H^{2j-1}(SU(n);\Q)$. We identify $S^1\cong \R/\Z$.
 On $S^{1}\times SU(n)$ we consider the suspension bundle $V_{susp}\to S^{1}\times SU(n)$ given as a quotient of the $\Z$-equivariant trivial bundle
 $$ \R\times SU(n)\times \C^{n}  \to  \R\times SU(n) \ ,\quad (t,g,v)\mapsto (t,g) \ ,$$ where the actions of
 $\Z$ are given by 
 $$n(t,g,v):=(t+n,g,g^{n}v)\ , \quad n(t,g):=(t+n,g)\ .$$
Then
$$\ch_{2j-1}:=\int_{S^{1}\times SU(n)/SU(n)} \ch_{2j}(V_{susp})\in H^{2j-1}(SU(n);\Q)$$
is primitive. The normalization of the Chern character is such that
$$\ch_{2j-1}(\pi_{2j-1}(SU(n)))=\Z\ .$$
We conclude that we can take
\begin{equation}\label{jul0210}\beta(n)_{2j-1}:= (2\pi i)^{j} \ch_{2j-1}\ .\end{equation}

We now calculate a differential form representative.
We choose a function $\chi\in C^{\infty}([0,1])$ which is constant near the end points of the interval and satisfies 
$\chi(0)=0$ and $\chi(1)=1$.
We define a connection on the trivial bundle $[0,1]\times SU(n)\times \C^{n}\to [0,1]\times SU(n)$ by
$$\tilde \nabla:=d+{\chi(t) g^{-1} dg} \ .$$  It can be extended to a $\Z$-invariant connection on 
$\R\times SU(n)\times \C^{n}\to \R\times SU(n)$. This $\Z$-invariant connection induces a connection $\nabla^{V_{susp}}$ on
$V_{susp}$ which we use to calculate a form representative 
$$\ch_{2j}(\nabla^{V_{susp}})=\frac{(-1)^{j}}{j!(2\pi i)^{j}} \Tr (R^{\nabla^{V_{susp}}})^{j}\ .$$
We have on $[0,1]\times SU(n)$ that
$$R^{\nabla^{V_{susp}}}={\chi^{\prime}(t) dt\wedge g^{-1} dg + (\chi(t)^{2}  - \chi(t))g^{-1}dg \wedge g^{-1}dg }\ .$$
The form $\ch_{2j}(\nabla^{V_{susp}})$ is thus bi-invariant under the action of $SU(n)\times SU(n)$ on $S^{1}\times SU(n)$. For further calculation we therefore restrict to the origin $g=1$, where
$$R^{\nabla^{V_{susp}}}(t,1)=\chi^{\prime}dt\wedge  dg + (\chi(t)^{2}-\chi(t))dg \wedge dg \ .$$
We have  for $t\in [0,1]$ that
$$\Tr\  R^{\nabla^{V_{susp}}}(t,1)^{2j}= j\chi^{\prime}(t) (\chi(t)^{2}-\chi(t))^{j-1} dt \wedge \Tr((dg)^{2j-1})+\dots\ ,$$
where the terms not written do not contain $dt$.
Using \begin{eqnarray*}
\int_{0}^{1} (x^{2}-x)^{j-1} dx&=&(-1)^{j-1}\int_{0}^{1} x^{j-1}(1-x)^{j-1}dx\\&=&(-1)^{j-1} B(j,j)\\&=&\frac{(-1)^{j-1} \Gamma(j)^{2}}{\Gamma(2j)}\\&=&\frac{(-1)^{j-1}((j-1)!)^{2}}{(2j-1)!}\end{eqnarray*}
we get
 \begin{eqnarray*}\int_{S^{1}\times SU(n)} \ch_{2j}(\nabla^{V_{susp}})&=& \frac{j}{j!(2\pi i)^{j}} \Tr (dg)^{2j-1}  \int_{0}^{1} \chi^{\prime}(t) (\chi(t)^{2}-\chi(t))^{j-1}\\&=&\frac{1}{(2\pi i)^{j}}\frac{(-1)^{j-1} (j-1)!}{(2j-1)!} \Tr (dg)^{2j-1} \ .\end{eqnarray*}
This gives 
$$\hat \beta_{2j-1}(X_{1},\dots,X_{2j-1})= \frac{(-1)^{j-1}(j-1)!}{(2j-1)!}\sum_{s\in \Sigma^{2j-1}} \sign(s) \Tr (X_{s(1)}\dots X_{s(2j-1)})\ .$$
This implies that
$$ \frac{(-1)^{j-1}(2j-1)!}{(j-1)! (2\pi i)^{j-1}}\hat \beta_{2j-1}=\hat \omega_{2j-1}\ .$$
 \hB

\subsection{More normalizations}\label{jul1705}

In order to compare calculations of higher torsion invariants one must be careful with normalizations of characteristic forms for flat bundles and the corresponding normalizations of torsion forms.
In this subsection we compare various normalizations occuring in the literature.
We give explicit formulas for the renormalizing factors needed to transfer between different normalizations.

The Bismut-Lott normalization is related to the choice (\ref{jun102})
$$\omega_{2j+1}(h^{V_{\sigma}})=\frac{1}{(2\pi i)^{j}2^{2j+1}} \Tr \ \omega(h^{V_{\sigma}})^{2j+1}\ .$$ This is the standard choice adopted in the present paper.
Another normalization, the Chern normalization, is fixed by taking
$$\omega^{Chern}_{2j+1}(h^{V_{\sigma}})_{2j+1}:=\Imm\left(\tilde\ch(\nabla^{V_{\sigma},*},\nabla^{V_{\sigma}})_{2j+1}\right)\ .$$
Here we use the notation introduced in Subsection \ref{klejflkewfjwlfjlwefwefwef} and
$\tilde \ch(\nabla_{0},\nabla_{1})$ is the transgression Chern form such that
$d\tilde \ch(\nabla_{0},\nabla_{1})=\ch(\nabla_{0})-\ch(\nabla_{1})$.
{The proof of Lemma \ref{mar1001} gives}
$$\tilde\ch(\nabla^{V_{\sigma},*},\nabla^{V_{\sigma}})_{2j+1}=(-1)^{j}\frac{j!}{(2\pi i)^{j+1}(2j+1)!} \Tr\   \omega(h^{V_{\sigma}})^{2j+1}$$
and therefore
\begin{equation}\label{jul1501}\omega^{Chern}_{2j+1}(h^{V_{\sigma}})_{2j+1}=(-1)^{j}\frac{2^{2j+1}j!}{2\pi(2j+1)!}\omega_{2j+1}(h^{V_{\sigma}})\ .\end{equation}
We define
the factor
\begin{equation}\label{jul1604}N_{Chern}(2j+1):=(-1)^{j}\frac{2\pi(2j+1)!}{2^{2j+1}j!}\ .\end{equation}

The Chern normalization is used e.g. in the work of Goette \cite{MR2674876}.
The normalization of Igusa \cite{MR1945530} is
$$\omega^{Igusa}_{2j+1}(h^{V_{\sigma}})=\frac{1}{(2j+1)! \ 2 \ i^{j}} \Tr  \: \omega(h^{V_{\sigma}})^{2j+1}\ .$$
It follows that
$$\omega^{Igusa}_{2j+1}(h^{V_{\sigma}})=\frac{(2\pi)^{j}2^{2j}}{(2j+1)!} \ \omega_{2j+1}(h^{V_{\sigma}})\ .$$
The Igusa normalization is used for the topological version of higher torsion, the higher Reidemeister-Franz torsion, which is also called Igusa-Klein torsion.
We define the factor
\begin{equation}\label{jul1605}N_{Igusa}(2j+1):=\frac{(2j+1)!}{(2\pi)^{j}2^{2j}} \ .\end{equation}
Finally we have the Borel normalization fixed by the normalization of the Borel regulator
\cite{MR1869655}.
We have by Proposition \ref{jul0231}
$$\omega^{Borel}_{2j+1}(h^{V_{\sigma}}):=\frac{(-1)^{j}(2\pi i)^{j}j!}{(2j+1)!} \ \omega_{2j+1}(h^{V_{\sigma}})\ .$$
Hence we define
\begin{equation}\label{jul1606}N_{Borel}(2j+1):=\frac{(-1)^{j}(2j+1)!}{(2\pi i)^{j} j!}\end{equation}

\bibliographystyle{halpha}
\bibliography{differential}

\begin{thebibliography}{EKMM97}

\bibitem[Ada63]{MR0159336}
J.F. Adams.
\newblock On the groups {$J(X)$}. {I}.
\newblock {\em Topology}, 2:181--195, 1963.

\bibitem[AHW15]{Asok:2015uq}
Aravind Asok, Marc Hoyois, and Matthias Wendt.
\newblock Affine representability results in $\mathbf{A}^1$-homotopy theory i:
  vector bundles, 06 2015, 1506.07093.

\bibitem[Bas62]{Bass62}
H.~Bass.
\newblock Torsion free and projective modules.
\newblock {\em Trans. AMS}, 102(2):319--327, 1962.

\bibitem[BDKW11]{MR2739777}
B.~Badzioch, W.~Dorabia{\l}a, J.R. Klein, and B.~Williams.
\newblock Equivalence of higher torsion invariants.
\newblock {\em Adv. Math.}, 226(3):2192--2232, 2011.

\bibitem[Be{\u\i}86]{MR862627}
A.A. Be{\u\i}linson.
\newblock Higher regulators of modular curves.
\newblock In {\em Applications of algebraic {$K$}-theory to algebraic geometry
  and number theory, {P}art {I}, {II} ({B}oulder, {C}olo., 1983)}, volume~55 of
  {\em Contemp. Math.}, pages 1--34. Amer. Math. Soc., Providence, RI, 1986.

\bibitem[BF78]{MR513569}
A.K. Bousfield and E.M. Friedlander.
\newblock Homotopy theory of {$\Gamma $}-spaces, spectra, and bisimplicial
  sets.
\newblock In {\em Geometric applications of homotopy theory ({P}roc. {C}onf.,
  {E}vanston, {I}ll., 1977), {II}}, volume 658 of {\em Lecture Notes in Math.},
  pages 80--130. Springer, Berlin, 1978.

\bibitem[BG75]{MR0377873}
J.C. Becker and D.H. Gottlieb.
\newblock The transfer map and fiber bundles.
\newblock {\em Topology}, 14:1--12, 1975.

\bibitem[BG01]{MR1867006}
J.-M. Bismut and S.~Goette.
\newblock Families torsion and {M}orse functions.
\newblock {\em Ast\'erisque}, (275):x+293, 2001.

\bibitem[BG02]{MR1869655}
J.I. Burgos~Gil.
\newblock {\em The regulators of {B}eilinson and {B}orel}, volume~15 of {\em
  CRM Monograph Series}.
\newblock American Mathematical Society, Providence, RI, 2002.

\bibitem[BGV92]{MR1215720}
N.~Berline, E.~Getzler, and M.~Vergne.
\newblock {\em Heat kernels and {D}irac operators}, volume 298 of {\em
  Grundlehren der Mathematischen Wissenschaften [Fundamental Principles of
  Mathematical Sciences]}.
\newblock Springer-Verlag, Berlin, 1992.

\bibitem[BKS10]{MR2674652}
U.~Bunke, M.~Kreck, and T.~Schick.
\newblock A geometric description of differential cohomology.
\newblock {\em Ann. Math. Blaise Pascal}, 17(1):1--16, 2010.

\bibitem[BL95]{MR1303026}
J.-M. Bismut and J.~Lott.
\newblock Flat vector bundles, direct images and higher real analytic torsion.
\newblock {\em J. Amer. Math. Soc.}, 8(2):291--363, 1995.

\bibitem[BM06]{MR2231056}
M.-T. Benameur and M.~Maghfoul.
\newblock Differential characters in {$K$}-theory.
\newblock {\em Differential Geom. Appl.}, 24(4):417--432, 2006.

\bibitem[BN14]{Bunke:2014fk}
Ulrich Bunke and Thomas Nikolaus.
\newblock Twisted differential cohomology, 06 2014,
  \url{http://arxiv.org/abs/1406.3231}.

\bibitem[BNV14]{Bunke:2013uq}
Ulrich Bunke, Thomas Nikolaus, and Michael V{\"o}lkl.
\newblock Differential cohomology theories as sheaves of spectra.
\newblock {\em Journal of Homotopy and related structures}, 2014,
  \url{http://arxiv.org/abs/1311.3188}.

\bibitem[Bor74]{MR0387496}
A.~Borel.
\newblock Stable real cohomology of arithmetic groups.
\newblock {\em Ann. Sci. \'Ecole Norm. Sup. (4)}, 7:235--272 (1975), 1974.

\bibitem[Bor77]{MR0506168}
A.~Borel.
\newblock Cohomologie de {${\rm SL}_{n}$} et valeurs de fonctions zeta aux
  points entiers.
\newblock {\em Ann. Scuola Norm. Sup. Pisa Cl. Sci. (4)}, 4(4):613--636, 1977.

\bibitem[Bou79]{MR551009}
A.K. Bousfield.
\newblock The localization of spectra with respect to homology.
\newblock {\em Topology}, 18(4):257--281, 1979.

\bibitem[BS98]{MR1621939}
J.C. Becker and R.E. Schultz.
\newblock Axioms for bundle transfers and traces.
\newblock {\em Math. Z.}, 227(4):583--605, 1998.

\bibitem[BS09]{MR2664467}
U.~Bunke and T.~Schick.
\newblock Smooth {$K$}-theory.
\newblock {\em Ast\'erisque}, (328):45--135 (2010), 2009.

\bibitem[BS10]{MR2608479}
U.~Bunke and T.~Schick.
\newblock Uniqueness of smooth extensions of generalized cohomology theories.
\newblock {\em J. Topol.}, 3(1):110--156, 2010.

\bibitem[BSGT]{bst}
U.~Bunke, F.~Strunk, and G-Tamme.
\newblock The glueing lemma revisited.
\newblock in preparation, 2016.

\bibitem[BSSW09]{MR2550094}
U.~Bunke, T.~Schick, I.~Schr{\"o}der, and M.~Wiethaup.
\newblock Landweber exact formal group laws and smooth cohomology theories.
\newblock {\em Algebr. Geom. Topol.}, 9(3):1751--1790, 2009.

\bibitem[BT13]{Bunke:2013kx}
Ulrich Bunke and Georg Tamme.
\newblock Multiplicative differential algebraic k-theory and applications, 11
  2013, \url{http://arxiv.org/abs/1311.1421}.

\bibitem[BT15]{buta}
Ulrich Bunke and Georg Tamme.
\newblock Regulators and cycle maps in higher-dimensional differential
  algebraic {$K$}-theory.
\newblock {\em Adv. Math.}, 285:1853--1969, 2015.

\bibitem[Bun]{skript}
U.~Bunke.
\newblock Differential cohomology.
\newblock Course notes, Universit\"at Regensburg, 2012.
  \url{http://arxiv.org/abs/1208.3961}.

\bibitem[Bun09]{MR2191484}
U.~Bunke.
\newblock Index theory, eta forms, and {D}eligne cohomology.
\newblock {\em Mem. Amer. Math. Soc.}, 198(928):vi+120, 2009.

\bibitem[Bun10a]{MR2740650}
U.~Bunke.
\newblock Adams operations in smooth {$K$}-theory.
\newblock {\em Geom. Topol.}, 14(4):2349--2381, 2010.

\bibitem[Bun10b]{MR2734150}
U.~Bunke.
\newblock Chern classes on differential {$K$}-theory.
\newblock {\em Pacific J. Math.}, 247(2):313--322, 2010.

\bibitem[BW87]{MR921487}
M.~B{\"o}kstedt and F.~Waldhausen.
\newblock The map {$BSG\to A(*)\to QS^0$}.
\newblock In {\em Algebraic topology and algebraic {$K$}-theory ({P}rinceton,
  {N}.{J}., 1983)}, volume 113 of {\em Ann. of Math. Stud.}, pages 418--431.
  Princeton Univ. Press, Princeton, NJ, 1987.

\bibitem[BZ92]{MR1185803}
J.-M. Bismut and W.~Zhang.
\newblock An extension of a theorem by {C}heeger and {M}\"uller.
\newblock {\em Ast\'erisque}, (205):235, 1992.
\newblock With an appendix by Fran{\c{c}}ois Laudenbach.

\bibitem[Che79]{MR528965}
J.~Cheeger.
\newblock Analytic torsion and the heat equation.
\newblock {\em Ann. of Math. (2)}, 109(2):259--322, 1979.

\bibitem[CS1v]{MR827262}
J.~Cheeger and J.~Simons.
\newblock Differential characters and geometric invariants.
\newblock In {\em Geometry and topology (College Park, Md., 1983/84)}, volume
  1167 of {\em Lecture Notes in Math.}, pages 50--80. Springer, Berlin, 1v.

\bibitem[DHI04]{MR2034012}
Daniel Dugger, Sharon Hollander, and Daniel~C. Isaksen.
\newblock Hypercovers and simplicial presheaves.
\newblock {\em Math. Proc. Cambridge Philos. Soc.}, 136(1):9--51, 2004.

\bibitem[Dou06]{MR2248973}
C.L. Douglas.
\newblock On the fibrewise {P}oincare-{H}opf theorem.
\newblock In {\em Recent developments in algebraic topology}, volume 407 of
  {\em Contemp. Math.}, pages 101--111. Amer. Math. Soc., Providence, RI, 2006.

\bibitem[Dug01]{MR1870515}
D.~Dugger.
\newblock Universal homotopy theories.
\newblock {\em Adv. Math.}, 164(1):144--176, 2001.

\bibitem[DWW03]{MR1982793}
W.~Dwyer, M.~Weiss, and B.~Williams.
\newblock A parametrized index theorem for the algebraic {$K$}-theory {E}uler
  class.
\newblock {\em Acta Math.}, 190(1):1--104, 2003.

\bibitem[EKMM97]{MR1417719}
A.D. Elmendorf, I.~Kriz, M.A. Mandell, and J.P. May.
\newblock {\em Rings, modules, and algebras in stable homotopy theory},
  volume~47 of {\em Mathematical Surveys and Monographs}.
\newblock American Mathematical Society, Providence, RI, 1997.
\newblock With an appendix by M. Cole.

\bibitem[Esn89]{MR1014822}
H.~Esnault.
\newblock On the {L}oday symbol in the {D}eligne-{B}e\u\i linson cohomology.
\newblock {\em $K$-Theory}, 3(1):1--28, 1989.

\bibitem[FH00]{MR1769477}
D.S. Freed and M.~Hopkins.
\newblock On {R}amond-{R}amond fields and {$K$}-theory.
\newblock {\em J. High Energy Phys.}, (5):Paper 44, 14, 2000.

\bibitem[FL10]{MR2602854}
D.S. Freed and J.~Lott.
\newblock An index theorem in differential {$K$}-theory.
\newblock {\em Geom. Topol.}, 14(2):903--966, 2010.

\bibitem[FMS07]{MR2286784}
D.S. Freed, G.W. Moore, and G.~Segal.
\newblock Heisenberg groups and noncommutative fluxes.
\newblock {\em Ann. Physics}, 322(1):236--285, 2007.

\bibitem[Fre00]{MR1919425}
D.S. Freed.
\newblock Dirac charge quantization and generalized differential cohomology.
\newblock In {\em Surveys in differential geometry}, Surv. Differ. Geom., VII,
  pages 129--194. Int. Press, Somerville, MA, 2000.

\bibitem[GGN13]{2013arXiv1305.4550G}
D.~{Gepner}, M.~{Groth}, and T.~{Nikolaus}.
\newblock {Universality of multiplicative infinite loop space machines}.
\newblock \href{http://arxiv.org/abs/1305.4550}{arXiv:1305.4550}, 2013.

\bibitem[GI14]{2010arXiv1011.4653G}
Sebastian Goette and Kiyoshi Igusa.
\newblock Exotic smooth structures on topological fiber bundles {II}.
\newblock {\em Trans. Amer. Math. Soc.}, 366(2):791--832, 2014.

\bibitem[Goe09]{MR2674876}
S.~Goette.
\newblock Torsion invariants for families.
\newblock {\em Ast\'erisque}, (328):161--206 (2010), 2009.

\bibitem[HK03]{MR2002643}
A.~Huber and G.~Kings.
\newblock Bloch-{K}ato conjecture and {M}ain {C}onjecture of {I}wasawa theory
  for {D}irichlet characters.
\newblock {\em Duke Math. J.}, 119(3):393--464, 2003.

\bibitem[HS05]{MR2192936}
M.J. Hopkins and I.M. Singer.
\newblock Quadratic functions in geometry, topology, and {M}-theory.
\newblock {\em J. Differential Geom.}, 70(3):329--452, 2005.

\bibitem[Igu02]{MR1945530}
K.~Igusa.
\newblock {\em Higher {F}ranz-{R}eidemeister torsion}, volume~31 of {\em AMS/IP
  Studies in Advanced Mathematics}.
\newblock American Mathematical Society, Providence, RI, 2002.

\bibitem[Igu08]{MR2365656}
K.~Igusa.
\newblock Axioms for higher torsion invariants of smooth bundles.
\newblock {\em J. Topol.}, 1:159--186, 2008.

\bibitem[Kim14]{Kim:2014fk}
Youngsoo Kim.
\newblock Standard vector bundles, 04 2014,
  \url{http://arxiv.org/abs/1404.1485}.

\bibitem[KS10]{MR2775352}
M.~Kreck and W.~Singhof.
\newblock Homology and cohomology theories on manifolds.
\newblock {\em M\"unster J. Math.}, 3:1--9, 2010.

\bibitem[Lot94]{MR1297676}
J.~Lott.
\newblock Equivariant analytic torsion for compact {L}ie group actions.
\newblock {\em J. Funct. Anal.}, 125(2):438--451, 1994.

\bibitem[Lot00]{MR1724894}
J.~Lott.
\newblock Secondary analytic indices.
\newblock In {\em Regulators in analysis, geometry and number theory}, volume
  171 of {\em Progr. Math.}, pages 231--293. Birkh\"auser Boston, Boston, MA,
  2000.

\bibitem[Lur]{highalg}
J.~Lurie.
\newblock Higher algebra.
\newblock \url{http://www.math.harvard.edu/~lurie/papers/HigherAlgebra.pdf}.

\bibitem[Lur09]{MR2522659}
J.~Lurie.
\newblock {\em Higher topos theory}, volume 170 of {\em Annals of Mathematics
  Studies}.
\newblock Princeton University Press, Princeton, NJ, 2009.

\bibitem[Mil71]{MR0349811}
J.~Milnor.
\newblock {\em Introduction to algebraic {$K$}-theory}.
\newblock Princeton University Press, Princeton, N.J., 1971.
\newblock Annals of Mathematics Studies, No. 72.

\bibitem[M{\"u}l78]{MR498252}
W.~M{\"u}ller.
\newblock Analytic torsion and {$R$}-torsion of {R}iemannian manifolds.
\newblock {\em Adv. in Math.}, 28(3):233--305, 1978.

\bibitem[M{\"u}l93]{MR1189689}
W.~M{\"u}ller.
\newblock Analytic torsion and {$R$}-torsion for unimodular representations.
\newblock {\em J. Amer. Math. Soc.}, 6(3):721--753, 1993.

\bibitem[Neu88]{MR944995}
J.~Neukirch.
\newblock The {B}e\u\i linson conjecture for algebraic number fields.
\newblock In {\em Be\u\i linson's conjectures on special values of
  {$L$}-functions}, volume~4 of {\em Perspect. Math.}, pages 193--247. Academic
  Press, Boston, MA, 1988.

\bibitem[Neu99]{MR1697859}
J.~Neukirch.
\newblock {\em Algebraic number theory}, volume 322 of {\em Grundlehren der
  Mathematischen Wissenschaften [Fundamental Principles of Mathematical
  Sciences]}.
\newblock Springer-Verlag, Berlin, 1999.
\newblock Translated from the 1992 German original and with a note by Norbert
  Schappacher, With a foreword by G. Harder.

\bibitem[Qui73]{MR0338129}
D.~Quillen.
\newblock Higher algebraic {$K$}-theory. {I}.
\newblock In {\em Algebraic {$K$}-theory, {I}: {H}igher {$K$}-theories ({P}roc.
  {C}onf., {B}attelle {M}emorial {I}nst., {S}eattle, {W}ash., 1972)}, pages
  85--147. Lecture Notes in Math., Vol. 341. Springer, Berlin, 1973.

\bibitem[RS71]{MR0295381}
D.B. Ray and I.M. Singer.
\newblock {$R$}-torsion and the {L}aplacian on {R}iemannian manifolds.
\newblock {\em Advances in Math.}, 7:145--210, 1971.

\bibitem[Seg68]{MR0232393}
G.~Segal.
\newblock Classifying spaces and spectral sequences.
\newblock {\em Inst. Hautes \'Etudes Sci. Publ. Math.}, (34):105--112, 1968.

\bibitem[Shi07]{MR2306038}
B.~Shipley.
\newblock {$H\Bbb Z$}-algebra spectra are differential graded algebras.
\newblock {\em Amer. J. Math.}, 129(2):351--379, 2007.

\bibitem[SS08]{MR2365651}
J.~Simons and D.~Sullivan.
\newblock Axiomatic characterization of ordinary differential cohomology.
\newblock {\em J. Topol.}, 1(1):45--56, 2008.

\bibitem[SS10]{MR2732065}
J.~Simons and D.~Sullivan.
\newblock Structured vector bundles define differential {$K$}-theory.
\newblock In {\em Quanta of maths}, volume~11 of {\em Clay Math. Proc.}, pages
  579--599. Amer. Math. Soc., Providence, RI, 2010.

\bibitem[Swi02]{MR1886843}
R.M. Switzer.
\newblock {\em Algebraic topology---homotopy and homology}.
\newblock Classics in Mathematics. Springer-Verlag, Berlin, 2002.
\newblock Reprint of the 1975 original [Springer, New York; MR0385836 (52
  \#6695)].

\bibitem[Wal85]{MR802796}
F.~Waldhausen.
\newblock Algebraic {$K$}-theory of spaces.
\newblock In {\em Algebraic and geometric topology ({N}ew {B}runswick,
  {N}.{J}., 1983)}, volume 1126 of {\em Lecture Notes in Math.}, pages
  318--419. Springer, Berlin, 1985.

\bibitem[Wei13]{weibel}
Charles~A. Weibel.
\newblock {\em The {$K$}-book}, volume 145 of {\em Graduate Studies in
  Mathematics}.
\newblock American Mathematical Society, Providence, RI, 2013.
\newblock An introduction to algebraic $K$-theory.

\end{thebibliography}
\end{document}